\documentclass[12pt]{article}
\usepackage{color}
\usepackage{amsbsy}
\usepackage{amsmath}
\usepackage{amssymb}
\usepackage{latexsym}
\usepackage{comment}
\usepackage{graphicx}
\usepackage{amsfonts}
\usepackage{url}
\usepackage{hyperref}
\usepackage[hyphenbreaks]{breakurl}
\usepackage{enumitem}
\usepackage{bbm}
\usepackage{enumerate}
\usepackage{yhmath}
\usepackage{tikz}
\usetikzlibrary{matrix}
\def\singlespace{\def\baselinestretch{1}\@normalsize}

\usepackage[square,numbers]{natbib}

\usepackage{geometry}

\renewcommand{\baselinestretch} {1.2}
\makeatletter 
\def\singlespace{\def\baselinestretch{1}\@normalsize}

\@addtoreset{equation}{section}
\renewcommand{\theequation} {\arabic{section}.\arabic{equation}}

\renewcommand{\thefootnote}{\fnsymbol{footnote}}

\newtheorem{theorem}{Theorem}
\newtheorem{corollary}{Corollary}
\newtheorem{lemma}{Lemma}
\newtheorem{remark}{Remark}
\newtheorem{proposition}{Proposition}
\def\graphheight{6in}                                    
\def\graphwidth{6in}                                    

\def\conP{\stackrel{\cal {P}} {\to}}                     
\def\conD{\stackrel{\cal D} {\to}}                         

\def\pf{\textsl{Proof}:\ }                          
\def\endpf{$\blacksquare$}                                 

\def\ni{\noindent}
\def\tilde{\widetilde}
\def\hat{\widehat}
\def\Widehat{\widehat}
\def\bfa{\mathbf{a}}

\def\bfc{\boldsymbol{c}}
\def\bfe{\mathbf{e}}

\def\bfI{\mathbf{I}}

\def\E{\mathrm{E}}
\def\bfW{\mathbf{W}}

\def\bfY{\boldsymbol{Y}}
\def\bfv{\boldsymbol{v}}
\def\bfW{\mathbf{W}}
\def\bfV{\mathbf{V}}
\def\bfZ{\boldsymbol{Z}}

\def\bfzero{\mathbf{0}}
\def\bfone{\mathbf{1}}
\def\beps{\boldsymbol{\epsilon}}
\def\bbeta{\boldsymbol{\beta}}
\def\bgamma{\boldsymbol{\gamma}}
\def\bfb{\boldsymbol{b}}
\def\bfa{\boldsymbol{a}}
\def\bdelta{\boldsymbol{\delta}}

\def\bfx{\boldsymbol{x}}

\def\bfX{\boldsymbol{X}}
\def\bmu{\boldsymbol{\mu}}

\def\SNP{\mathrm{SNP}}

\def\bfbeta{\boldsymbol{\beta}}

\def\ID{\,\mathrm {I}}                                 

\def\KS{\mathrm{KS}}

\def\tr{\mathrm{tr}}                                   
\def\diag{\mathrm{diag}}

\def\MM{\mathrm{MM}}
\def\cov{\mathrm{Cov}}
\def\var{\mathrm{Var}}

\def\SNR{\mathrm{SNR}}
\def\nullH0{\mathrm{null}}

\def\iid{\mathrm{i.i.d.}}
\def\I{\mathrm{I}}
\def\II{\mathrm{II}}
\def\III{\mathrm{III}}
\def\IV{\mathrm{IV}}
\def\pr{\mathrm{P}}

\def\KL{\mathrm{KL}}
\def\OLS{\mathrm{OLS}}
\def\MLE{\mathrm{MLE}}

\def\mathbbR{\mathbb{R}}

\begin{document}

\title
{
Moderate-Dimensional Inferences on Quadratic Functionals in Ordinary Least Squares
}

\author{Xiao Guo\footnotemark[1] \footnotemark[3]
\ and    Guang Cheng\footnotemark[2] \footnotemark[3] }
\renewcommand{\thefootnote}{\fnsymbol{footnote}}

\footnotetext[1]{International Institute of Finance, School of Management,
University of Science and Technology of China,
Hefei, Anhui 230026, People's Republic of China. Research Sponsored by the National Natural Science Foundation of China, grants 11601500, 11671374 and 11771418, and the Fundamental Research Funds for the Central Universities. Email: xiaoguo@ustc.edu.cn.}

\footnotetext[2]{Department of Statistics,
Purdue University, West Lafayette, IN 47906. Research Sponsored by NSF DMS-1712907, DMS-1811812, DMS-1821183, and Office of Naval Research, (ONR N00014-18-2759). Email: chengg@purdue.edu.}

\footnotetext[3]{We thank Mr. Ching-Wei Cheng for implementing our codes in Purdue supercomputer.}

\date{}

\maketitle

\begin{abstract}
Statistical inferences for quadratic functionals of linear regression parameter have found wide applications including signal detection, global testing, inferences of error variance and fraction of variance explained. Classical theory based on ordinary least squares estimator works perfectly in the low-dimensional regime, but fails when the parameter dimension $p_n$ grows proportionally to the sample size $n$. In some cases, its performance is not satisfactory even when $n\ge 5p_n$.

The main contribution of this paper is to develop {\em dimension-adaptive} inferences for quadratic functionals when $\lim_{n\to \infty} p_n/n=\tau\in[0,1)$. We propose a bias-and-variance-corrected test statistic and demonstrate that its theoretical validity (such as consistency and asymptotic normality) is adaptive to both low dimension with $\tau = 0$ and moderate dimension with $\tau \in(0, 1)$. Our general theory holds, in particular, without Gaussian design/error or structural parameter assumption, and applies to a broad class of quadratic functionals covering all aforementioned applications. As a by-product, we find that the classical fixed-dimensional results continue to hold {\em if and only if} the signal-to-noise ratio is large enough, say when $p_n$ diverges but slower than $n$. Extensive numerical results demonstrate the satisfactory performance of the proposed methodology even when $p_n\ge 0.9n$ in some extreme cases. The mathematical arguments are based on the random matrix theory and leave-one-observation-out method.
\end{abstract}

\section{Introduction} \label{Section-1}
The linear regression model is one of the most widely used statistical tools to discover the relation  between a continuous response and a class of explanatory variables in different scientific areas. Specifically, we
consider
\begin{equation} \label{1.1}
Y_i = \bfX_i^T \bfbeta_0 + \epsilon_i, \quad \text{for $i = 1, \ldots, n$,}
\end{equation}
where $\bfbeta_0 = (\beta_{0,1}, \ldots, \beta_{0, p_n})^T \in \mathbbR^{p_n}$ is an unknown vector of parameters, and $\{\epsilon_i\}_{i=1}^n$ are $\iid$ errors independent of $\{\bfX_i\}_{i=1}^{n}$ with $\E(\epsilon_i) = 0$ and $\var(\epsilon_i)=\sigma_{\epsilon}^2$. We assume $\{Y_i, \bfX_i\}_{i=1}^n$ are $\iid$ observations with $\E(\bfX_i) = \bfzero_{p_n}$ and $\cov(\bfX_i)= \Sigma$, without imposing any specific distributional assumption on either $\bfX_i$ or $\epsilon_i$ throughout this paper. Denoting $\bfY = (Y_1, \ldots, Y_n)^T$, $X = (\bfX_1, \ldots, \bfX_n)^T$ and $\beps = (\epsilon_1, \ldots, \epsilon_n)^T$, (\ref{1.1}) can be re-expressed as $$\bfY = X \bbeta_0 + \beps.$$

For fixed dimension, statistical estimation and inference for $\bbeta_0$ and $\sigma_{\epsilon}^2$ have been well studied based on the ordinary least squares ($\OLS$) estimator,
$$\hat\bbeta = (X^TX)^{-1}X^T \bfY.$$ In the modern high-dimensional regime, the parameter dimension $p_n$ is allowed to be much larger than $n$, e.g. $(\log p_n) / n = o(1)$, but in most cases the number of non-zero elements in $\bbeta_0$ is a vanishing fraction of $n$. Such a sparsity condition is commonly assumed in the high-dimensional literature, e.g., \cite{Meinshausen_Yu_2009, van_de_Geer_2008, Zhang_Huang_2008} on oracle inequality and parameter estimation; \cite{Tibshirani_1996, Fan_Lv_2008, Meinshausen_Buhlmann_2006} on variable selection and \cite{Javanmard_Montanari_2014, van_de_Geer_etal_2014, Zhang_Zhang_2014} on statistical inference. However, in reality, $p_n$ may be moderately large, i.e., of  the same magnitude as $n$, and $\bbeta_0$ is not necessarily sparse. One example is the genomic study, where the number of {\em significantly identified} genes with association in {\em trans}, i.e., $p_n = 108$, is moderately large compared with $n=270$; see \cite{Stranger_etal_2007}.

For moderate dimension with $\lim_{n\to\infty} p_n/n = \tau \in(0,1)$, which is of major concern in this paper, some classical statistical inference procedures developed for fixed-dimensional data are no longer valid. For example, when $p_n$ is fixed, we can test
\begin{eqnarray} \label{1.2}
H_0: \|\bbeta_0\|_2=c_0 \quad \text{versus} \quad H_1: \|\bbeta_0\|_2\neq c_0,
\end{eqnarray}
for a known constant $c_0 \geq 0$, by calculating the Z-score
\begin{eqnarray}\label{1.3}
\mathbb{Z}_0 = \frac{\|\hat\bfbeta\|_2^2 - c_0^2}
{\hat\zeta_0},
\end{eqnarray}
 where
\begin{eqnarray} \label{1.4}
\hat \zeta_0^2=  4 \hat \sigma_{\epsilon}^2  \hat\bbeta^T (X^TX)^{-1} \hat\bbeta\quad \mbox{and}\quad \hat \sigma_{\epsilon}^2 = \frac{\|\bfY- X \hat\bbeta\|_2^2}{n-p_n}.
\end{eqnarray}
Under the null hypothesis, $\mathbb{Z}_0\conD N(0, 1)$; see Theorem \ref{Theorem4}. Hence, the P-value for testing \eqref{1.2} is $2\Phi(-|\mathbbm{z}_0|)$, where $\mathbbm{z}_0$ is a realization of $\mathbb{Z}_0$ and $\Phi(\cdot)$ is the cumulative distribution function of the standard normal distribution.

We next examine the empirical performance of the conventional Z-test by setting $n = 1000$ with $p_n = 10$ for fixed dimension and  $p_n=200$, $500$ and $900$ for moderate dimension. Consider $\bfX_i \stackrel{\iid}{\sim}  N(\bfzero_{p_n}, \bfI_{p_n})$ and $\epsilon_i \stackrel{\iid}{\sim} N(0, 1 )$, where $\bfI_{p_n}$ denotes the $p_n\times p_n$ identity matrix. The true parameter $\beta_{0,j}$'s were generated independently from $\text{Unif}(0, 1)$, and 40000 replications were conducted in each setup. The plots of the P-values under the valid null hypothesis are given in the top panels of Figure \ref{Figure_1} below. The uniform distribution of the P-values when $p_n=10$ is consistent with the classical fixed-dimensional theory. But for $p_n=200$, $ 500$ and $900$, P-values are relatively concentrated around $0$. We further test the uniform distribution of the P-values by the formal Kolmogorov-Smirnov ($\KS$) test \citep{Kolmogorov_1933, Smirnov_1939}, and find that the P-values for $p_n=10$, $200$, $500$ and $900$ are $0.2518$,  $8.05\times 10^{-68}$, $0$ and $0$, respectively. Hence, the naive Z-score does not work under moderate dimension, say even when $n\ge 5p_n$.
\begin{figure}[!ht]
\centering
\includegraphics[height = 2.6in, width = 6.8in]{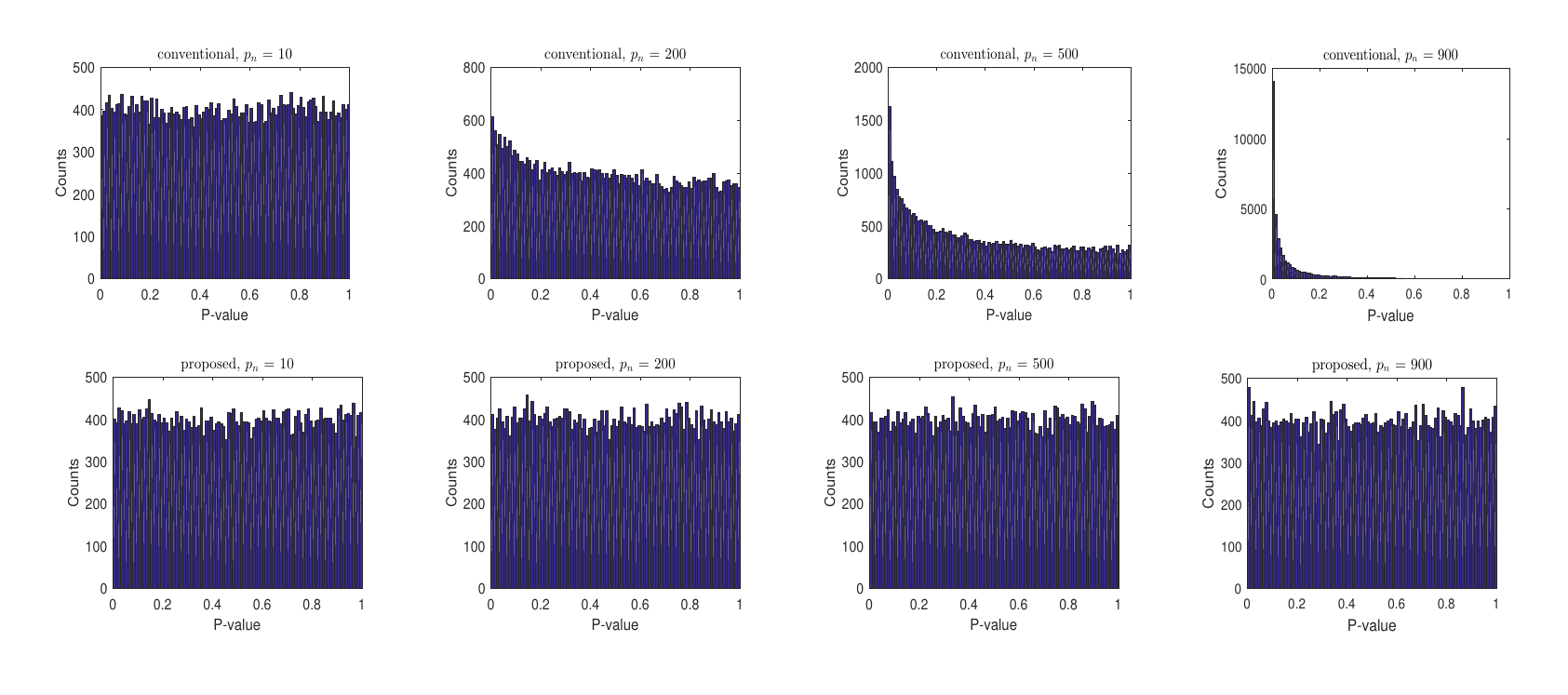}
\caption{\textsl{P-values of $\mathbb{Z}_0$ (top panels) and $\mathbb{Z}_n$ (bottom panels). The  panels from left to right are for $p_n = 10/200/500/900$.}}
\label{Figure_1}
\end{figure}

The main focus of this paper is on the moderate-dimensional inference without imposing any type of structural conditions on $\bbeta_0$ and $\Sigma$, while our results are also adaptive to the low-dimensional case with $\tau = 0$\footnote{We call it low-dimensional regime  when $p_n\to\infty$ but $p_n/n\to 0$. Hence, both fixed- and low-dimensional regimes correspond to that $\tau=0$.}.  Specifically, we conduct statistical inferences for a class of quadratic functionals such as $\|\bbeta_0\|_2^2$ and $\sigma_{\epsilon}^2$, which cover a wide range of applications including signal detection and global testing.
A related line of work is the study of the signal strength $\bbeta_0^T \Sigma \bbeta_0$ by \cite{Dicker_2014} and \cite{Janson_etal_2017}. However, their procedures crucially rely on the fact that $Y_i\sim N(0, \bbeta_0^T \Sigma \bbeta_0+\sigma_\epsilon^2)$, and their theoretical results hold only when $\bfX_i$ and $\epsilon_i$ are both Gaussian. Hence, their results are not readily carried over into our case, e.g., two-sample inference. Additionally, different tools such as leave-one-observation-out method \cite{El_Karoui_2013, El_Karoui_2017} are used in our paper. Please see more discussions in the end of Section \ref{Section-3.4}. As a side remark, we point out that the classical fixed-dimensional inference may still be applied to the low-dimensional regime {\em if and only if} the signal-to-noise ratio $\SNR:={\var(\bfX_i^T \bbeta_0)}/{\var(\epsilon_i)}= {\bbeta_0^T \Sigma \bbeta_0}/{\sigma_{\epsilon}^2}$ is large. However, the strength of the $\SNR$ cannot be directly examined in practice. Hence, the adaptiveness of our proposed method (without relying on $\SNR$) is practically important. In case of interest, readers may refer to Figure \ref{Figure_11} in Section \ref{app-A2} for the precise relation between $\tau$ and $\SNR$.

\begin{figure}[!ht]
\centering
\includegraphics[height = 2in, width = 6.5 in]{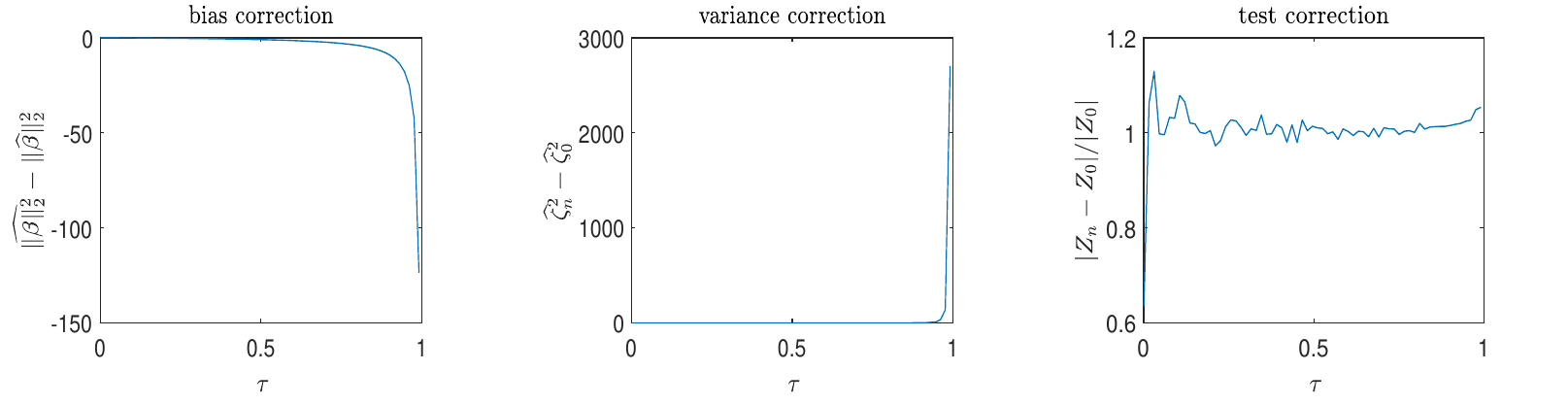}
\caption{\textsl{Amount of empirical corrections of bias (left penal) and variance (middle penal) versus $\tau$ for $\hat{\|\bbeta\|_2^2}$ compared with $\|\hat\bbeta\|_2^2$. The right panel plots  $|\mathbb{Z}_n -\mathbb{Z}_0|/|\mathbb{Z}_0| $ versus $\tau$.
}}
\label{Figure_2}
\end{figure}

Our primary contribution is to propose a bias-corrected estimator $\hat{\|\bbeta\|_2^2}$ for $\|\bbeta_0\|_2^2$  in \eqref{2.3}, based on which a bias-and-variance-corrected test statistic $\mathbb{Z}_n$ is developed in \eqref{2.5}. The bottom panels of Figure \ref{Figure_1} plot the P-values of $\mathbb{Z}_n$ for $p_n=10$, $200$, $500$ and $900$. The P-values of the $\KS$ test for the uniformity are 0.4755, 0.1175, 0.8972 and 0.2672 correspondingly.
Figure \ref{Figure_2} plots the amount of empirical corrections of bias and variance needed in $\hat{\|\bbeta\|_2^2}$ under the same setting. It reveals that the bias correction tends to $-\infty$ as $\tau \to 1$, while the  variance correction diverges to $\infty$. The right panel of Figure \ref{Figure_2} plots the relative difference between $\mathbb{Z}_n$ and $\mathbb{Z}_0$ versus $\tau$. As $\tau$ deviates from zero, the amount of correction rapidly increases to its largest value, and then decreases and stabilizes around $1$.
As an immediate application, global testing
\begin{eqnarray} \label{1.5}
H_0: \bfbeta_0 =  \bbeta_0^{\nullH0}\quad \mbox{versus}\quad H_1: \bfbeta_0 \neq  \bbeta_0^{\nullH0},
 \end{eqnarray}
 can also be performed
with a bias-and-variance-corrected version of $\|\hat \bbeta - \bbeta_0^{\nullH0}\|_2^2$ as the test statistic. Please see \cite{Portnoy_1985, Arias-Castro_etal_2011, Zhong_Chen_2011, Zhang_Cheng_2017} for low- and high-dimensional results, respectively.

Our general moderate-dimensional theory can also be applied to other statistical inference problems. For example, we can detect the existence of signal by setting $c_0=0$ in \eqref{1.2}. By formulating a sequence of alternatives  $H_{1n}: \|\bbeta_0\|_2^2 = \delta_n$, we further show that $\delta_n^*:=\sigma_{\epsilon}^2\sqrt{p_n }n^{-1}$ is the smallest separation rate such that successful detection of $H_{1n}$ is still possible, which matches with the minimax detection rate in \cite{Ingster_etal_2010}. As far as we are aware, the existing results concerned with detection boundary only focus on either Gaussian mean models with $p_n = n$, e.g., \cite{Cai_etal_2007, Donoho_Jin_2004, Hall_Jin_2010}, or high-dimensional data, e.g., \cite{Ingster_etal_2010, Arias-Castro_etal_2011}.

New results of inference on the error variance will also be established for moderate dimension. We still use the estimator $\hat \sigma_{\epsilon}^2$ defined in \eqref{1.4} for low-dimensional data, but modify its asymptotic variance as $\zeta_{\epsilon}^2 =  \{\nu_4  + \sigma_{\epsilon}^4  ( 3\tau - 1)/(1-\tau) \} /n$ with $\nu_4 = \E ( \epsilon_i^4)$ to derive that $$\frac{\hat \sigma_{\epsilon}^2  - \sigma_{\epsilon}^2}{\zeta_{\epsilon}}\conD N(0, 1).$$ One related result is concerned with the fraction of variance explained (and also $\SNR$), defined as
\begin{eqnarray}\label{1.6}
\rho_0:=\frac{\bbeta_0^T \Sigma \bbeta_0}{\bbeta_0^T \Sigma \bbeta_0 + \sigma_{\epsilon}^2}=\frac{\SNR}{\SNR+1},
\end{eqnarray}
which describes the proportion of the variance in the dependent variable that is predictable from the independent variable.
The high-dimensional estimation of $\sigma_{\epsilon}^2$, $\rho_0$ and $\SNR$ can be found in \cite{Sun_Zhang_2012, Fan_etal_2012, Verzelen_Gassiat_2018}.

Our results can be naturally extended to two-sample inference. Here, we give two examples in Section S.4
in the supplementary material. Let $\bgamma_0 \in \mathbb{R}^{p_n}$ be the regression parameter in another linear regression model independent of \eqref{1.1}. The first issue is to test the equality of $\bgamma_0 $ and $\bbeta_0$, while the second is concerned with the co-heritability, defined as
\begin{eqnarray} \label{1.7}
\theta_0 = \frac{\bgamma_0 ^T \bbeta_0}{\|\bgamma_0\|_2 \|\bbeta_0\|_2 }.
 \end{eqnarray}
The measure $\theta_0$ is an important concept that characterizes the genetic associations within pairs of quantitative traits, whose high-dimensional estimation has recently been studied in \cite{Guo_etal_2016}. Besides, an immediate application of our arguments is the inference for the linear functionals as discussed in Section \ref{app-A1} of the Appendix.

As a summary, a list of hypotheses in consideration together with potential applications is given below:
\begin{singlespace}
\begin{itemize}
\item Testing the quadratic functional: hypotheses in \eqref{1.2};
\item Signal detection: hypotheses in \eqref{1.2} with $c_0=0$;
\item Global testing: hypotheses in \eqref{1.5};
\item Inference for the error variance $\sigma_{\epsilon}^2$ using Proposition \ref{Proposition-1};
\item Testing the fraction of variance explained (or $\SNR$): hypotheses in \eqref{1.6};
\item Inference for the signal strength $\bbeta_0^T \Sigma \bbeta_0$ using \eqref{3.4};
\item Two-sample inferences: hypotheses in (S.4.2)
and (S.4.3).
\end{itemize}
\end{singlespace}

Our asymptotic normality result relies on the application of the martingale difference central limit theorem (CLT) \cite{Heyde_Brown_1970} to linear-quadratic forms, i.e., $\beps^T A_n \beps + \bfb_n^T \beps$, where $A_n \in \mathbb{R}^{n\times n}$ ($\bfb_n\in \mathbb{R}^n$) is some random matrix (vector) independent of $\beps$.
Although CLT has been studied for quadratic and linear-quadratic forms, to the best of our knowledge, those results can not be directly applied to our problem. For example, \cite{Dicker_Erdogdu_2017} provides the concentration bounds and finite sample
multivariate normal approximation for quadratic forms, but these results are not applicable to the linear-quadratic forms. \cite{deJong_1987} developed CLT for ``clean'' quadratic forms requiring zero elements on the diagonal (see definition 2.1 therein), which however is not satisfied by the linear-quadratic form in our paper. A more related example is the CLT for the linear-quadratic form  in \cite{Kelejian_Prucha_2001} with $A_n $ and $\bfb_n$ being deterministic. An important assumption in \cite{Kelejian_Prucha_2001} is $\|A_n\|_{1} \leq C <\infty$ which   is violated in our work (see Section S.1
in the supplementary material for detailed explanations).
Besides, two technical tools have been used in our paper: random matrix theory \cite{Bai_Silverstein_2010} and leave-one-observation-out method \cite{El_Karoui_2013, El_Karoui_2017}. The former contributes to bounding the eigenvalues of $X^TX/n$ from 0 and $\infty$ as in Lemma \ref{Lemma-1}, while the latter is employed here to demonstrate the consistency of terms like $\tr\{(X^TX)^{-1}\}$ as in Lemma \ref{Lemma-2}. Note that no sparsity assumption on $\Sigma$ is needed in our technical analysis. It is worth pointing out that the theoretical results above are adaptive to the low-dimensional regime, which makes our proposed method concretely applicable in practice.

In the end, we conduct a real data analysis on the relationship between gene expression and single nucleotide polymorphism ($\SNP$) with $n = 377$ and $p_n$ ranging from 33 to 87. Specifically, confidence intervals for the fraction of variance explained, i.e. $\rho_0$, are constructed using our proposed method, the conventional method and the one in \cite{Dicker_2014} based on the method of moment. We find that the conventional method may falsely discover non-zero $\rho_0$ for some genes due to the moderate dimension  and insufficient $\SNR$, and that our confidence interval is mostly narrower than that by \cite{Dicker_2014}.

{\bf Related Works}. Some earlier studies, e.g., \cite{Portnoy_1984, Portnoy_1985}, focused on the quadratic functional $(\hat\bbeta - \bbeta_0)^T X^TX (\hat\bbeta - \bbeta_0)$, under the low-dimensional regime, i.e., $\tau=0$. In the moderate-dimensional regime, \cite{El_Karoui_2013, El_Karoui_2017, El_Karouia_etal_2013, Donoho_Montanari_2016} studied the consistency of $\|\hat\bbeta - \bbeta_0\|_2$ for a general M-estimator $\hat\bbeta$. As far as we are aware, these techniques and results for consistency are not ready for deriving the asymptotic distributions of the quadratic functionals, which is the main contribution of our work. Another line of research is the element-wise inference \cite{Bai_etal_2013, Lei_etal_2018, Dobriban_Su_2018, Sur_etal_2017} whose strategies for analyzing single-element estimation error cannot be easily adapted for the analysis of aggregated estimation errors, e.g., quadratic functionals. To elucidate the difference between the two types of inferences, we plot $\sqrt n (\hat \beta_{j}^2 - \beta_{0,j}^2)$ versus $j$ and $\|\hat\bbeta\|_2^2 - \|\bbeta_0\|_2^2$ in Figure \ref{Figure_10}. In the high-dimensional regime, a more recent result is \cite{Cai_Guo_2017} who studied the point and interval estimations of $\|\hat\bbeta - \bbeta_0\|_q^2$ with $1 \leq q\leq 2$.

The rest of the paper is organized as follows. Section \ref{Section-2} develops the bias-and-variance-corrected inference for $\|\bbeta_0\|_2^2$ and demonstrates that the conventional procedure works if and only if the $\SNR$ is large. Section \ref{Section-3} consists of important applications of inferences for the quadratic functionals, including signal detection, global testing, inferences for the error variance and fraction of variance explained. Simulations are conducted in Section \ref{Section-4} and a real data analysis is performed in Section \ref{Section-5}. The proofs of some main theoretical results are included in the Appendix while the remaining proofs are relegated to the supplementary material.

{\bf Notation}. Let $\lfloor \cdot \rfloor$ be the floor function.
For any set $G$, denote by $\bar G$ the complement of $G$. Let $\ID(\cdot)$ be the indicator function.
Denote by $\bfI_m$ the $m \times m$ identity matrix and by $\bfe_{j,m}$ ($j=1,\ldots, m$) the $j$th column of $\bfI_{m}$. Let $\bfzero_m \in \mathbbR^m$ and $\bfone_m \in \mathbbR^m$ be the vectors of zeros and ones respectively.
For a vector $\boldsymbol{v}=(v_1,\ldots, v_m)^T$, the $L_1$, $L_2$ and $L_{\infty}$ norms are $\|\boldsymbol{v}\|_1=\sum_{i=1}^m|v_i|$, $\|\boldsymbol{v}\|_2 = (\sum_{i=1}^m v_i^2)^{1/2}$ and $\|\boldsymbol{v}\|_{\infty}=\max_{i \leq m}|v_i|$, respectively.
For an $m \times m$ matrix $A = \{a_{ij}\}_{1\leq i,j\leq m} $, denote by $\lambda_{\max}(A)$ and $\lambda_{\min}(A)$ the maximum and minimum eigenvalues of $A$, respectively. Let $|A|$ be the determinant of $A$.
The $L_1$, $L_2$ and $L_\infty$ norms of $A$ are defined as $\|A\|_1 = \max_{1\leq j \leq m} \sum_{i = 1}^{m} |a_{ij}|$, $\|A\|_2 = \{\lambda_{\max}(A' A)\}^{1/2}$ and $\|A\|_\infty = \max_{1\leq i \leq m} \sum_{j = 1}^{m} |a_{ij}|$, respectively.
For sequences $\{a_n\}_{n\geq 1}$ and $\{b_n \}_{n\geq 1}$, we write $a_n \lesssim b_n$ ($a_n \gtrsim b_n$) if there exists a constant $C>0$ independent with $n$ such that $|a_n| \leq C |b_n|$ ($|a_n| \geq C |b_n|$). Denote $a_n = \Omega(b_n)$ if $a_n = O(b_n)$ and $b_n = O(a_n)$. If $\{U_n\}_{n\geq 1}$ and $\{V_n\}_{n\geq 1}$ are random sequences, then  $U_n = \Omega_{\pr}(V_n)$ denotes that $U_n = O_{\pr}(V_n)$ and $V_n = O_{\pr}(U_n)$.
Notation ``$S_1  \Longleftrightarrow S_2$'' means that statements $S_1$ and $S_2$ are equivalent, while ``$S_1 \Longrightarrow S_2$'' denotes that $S_1$ implies $S_2$.
In the following, $C$ and $c$ are generic finite constants which may vary from place to place and do not depend on sample size $n$.

\section{Statistical Inference for Quadratic Functionals} \label{Section-2}

This section  establishes the dimension-adaptive inference  for $\|\bfbeta_0\|_2^2$, which is the main theoretical result of this paper. As a by-product, we discover that the classical (fixed-dimensional) statistical inference procedure continues to work in the low-dimensional regime if and only if the signal-to-noise ratio is large. As far as we are aware, this finding is new for quadratic functional $\|\bfbeta_0\|_2^2$.

We start with an examination of the plug-in estimator $\|\hat\bbeta\|_2^2$ for $\|\bfbeta_0\|_2^2$.
The estimation error $\|\hat\bbeta\|_2^2 - \|\bfbeta_0\|_2^2$ can be expressed as a linear-quadratic form, i.e.,  sum of a quadratic term and  a linear term, as follows
\begin{equation}\label{2.1}
\|\hat\bbeta\|_2^2 - \|\bfbeta_0\|_2^2 = \beps^T X(X^TX)^{-2} X^T \beps + 2 \bbeta_0^T (X^TX)^{-1} X^T \beps.
\end{equation}
The linear term has  zero mean and hence
 the bias of $\|\hat\bbeta\|_2^2$ is
\begin{eqnarray}\label{2.2}
 \E (\|\hat\bfbeta\|_2^2) - \|\bfbeta_0\|_2^2 = \E\{\beps^T X(X^TX)^{-2} X^T \beps\}  = \E\tr\{(X^TX)^{-1}\} \sigma_{\epsilon}^2 > 0.
\end{eqnarray}
For a special case that $ \bfX_i \stackrel{\iid}{\sim} N(\bfzero_{p_n}, \bfI_{p_n})$, $(X^TX)^{-1} $ follows the inverse Wishart distribution and hence
$$\E \tr\{(X^TX)^{-1}\} = p_n / (n-p_n-1)\to \tau/(1-\tau), \quad \text{as $n\to \infty$}.$$
For low dimension with $\tau = 0$, $\|\hat\bfbeta\|_2^2$ is asymptotically unbiased and its asymptotic distribution can be established based on the dominating linear term as in Theorem \ref{Theorem4} to be introduced later.
However, when $\tau>0$, the bias \eqref{2.2} is non-ignorable, leading to failure of the conventional low-dimensional results.

The  analysis above suggests a bias-corrected estimator for $\|\bfbeta_0\|_2^2$:
\begin{eqnarray} \label{2.3}
\hat{\|\bbeta\|_2^2} = \|\hat\bfbeta\|_2^2 - \tr\{(X^TX)^{-1}\} \hat\sigma_{\epsilon}^2,
\end{eqnarray}
where $\hat\sigma_{\epsilon}^2$ is defined in  \eqref{1.4}.
Since $X$ and  $\beps$ are independent,
$\hat{\|\bbeta\|_2^2}$ is unbiased  for $\|\bfbeta_0\|_2^2$.
Before presenting the asymptotic properties of $\hat{\|\bbeta\|_2^2} $, we first provide our assumptions below.

\ni{\bf Condition A.}
\begin{itemize}
\item[A1.] Assume $\{\bfX_i\}_{i=1}^n$ are $\iid$, $\bfX_i = \Sigma^{1/2} \bfZ_i$ where $\bfZ_i = (z_{i1}, \ldots, z_{ip_n})^T $, $\{z_{ij}\}_{j=1}^{p_n}$ are independent for each $i \leq n$, $\E(z_{ij}) = 0$, $\E(z_{ij}^2) = 1$ and there exists a constant $c^*>0$ such that for any $n\geq 1$, $i\leq n$, $j\leq p_n$ and $t>0$, $ \pr(|z_{ij}| \geq t) \leq 2 \exp( -c^*t^2)$.
\item[A2.] Suppose $\{\epsilon_i\}_{i=1}^n$ are $\iid$ and independent of $\{\bfX_i\}_{i=1}^n$, $\E(\epsilon_i)=0$, $\E(\epsilon_i^2)=\sigma_{\epsilon}^2\ge c>0$ and $\E(\epsilon_i^8)= O(\sigma_{\epsilon}^8)$.
\item[A3.] There exist constants $c$ and $C$, such that $0<c < \lambda_{\min}(\Sigma) \leq \lambda_{\max}(\Sigma) <C <\infty$.
\item[A4.] There exists a  constant $C$, such that $\|\bbeta_0\|_{\infty} \leq C < \infty$.
\end{itemize}
Conditions A1 and A2 only require sub-Gaussian tail for $\bfX_i$ and moment conditions on $\beps$, rather than impose any specific distributional restriction. The independence between  $\{\epsilon_i\}_{i=1}^n$ and $\{\bfX_i\}_{i=1}^n$ is crucial for applying the martingale difference CLT, and is a standard assumption for inference of the quadratic functionals in, e.g., \cite{Dicker_2014, Janson_etal_2017}.
The error variance $\sigma_{\epsilon}^2$ could either be bounded or diverging with $n$. Under Condition A4, $\|\bbeta_0\|_2 = O(\sqrt {p_n})$, and will reach $\Omega(\sqrt {p_n})$ when $\bbeta_0$ is not sparse. Throughout this paper, both $\bbeta_0$ and $\sigma_{\epsilon}^2$ are allowed to vary with $n$, except when $p_n$ is fixed.

Under Condition  $\mathrm{A}$,   $\hat\bbeta$ is well defined, i.e., the $p_n\times p_n$ matrix $(X^TX)^{-1}$ exists with probability tending to 1.
Lemma \ref{Lemma-1} below shows that the eigenvalues of $X^TX /n$ are bounded away from $0$ and $\infty$ with probability tending to 1 based on the random matrix theory in \cite{Bai_Silverstein_2010}. The proof is given in Section \ref{app-A3} of the Appendix.
\begin{lemma} \label{Lemma-1} If $ \tau \in [0, 1)$ and Conditions $\mathrm{A1}$ and $\mathrm{A3}$ hold, then for any $\ell \in \mathbb{N}$, we have $$\pr(\|X^T X / n\|_2 \geq x_1)= o(n^{ - \ell}), \quad \pr ( \|(X^TX / n)^{-1} \|_2 \geq  x_2^{-1}) = o(n^{ - \ell}),$$ where $x_1 = 4(1 + \sqrt{\tau})^2\|\Sigma\|_2$ and $x_2 = (1 - \sqrt{\tau})^2/(4\|\Sigma^{-1}\|_2 )$.
\end{lemma}
Define event
$K= H \cap J$, where $H$ and $J$ denote the events
$\|(X^TX /n )^{-1}\|_2 < x_2^{-1}$ and
$ \|X^TX /n \|_2 <x_1$, respectively.
Event $K$ is introduced to truncate the eigenvalues of $X^TX /n$.
Constants $x_1$ and $x_2^{-1}$ may not be the smallest for our analysis, and can be replaced by any constants larger than them.
From Lemma \ref{Lemma-1}, for any $\ell  \in \mathbb{N}$, we have $\pr(\bar K) = o(n^{-\ell})$.

We now present our main result: the asymptotic normality and ratio consistency of $\hat{\|\bbeta\|_2^2}$.

\begin{theorem} \label{Theorem1}
$\mathrm{(a)}$ 	Assume   $ \tau \in [0,1)$ and Condition $\mathrm{A}$ for \eqref{1.1}. If either of the following conditions hold: $\mathrm{(1)}$ $\lim_{n\to\infty}p_n = \infty$; $\mathrm{(2)}$ $p_n$ is fixed and $\bfbeta_0 \neq \bfzero_{p_n}$, then,
 $$     \zeta_n^{-1} (\hat{\|\bbeta\|_2^2} - \|\bfbeta_0\|_2^2 )  \conD N(0, 1),$$
 where
 \begin{eqnarray*}
  \zeta_n^2 &=&  4 \sigma_{\epsilon}^2  \bbeta_0^T \E\{(X^TX)^{-1} \ID(K)\} \bbeta_0 + 2 \sigma_{\epsilon}^4 \E \tr\{ (X^TX)^{-2}\ID(K)\} \cr
  &&+  2  \sigma_{\epsilon}^4  [\E\tr\{(X^TX)^{-1}\ID(K)\}]^2 /(n-p_n).
  \end{eqnarray*}
$\mathrm{(b)}$ Additionally, if $p_n^{1/2} / n = o(\SNR)$,
then
\begin{eqnarray} \label{2.4}
  \frac{\hat{\|\bbeta\|_2^2}}{\|\bfbeta_0\|_2^2} \conP 1.
  \end{eqnarray}
\end{theorem}

The proof of Theorem \ref{Theorem1} relies on the martingale difference CLT  \cite{Heyde_Brown_1970} and is provided in Section \ref{app-A3} of the Appendix.

A few remarks are in order: (i) After bias correction, $\hat{\|\bbeta\|_2^2}$ is asymptotically normal under fixed, low or moderate dimension. Hence, the proposed method is adaptive to dimension and generally applicable for $p_n < n$ in practice.
(ii) As $\|\bbeta_0\|_2$ may vary with $n$, the ratio consistency of $\hat{\|\bbeta\|_2^2}$ in \eqref{2.4}  is not automatically implied by the asymptotic normality but requires an additional assumption $p_n^{1/2} / n = o(\SNR)$.
(iii) Random design is assumed in Theorem \ref{Theorem1}, but the result also holds for fixed design, i.e., conditioning on $X$. (The $\SNR$ is not well defined for fixed design, and needs to be replaced by ${\|\bbeta_0\|_2^2}/{\sigma_{\epsilon}^2}$ in the condition of part (b).) As discussed in the end of the proof of Theorem \ref{Theorem1},  there exists a set $\mathcal{X}_n \subseteq \mathbb{R}^{n\times p_n}$ as in \eqref{A.11} satisfying $\pr(X \in \mathcal{X}_n ) \to 1$, such that for any $x \in \mathcal{X}_n$, $t\in \mathbb{R}$ and $\varepsilon>0$, $ \pr(\zeta_n^{-1} (\hat{\|\bbeta\|_2^2} - \|\bfbeta_0\|_2^2 )\leq t |X=x) \to \Phi(t)$ and $ \pr(|{\hat{\|\bbeta\|_2^2}} / {\|\bfbeta_0\|_2^2} -1 | \geq \varepsilon|X=x) \to 0$. In the following, all theoretical results (Theorems, Corollaries and Propositions) are applicable to fixed design unless otherwise specified.

\begin{remark} \label{Remark1}
In Theorem \ref{Theorem1}, we assume homoskedasticity for the error. Here, we consider heteroskedasticity, i.e., $\epsilon_i$ are independent with different variances $\sigma_i^2 = \E(\epsilon_i^2)$ for $i=1,\ldots, n$. We first introduce a general result. For $k \in \mathbb{N} $, denoting $D = \diag(\sigma_1^2, \ldots, \sigma_n^2)$, then
\begin{eqnarray*}
&&\E \{\beps^T X (X^TX)^{-k}X^T \beps\}   = \sum_{i=1}^n \E \{\bfX_i^T (X^TX)^{-k}\bfX_i\}\sigma_i^2 \cr
&=& \E\tr\{X (X^TX)^{-k}X^T/n\} \tr(D) =\E\tr\{  (X^TX)^{-k+1} \} \tr(D)/n,
\end{eqnarray*}
since $\E \{\bfX_i^T (X^TX)^{-k}\bfX_i\}$ are identical and equal $\E\tr\{X (X^TX)^{-k}X^T/n\}$ for $i=1,\ldots, n$.
From \eqref{2.1}, the bias of $\|\hat\bfbeta\|_2^2$ becomes $\E \{\beps^T X (X^TX)^{-2}X^T \beps\}=  \E\tr\{  (X^TX)^{-1} \} \tr(D)/n$. In Lemma \ref{Lemma-2} below, $\tr\{  (X^TX)^{-1} \}$ is ratio consistent for $\E\tr\{  (X^TX)^{-1} \}$ given event $K$, while $\hat \sigma_{\epsilon}^2$ is unbiased for $\tr(D)/n$ because $\E(\hat \sigma_{\epsilon}^2) = \E [ \beps^T \{\bfI_{n} - X(X^TX)^{-1} X^T\} \beps / (n-p_n) ] = \{\tr(D) - p_n/n \tr(D)\} / (n-p_n) = \tr(D)/n $. Hence, $\hat{\|\bbeta\|_2^2}$  is still a bias-corrected estimator for $\|\bfbeta_0\|_2^2$ under heteroskedasticity.
It will be an interesting future work to derive the asymptotic distribution of $\hat{\|\bbeta\|_2^2}$ under heteroskedasticity, and derive consistent estimators for the parameters in the limiting distribution.
\end{remark}
\begin{remark} \label{Remark2}
For ease of presenting the proofs, we assume $\E(\bfX_i) = \bfzero_{p_n}$ in Condition $\mathrm{A1}$ and hence $\E(Y_i) = 0$.
For the general form of the linear model
$$Y_i = \alpha_0+ \bfX_i^T \bfbeta_0 + \epsilon_i, \quad \text{for $i = 1, \ldots, n$,}$$
where  $\alpha_0$ is the intercept and $\E(\bfX_i) = \bmu $, our method is still applicable to the centralized data  $\{Y_i - \bar Y, \bfX_i - \bar \bfX\}_{i=1}^n$ with $\bar Y = n^{-1} \sum_{i=1}^n Y_i$ and $\bar \bfX = n^{-1} \sum_{i=1}^n \bfX_i$.
Specifically, Theorem \ref{Theorem1} together with the Theorems, Corollaries and Propositions below still holds after data centralization. Please see Section S.2
in the supplementary material for a brief explanation.
\end{remark}
\begin{remark} \label{Remark3}
From \eqref{1.4}, the variance term of the conventional inference procedure
$\hat \zeta_0^2=  4 \hat \sigma_{\epsilon}^2  \hat\bbeta^T (X^TX)^{-1} \hat\bbeta$ is ratio consistent for
$\zeta_0^2= 4 \sigma_{\epsilon}^2  \bbeta_0^T \E\{(X^TX)^{-1} \ID(K)\} \bbeta_0$
(see Theorem \ref{Theorem4}   to be introduced later).
Hence, the removal of bias in $\|\hat\bbeta\|_2^2$ leads to a larger variance $\zeta_n^2$ than $\zeta_0^2$. To see that more clearly, we consider a special case that $\bfX_i \stackrel{\iid}{\sim} N(\bfzero_{p_n}, \bfI_{p_n})$. In this case, $\E\{(X^TX)^{-1}\} = \bfI_{p_n}/(n-p_n-1)$ and $n \E \tr\{ (X^TX)^{-2}\} \to  \tau/(1 - \tau)^3  $ based on \cite{Letac_Massam_2004}. From \eqref{2.2}, we know that the amount of theoretical correction of bias for $\hat{\|\bbeta\|_2^2}$ compared with $ \|\hat\bbeta\|_2^2$ is $-\tau\sigma_{\epsilon}^2/(1-\tau)$. Also,
\begin{eqnarray}
\zeta_n^2
= \zeta_0^2 + \frac{2\sigma_{\epsilon}^4} n \frac{\tau(1+\tau)}{(1-\tau)^3} \{1 + o(1)\}\nonumber
\end{eqnarray}
for $\tau \in (0,1)$.
Both bias and variance corrections deviate from zero significantly as $\tau\to 1$; see Figure~\ref{Figure_3} for $n=100$ and $\sigma_\epsilon^2=1$. The  patterns in Figure~\ref{Figure_3} are consistent with the empirical ones observed in Figure \ref{Figure_2}.
\begin{figure}[!ht]
	\centering
	\includegraphics[height = 2in, width = 5in]{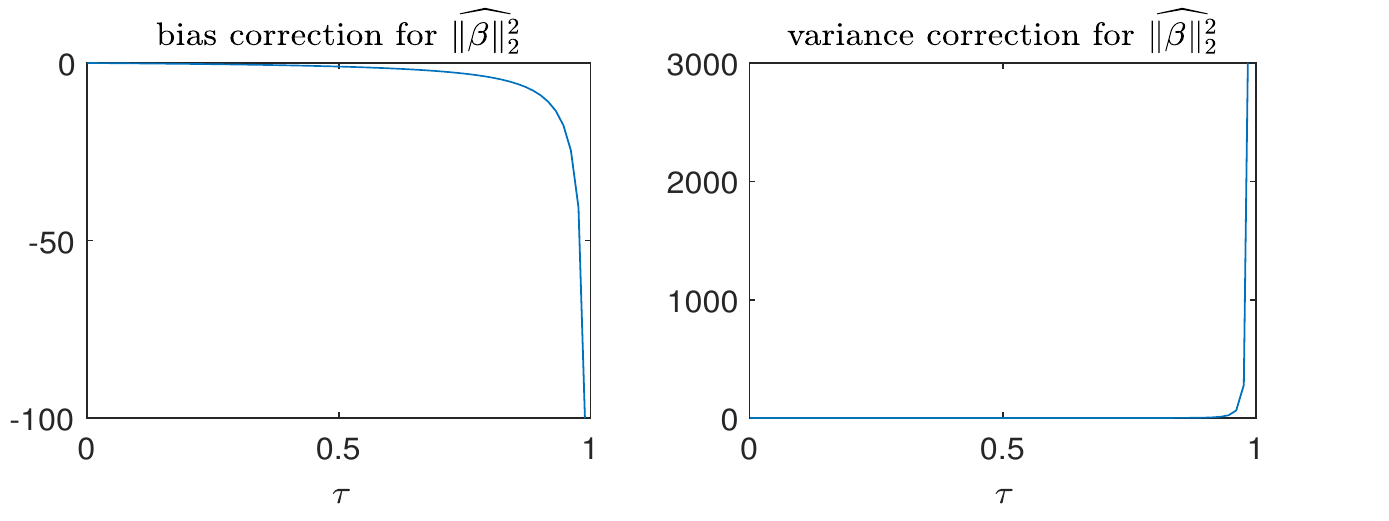}
	\caption{\textsl{Amount of theoretical corrections of bias (left penal) and variance (right penal) versus $\tau$ for $\hat{\|\bbeta\|_2^2}$ compared with $\|\hat\bbeta\|_2^2$. }}	
\label{Figure_3}
\end{figure}
\end{remark}

\begin{remark} \label{Remark4}
We now discuss that  $\hat{\|\bbeta\|_2^2}$ is not a uniformly minimum variance unbiased estimator (UMVUE). Denote by $T(\bfY, X)$ a generic unbiased estimator of $\|\bbeta_0\|_2^2$. If we assume $\bfX_i \stackrel{\iid}{\sim} N(\bfzero_{p_n}, \Sigma)$ and $\beps \sim N(\bfzero_{n}, \sigma_{\epsilon}^2 \bfI_{n})$, then the joint probability density function of $(Y_i, \bfX_i)$ is $f(y_i, \bfx_i) =   \{({2\pi})^{p_n+1} \sigma_{\epsilon}^2 |\Sigma|\}^{-1/2}
\exp \{-(y_i - \bfx_i^T \bbeta_0)^2 / (2\sigma_{\epsilon}^2)- \bfx_i^T \Sigma^{-1} \bfx_i/2\}  $ implying that the Fisher information matrix with the full data $(\bfY,X)$ for $\bbeta_0$ is $I(\bbeta_0) = n \Sigma / \sigma_{\epsilon}^2$.
Using the Cram\'{e}r-Rao lower bound (see, e.g., Theorem 3.3 in \cite{Shao_2003}), we have $\var\{T(\bfY, X)\} \geq 4 \sigma_{\epsilon}^2 \bbeta_0^T \Sigma^{-1} \bbeta_0 /n$. From the fact that $\E\{(X^TX)^{-1}\} = \Sigma^{-1}/(n-p_n-1)$, we know $\zeta_n^2 > 4 \sigma_{\epsilon}^2 \bbeta_0^T \Sigma^{-1} \bbeta_0 /n$ and hence $\hat{\|\bbeta\|_2^2}$ may not be a UMVUE of $\|\bbeta_0\|_2^2$.
\end{remark}

As discussed in Remarks \ref{Remark3} and \ref{Remark4},
 the bias correction for $\hat{\|\bbeta\|_2^2}$ leads to larger asymptotic variance and $\hat{\|\bbeta\|_2^2}$ is not a UMVUE for $\|\bbeta_0\|_2^2$. However, we can show that $\hat{\|\bbeta\|_2^2}$ achieves the optimal rate of convergence in terms of the quadratic loss.
\begin{theorem} \label{Theorem2}
Assume model \eqref{1.1} with $\bfX_i \stackrel{\iid}{\sim} N(\bfzero_{p_n}, \Sigma)$, $\beps \sim N(\bfzero_{n}, \sigma_{\epsilon}^2 \bfI_{n})$, $\sigma_{\epsilon}^2 = O(n)$, that $\{\epsilon_i\}_{i=1}^n$ are independent of $\{\bfX_i\}_{i=1}^n$, and Conditions $\mathrm{A3}$ and $\mathrm{A4}$ hold.
For any estimator $T$ of $\|\bbeta_0\|_2^2$,
we have
$$\inf_{T} \sup_{\bbeta \in \mathcal{G}_{\bbeta_0}(c) } \E_{(\bfY,X)|(\bbeta, \sigma_{\epsilon}^2, \Sigma)} (T - \|\bbeta\|_2^2)^2 = \Omega(\zeta_n^2),$$
where $\mathcal{G}_{\bbeta_0}(c) = \{\bbeta \in \mathbb{R}^{p_n}: \|\bbeta\|_{\infty} \leq C < \infty, \|\bbeta\|_2 \leq c(\|\bbeta_0\|_2 + \sigma_{\epsilon}\sqrt{p_n} /\sqrt n)\}$, $c>1$ is a generic constant and $\E_{(\bfY,X)|(\bbeta, \sigma_{\epsilon}^2, \Sigma)} (\cdot)$ denotes taking expectation with respect to $(\bfY,X)$ given parameters $(\bbeta, \sigma_{\epsilon}^2, \Sigma)$.
\end{theorem}
Theorem \ref{Theorem2} implies that,
$\hat{\|\bbeta\|_2^2}$ achieves the optimal convergence rate over all $\bbeta \in \mathcal{G}_{\bbeta_0}(c)$ under the quadratic loss.
Since $\zeta_n^2$ involves the true parameter $\bbeta_0$, the set of parameters $\mathcal{G}_{\bbeta_0}(c)$ also depends on $\bbeta_0$, which covers a wide range of $p_n$-dimensional vectors including $\bbeta_0$.

To estimate the variance term $\zeta_n^2$, we need to estimate the following four terms: $\sigma_{\epsilon}^2 $, $\II := \bbeta_0^T \E\{(X^TX)^{-1} \ID(K)\} \bbeta_0$, $\III := \E \tr\{ (X^TX)^{-2}\ID(K)\}$ and $\IV := \E\tr\{(X^TX)^{-1}\ID(K)\}$. The error variance $\sigma_{\epsilon}^2 $ can be consistently estimated by $\hat\sigma_{\epsilon}^2 $ as in Proposition \ref{Proposition-1} to be introduced in Section \ref{section-3.3}. For the other three terms, we need to utilize the following general result.
\begin{lemma} \label{Lemma-2} Assume  $ \tau \in [0,1)$ and Conditions $\mathrm{A1}$ and $\mathrm{A3}$ for \eqref{1.1}. For any $k \in \mathbb{N}$,
\begin{eqnarray*}
\var[   \tr\{(X^T X/n)^{-k}\} \ID(K)] = o(p_n^2), \quad
\var\{ \bbeta_0^T (X^TX/n)^{-k} \bbeta_0  \ID(K)\} = o(\|\bbeta_0\|_2^4).
\end{eqnarray*}
\end{lemma}
The key strategy to prove Lemma \ref{Lemma-2} is the leave-one-observation-out method. See Section \ref{app-A3} of the Appendix for the detailed proof. According to Lemma \ref{Lemma-2},
 $\tr\{(X^TX)^{-2}\}$ and $ \tr\{(X^TX)^{-1}\}$ are ratio consistent for terms $ \III$ and $\IV$, respectively. It's not necessary to include $\ID(K)$ in the estimators, since $\ID(K)\conP 1$ due to $\pr(K) \to 1$.
Lemma \ref{Lemma-2} further induces  Lemma S.12
in the supplementary material that,
$$   \hat\bbeta^T (X^TX)^{-1} \hat\bbeta -  \hat\sigma_{\epsilon}^2   \tr\{ (X^TX)^{-2}  \} -   \II = o_{\pr}(\zeta_n^2 / \sigma_{\epsilon}^2).$$
Subsequently, the plug-in estimator of $\zeta_n^2$ is
$$ \hat \zeta_n^2 =  4 \hat \sigma_{\epsilon}^2  \hat\bbeta^T (X^TX)^{-1} \hat\bbeta - 2 \hat \sigma_{\epsilon}^4    \tr\{ (X^TX)^{-2}\} +  2 \hat \sigma_{\epsilon}^4   [\tr\{(X^TX)^{-1}\}]^2/(n-p_n).$$
 We summarize the above discussion into  Theorem \ref{Theorem3} below.

\begin{theorem} \label{Theorem3}
Under the conditions in part $\mathrm{(a)}$ of Theorem \ref{Theorem1}, we have
 $$  {\hat \zeta_n^2 }  / { \zeta_n^2 } \conP 1.$$
\end{theorem}
The proof of Theorem \ref{Theorem3} is provided in Section \ref{app-A3} of the Appendix.

We are now ready to test the hypothesis in \eqref{1.2} by proposing the following test statistic
\begin{eqnarray} \label{2.5}
\mathbb{Z}_n = \frac{\hat{\|\bbeta\|_2^2} - c_0^2}{\hat \zeta_n }.
\end{eqnarray}
Theorems \ref{Theorem1} and \ref{Theorem3} directly imply that
  the null limiting distribution of $\mathbb{Z}_n$ is standard normal, the P-value for testing \eqref{1.2} is
$2\Phi(-|\mathbbm{z}_n|)$, where $\mathbbm{z}_n$ is a realization of $\mathbb{Z}_n$,
   and the asymptotic  power function is $1 - \Phi\{ -(c_1^2 - c_0^2) / {\hat \zeta_n } + \Phi^{-1}(1-\alpha / 2)\} + \Phi\{ -(c_1^2 - c_0^2) / {\hat \zeta_n } - \Phi^{-1}(1-\alpha / 2)\} $ under the fixed alternative $ H_1: \|\bbeta_0\|_2 = c_1\neq c_0$.

In the end, we point out that as long as the $\SNR$ is large enough, the conventional estimator $\|\hat\bbeta\|_2^2$ is still ratio consistent and asymptotically normal. However, the strength of $\SNR$ is usually unknown in practice. Hence, this result highlights the importance of our proposed {\em adaptive} method that works for both moderate and low dimensions, regardless of weak or strong signals.
\begin{theorem} \label{Theorem4}
Assume $  \tau \in [0,1)$ and Condition $\mathrm{A}$ for \eqref{1.1}. Then,
 \begin{eqnarray}
\frac{\|\hat\bbeta\|_2^2}{\|\bfbeta_0\|_2^2} \conP 1  &\Longleftrightarrow& p_n/n = o(\SNR )\Longleftrightarrow
\frac{\hat \zeta_0^2 }{ \zeta_0^2 } \conP 1, \nonumber \\
\frac{\|\hat\bbeta\|_2^2 - \|\bfbeta_0\|_2^2}{\zeta_0}\conD N(0, 1)  &\Longleftrightarrow& p_n^2/n = o(\SNR) \Longrightarrow \tau = 0, \nonumber
\end{eqnarray}
where $\zeta_0^2=  4 \sigma_{\epsilon}^2   \bbeta_0^T \E\{(X^TX)^{-1} \ID(K)\} \bbeta_0$  and $\hat\zeta_0^2$ is defined in \eqref{1.4}.
\end{theorem}
Theorem \ref{Theorem4} verifies the asymptotic normality of $\mathbb{Z}_0$ in \eqref{1.3} under  valid null hypothesis. The proof
is provided in the supplementary material.

\section{Applications} \label{Section-3}
This section consists of four important applications of our general theory: signal detection, global testing, inferences for the error variance and fraction of variance explained. The results on two-sample inference are postponed to supplementary.

\subsection{Detection boundary for  $\|\bbeta_0\|_2^2$}

Hypothesis \eqref{1.2} can be used to perform signal detection by setting $c_0= 0$. In this problem, the detection boundary is often of interest, which is the smallest separation rate between the null and a sequence of contiguous alternatives $H_{1n}$ indexed by $\delta_n\to 0$, i.e.,
\begin{eqnarray*}
H_{1n}: \|\bbeta_0\|_2^2  =  \delta_n,
\end{eqnarray*}
such that successful detection is still possible. From Theorem \ref{Theorem1}, we propose the following test statistic for hypothesis \eqref{1.2} with $c_0= 0$
\begin{eqnarray*}
\mathbb{Z}_n ^* = \frac{\hat{\|\bbeta\|_2^2} }{\hat \zeta_*},
\end{eqnarray*}
where $ \hat \zeta_*^2 =   2 \hat \sigma_{\epsilon}^4    \tr\{ (X^TX)^{-2}\}   +  2 \hat \sigma_{\epsilon}^4   [\tr\{(X^TX)^{-1}\}]^2/(n-p_n)$. The difference between $\mathbb{Z}_n ^*$ and $\mathbb{Z}_n$ lies in the variance term $\hat \zeta_*^2 $.
Using Lemma \ref{Lemma-2},  $ \hat \zeta_*^2$ is ratio consistent for $  \zeta_*^2  =2 \sigma_{\epsilon}^4 \E \tr\{ (X^TX)^{-2}\ID(K)\}  +  2  \sigma_{\epsilon}^4  [\E\tr\{(X^TX)^{-1}\ID(K)\}]^2 /(n-p_n)$ which equals
$\zeta_n^2$ when $\bbeta_0 = \bfzero_{p_n}$. In other words, $\hat \zeta_*^2 $ is a refined estimator of $\zeta_n^2$ under the null hypothesis \eqref{1.2} with $c_0=0$.
Then, the asymptotic standard normality of $\mathbb{Z}_n ^*$ under the null  follows directly from Theorem \ref{Theorem1} for diverging $p_n$. Corollary~\ref{Corollary-1} below presents the detection boundary using $\mathbb{Z}_n ^*$.

\begin{corollary} \label{Corollary-1}
Assume  that $ \tau \in [0,1)$, $\lim_{n\to\infty} p_n = \infty$, Condition $\mathrm{A}$ holds for \eqref{1.1}.
If $\delta_n = \Omega( \sigma_{\epsilon}^2 p_n^{1/2} /n )$, then
$$\mathbb{Z}_n^* -  \hat\zeta_*^{-1} \delta_n \conD N(0, 1)$$ where $\hat\zeta_*^{-1} \delta_n = \Omega_{\pr}(1)$. If $\delta_n = o( \sigma_{\epsilon}^2 p_n^{1/2} /n )$, then $$\mathbb{Z}_n^*  \conD N(0, 1).$$
\end{corollary}
Therefore, the detection boundary is $\sigma_{\epsilon}^2 p_n^{1/2} /n$, which matches with the minimax detection rate in \cite{Ingster_etal_2010} (see (1.2) therein). It is worth mentioning that Corollary \ref{Corollary-1} requires diverging $p_n$.

\subsection{Global inference for $\bbeta_0$}
This section is concerned with the global hypothesis \eqref{1.5}
\begin{eqnarray*}
H_0: \bfbeta_0 =  \bbeta_0^{\nullH0}\quad \mbox{versus}\quad H_1: \bfbeta_0 \neq  \bbeta_0^{\nullH0},
\end{eqnarray*}
by proposing a bias-and-variance-corrected test statistic  based on $\|\hat \bbeta - \bbeta_0^{\nullH0} \|_2^2$ as follows
\begin{eqnarray} \label{3.1}
\mathbb{G}_n =   \frac{\|\hat \bbeta - \bbeta_0^{\nullH0} \|_2^2 - \tr\{(X^TX)^{-1}\} \hat\sigma_{\epsilon}^2} {\hat \zeta_*}.
\end{eqnarray}
The  construction of $\mathbb{G}_n$ is based on the fact that the distribution of $ \hat \bbeta - \bbeta_0^{\nullH0} $ under $H_0$ in \eqref{1.5} is the same as that of $ \hat \bbeta  $ with $\bbeta_0 = \bfzero_{p_n}$. Thus, the amount of bias correction  $- \tr\{(X^TX)^{-1}\} \hat\sigma_{\epsilon}^2$ and the variance term  $ \hat \zeta_*^2 $ for $\mathbb{G}_n$ are the same as those for $\mathbb{Z}_n ^*$.
From the asymptotic results of  $\mathbb{Z}_n ^*$,  $\mathbb{G}_n$ is also asymptotically standard normal under the null for diverging $p_n$, and
the smallest separation rate for $H_{1n}:  \|\bfbeta_0- \bbeta_0^{\nullH0}\|_2 ^2 = \delta_n$ is  $\delta_n^* = \sigma_{\epsilon}^2 p_n^{1/2} /n$,  the same as that identified by Corollary \ref{Corollary-1}.

From \eqref{3.1}, we can construct  $1-\alpha$ confidence regions for $\bbeta_0$ using one-sided and two-sided strategies as
\begin{equation}\label{3.2}
\begin{aligned}
\mathrm{CR_1}&=	\{ \bbeta:   \|\hat\bbeta - \bbeta \|_2^2 -\tr\{(X^TX)^{-1}\} \hat\sigma_{\epsilon}^2 \leq   \Phi^{-1}(1 - \alpha ) \hat \zeta_* \};\\
\mathrm{CR_2}&=	\{ \bbeta:  | \|\hat\bbeta - \bbeta \|_2^2 -\tr\{(X^TX)^{-1}\} \hat\sigma_{\epsilon}^2| \leq   \Phi^{-1}(1 - \alpha /2) \hat \zeta_* \}.
\end{aligned}
\end{equation}
Confidence region for high-dimensional sparse $\bbeta_0$ was studied in \cite{Cai_Guo_2018, Nickl_van_de_Geer_2013}.

\subsection{Inference for $\sigma_{\epsilon}^2$} \label{section-3.3}

This section is concerned with moderate-dimensional inference for the error variance $\sigma_{\epsilon}^2$.

\begin{proposition} \label{Proposition-1}
Assume   $\tau \in [0,1)$ and Conditions $\mathrm{A1}$--$\mathrm{A3}$ for \eqref{1.1}.
Then $$ \frac{\hat \sigma_{\epsilon}^2 }{ \sigma_{\epsilon}^2}\conP 1 \quad \text{and}\quad \frac{\hat \sigma_{\epsilon}^2  - \sigma_{\epsilon}^2 }{ \zeta_{\epsilon}  } \conD N(0, 1),$$
 where
 $  \zeta_{\epsilon}^2 = n^{-1} \{\nu_4  + \sigma_{\epsilon}^4  ( 3\tau - 1)/(1-\tau) \}  $ and $\nu_4 = \E (\epsilon_i^4)$.
\end{proposition}
Our result is adaptive to data dimension  by incorporating $\tau$ in
the variance term $ \zeta_{\epsilon}^2 $ for both fixed and diverging $p_n$. Specifically, $ \zeta_{\epsilon}^2 $ increases with $\tau$. For a special case that $\epsilon_i \sim N(0, \sigma_{\epsilon}^2)$, $ \zeta_{\epsilon}^2 = 2 \sigma_{\epsilon}^4 / \{n(1-\tau)\}$.

To estimate $ \zeta_{\epsilon}^2 $, it suffices to provide a ratio consistent estimator for $\nu_4$. We first examine a straightforward estimator $ 1/n \sum_{i = 1}^n  \hat\epsilon_i^4$  where $(\hat\epsilon_1, \ldots, \hat\epsilon_{n})^T = \bfY - X\hat\bbeta$. However, as in the proof of Lemma S.15
in the supplementary material,
 $ 1/n \E(\sum_{i = 1}^n  \hat\epsilon_i^4) =(1  - \tau )^4\nu_4 +3\sigma_{\epsilon} ^4 \tau (1-\tau)^2 (2-\tau) +o(\sigma_{\epsilon}^4)$.
Although the naive estimator is biased, it induces an estimator for  $\nu_4$ after centering and re-scaling
$$\hat\nu_4 = (1-p_n/n)^{-4} \Big\{1/n \sum_{i = 1}^n  \hat\epsilon_i^4-3 \hat \sigma_{\epsilon}^4 (p_n/n) (1-p_n/n)^2(2-p_n/n) \Big\}.$$
Lemma S.15
demonstrates that $\hat\nu_4$ is ratio consistent  for $\nu_4$.
 Hence, the
plug-in estimator $\hat\zeta_{\epsilon}^2 = n^{-1} \{\hat\nu_4  + \hat\sigma_{\epsilon}^4  ( 3p_n/n - 1)/(1-p_n/n) \}$ is ratio consistent for $\zeta_{\epsilon}^2$. From   Proposition \ref{Proposition-1}, we have
$  ({\hat \sigma_{\epsilon}^2  - \sigma_{\epsilon}^2 })/{\hat \zeta_{\epsilon}  } \conD N(0, 1)$, which can be utilized to conduct inference for $\sigma_{\epsilon}^2$.

\subsection{Inference for $\rho_0$} \label{Section-3.4}
 Consider the hypotheses
\begin{eqnarray} \label{3.3}
H_0: \rho_0 \geq \rho_0^{\nullH0} \quad \text{versus} \quad H_1:  \rho_0 < \rho_0^{\nullH0},
\end{eqnarray}
where $0 < \rho_0^{\nullH0} < 1$ is a given constant. Recalling the definition of $\rho_0$ in \eqref{1.6}, its conventional plug-in estimator can be obtained by replacing  $\eta_0:=\bbeta_0^T \Sigma \bbeta_0$ and $ \sigma_{\epsilon}^2$ with $\hat\bbeta^T (X^TX/n)\hat \bbeta$  and $\hat\sigma_{\epsilon}^2$, respectively.
The asymptotic normality of this estimator
is studied under $\tau = 0$  in  Theorem S.1
in the supplementary material followed by the low-dimensional inference for $\rho_0$.

However, in the moderate-dimensional regime, the bias of $\hat\bbeta^T (X^TX/n)\hat \bbeta$ for $\eta_0$ is non-ignorable, i.e.,
\begin{eqnarray*}
&& \E\{ \hat\bbeta^T (X^TX/n)\hat \bbeta\} -\eta_0 \cr
  &=& \E\{\bbeta_0^T (X^TX/n) \bbeta_0\} -\eta_0 + 2n^{-1} \E(\bbeta_0^T  X^T\beps) + n^{-1}\E\{\beps^T X(X^TX)^{-1} X^T\beps\} \cr
  &=& n^{-1}\E\{\beps^T X(X^TX)^{-1} X^T\beps\} \cr
  &=&\sigma_{\epsilon}^2 p_n/n \to \sigma_{\epsilon}^2 \tau > 0.
\end{eqnarray*}
Consequently, we propose an unbiased estimator for $\eta_0$ as
$$\hat \eta = \hat\bbeta^T (X^TX/n)\hat \bbeta - \hat\sigma_{\epsilon}^2 p_n/n.$$ Hence, a new plug-in estimator for $\rho_0$ is
$$\hat \rho = \frac {\hat \eta}{\hat \eta + \hat \sigma_{\epsilon}^2} = \frac{\hat\bbeta^T (X^TX/n)\hat \bbeta - \hat\sigma_{\epsilon}^2 p_n/n} { \hat\bbeta^T  (X^TX/n)\hat \bbeta + \hat\sigma_{\epsilon}^2 (1-p_n/n )}$$ with the following asymptotic distribution.

\begin{theorem} \label{Theorem5}
Assume   $ \tau \in [0,1)$, $\rho_0 \in[C_1, C_2] $ for some constants $0<C_1\leq C_2<1$ and Condition $\mathrm{A}$ holds for \eqref{1.1}. Then, $\hat \rho-\rho_0 = o_{\pr}(1)$ and
$$  \frac{\hat \rho - \rho_0}{ \sigma_{\hat \rho} } \conD N(0, 1),$$
where $ \sigma^{2}_{\hat \rho} = n^{-1}(\eta_0 + \sigma_{\epsilon}^2)^{-4} [2 \sigma_{\epsilon}^8 \tau /(1-\tau) -\{2 + 4 \tau/(\tau-1)\} \sigma_{\epsilon}^6 \eta_0  + \sigma_{\epsilon}^4 \{\E( Y_1^4) - \nu_4   +\eta_0^2  ( 4\tau - 2)/(1-\tau)\}
 + \eta_0^2 \nu_4  ] $.
\end{theorem}
Simple calculation implies that $\hat\rho = 1- (1-p_n/n)^{-1}\|\bfY- X \hat\bbeta\|_2^2 / \|\bfY\|_2^2 = 1- (1-p_n/n)^{-1}(1-R^2)$, where $R^2$ is the coefficient of determination. Therefore, if $\tau = 0$, then $R^2$ is asymptotically unbiased for $\rho_0$, but when $\tau > 0$, a re-scaled $R^2$, i.e. $\hat\rho$, is required for the inference of $\rho_0$.
The definition of $\SNR$ is not applicable to fixed design, and hence the results of Theorem \ref{Theorem5} and Proposition \ref{Proposition2} below are not available for fixed design.

For the variance term $ \sigma^{2}_{\hat \rho}$,
the plug-in estimator $\hat\sigma^{2}_{\hat \rho}$ is obtained by replacing $\E( Y_1^4)$, $\eta_0$, $\sigma_\epsilon^2$, $\nu_4$ and $\tau$ in $\sigma^{2}_{\hat \rho}$ with $ n^{-1} \sum_{i=1}^n Y_i^4$, $\hat\eta$, $\hat{\sigma}_\epsilon^2$, $\hat\nu_4$ and $p_n/n$, respectively, and its consistency is demonstrated below.
\begin{proposition} \label{Proposition2}
Assume the conditions in Theorem \ref{Theorem5}. Then, $ \sigma^{2}_{\hat \rho}  = \Omega(1/n)$ and
$$ \hat \sigma^{2}_{\hat \rho} - \sigma^{2}_{\hat\rho} = o_{\pr}(1/n).$$
\end{proposition}
Hence, \eqref{3.3} can be tested by $\hat \sigma^{-1}_{\hat \rho}  (\hat \rho - \rho_0^{\nullH0}),$
 whose null limiting distribution is standard normal. Also the smallest separation rate for contiguous alternative is $n^{-1/2}$.

For the inference of $\eta_0$, the proof of Theorem \ref{Theorem5} immediately implies that
 \begin{eqnarray} \label{3.4}
 \sigma_{\hat\eta}^{-1}( \hat\eta-  \eta_0) \conD N(0, 1),
  \end{eqnarray}
  with $\sigma_{\hat\eta}^2=n^{-1}\{\E  ( Y_1^4 ) - \nu_4 - 2 \sigma_{\epsilon}^2 \eta_0 - \eta_0^2 + 2 \sigma_{\epsilon}^4 p_n /(n-p_n)\}$ which is consistently estimated by the plug-in estimator following the proof of Proposition \ref{Proposition2}.

 In the end, we comment on  related works concerned with signal strength, i.e., \cite{Dicker_Erdogdu_2016, Janson_etal_2017, Dicker_2014, Dicker_Erdogdu_2017, Verzelen_Gassiat_2018}.
The first three works, i.e., \cite{Dicker_Erdogdu_2016, Janson_etal_2017, Dicker_2014}, conduct statistical inference for moderate-dimensional fixed effect models. However, our OLS-based methods are essentially different from their methods in the following aspects:  parameters of interest and weak assumptions.

First, the parameter of interest in the aforementioned three works is $\bbeta_0^T \Sigma \bbeta_0$, and their procedures crucially rely on the fact that $Y_i \sim N(0, \bbeta_0^T \Sigma \bbeta_0 + \sigma_{\epsilon}^2)$ and that $\bbeta_0^T \Sigma \bbeta_0$ is a part of the variance term. Therefore, their results are not readily translated into inference for our parameter of interest, i.e., $\|\bbeta_0\|_2^2$, unless $\Sigma$ is identity. In contrast, our strategy depends on the OLS estimator. This flexibility also allows us to conduct inference for $\bbeta_0^T \Sigma \bbeta_0$ based on a bias-corrected version of $\hat\bbeta^T (X^TX/n) \hat\bbeta$ as in \eqref{3.4}, and even two-sample inferences, e.g., the co-heritability \eqref{1.7}, to which it is unclear how their methods can be applied.

Second, some assumptions of our paper are weaker due to the use of different technical tools. Specifically, the proofs as well as the development of the estimation and inference procedures in the aforementioned three works rely heavily on the Gaussian assumption of the design matrix $X$ and error $\beps$. For example, among other implications, the Gaussian design is important in deriving the invariant distribution of $X$ under orthogonal transformations in \cite{Dicker_Erdogdu_2016}, the Haar distribution of the right-singular vectors from the singular value decomposition of $X$ in \cite{Janson_etal_2017} and  the Wishart distribution of $X^T X$ in \cite{Dicker_2014}. However, our OLS-based result is derived using the martingale difference CLT without requiring any specific distributional assumption of $X$ or $\beps$.  Also, our results can be easily extended to fixed design. Besides, the three works above need $\Sigma = \bfI_{p_n}$ to conduct inference for $\|\bbeta_0\|_2^2$.
Even for the inference of $\bbeta_0^T \Sigma \bbeta_0$, they still require strong conditions on $\Sigma$, e.g., known or consistently estimable $\Sigma$.
And, some further sparsity assumptions need to be imposed if $\Sigma$ will be estimated.  Our inference methods for $\|\bbeta_0\|_2^2$ and $\bbeta_0^T \Sigma \bbeta_0$ neither need known or a consistent estimator of $\Sigma$ nor require any sparsity assumption on $\Sigma$. From the simulations in the end of Section \ref{Section-4}, our method performs better than or at least as well as those in \cite{Dicker_2014, Dicker_Erdogdu_2016}.

As for \cite{Dicker_Erdogdu_2017, Verzelen_Gassiat_2018}, the former conducted inference for the variance of the regression parameter by considering the ``random effect'' model conditioning on the design matrix, and hence is different from the setup of fixed effect model in our paper;
the latter derives the minimax estimators of $\rho_0$ under Gaussian design and error, but they didn't derive the asymptotic distribution of the estimators and hence their results can not be applied to the inference for $\rho_0$.

\section{Simulations} \label{Section-4}

Numerical studies are conducted to support the proposed statistical inference procedures. Set $n\in\{400,\, 800\}$ and
  $p_n =  4, \lfloor n/6 \rfloor, n/4, n/2.5$ corresponding to fixed dimension ($p_n=4$), low dimension ($p_n = \lfloor n/6 \rfloor$) and moderate dimension ($p_n=n/4, n/2.5$), unless otherwise specified. In the simulations, we consider a general form of the linear model $ Y_i =1 + \bfX_i^T \bfbeta_0 + \epsilon_i$ with $\E(\bfX_i)=\bmu = (\mu_1, \ldots, \mu_{p_n})^T$ generated by $\{\mu_i\}_{i=1}^{p_n}\stackrel{\iid}{\sim} \text{Unif}[1,2]$. Both $\bfY$ and $\{\bfX_i\}_{i=1}^n$ are centralized before conducting inference for the quadratic functionals.
  To generate data, we consider the following four cases representing various situations in reality:
  \begin{itemize}
\item [I.] Gaussian design with $\Sigma = \bfI_{p_n}$:
 $\bfX_i \stackrel{\iid}{\sim} N(\bmu, \bfI_{p_n})$ for $i=1,\ldots, n$, $\beps \sim N(\bfzero_n,  \bfI_{n})$ and $\bbeta_0 = \tilde\bbeta / \|\tilde\bbeta\|_2$ with $\{\tilde\beta_{j}\}_{j=1}^{p_n}\stackrel{\iid}{\sim} \text{Unif}[1, 2]$;
 \item[II.] Gaussian design with general $\Sigma$:
 $\bfX_i \stackrel{\iid}{\sim} N(\bmu, \Sigma)$ for $i=1,\ldots, n$ and $\beps \sim N(\bfzero_n,  \bfI_{n})$, where $\Sigma = \Sigma^{*T} \Sigma^* / \lambda_{\max} (\Sigma^{*T} \Sigma^*) + \diag(d_1, \ldots, d_{p_n})$,  $\Sigma_{ij}^* \stackrel{\iid}{\sim} \text{Unif}[-0.5, 0.5]$ for $1\leq i,j\leq p_n$ and $\{d_i\}_{i=1}^{p_n}\stackrel{\iid}{\sim}\text{Unif}[0.4, 1]$. Here $\bbeta_0 = c_{p_n} \bbeta^*$ where $\bbeta^*$ is the normalized eigenvector of $\Sigma$ corresponding to the smallest eigenvalue and $c_{p_n} = 1$ for $p_n=4$, $c_{p_n} = 2$ for $p_n=\lfloor n/6 \rfloor, n/4$ and $c_{p_n} = 5$ for $p_n=n/2.5$;
 \item [III.] t-distributed design:  $X_{ij} - \mu_{j} \stackrel{\iid}{\sim} t_5 / \sqrt {5/3}$ for $i=1,\ldots, n; j=1,\ldots, p_n$, $\epsilon_i \stackrel{\iid}{\sim} t_{16} / \sqrt{8/7}$ and $\bbeta_0 = (1,1,1,0,\ldots, 0)^T$;
 \item[IV.] fixed design: $X$ is identical  for all replications and generated by $\bfX_i \stackrel{\iid}{\sim} N(\bmu, \bfI_{p_n})$, while $\beps \sim N(\bfzero_n,  \bfI_{n})$ are independently generated for each replication, and $\bbeta_0$ is the normalized eigenvector of $\sum_{i=1}^n (\bfX_i - \bar\bfX)(\bfX_i - \bar\bfX)^T$ corresponding to the largest eigenvalue.
  \end{itemize}
  In case II, $\Sigma$ is not necessarily sparse and allows different diagonal elements. Case III is for non-Gaussian design and case IV corresponds to fixed design matrix.
In what follows, QQ plots (for the 1st to 99th percentiles) under the valid null hypotheses and confidence
intervals were obtained with 1000 replications, while the power function was computed using 500 replications for each setup.

\begin{figure}[!ht]
\centering
\includegraphics[height = \graphheight, width = \graphwidth]{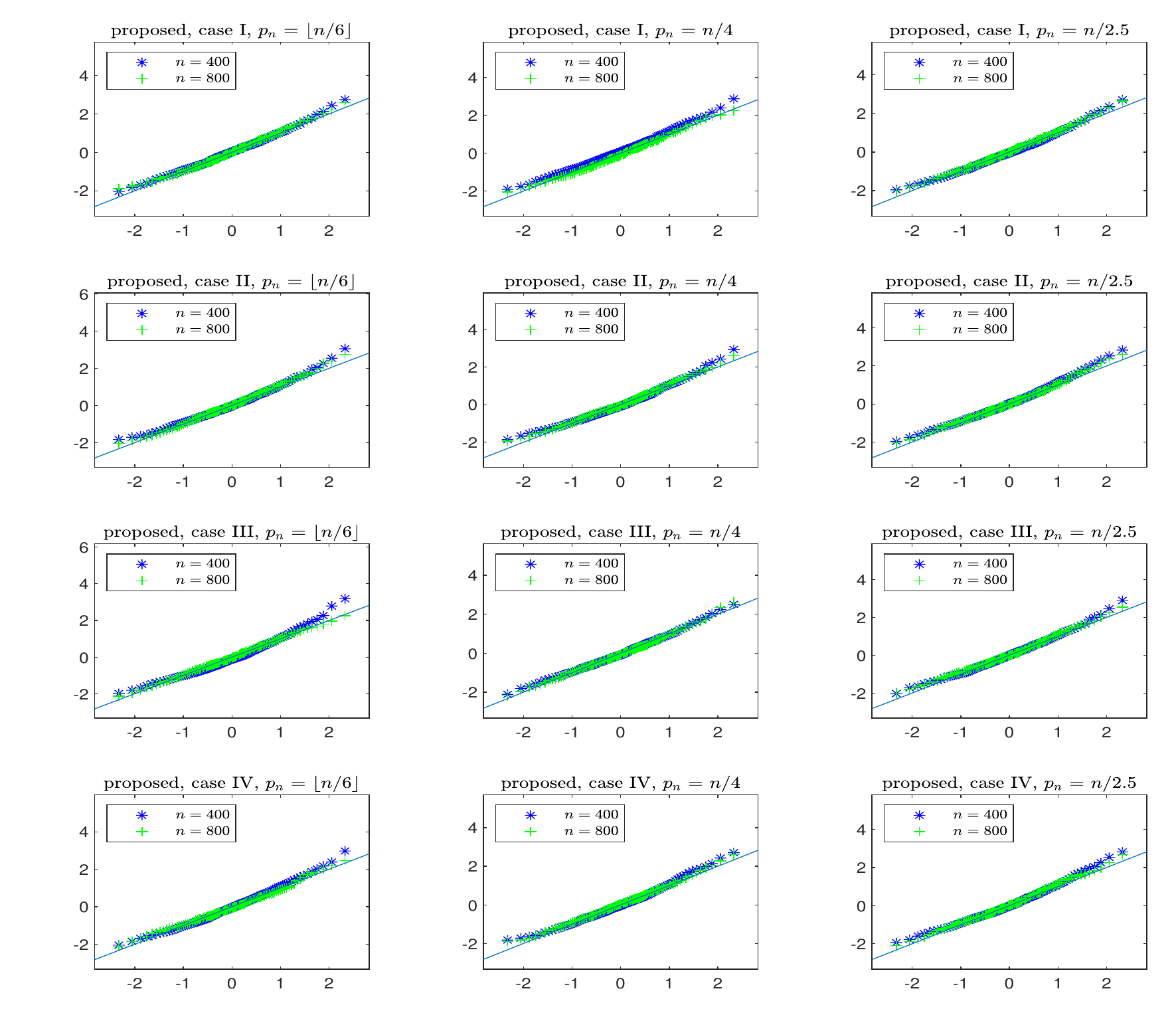}
\caption{\textsl{ QQ plots for testing $H_0$ in \eqref{1.2} with $c_0=0$ using $\mathbb{Z}_n^*$. Panels from top to bottom are for cases I--IV, respectively, while panels from left to right are for $p_n=\lfloor n/6\rfloor,n/4, n/2.5$, respectively.}}
\label{Figure_4}
\end{figure}

\begin{figure}[!ht]
\centering
\includegraphics[height = \graphheight, width = \graphwidth]{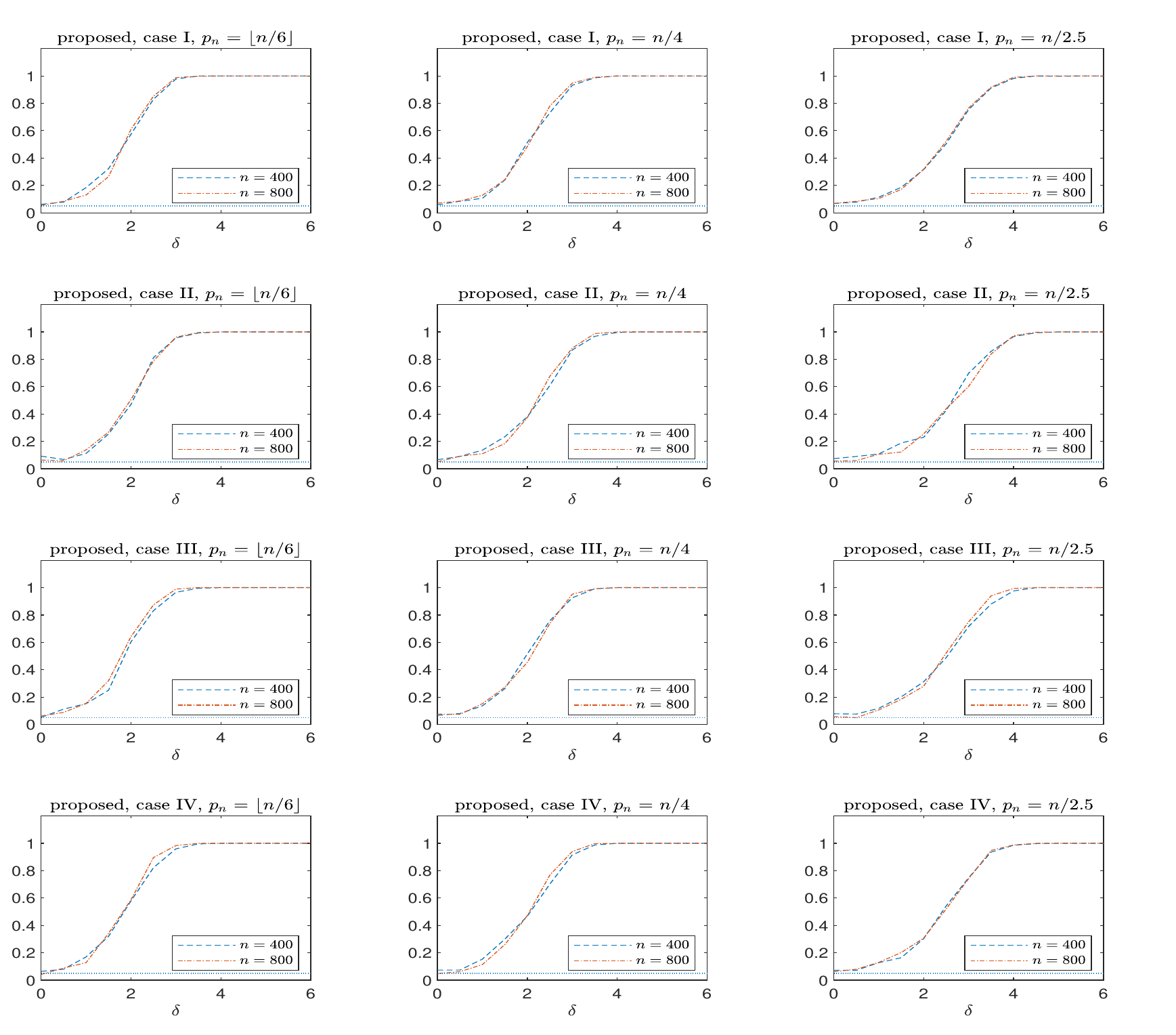}
\caption{\textsl{Empirical rejection rates versus $\delta$ for testing $H_1$ in \eqref{1.2} using $\mathbb{Z}_n^*$. Panels from top to bottom are for cases I--IV, respectively, while panels from left to right are for $p_n=\lfloor n/6\rfloor,n/4, n/2.5$, respectively. The dotted line indicates the true significance level $\alpha=0.05$.}}
\label{Figure_5}
\end{figure}

First, consider hypothesis \eqref{1.2} with $c_0 = 0$. Data are generated by cases I--IV with $\bbeta_0 = \bfzero_{p_n}$.
From Figure \ref{Figure_4}, under low and moderate dimensions, $\mathbb{Z}_n^*$ follows standard normal distribution under the null hypothesis.
The fixed-dimensional results are not reported as the signal detection is only conducted for diverging $p_n$. The empirical power of $\mathbb{Z}_n^*$ is given in Figure  \ref{Figure_5}  by varying $\bbeta_0 =  \bfone_{p_n} \delta \sigma_{\epsilon} / (n^{1/2} p_n^{1/4})$ with $\delta = 0, 0.5, 1, 1.5, \ldots, 6$. This choice of alternative values is supported by the derived detection boundary $\delta_n^* = \sigma_{\epsilon}^2 p_n^{1/2}/n$ for signal detection. From  Figure \ref{Figure_5}, we can tell that the empirical rejection rate grows from the nominal level to one as $\delta$ increases from zero.

We also check the coverage probability of the two-sided ($\mathrm{CR}_2$) and one-sided ($\mathrm{CR}_1$) confidence regions of $\bbeta_0$ based on \eqref{3.2}, with 1000 replications  at $\alpha = 0.05$.
  Table \ref{Table-1} reveals that both $\mathrm{CR}_1$ and $\mathrm{CR}_2$ are satisfactory while the latter  slightly outperforms the former.
 The coverage probabilities of $\mathrm{CR}_2$ are around 0.95 while those of  $\mathrm{CR}_1$ are generally below 0.95. Hence, we suggest to use $\mathrm{CR}_2$ in practice.
 Note that our proposed method particularly works for diverging $p_n$, but when $p_n$ is fixed, the finite-sample performance is still satisfactory.

\begin{table}[!ht]
\caption{\textsl{{Coverage probability of $95\%$ confidence regions for $\bbeta_0$ }}}
\label{Table-1}
\begin{center}
\begin{tabular}{ l  cccc}
\hline
      $ (n, p_n)$    &  {case I}  &     {case II}  &  {case III}  &     {case IV}   \\ \hline
      \multicolumn{5}{c}{ $\mathrm{CR}_2$}\\
   $ (400, \ \ \ 4)$&$0.943$&$0.946$ &0.950 & 0.932  \\
   $  (400, \ \, 66)$&$0.945$&$0.946$ &0.941& 0.945\\
   $ (400, 100)$ &$0.949$ &$0.948$&0.951&  0.953\\
   $ (400, 160)$&$0.950$& $0.940$ &0.946&0.942 \\
   $ (800,\ \ \ 4)$ &$0.944$&$0.960$&0.964&0.948  \\
   $  (800, 133)$&$0.958$&$0.945$&0.959&0.951  \\
   $ (800, 200)$&$0.944$&$0.961$&0.947&0.957  \\
   $ (800, 320)$&$0.952$&$0.943$&0.957& 0.954 \\
\hline
\multicolumn{5}{c}{ $\mathrm{CR}_1$} \\
   $ (400, \ \ \ 4)$&$0.913$&$0.921$ &0.936 & 0.908  \\
   $  (400, \ \, 66)$&$0.925$&$0.931$ &0.920& 0.928\\
   $ (400, 100)$ &$0.933$ &$0.933$&0.944&  0.916\\
   $ (400, 160)$&$0.941$& $0.923$ &0.936&0.911 \\
   $ (800,\ \ \ 4)$ &$0.929$&$0.938$&0.942&0.923  \\
   $  (800, 133)$&$0.951$&$0.935$&0.953&0.937  \\
   $ (800, 200)$&$0.936$&$0.940$&0.953&0.941  \\
   $ (800, 320)$&$0.922$&$0.928$&0.942& 0.934 \\
\hline
\end{tabular}
\end{center}
\end{table}

\begin{figure}[!ht]
\centering
\includegraphics[height = 6 in, width = 6 in]{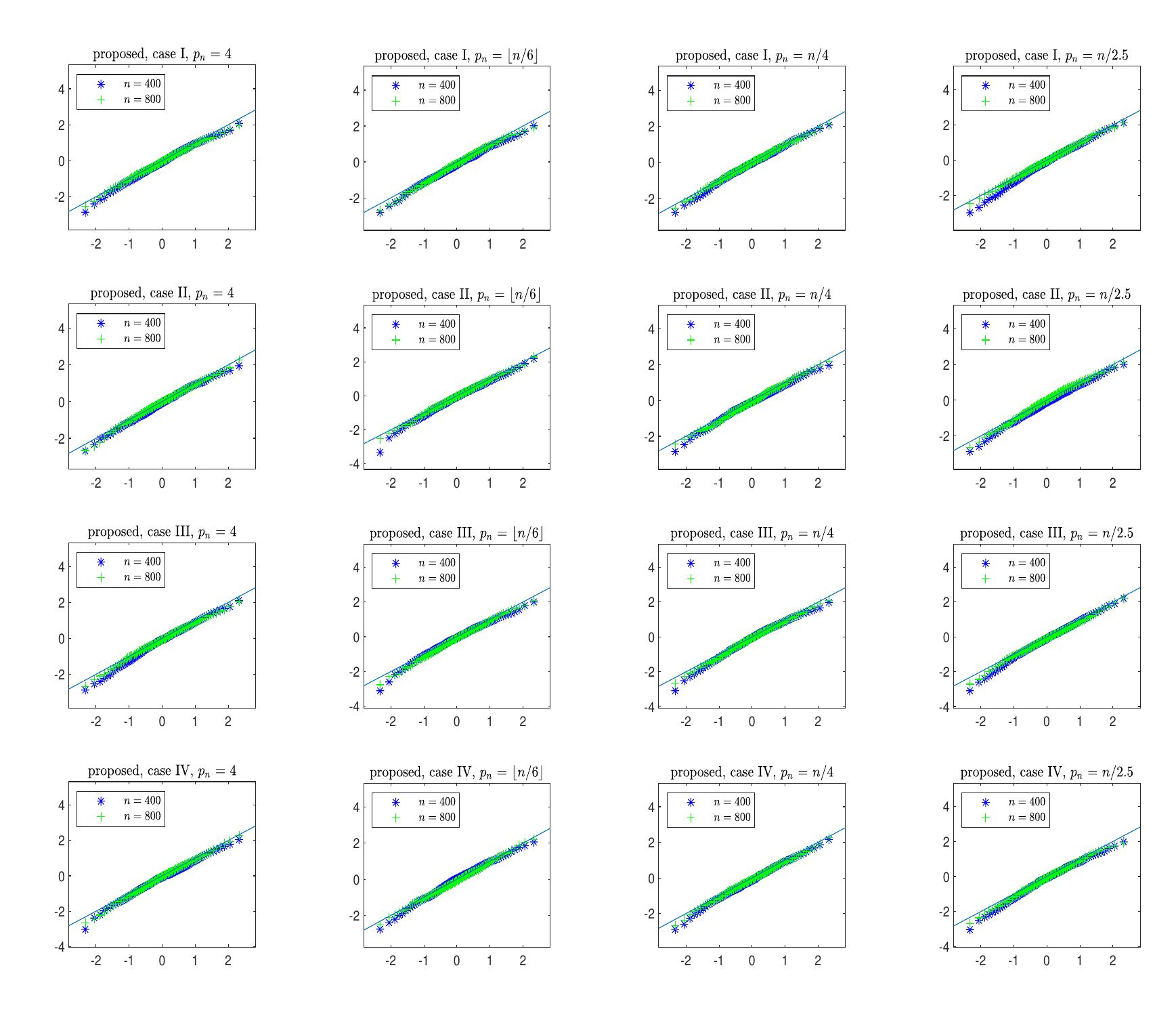}
\caption{\textsl{ QQ plots for testing $H_0$ in \eqref{4.1} using the proposed test statistics. Panels from top to bottom are for cases I--IV, respectively, while panels from left to right are for $p_n=4, \lfloor n/6\rfloor,n/4, n/2.5$, respectively.}}
\label{Figure_6}
\end{figure}

\begin{figure}[!ht]
\centering
\includegraphics[height = 6 in, width = 6 in]{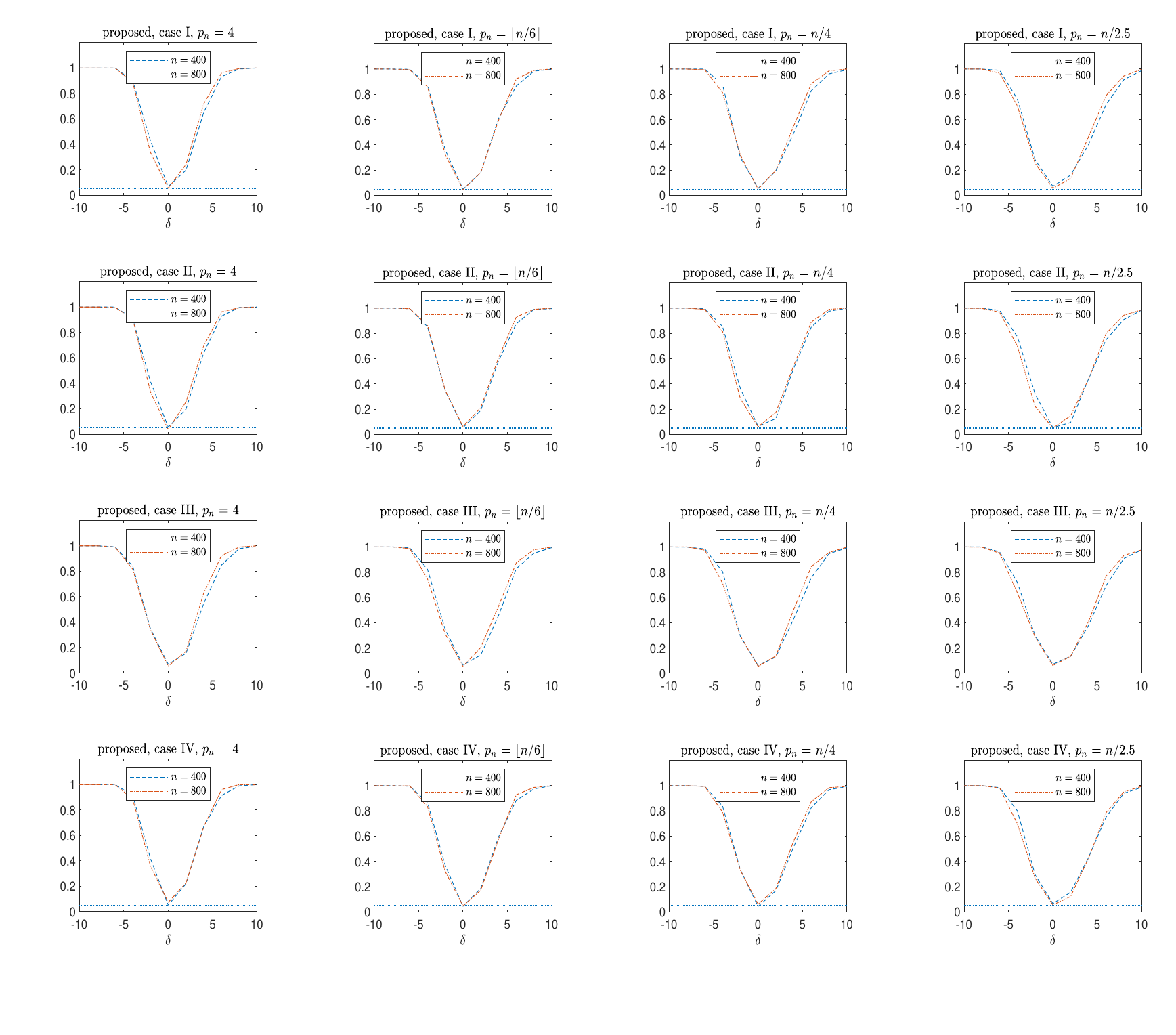}
\caption{\textsl{Empirical rejection rates versus $\delta$ for testing $H_1$ in \eqref{4.1} using the proposed test statistics. Panels from top to bottom are for cases I--IV, respectively, while panels from left to right are for $p_n=4, \lfloor n/6\rfloor,n/4, n/2.5$, respectively. The dotted line indicates the true significance level $\alpha=0.05$.}}
\label{Figure_7}
\end{figure}

Testing error variance:
\begin{eqnarray} \label{4.1}
H_0: \sigma_{\epsilon}^2 =  1  \quad &\text{versus}& \quad H_1:  \sigma_{\epsilon}^2 \neq 1
\end{eqnarray}
is performed by test statistic $ ( \hat \sigma_{\epsilon}^2  - 1 ) / {\hat \zeta_{\epsilon}  }$.
Figure \ref{Figure_6} provides the QQ plots of the test statistic under the null hypothesis. Clearly, the proposed test statistic well adapts to  fixed-, low- and moderate-dimensional regimes. The empirical powers under $\sigma_{\epsilon}^2 - 1 = \delta / n^{1/2}$ are provided in Figure  \ref{Figure_7}  with $\delta = -10, -8, \ldots,0, \ldots, 8, 10$. Again, the power behaviors are satisfactory.

\begin{figure}[!ht]
\centering
\includegraphics[height = 5 in, width = 6.5 in]{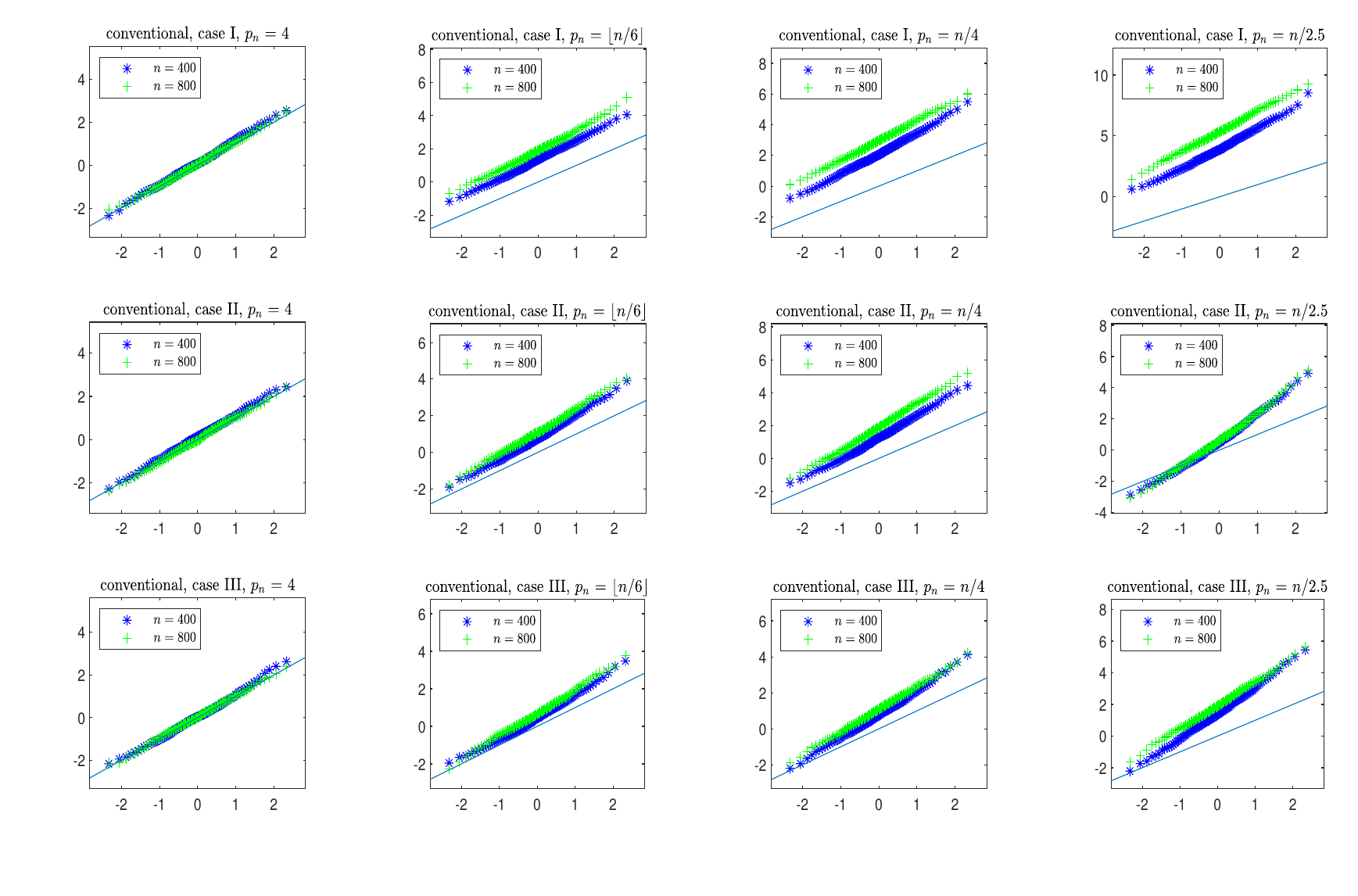}
\caption{\textsl{QQ plots for testing $H_0$ in \eqref{3.3} using the conventional  method. Panels from top to bottom are for cases I--III, respectively, while panels from left to right are for $p_n=4, \lfloor n/6\rfloor,n/4, n/2.5$, respectively. }}
\label{Figure_8}
\end{figure}

\begin{figure}[!ht]
\centering
\includegraphics[height = 5 in, width = 6.5 in]{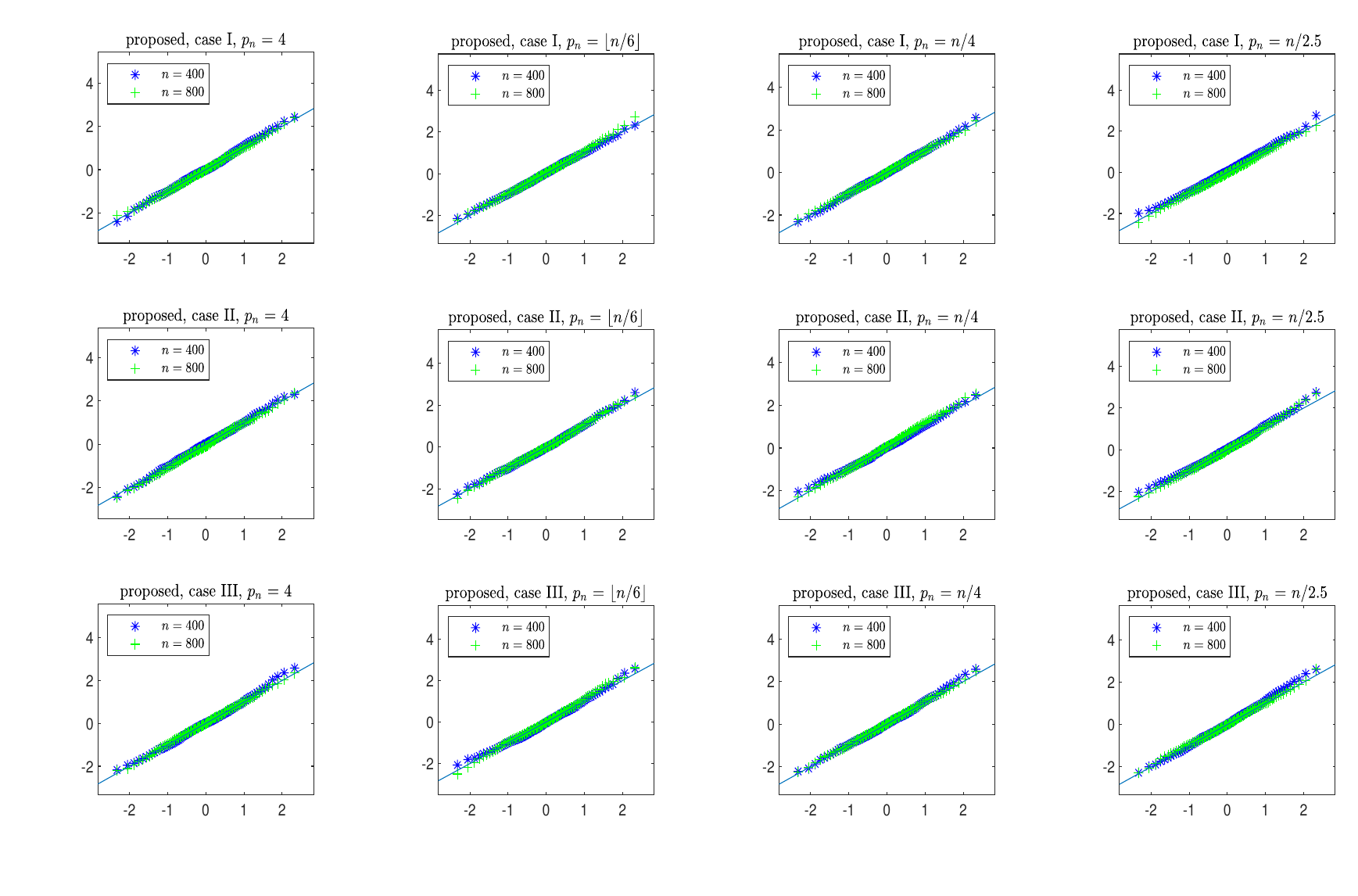}
\caption{\textsl{QQ plots for testing $H_0$ in \eqref{3.3} using the proposed method. Panels from top to bottom are for cases I--III, respectively, while panels from left to right are for $p_n=4, \lfloor n/6\rfloor,n/4, n/2.5$, respectively. }}
\label{Figure_9}
\end{figure}

We compare the performance of the conventional (in Section S.3)
and proposed test statistics for testing $H_0: \rho_0 \geq \rho_0^{\nullH0}$, i.e., \eqref{3.3}.
Figures \ref{Figure_8} and \ref{Figure_9}  provide the QQ   plots of the conventional and proposed test statistics, respectively. We find that both the conventional and proposed tests perform well for the fixed dimension. Under low and moderate dimensions, the conventional method fails but the proposed test continues to perform satisfactorily.

In the end, we consider the performance of the confidence intervals for three parameters $\bbeta_0^T \Sigma \bbeta_0$, $\sigma_{\epsilon}^2$ and $\rho_0$ in Tables \ref{Table-3}--\ref{Table-5}, respectively. The averaged coverage probability  and length of the confidence intervals are calculated using the proposed method, the MLE method in \cite{Dicker_Erdogdu_2016}, the method of moment for $\Sigma = \bfI_{p_n}$ ($\MM_1$) in \cite{Dicker_2014} (see Corollary 1 therein) and its alternative version for unknown $\Sigma$ ($\MM_2$) (see Proposition 2 therein). The  $\MLE$ and $\MM_1$ methods are applied by assuming that $\Sigma=\bfI_{p_n}$, i.e., $\Sigma$ is correctly identified in cases I and III but misidentified in cases II and IV.
Since the EigenPrism method in \cite{Janson_etal_2017} is particularly proposed for $p_n>n$, it is not compared with our methods.
For all three parameters, in case I with Gaussian design and $\Sigma =\bfI_{p_n}$, all methods are satisfactory with coverage probability close to the nominal level $95\%$.
For cases II--IV, our method still performs well with the coverage probability close to $95\%$ for all three parameters.
But, due to the misidentified $\Sigma$ or non-Gaussian design,
the coverage probabilities of the $\MLE$ method are away from $95\%$ for $ \bbeta_0^T \Sigma \bbeta_0$ and $\rho_0$ in case II and for $\sigma_{\epsilon}^2$ in cases III and IV, while the $\MM_1$ and $\MM_2$ methods  result in invalid confidence intervals in most situations.
Since $ \bbeta_0^T \Sigma \bbeta_0$ and $\rho_0$ are not defined for fixed design, the corresponding results for case IV are not reported.

\section{Real Data} \label{Section-5}
We study a data set from the International HapMap Project
 to investigate the relationship between gene expression
and single nucleotide polymorphism ($\SNP$).
In \cite{Stranger_etal_2007}, it's revealed that the expression levels of certain genes are associated with its nearby $\SNP$s.
Specifically, they identified 803 genes that were significantly associated
with certain $\SNP$s located within 1-Mbp of the gene midpoint  using
30
Caucasian trios of northern and western European origin (CEU),
45 unrelated Chinese individuals from Beijing University (CHB), 45
unrelated Japanese individuals from Tokyo (JPT), and 30 Yoruba trios
from Ibadan, Nigeria (YRI).
We select 9 genes among these 803 genes and investigate the relationship between each  gene  and its nearby $\SNP$s from $n=377$ individuals (80 individuals in CHB population, 82 from JPT, 107 from CEU and 108 from YRI). We use the gene expression dataset
from \url{https://www.ebi.ac.uk/arrayexpress/experiments/E-MTAB-264/} and \url{https://www.ebi.ac.uk/arrayexpress/experiments/E-MTAB-198/}. The $\SNP$ data were obtained from the International HapMap Project (\url{https://www.ncbi.nlm.nih.gov/variation/news/NCBI_retiring_HapMap/}),
release 28.

For each selected gene, we only focus on the $\SNP$s that are  significantly associated with this gene. A list of these significant $\SNP$s  for each identified gene is provided in Supplementary Table S2 of \cite{Stranger_etal_2007}.
Furthermore, among these significant $\SNP$s, we only choose those with minor allele frequency greater than $5\%$ and missingness no larger than $  20\%$. For the selected $\SNP$s included in the analysis, we impute the missing
values by the marginal mean. The minor allele counts are assigned as the numerical values for the $\SNP$s.

For each gene, we aim to regress the gene expression levels on the minor allele counts of its related $\SNP$s. We first center the gene expression levels
and $\SNP$ minor allele counts, and hence each variable has mean
0. Denote by $Y_{ik}$ the centered expression level for the $k$th gene $(k = 1, \ldots, 9)$ and
$i$th individual $(i = 1, \ldots, n = 377)$, and by $X_{ijk}$ the centered minor allele count
for the $j$th $\SNP$ ($j=1,\ldots, p_k$) corresponding to the $k$th gene and $i$th individual. For the $k$th gene, if the design matrix $X_k = \{X_{ijk}\}_{i=1,\ldots, n, j=1,\ldots,p_k}$ does not have full column rank, then we will randomly  delete one column  which is linearly correlated with the others and repeat this procedure until the design matrix is of full column rank.
With a slight abuse of notation, denote by $p_k$  the number of eventually selected $\SNP$s for the $k$th gene. The list of the selected genes and the corresponding $p_k$ are provided in Table \ref{Table-2} below. Our strategy for selecting the 9 genes in this study is that their corresponding $p_k$ is large enough, i.e., at least 33, such that $\tau$ ranges from 0.088 to 0.231.

The linear model for regressing  $\bfY_k=(Y_{1k},\ldots, Y_{nk})^T$ on $ X_k$ is fitted for each gene and the confidence intervals for $\rho_0$ are given in Table \ref{Table-2} using the conventional method in Section S.3
in the supplementary material, our proposed method and the $\MM_2$ method in \cite{Dicker_2014} (see Proposition 2 therein) proposed for general covariance matrix $\Sigma$. The $\MLE$ method in \cite{Dicker_Erdogdu_2016} and $\MM_1$ in \cite{Dicker_2014} are not designed for general $\Sigma$, and hence are not compared here.
We observe that the upper bounds of the confidence intervals of $\rho_0$ by our method are {bounded away from 1}, which indicates that the $\SNR$s are not exploding. Specifically, using \eqref{1.6}, we can calculate the confidence intervals for $\SNR$ directly  by those of $\rho_0$, and  we find that  the largest value of the upper bounds of the confidence intervals for $\SNR$ is 2.9526. Also, the values of $p_k^2/n$ range from 2.8886 to 20.0769, and hence, the strong signal condition (as in Theorem \ref{Theorem4}) that $p_k^2/n = o(\SNR) $ is not satisfied by this data.
From Table \ref{Table-2}, for most genes, the  confidence intervals of $\rho_0$ do not cover 0, indicating that these selected $\SNP$s are indeed significantly associated with the genes, which supports the findings in \cite{Stranger_etal_2007}.
However, for  gene AKAP10, 0 is covered by the confidence interval using our proposed method but excluded by the conventional method. This discrepancy may be due to the failure of the conventional method under moderate dimension and insufficient $\SNR$.
More importantly, we can see that the confidence intervals for  $\rho_0$ using our method are narrower than those using $\MM_2$ in \cite{Dicker_2014} for most genes, which means that our method is more accurate in the moderate-dimensional case with $\tau \in [0,1)$.

\begin{table}[!ht]
\caption{\textsl{ $90\%$ confidence intervals of  $\rho_0$ for gene data}}
\label{Table-2}
\begin{center}
\begin{tabular}{ l  cc ccc}
\hline
  Gene & Probe  & $p_k$ & conventional & proposed& $\MM_2$ \cite{Dicker_2014}  \\ \hline
AKAP10 & ILMN\_1718808 & 33 & [0.071, 0.161] & [0.000, 0.092] & [0.000, 0.040] \\
CPNE1 & ILMN\_1670841 & 35 & [0.293, 0.524] & [0.244, 0.504] & [0.259, 0.474] \\
NUDT13 & ILMN\_1680420 & 59 & [0.366, 0.467] & [0.294, 0.422] & [0.274, 0.482] \\
PIGN & ILMN\_1691112 & 36 & [0.225, 0.325] & [0.162, 0.280] & [0.183, 0.376] \\
PKHD1L1 & ILMN\_1717886 & 87 & [0.680, 0.747] & [0.640, 0.747] & [0.831, 1.000] \\
SPG7 & ILMN\_1675583 & 38 & [0.265, 0.375] & [0.212, 0.328] & [0.208, 0.406] \\
ST7L & ILMN\_1659926 & 40 & [0.548, 0.637] & [0.521, 0.627] & [0.522, 0.796] \\
TGM5 & ILMN\_1699925 & 39 & [0.298, 0.406] & [0.243, 0.368] & [0.232, 0.441] \\
TSGA10 & ILMN\_1674645 & 44 & [0.512, 0.613] & [0.479, 0.599] & [0.496, 0.761] \\
 \hline
\end{tabular}
\end{center}
\end{table}

\newpage

\newpage

\newpage

\newpage

\newpage

\newpage

\newpage

\begin{table}[!ht]
\caption{\textsl{{Coverage probabilities (Cov) and length (Len) of $95\%$ confidence intervals for $ \bbeta_0^T \Sigma \bbeta_0$ }}}
\label{Table-3}
\begin{center}
\small
\begin{tabular}{ c r ccccccccccc}
\hline
&& \multicolumn{2}{c}{ Proposed}  && \multicolumn{2}{c}{   $\MLE$} &&\multicolumn{2}{c}{ $\MM_1$} &&\multicolumn{2}{c}{$\MM_2$}   \\
\cline{3-4}\cline{6-7}\cline{9-10}\cline{12-13}
$n$ &$p_n$ & Cov &Len &&Cov & Len && Cov & Len && Cov & Len\\ \hline
\multicolumn{13}{c}{Case I} \\
400 & 4 & 0.941 & 0.478  && 0.940  & 0.395  && 0.943 & 0.681 && 0.959 & 0.682 \\
    & 66 & 0.939 & 0.489  && 0.947  & 0.440  && 0.949 & 0.712 && 0.951 & 0.712 \\
    & 100 & 0.946 & 0.504  && 0.950  & 0.468  && 0.933 & 0.729 && 0.934 & 0.730 \\
    & 160 & 0.957 & 0.529  && 0.948  & 0.514  && 0.960 & 0.762 && 0.960 & 0.763 \\ \cline{2-13}
800 & 4 & 0.950 & 0.338  && 0.949  & 0.278  && 0.939 & 0.479 && 0.953 & 0.479 \\
    & 133 & 0.940 & 0.350  && 0.940  & 0.312  && 0.951 & 0.506 && 0.954 & 0.506 \\
    & 200 & 0.947 & 0.356  && 0.954  & 0.332  && 0.955 & 0.517 && 0.952 & 0.517 \\
    & 320 & 0.948 & 0.373  && 0.948  & 0.363  && 0.932 & 0.536 && 0.933 & 0.536 \\ \hline
\multicolumn{13}{c}{Case II} \\
400 & 4 & 0.946 & 0.317  && 0.015  & 0.456  && 0.070 & 0.275 && 0.024 & 0.275 \\
    & 66 & 0.942 & 0.767  && 0.000  & 0.887  && 0.008 & 0.808 && 0.004 & 0.811 \\
    & 100 & 0.954 & 0.812  && 0.000  & 0.895  && 0.023 & 0.923 && 0.005 & 0.888 \\
    & 160 & 0.927 & 3.625  && 0.000  & 2.658  && 0.000 & 4.690 && 0.000 & 4.689 \\ \cline{2-13}
800 & 4 & 0.938 & 0.284  && 0.117  & 0.263  && 0.405 & 0.327 && 0.002 & 0.279 \\
    & 133 & 0.928 & 0.546  && 0.000  & 0.640  && 0.000 & 0.568 && 0.000 & 0.579 \\
    & 200 & 0.932 & 0.556  && 0.000  & 0.646  && 0.000 & 0.600 && 0.000 & 0.599 \\
    & 320 & 0.946 & 2.563  && 0.000  & 1.890  && 0.000 & 3.283 && 0.000 & 3.298 \\ \hline
    \multicolumn{13}{c}{Case III} \\
400 & 4 & 0.938 & 1.265  && 0.945  & 0.686  && 0.871 & 1.804 && 0.936 & 1.796 \\
    & 66 & 0.938 & 1.260  && 0.948  & 0.756  && 0.872 & 1.852 && 0.868 & 1.844 \\
    & 100 & 0.944 & 1.285  && 0.951  & 0.796  && 0.870 & 1.888 && 0.876 & 1.880 \\
    & 160 & 0.933 & 1.290  && 0.937  & 0.875  && 0.867 & 1.943 && 0.860 & 1.935 \\ \cline{2-13}
800 & 4 & 0.931 & 0.912  && 0.939  & 0.483  && 0.848 & 1.274 && 0.913 & 1.270 \\
    & 133 & 0.933 & 0.914  && 0.954  & 0.532  && 0.850 & 1.315 && 0.850 & 1.312 \\
    & 200 & 0.947 & 0.919  && 0.955  & 0.561  && 0.870 & 1.334 && 0.868 & 1.331 \\
    & 320 & 0.942 & 0.925  && 0.937  & 0.617  && 0.862 & 1.364 && 0.855 & 1.361 \\ \hline
\end{tabular}
\end{center}
\end{table}

\newpage

\begin{table}[!ht]
\caption{\textsl{{Coverage probabilities (Cov) and length (Len) of  $95\%$ confidence intervals for $\sigma_{\epsilon}^2$ }}}
\label{Table-4}
\begin{center}
\small
\begin{tabular}{ c r ccccccccccc}
\hline
&& \multicolumn{2}{c}{ Proposed}  && \multicolumn{2}{c}{   $\MLE$} &&\multicolumn{2}{c}{ $\MM_1$} &&\multicolumn{2}{c}{$\MM_2$}   \\
\cline{3-4}\cline{6-7}\cline{9-10}\cline{12-13}
$n$ &$p_n$ & Cov &Len &&Cov & Len && Cov & Len && Cov & Len\\ \hline
\multicolumn{13}{c}{Case I} \\
400 & 4 & 0.940 & 0.276  && 0.944  & 0.278  && 0.953 & 0.398 && 0.973 & 0.398 \\
    & 66 & 0.938 & 0.299  && 0.944  & 0.301  && 0.951 & 0.452 && 0.953 & 0.453 \\
    & 100 & 0.938 & 0.318  && 0.938  & 0.319  && 0.951 & 0.483 && 0.953 & 0.483 \\
    & 160 & 0.913 & 0.354  && 0.921  & 0.350  && 0.967 & 0.529 && 0.966 & 0.529 \\ \cline{2-13}
800 & 4 & 0.952 & 0.196  && 0.956  & 0.196  && 0.962 & 0.279 && 0.978 & 0.279 \\
    & 133 & 0.941 & 0.213  && 0.941  & 0.214  && 0.953 & 0.321 && 0.956 & 0.321 \\
    & 200 & 0.946 & 0.226  && 0.946  & 0.226  && 0.948 & 0.341 && 0.949 & 0.341 \\
    & 320 & 0.956 & 0.252  && 0.959  & 0.249  && 0.951 & 0.372 && 0.953 & 0.372 \\ \hline
\multicolumn{13}{c}{Case II} \\
400 & 4 & 0.952 & 0.276  && 0.956  & 0.278  && 0.144 & 0.357 && 0.088 & 0.379 \\
    & 66 & 0.936 & 0.299  && 0.946  & 0.305  && 0.000 & 0.698 && 0.000 & 0.721 \\
    & 100 & 0.951 & 0.317  && 0.949  & 0.326  && 0.000 & 0.759 && 0.000 & 0.786 \\
    & 160 & 0.935 & 0.352  && 0.938  & 0.361  && 0.000 & 3.489 && 0.000 & 3.574 \\ \cline{2-13}
800 & 4 & 0.945 & 0.196  && 0.944  & 0.196  && 0.209 & 0.259 && 0.001 & 0.285 \\
    & 133 & 0.951 & 0.214  && 0.954  & 0.216  && 0.000 & 0.495 && 0.000 & 0.510 \\
    & 200 & 0.942 & 0.224  && 0.936  & 0.231  && 0.000 & 0.525 && 0.000 & 0.538 \\
    & 320 & 0.948 & 0.253  && 0.937  & 0.258  && 0.000 & 2.456 && 0.000 & 2.513 \\ \hline
    \multicolumn{13}{c}{Case III} \\
400 & 4 & 0.926 & 0.304  && 0.906  & 0.277  && 0.883 & 0.868 && 0.968 & 0.894 \\
    & 66 & 0.942 & 0.328  && 0.925  & 0.303  && 0.879 & 0.958 && 0.889 & 1.005 \\
    & 100 & 0.940 & 0.345  && 0.919  & 0.319  && 0.878 & 1.004 && 0.896 & 1.060 \\
    & 160 & 0.921 & 0.377  && 0.912  & 0.355  && 0.892 & 1.076 && 0.891 & 1.147 \\ \cline{2-13}
800 & 4 & 0.943 & 0.218  && 0.920  & 0.196  && 0.853 & 0.622 && 0.957 & 0.627 \\
    & 133 & 0.938 & 0.233  && 0.920  & 0.214  && 0.859 & 0.697 && 0.869 & 0.707 \\
    & 200 & 0.944 & 0.245  && 0.929  & 0.226  && 0.890 & 0.734 && 0.892 & 0.742 \\
    & 320 & 0.936 & 0.271  && 0.927  & 0.251  && 0.892 & 0.785 && 0.895 & 0.800 \\ \hline
    \multicolumn{13}{c}{Case IV} \\
400 & 4 & 0.945 & 0.277  && 0.951  & 0.277  && 0.978 & 0.407 && 0.877 & 0.411 \\
    & 66 & 0.941 & 0.301  && 0.940  & 0.300  && 0.000 & 0.004 && 0.000 & 1.141 \\
    & 100 & 0.936 & 0.315  && 0.912  & 0.309  && 0.000 & 0.001 && 0.000 & 1.212 \\
    & 160 & 0.931 & 0.351  && 0.812  & 0.326  && 0.000 & 0.000 && 0.000 & 1.773 \\ \cline{2-13}
800 & 4 & 0.932 & 0.195  && 0.935  & 0.196  && 0.960 & 0.286 && 0.951 & 0.285 \\
    & 133 & 0.946 & 0.214  && 0.936  & 0.212  && 0.000 & 0.000 && 0.000 & 0.748 \\
    & 200 & 0.945 & 0.225  && 0.900  & 0.219  && 0.000 & 0.000 && 0.000 & 0.924 \\
    & 320 & 0.948 & 0.251  && 0.725  & 0.231  && 0.000 & 0.000 && 0.000 & 1.383 \\ \hline
\end{tabular}
\end{center}
\end{table}

\newpage

\begin{table}[!ht]
\caption{\textsl{{Coverage probabilities (Cov) and length (Len) of  $95\%$ confidence intervals for $\rho_0$ }}}
\label{Table-5}
\begin{center}
\small
\begin{tabular}{ c r ccccccccccc}
\hline
&& \multicolumn{2}{c}{ Proposed}  && \multicolumn{2}{c}{   $\MLE$} &&\multicolumn{2}{c}{ $\MM_1$} &&\multicolumn{2}{c}{$\MM_2$}   \\
\cline{3-4}\cline{6-7}\cline{9-10}\cline{12-13}
$n$ &$p_n$ & Cov &Len &&Cov & Len && Cov & Len && Cov & Len\\ \hline
\multicolumn{13}{c}{Case I} \\
400 & 4 & 0.944 & 0.138  && 0.940  & 0.121  && 0.949 & 0.241 && 0.970 & 0.241 \\
    & 66 & 0.955 & 0.150  && 0.958  & 0.142  && 0.945 & 0.264 && 0.950 & 0.264 \\
    & 100 & 0.942 & 0.160  && 0.950  & 0.154  && 0.938 & 0.275 && 0.941 & 0.276 \\
    & 160 & 0.959 & 0.178  && 0.948  & 0.176  && 0.966 & 0.296 && 0.966 & 0.296 \\ \cline{2-13}
800 & 4 & 0.955 & 0.098  && 0.942  & 0.085  && 0.951 & 0.170 && 0.970 & 0.170 \\
    & 133 & 0.945 & 0.107  && 0.939  & 0.100  && 0.964 & 0.187 && 0.963 & 0.187 \\
    & 200 & 0.961 & 0.113  && 0.951  & 0.109  && 0.951 & 0.195 && 0.952 & 0.195 \\
    & 320 & 0.948 & 0.127  && 0.947  & 0.125  && 0.937 & 0.209 && 0.942 & 0.209 \\ \hline
\multicolumn{13}{c}{Case II} \\
400 & 4 & 0.935 & 0.150  && 0.026  & 0.135  && 0.008 & 0.159 && 0.000 & 0.164 \\
    & 66 & 0.951 & 0.115  && 0.001  & 0.069  && 0.000 & 0.230 && 0.000 & 0.234 \\
    & 100 & 0.951 & 0.117  && 0.023  & 0.077  && 0.000 & 0.251 && 0.000 & 0.248 \\
    & 160 & 0.947 & 0.033  && 0.000  & 0.015  && 0.000 & 0.287 && 0.000 & 0.291 \\ \cline{2-13}
800 & 4 & 0.936 & 0.103  && 0.203  & 0.083  && 0.211 & 0.146 && 0.000 & 0.134 \\
    & 133 & 0.938 & 0.082  && 0.000  & 0.049  && 0.000 & 0.162 && 0.000 & 0.167 \\
    & 200 & 0.939 & 0.085  && 0.001  & 0.055  && 0.000 & 0.173 && 0.000 & 0.175 \\
    & 320 & 0.933 & 0.024  && 0.000  & 0.011  && 0.000 & 0.203 && 0.000 & 0.206 \\ \hline
    \multicolumn{13}{c}{Case III} \\
400 & 4 & 0.938 & 0.098  && 0.919  & 0.068  && 0.876 & 0.268 && 0.960 & 0.270 \\
    & 66 & 0.954 & 0.102  && 0.934  & 0.076  && 0.868 & 0.289 && 0.873 & 0.285 \\
    & 100 & 0.940 & 0.106  && 0.931  & 0.083  && 0.868 & 0.300 && 0.881 & 0.294 \\
    & 160 & 0.930 & 0.113  && 0.914  & 0.095  && 0.874 & 0.319 && 0.874 & 0.307 \\ \cline{2-13}
800 & 4 & 0.946 & 0.070  && 0.939  & 0.048  && 0.842 & 0.190 && 0.950 & 0.191 \\
    & 133 & 0.930 & 0.073  && 0.937  & 0.054  && 0.848 & 0.205 && 0.857 & 0.205 \\
    & 200 & 0.948 & 0.076  && 0.939  & 0.058  && 0.872 & 0.213 && 0.870 & 0.213 \\
    & 320 & 0.954 & 0.081  && 0.941  & 0.067  && 0.880 & 0.226 && 0.885 & 0.224 \\ \hline
\end{tabular}
\end{center}
\end{table}

\newpage

\setcounter{equation}{0}
\renewcommand{\theequation} {A.\arabic{equation}}
\setcounter{theorem}{0}
\renewcommand{\thetheorem} {A.\arabic{theorem}}
\setcounter{lemma}{0}
\renewcommand{\thelemma} {A.\arabic{lemma}}
\setcounter{proposition}{0}
\renewcommand{\theproposition} {A.\arabic{proposition}}
\setcounter{section}{0}
\renewcommand{\thesection} {A.\arabic{section}}
\begin{center}
\bf \Large Appendix
\end{center}
\section{Inference for single element and linear functional of $\bbeta_0$
}
\label{app-A1}
 We provide a brief discussion of the
element-wise inference  for $\beta_{0,j}$ and $\beta_{0,j}^2$ $(j=1, \ldots, p_n)$. The estimator for $\beta_{0,j}$ is $$\hat \beta_{j} = \bfe_{j,p_n}^T \hat\bbeta = \beta_{0,j} + \bfe_{j,p_n}^T (X^TX)^{-1}X^T \beps.$$ If $\sigma^2_{\epsilon} = o(n)$, then  $\sigma^{-1}_{\hat \beta_{j} }(\hat \beta_{j} - \beta_{0,j})  \conD N(0,1),$ where $\sigma^{2}_{\hat \beta_{j} } = \sigma^2_{\epsilon} \bfe_{j,p_n}^T  \E\{(X^TX)^{-1} \ID(K)\}\bfe_{j,p_n} = \Omega( \sigma^2_{\epsilon} / n )$. If $\sigma^2_{\epsilon} = \Omega(1)$ and $ \bfX_i \sim N(\bfzero_{p_n}, \bfI_{p_n})$, the bias of $\hat \beta_{j}^2$ is
$$\E (\hat \beta_{j}^2) - \beta_{0,j}^2 = \sigma^2_{\epsilon} \bfe_{j,p_n}^T \E\{(X^TX)^{-1}  \}\bfe_{j,p_n}  = \Omega(1/n).$$
 Therefore, the bias of $\hat \beta_{j}^2$ is ignorable if $n^{-1} = o(\beta_{0,j}^2) $. Inference for $ \beta_{0,j}^2$ can be conducted using $ (2 \beta_{0,j} \sigma_{\hat \beta_{j} })^{-1}(\hat \beta_{j}^2 - \beta_{0,j}^2)  \conD N(0,1)$.
The $\sqrt n (\hat \beta_{j}^2 - \beta_{0,j}^2)$ versus $j$ and the bias for $\|\hat\bbeta\|_2^2$ are plotted in Figure \ref{Figure_10} below.

\vspace{0.7cm}

\begin{figure}[!htb]
\centering
\includegraphics[height = 2 in, width = 2.5 in]{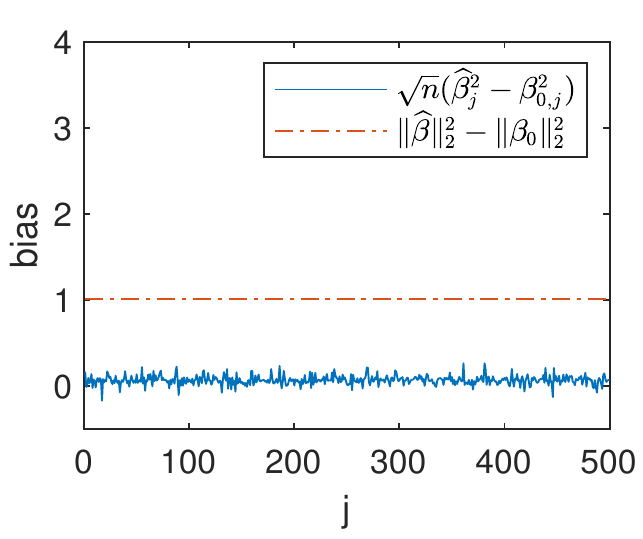}
\caption{\textsl{Plots of $\sqrt n (\hat \beta_{j}^2 - \beta_{0,j}^2)$ versus $j$ (solid line) and bias for $\|\hat\bbeta\|_2^2$ (dash-dotted line) under the same setting as in Figure \ref{Figure_1} with $p_n = 500$. }}
\label{Figure_10}
\end{figure}

We then discuss the inference for linear functionals $\bfc^T \bbeta_0$ where $\bfc \in \mathbb{R}^{p_n}$ is deterministic. The estimator for $\bfc^T \bbeta_0$ is $\bfc^T \hat\bbeta=\bfc^T\bbeta_0 + \bfc^T(X^TX)^{-1}X^T \beps$. Following the proof of Theorem \ref{Theorem1}, its limiting distribution is $\sigma^{-1}_{L }(\bfc^T \hat\bbeta- \bfc^T \bbeta_0)  \conD N(0,1)$ where $\sigma^{2}_{L } = \sigma^2_{\epsilon} \bfc^T  \E\{(X^TX)^{-1} \ID(K)\}\bfc $ with a ratio consistent estimator $\hat\sigma^{2}_{L } = \hat\sigma^2_{\epsilon} \bfc^T (X^TX)^{-1} \bfc $.

\section{Relation between $\tau$ and $\SNR$} \label{app-A2}
\begin{figure}[!ht]
\begin{center}
\begin{tikzpicture}
  \matrix (m) [matrix of math nodes,row sep=7em,column sep=16em,minimum width=1.8em]
  {
     \text{$\tau=0$}&
     \text{$\text{\small Strong Signal}$ }\\
     \text{$\tau > 0$  } &
      \text{$\text{\small Weak Signal} $}\\};
  \path[-stealth]
    (m-1-1) edge node [above] {$\sigma_{\epsilon}^2 = \Omega(1)$, $\|\bbeta_0\|_2^2 = \Omega(p_n)$   } (m-1-2)
            edge node [pos= 0.3,above,sloped] {$\sigma_{\epsilon}^2\geq C>0$} node [pos = 0.8, above, sloped] { $\|\bbeta_0\|_2^2 = O(1/n)$}   (m-2-2)
    (m-2-1) edge node [pos= 0.25,above, sloped]{$\sigma_{\epsilon}^2=o(1)$  }  node [pos= 0.7,above, sloped] {$\|\bbeta_0\|_2^2 = \Omega(p_n)$}  (m-1-2)
    edge node [above]{$\sigma_{\epsilon}^2\geq C>0$, $\|\bbeta_0\|_2^2 = O(p_n)$} (m-2-2);
\end{tikzpicture}
\end{center}
\caption{\textsl{Relation between $\tau = 0$/$\tau > 0$ and Strong/Weak Signal.}}
\label{Figure_11}
\end{figure}
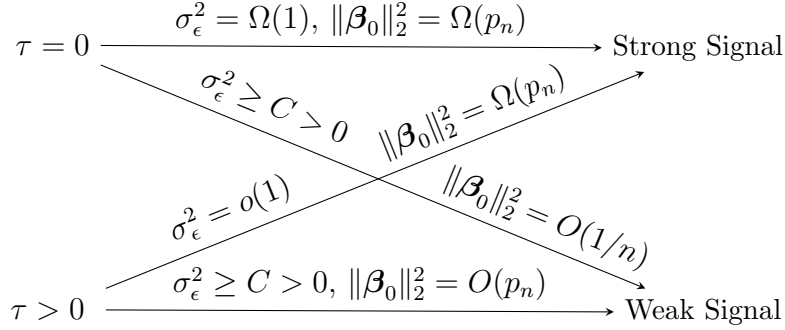

According to Theorem \ref{Theorem4},
 define
	\begin{itemize}
		\item strong $\SNR$: $ p_n^2/n = o(\SNR)$;
		\item weak $\SNR$: $\SNR \lesssim p_n^2/n$.
	\end{itemize}
Figure \ref{Figure_11} describes the precise relation between $\tau$ and the signal strength under mild conditions.
	In particular, $\tau=0/\tau>0$ may imply strong/weak signals unless we allow $\|\bbeta_0\|_2$ or $\sigma_{\epsilon}^2$ to diminish.

\section{Proofs of main theoretical results} \label{app-A3}
This section includes the proofs of Lemmas  \ref{Lemma-1} and \ref{Lemma-2} and Theorems \ref{Theorem1} and \ref{Theorem3}. In all the proofs, we only consider the case that $\sigma_{\epsilon}^2$ is fixed. The results for diverging $\sigma_{\epsilon}^2 $ can be simply obtained by replacing $Y_i$, $\beps $ and $\bbeta_0 $ with $Y_i/\sigma_{\epsilon}$, $\beps /\sigma_{\epsilon}$ and $\bbeta_0 /\sigma_{\epsilon} $ respectively in the proofs.

We introduce some notations and equations.
Let $X_{(i)} = (\bfX_1, \ldots,\bfX_{i-1},\bfX_{i+1}, \ldots, \bfX_n)^T$ for $i=1,\ldots, n$, i.e. the design matrix without the $i$th observation. Similarly, $X_{(i,j)}$ denotes the design matrix without the $i$th  and $j$th observations for $1 \leq i\neq j \leq n$. From the Sherman-Morrison formula \citep{Sherman_Morrison_1950},
\begin{eqnarray} \label{A.1}
 (X^T X)^{-1}   =   (X_{(1)}^T X_{(1)} + \bfX_1 \bfX_1^T)^{-1}=(X_{(1)}^T X_{(1)})^{-1}   -\frac{ (X_{(1)}^T X_{(1)})^{-1} \bfX_1 \bfX_1^T (X_{(1)}^T X_{(1)})^{-1}}{1 + \bfX_1^T (X_{(1)}^T X_{(1)})^{-1}\bfX_1}, \quad
\end{eqnarray}
and hence,
\begin{eqnarray} \label{A.2}
  (X^T X)^{-2}
  &=& (X_{(1)}^T X_{(1)})^{-2}  - (X_{(1)}^T X_{(1)})^{-2} \bfX_1 \bfX_1^T (X_{(1)}^T X_{(1)})^{-1} /\{1 + \bfX_1^T (X_{(1)}^T X_{(1)})^{-1}\bfX_1\}\cr
&&                                         - (X_{(1)}^T X_{(1)})^{-1} \bfX_1 \bfX_1^T (X_{(1)}^T X_{(1)})^{-2} /\{1 + \bfX_1^T (X_{(1)}^T X_{(1)})^{-1}\bfX_1\}\cr
&&  + \{(X_{(1)}^T X_{(1)})^{-1} \bfX_1 \bfX_1^T (X_{(1)}^T X_{(1)})^{-1}\}^2 /\{1 + \bfX_1^T (X_{(1)}^T X_{(1)})^{-1}\bfX_1\}^2.
\end{eqnarray}
Therefore,
\begin{eqnarray}
(X^T X)^{-1} \bfX_1 &=&  \frac{(X_{(1)}^T X_{(1)})^{-1}\bfX_1 }{ 1 + \bfX_1^T (X_{(1)}^T X_{(1)})^{-1}\bfX_1 }, \label{A.3} \\
(X^T X)^{-2} \bfX_1&=& \frac{(X_{(1)}^T X_{(1)})^{-2} \bfX_1 } { 1 + \bfX_1^T (X_{(1)}^T X_{(1)})^{-1}\bfX_1 }
 - \frac{ (X_{(1)}^T X_{(1)})^{-1}  \bfX_1 \bfX_1^T (X_{(1)}^T X_{(1)})^{-2} \bfX_1 }{\{1 + \bfX_1^T (X_{(1)}^T X_{(1)})^{-1}\bfX_1\}^2}, \qquad \label{A.4} \\
\bfX_1^T(X^T X)^{-2} \bfX_1&=& \frac{\bfX_1^T(X_{(1)}^T X_{(1)})^{-2}\bfX_1 }{\{1 + \bfX_1^T (X_{(1)}^T X_{(1)})^{-1}\bfX_1\}^2}. \label{A.5}
\end{eqnarray}

The following are the proofs of the main results in this paper.

\ni\textbf{Proof of Lemma \ref{Lemma-1}: } For $\bfZ_i = (z_{i1}, \ldots, z_{ip_n})^T $ defined in Condition $\mathrm{A1}$,
let $z_{ij}^* = z_{ij} \ID(|z_{ij}| \leq \sqrt n / \sqrt{\log n}) - \E\{z_{ij} \ID(|z_{ij}| \leq \sqrt n / \sqrt{\log n})\}$, $\tilde z_{ij} = z_{ij} -  z_{ij}^* = z_{ij} \ID(|z_{ij}| > \sqrt n / \sqrt{\log n}) + \E\{z_{ij} \ID(|z_{ij}| \leq \sqrt n / \sqrt{\log n})\}$, $\bfZ_i^* = (z_{i1}^*, \ldots, z_{ip_n}^*)^T$, $\tilde \bfZ_i = (\tilde z_{i1}, \ldots, \tilde z_{ip_n})^T$, $ Z^* = ({\bfZ}_1^*, \ldots, {\bfZ}_n^*)^T = (z_{ij}^*) _{i\leq n, j\leq p_n}$ and $\tilde Z = (\tilde{\bfZ}_1, \ldots, \tilde{\bfZ}_n)^T = (\tilde z_{ij}) _{i\leq n, j\leq p_n}$.

Then, $\E(z_{ij}^*) = 0$ and by Cauchy-Schwarz inequality and Chebyshev's inequality,
\begin{eqnarray*}
  1 - \E(z_{ij}^{*2} )
&=& 1 - \E\{z_{ij}^2 \ID(|z_{ij}| \leq \sqrt n / \sqrt{\log n})\} +  [\E\{z_{ij} \ID(|z_{ij}| \leq \sqrt n / \sqrt{\log n})\}]^2 \cr
&=& 1 -1  + \E\{z_{ij}^2 \ID(|z_{ij}| > \sqrt n / \sqrt{\log n})\} +  [\E\{z_{ij} \ID(|z_{ij}| > \sqrt n / \sqrt{\log n})\}]^2 \cr
&\leq& 2 \E\{z_{ij}^2 \ID(|z_{ij}| > \sqrt n / \sqrt{\log n})\} \leq 2 \{ \E (z_{ij}^4) \pr (|z_{ij}| > \sqrt n / \sqrt{\log n})\}^{1/2} \cr
&\lesssim& \{ \pr (|z_{ij}| > \sqrt n / \sqrt{\log n}) \}^{1/2} \leq \{\E (z_{ij}^4) / (\sqrt n / \sqrt{\log n})^4\}^{1/2} \cr
&\lesssim&  (\log n) / n,
\end{eqnarray*}
which implies that $ \max_{j \leq p_n} \sum_{i=1}^n |1 - \E(z_{ij}^{*2} )| \lesssim \log n = o(n) $. Also,
 \begin{eqnarray*}
&& \sup_{i\leq n,j\leq p_n,n\geq 1} \E(z_{ij}^{*4}) \cr
&\lesssim& \sup_{i\leq n,j\leq p_n,n\geq 1} \big( \E\{ z_{ij} \ID(|z_{ij}| \leq \sqrt n / \sqrt{\log n})\}^4
+[ \E\{z_{ij} \ID(|z_{ij}| \leq \sqrt n / \sqrt{\log n})\}]^4\big) \cr
&\lesssim&  \sup_{i\leq n,j\leq p_n,n\geq 1} \E\{ z_{ij} \ID(|z_{ij}|  \leq \sqrt n / \sqrt{\log n})\}^4 +C \cr
&\leq& \sup_{i\leq n,j\leq p_n,n\geq 1} \E( z_{ij} ^4) +C \leq 2C < \infty.
\end{eqnarray*}
It's easy to see $|z_{ij}^* | \leq \sqrt n / \sqrt{\log n}$.
From Theorem 9.13 of \cite{Bai_Silverstein_2010}, for any $s_1 > (1 + \sqrt{\tau})^2$, $s_2 < (1 - \sqrt{\tau})^2$ and any $\ell > 0$, we have
$$\pr(\|Z^{*T} Z^* / n\|_2 > s_1)= o(n^{ - \ell}), \quad \pr ( \|(Z^{*T}Z^* / n)^{-1} \|_2 > 1 / s_2) = o(n^{ - \ell}).$$
Since $Z = Z^* + \tilde Z$, we have $ Z^TZ/n = Z^{*T} Z^* / n + \tilde Z ^T Z^* / n + Z^{*T} \tilde Z / n + \tilde Z^T \tilde Z /n$. We know that
\begin{eqnarray} \label{A.6}
&& \| \tilde Z^T \tilde Z /n \|_2 \leq \| \tilde Z^T \tilde Z /n \|_1 = \max_{i \leq p_n} \sum_{j=1}^{p_n}\Big| \sum_{t=1}^n \tilde z_{ti} \tilde z_{tj} /n \Big|.
\end{eqnarray}
Note  $\E(\tilde z_{ti} )=0$, $\E(\tilde z_{ti} \tilde z_{tj}) = 0$ for $i \neq j$, and, for any integer $k >0$, from Cauchy-Schwarz inequality,
\begin{eqnarray*}
 \E |\tilde z_{ti}|^k &=& \E |z_{ti} \ID(|z_{ti}| > \sqrt n / \sqrt{\log n}) + \E\{z_{ti} \ID(|z_{ti}| \leq \sqrt n / \sqrt{\log n})\}|^k \cr
&=& \E |z_{ti} \ID(|z_{ti}| > \sqrt n / \sqrt{\log n}) - \E\{z_{ti} \ID(|z_{ti}| > \sqrt n / \sqrt{\log n})\}|^k \cr
&\lesssim& \E |z_{ti} \ID(|z_{ti}| > \sqrt n / \sqrt{\log n}) |^k \leq \{\E  |z_{ti} |^{2k} \pr (|z_{ti}| > \sqrt n / \sqrt{\log n})\}^{1/2} \cr
&\lesssim& \{\pr (|z_{ti}| > \sqrt n / \sqrt{\log n})\}^{1/2}.
\end{eqnarray*}

Therefore,
for any integer $\ell >0$, taking $x = 1/\sqrt n$, from \eqref{A.6} and Markov's inequality,
\begin{eqnarray*}
&& \pr( \| \tilde Z^T \tilde Z /n \|_2 > x) \cr
&\leq& \pr\Big( \max_{i \leq p_n} \sum_{j=1}^{p_n}\Big| \sum_{t=1}^n \tilde z_{ti} \tilde z_{tj} /n \Big| >x \Big) \leq p_n  \pr\Big(  \sum_{j=1}^{p_n}\Big| \sum_{t=1}^n \tilde z_{ti} \tilde z_{tj} /n \Big| >x \Big) \cr
&\leq&  p_n ^2 \pr\Big(  \Big| \sum_{t=1}^n \tilde z_{ti} \tilde z_{tj} /n \Big| >x /p_n \Big)
\leq p_n ^2 \pr\Big(  \Big| \sum_{t=1}^n \tilde z_{ti} ^2 /n \Big| \Big| \sum_{t=1}^n \tilde z_{tj} ^2 /n \Big| >x^2 /p_n^2 \Big) \cr
&\lesssim&  p_n ^2 \pr\Big(  \Big| \sum_{t=1}^n \tilde z_{ti} ^2 /n \Big| >x /p_n \Big)
\leq p_n ^2 \E\Big( \Big| \sum_{t=1}^n \tilde z_{ti} ^2 /n \Big|^{2\ell} \Big) \Big / (x /p_n)^{2\ell}  \cr
&=& p_n ^{{2\ell} +2} x^{-{2\ell} } n^{-{2\ell} }  \E\Big( \Big| \sum_{t=1}^n \tilde z_{ti} ^2 \Big|^{2\ell}  \Big)\cr
&\lesssim&  p_n ^{{2\ell} +2} x^{-{2\ell} } n^{-{2\ell} }  n^{{2\ell} }  \E |\tilde z_{ti} ^2 |^{2\ell}
\lesssim  p_n ^{{2\ell} +2} x^{-{2\ell} }
 \E |\tilde z_{ti}|^{4\ell}    \cr
&\lesssim& p_n ^{{2\ell} +2} x^{-{2\ell} }  \{ \pr (|z_{ti}| > \sqrt n / \sqrt{\log n}) \}^{1/2 }\cr
&\lesssim&   p_n ^{{2\ell} +2} x^{-{2\ell} }   \{\E |z_{ti}|^{28\ell } / (\sqrt n / \sqrt{\log n})^{28\ell}\}^{1/2} \cr
&\lesssim&  p_n ^{{2\ell} +2} x^{-{2\ell} }  ( \sqrt{\log n}/ \sqrt n)^{{14\ell} } \cr
&\leq& (\log n)^{7\ell} n^{-5\ell+2} x^{-{2\ell} } = (\log n)^{7\ell} n^{-4\ell+2}  = o(n^{- \ell}) .
\end{eqnarray*}
Then, for $n$ large enough,
\begin{eqnarray*}
&& \pr( \| \tilde Z ^T Z^* / n \|_2 > 1/  {\log n}) \leq \pr( \| \tilde Z ^T\|_2 \| Z^*  \|_2 /n > 1/ {\log n}) \cr
&=& \pr( \| \tilde Z ^T  \tilde Z / n \|_2 \| Z^{*T} Z^* /n\|_2 > 1 /(\log n)^2 )\cr
&\leq& \pr( \| \tilde Z ^T  \tilde Z / n \|_2  > n^{-1/4} / \log n) + \pr(  \| Z^{*T} Z^* /n\|_2 > n^{1/4} / \log n ) = o(n^{-\ell}).
\end{eqnarray*}
Therefore, taking $\mu_1 = 4 (1 + \sqrt{\tau})^2$ and $\mu_2 = (1 - \sqrt{\tau})^2/4$, we have
\begin{eqnarray*}
  \pr( \|   Z ^T Z / n \|_2 \geq \mu_1 )
& \leq& \pr( \| Z^{*T} Z^* / n\|_2 > \mu_1 / 2) + \pr( \| \tilde Z ^T Z^* / n\|_2 > \mu_1 / 8) \cr
&&\quad  +\pr( \|Z^{*T} \tilde Z / n\|_2 > \mu_1 / 8)+\pr( \|\tilde Z^T \tilde Z /n\|_2 > \mu_1 / 8)  = o(n^{-\ell}),
\end{eqnarray*}
and
\begin{eqnarray*}
&& \pr( \|   (Z ^T Z / n)^{-1} \|_2 \geq 1/ \mu_2 )= \pr(\lambda_{\min} (Z ^T Z / n) \leq \mu_2 ) \cr
&\leq&  \pr(\lambda_{\min} (Z^{*T} Z^* / n) - \| \tilde Z ^T Z^* / n\|_2 - \|Z^{*T} \tilde Z / n\|_2 - \|\tilde Z^T \tilde Z /n\|_2 \leq \mu_2 ) \cr
&\leq&  \pr(\lambda_{\min} (Z^{*T} Z^* / n)  < 2 \mu_2) + \pr( \| \tilde Z ^T Z^* / n\|_2 \geq \mu_2 / 4) \cr
&&\quad + \pr(  \|Z^{*T} \tilde Z / n\|_2  \geq \mu_2/4)+ \pr(  \|\tilde Z^T \tilde Z /n\|_2  \geq \mu_2 /4 ) \cr
&=&  \pr(\|(Z^{*T} Z^* / n)^{-1}\|_2  > 1/( 2 \mu_2)) + o(n^{-\ell}) = o(n^{-\ell}).
\end{eqnarray*}

Then, taking $x_1 = \|\Sigma\|_2 \mu_1$ and $x_2 =  \mu_2 / \|\Sigma^{-1}\|_2$, we have
\begin{eqnarray*}
\pr(\|X^T X / n\|_2 \geq x_1) &\leq& \pr(\|\Sigma\|_2 \|Z^T Z / n\|_2 \geq x_1) \cr
&=& \pr(\|Z^T Z / n\|_2 \geq x_1/\|\Sigma\|_2) = o(n^{-\ell}), \cr
\pr ( \|(X^TX / n)^{-1} \|_2 \geq  x_2^{-1}) &\leq& \pr ( \|\Sigma^{-1}\|_2\|(Z^TZ / n)^{-1} \|_2 \geq  x_2^{-1}) \cr
&=& \pr ( \|(Z^TZ / n)^{-1} \|_2 \geq  (x_2\|\Sigma^{-1}\|_2)^{-1})= o(n^{-\ell}).
 \end{eqnarray*}
\endpf

\ni{\textbf{Proof of Theorem \ref{Theorem1}:}} We first consider the situation that $\lim _{n\to\infty} p_n = \infty$ for part (a).
From Lemma \ref{Lemma-2} that  $\tr\{(X^TX)^{-1}\} - \E \tr\{(X^TX)^{-1}\ID(K)\}= o_{\pr}(p_n/n)$. Under event $K$, the eigenvalues of $X^TX / n$ are bounded away from 0 and infinity. Hence, $\E \tr\{(X^TX)^{-1}\ID(K)\} = \Omega( p_n/n)$. Therefore, we have
 $\tr\{(X^TX)^{-1}\}=\E \tr\{(X^TX)^{-1}\ID(K)\}\{1 + o_{\pr}(1)\}$. From
\begin{eqnarray*}
&&  \|\hat\bfbeta\|_2^2 - \|\bfbeta_0\|_2^2 -  \tr\{(X^TX)^{-1}\} \hat\sigma_{\epsilon}^2  \cr
&=& \|\hat\bfbeta - \bfbeta_0\|_2^2   + 2 \bfbeta_0^T (\hat\bfbeta - \bfbeta_0)
 -  \tr\{(X^TX)^{-1}\} \hat\sigma_{\epsilon}^2   \cr
&=&   [\|\hat\bfbeta - \bfbeta_0\|_2^2  -\tr\{(X^TX)^{-1}\}  \sigma_{\epsilon}^2 ]
 -   \tr\{(X^TX)^{-1}\} (\hat\sigma_{\epsilon}^2 -  \sigma_{\epsilon}^2)  + 2 \bfbeta_0^T (\hat\bfbeta - \bfbeta_0)  \cr
 &\equiv& \I_1 - \E \tr\{(X^TX)^{-1}\ID(K)\} \{1 + o_{\pr}(1)\} \I_2 + 2 \I_3,
\end{eqnarray*}
we first demonstrate the asymptotic normality of $\zeta_n^{-1}(c_1\I_1 +c_2 \I_2 + c_3 \I_3) \ID(K)$ for any constants $c_1=\Omega( 1)$, $c_2 =\Omega( p_n / n)$ and $c_3=\Omega(1)$.

For notational simplicity, denote $M_1 = X (X^T X)^{-2} X^T$, $M_2 = \{\bfI_n - X(X^TX)^{-1}X^T\} / (n-p_n)$ and $\bfv^T = \bfbeta_0^T (X^T X)^{-1} X^T.$ Then,
 \begin{eqnarray*}
\I_1 &=&  \beps^T M_1 \beps - \tr\{(X^TX)^{-1}\}  \sigma_{\epsilon}^2 \cr
 &=& 2\sum_{1 \leq i<j \leq n} M_1(i,j)\epsilon_i \epsilon_j + \sum_{j=1}^n M_1(j,j)\epsilon_j ^2- \tr\{(X^TX)^{-1}\}  \sigma_{\epsilon}^2 \cr
 &=&   2\sum_{1 \leq i<j \leq n} M_1(i,j)\epsilon_i\epsilon_j + \sum_{j=1}^n M_1(j,j)(\epsilon_j ^2 - \sigma_{\epsilon}^2),\cr
\I_2 &=& \beps^T M_2 \beps -  \sigma_{\epsilon}^2
=   2\sum_{1 \leq i<j \leq n} M_2(i,j)\epsilon_i\epsilon_j + \sum_{j=1}^n M_2(j,j)(\epsilon_j ^2 - \sigma_{\epsilon}^2),\cr
\I_3 &=& \bfv^T\beps = \sum_{j=1}^n v_j \epsilon_j,
 \end{eqnarray*}
where $M_k(i,j)$ is the $(i,j)$th element of $M_k$ ($k=1,2$) and $\bfv = (v_1,\ldots,v_n)^T$. Hence,
  \begin{eqnarray*}
& & c_1\I_1 +c_2 \I_2 + c_3 \I_3  \cr
&=&   2\sum_{1\leq i<j \leq n}\{c_1 M_1(i,j) + c_2M_2(i,j)\} \epsilon_i \epsilon_j
 + \sum_{j=1}^n \{c_1 M_1(j,j) + c_2M_2(j,j)\}(\epsilon_j ^2 - \sigma_{\epsilon}^2)  +c_3 \sum_{j=1}^n v_j \epsilon_j \cr
&=&   \sum_{j = 1}^n \Big[\sum_{1 \leq i < j} 2\{c_1 M_1(i,j) + c_2M_2(i,j)\}  \epsilon_i\epsilon_j
 + \{c_1 M_1(j,j) + c_2M_2(j,j)\}(\epsilon_j ^2 - \sigma_{\epsilon}^2) + c_3 v_j \epsilon_j\Big]  \cr
&\equiv&   \sum_{j = 1}^n U_j.
 \end{eqnarray*}

Note that $U_j \ID(K) $, $j=1,2,\ldots,$ is a martingale difference, with $$\E(U_j  \ID(K)|X, \epsilon_1, \ldots, \epsilon_{j-1}) = 0$$
  and
    \begin{eqnarray*}
&& \sum_{j = 1}^n \E[\{ U_j  \ID(K)\} ^2 |X, \epsilon_1, \ldots, \epsilon_{j-1}] \cr
&=&   \ID(K)\sum_{j = 1}^n \E\Big(\Big[ \sum_{1 \leq i < j} 2\{c_1 M_1(i,j) + c_2M_2(i,j)\}  \epsilon_i\epsilon_j \cr
&&+ \{c_1 M_1(j,j) + c_2M_2(j,j)\}(\epsilon_j ^2 - \sigma_{\epsilon}^2) + c_3 v_j \epsilon_j\Big]^2 \Big|X, \epsilon_1, \ldots, \epsilon_{j-1}\Big) \cr
 &=& \ID(K)\sum_{j = 1}^n \E\Big( \Big[\sum_{1 \leq i < j} 2\{c_1 M_1(i,j) + c_2M_2(i,j)\} \epsilon_i \Big]^2 \epsilon_j ^2\cr
 &&+ \{c_1 M_1(j,j) + c_2M_2(j,j)\}^2(\epsilon_j ^2 - \sigma_{\epsilon}^2)^2 + c_3^2 v_j^2 \epsilon_j^2 \cr
 &&+ 2\sum_{1 \leq i < j} 2\{c_1 M_1(i,j) + c_2M_2(i,j)\} \epsilon_i \epsilon_j
  \{c_1 M_1(j,j) + c_2M_2(j,j)\}(\epsilon_j ^2 - \sigma_{\epsilon}^2) \cr
 &&+  2\sum_{1 \leq i < j} 2\{c_1 M_1(i,j) + c_2M_2(i,j)\} \epsilon_i \epsilon_j c_3 v_j \epsilon_j\cr
 && + 2\{c_1 M_1(j,j) + c_2M_2(j,j)\}(\epsilon_j ^2 - \sigma_{\epsilon}^2)c_3 v_j \epsilon_j   \Big|X, \epsilon_1, \ldots, \epsilon_{j-1}\Big)\cr
  &=& \ID(K)\sum_{j = 1}^n \Big( \Big[\sum_{1 \leq i < j} 2\{c_1 M_1(i,j) + c_2M_2(i,j)\} \epsilon_i\Big]^2 \sigma_{\epsilon} ^2\cr
 &&+ \{c_1 M_1(j,j) + c_2M_2(j,j)\}^2\var(\epsilon_j ^2) + c_3^2 v_j^2 \sigma_{\epsilon}^2 \cr
 && +2\sum_{1 \leq i < j} 2\{c_1 M_1(i,j) + c_2M_2(i,j)\} \{c_1 M_1(j,j) + c_2M_2(j,j)\}\E(\epsilon_j ^3 )\epsilon_i  \cr
 &&+  2\sum_{1 \leq i < j} 2\{c_1 M_1(i,j) + c_2M_2(i,j)\}  \epsilon_i c_3 v_j \sigma_{\epsilon}^2\cr
 && + 2\{c_1 M_1(j,j) + c_2M_2(j,j)\}\E(\epsilon_j ^3 )c_3 v_j \Big)\cr
 &\equiv& \ID(K) \sum_{j = 1}^n(\II_{1,j} +\II_{2,j} +\II_{3,j} +\II_{4,j} +\II_{5,j}+\II_{6,j} ).
 \end{eqnarray*}
 Denote $ t_n = \|\bbeta_0\|_2 /\sqrt n + \sqrt{p_n }/ n$.
 Lemmas S.6--S.11
 imply
 \begin{eqnarray} \label{A.7}
   \var\Big\{ \sum_{j = 1}^n \II_{k,j}  \ID(K)\Big\} = o(t_n^{4}), \quad \text{for $k = 1,\ldots, 6$.}
    \end{eqnarray}
Lemma S.5
indicates
     \begin{eqnarray} \label{A.8}
   \sum_{j=1}^n \E \{ U_j   \ID(K)\}^4 = o(t_n^{4}).
    \end{eqnarray}
 From  Lemmas S.12--S.14,
  we have
   \begin{eqnarray} \label{A.9}
\sum_{j = 1}^n \sum_{k = 1}^3 \E\{\II_{k,j} \ID(K)\} =O( t_n^{2}).
    \end{eqnarray}
    Lemmas S.9--S.11
    imply that
    \begin{eqnarray} \label{A.10}
 \sum_{j = 1}^n \sum_{k = 4}^6 \E\{\II_{k,j} \ID(K)\} =o( t_n^{2}).
    \end{eqnarray}
Checking conditions (2) and (4) with $\delta = 1$ in the theorem of \cite{Heyde_Brown_1970}, from \eqref{A.7}, \eqref{A.8}, \eqref{A.9}  and \eqref{A.10}, taking $c_1 = 1$, $c_2 = - \E\tr\{(X^TX)^{-1}\ID(K)\} $ and $c_3=2$,
 $$\zeta_n^{-1}  (c_1\I_1 +c_2 \I_2 + c_3 \I_3) \ID(K) \conD N(0,1),$$ where, from Lemmas S.12--S.14,
  \begin{eqnarray*}
  \zeta_n^2 &=&  4 \sigma_{\epsilon}^2  \bbeta_0^T \E\{(X^TX)^{-1} \ID(K)\} \bbeta_0 + 2 \sigma_{\epsilon}^4 \E \tr\{ (X^TX)^{-2}\ID(K)\}  + \frac{ 2  \sigma_{\epsilon}^4}{n-p_n} [\E\tr\{(X^TX)^{-1}\ID(K)\}]^2\cr
   &=& \Omega(\sigma_{\epsilon}^2\| \bbeta_0\|_2^2/n + \sigma_{\epsilon}^4p_n/n^2 + \sigma_{\epsilon}^4p_n^2/n^3  )
    = \Omega(\sigma_{\epsilon}^2\| \bbeta_0\|_2^2/n + \sigma_{\epsilon}^4p_n/n^2  )=\Omega(t_n^2).
  \end{eqnarray*}
  For part (b),
if  $p_n^{1/2} / n = o(\SNR)$, then $ \zeta_n^2 = o(\|\bbeta_0\|_2^4)$. Then, $\hat{\|\bbeta\|_2^2} - \|\bfbeta_0\|_2^2 = O_{\pr}(\zeta_n ) = o_{\pr}( \|\bbeta_0\|_2^2 )$.

Last, we consider the case for fixed $p_n$. Since $\bfbeta_0$ and $\sigma_{\epsilon}^2$ don't change with $n$ and $\bfbeta_0\neq \bfzero_{p_n}$, we have  $n^{-1} = o(\SNR)$. Note
\begin{eqnarray*}
&&  \|\hat\bfbeta\|_2^2 - \|\bfbeta_0\|_2^2 -  \tr\{(X^TX)^{-1}\} \hat\sigma_{\epsilon}^2  \cr
&=& \|\hat\bfbeta - \bfbeta_0\|_2^2   + 2 \bfbeta_0^T (\hat\bfbeta - \bfbeta_0)
 -  \tr\{(X^TX)^{-1}\} \hat\sigma_{\epsilon}^2\cr
 &=& \beps^T X(X^TX)^{-2} X^T \beps + 2 \bfbeta_0^T (\hat\bfbeta - \bfbeta_0)
 -  \tr\{(X^TX)^{-1}\} \hat\sigma_{\epsilon}^2.
\end{eqnarray*}
Following the proof for $\I_3$, we have $2 \bfbeta_0^T (\hat\bfbeta - \bfbeta_0)  / [4 \sigma_{\epsilon}^2  \bbeta_0^T \E\{(X^TX)^{-1} \ID(K)\} \bbeta_0]^{1/2} \conD N(0, 1)$ and hence $2 \bfbeta_0^T (\hat\bfbeta - \bfbeta_0) = \Omega_{\pr}(\sigma_{\epsilon} \|\bfbeta_0\|_2 /\sqrt n )$. From $\E \{\beps^T X(X^TX)^{-2} X^T \beps \ID(K)\} = \Omega( \sigma_{\epsilon}^2 p_n/  n ) $, $\tr\{(X^TX)^{-1}\} \hat\sigma_{\epsilon}^2 =\sigma_{\epsilon}^2 O_{\pr}(p_n/n )$ and $n^{-1} = o(\SNR)$, we have the following results $  \zeta_n^2 =\Omega [ 4 \sigma_{\epsilon}^2  \bbeta_0^T \E\{(X^TX)^{-1} \ID(K)\} \bbeta_0 ]$ and  $(\hat{\|\bbeta\|_2^2} - \|\bfbeta_0\|_2^2) / [4 \sigma_{\epsilon}^2  \bbeta_0^T \E\{(X^TX)^{-1} \ID(K)\} \bbeta_0]^{1/2} \conD N(0, 1)$. Hence, $(\hat{\|\bbeta\|_2^2} - \|\bfbeta_0\|_2^2) / \zeta_n \conD N(0, 1)$.
The proof for part(b) is similar to that for diverging $p_n$.

In the end, we will discuss the extension to fixed design.
From the proofs for \eqref{A.7}, \eqref{A.8}, \eqref{A.9} and \eqref{A.10} with $c_1 = 1$, $c_2 = - \E\tr\{(X^TX)^{-1}\ID(K)\} $ and $c_3=2$, there exists a sequence of positive real numbers $\{\omega_n\}_{n\geq 1}$ with $\omega_n = o(1)$ and constants $0<C_1<C_2<\infty$, such that $\pr(X \in \mathcal{X}_n) \to 1$ where  $\mathcal{X}_n \subseteq \mathbb{R}^{n\times p_n}$  is a collection of all  $x \in \mathbb{R}^{n\times p_n}$ satisfying
\begin{equation} \label{A.11}
\begin{aligned}
 \sum_{k=1}^6  \var\Big\{ \sum_{j = 1}^n \II_{k,j}  \ID(K)\Big | X=x\Big\} &\leq \omega_n t_n^{4}, \cr
   \sum_{j=1}^n \E [\{ U_j   \ID(K)\}^4 |X=x] &\leq \omega_n t_n^{4}, \cr
   \sum_{j = 1}^n \sum_{k = 1}^3 \E\{\II_{k,j} \ID(K)|X=x\} &\in[C_1 t_n^{2}, C_2 t_n^{2}], \cr
   \sum_{j = 1}^n \sum_{k = 4}^6 \E\{\II_{k,j} \ID(K)|X=x\} &\leq\omega_n  t_n^{2}\cr
   \Big|  4 \sigma_{\epsilon}^2  \bbeta_0^T (x^Tx)^{-1}  \bbeta_0 + 2 \sigma_{\epsilon}^4  \tr\{ (x^Tx)^{-2}\}  + \frac{ 2  \sigma_{\epsilon}^4}{n-p_n} [ \tr\{(x^Tx)^{-1}\}]^2-\zeta_n^2\Big|&\leq \omega_n  t_n^{2},
\end{aligned}
\end{equation}
while the last equation above is due to Lemma \ref{Lemma-2}.
Then, using the martingale difference CLT in \cite{Heyde_Brown_1970}, the asymptotic standard normality holds for $(\hat{\|\bbeta\|_2^2} - \|\bfbeta_0\|_2^2) / \zeta_n$ conditioning on $X=x$ for any $x \in \mathcal{X}_n$. The consistency result in part (b) can be derived using  similar arguments.
\endpf

\ni \textbf{Proof of Lemma \ref{Lemma-2}: } We provide the proof given event $H$. The results given event $K$  can be similarly derived.

From Efron-Stein inequality in \cite{Efron_Stein_1981}, if $W$ is a function of $n$ independent random variables and $W_{(i)}$ is any function of all those random variables except the $i$th, then
\begin{eqnarray} \label{A.12}
\var(W) \leq \sum_{i=1}^n \var(W - W_{(i)}) \leq \sum_{i=1}^n \E (W - W_{(i)})^2.
\end{eqnarray}
First, we use \eqref{A.12} with
$$W = n^{k}/p_n\tr\{(X^T X)^{-k}\} \ID(H), \quad W_{(i)} =  n^{k}/p_n  \tr\{(X_{(i)}^T X_{(i)})^{-k}\} \ID(H_{(i)})$$ where $H_{(i)}$ denotes the event that $\|(X_{(i)}^TX_{(i)} /n )^{-1}\|_2 \leq 1/x_2$.
Note
\begin{eqnarray*}
&& \sum_{i=1}^n \E (W - W_{(i)})^2 = n\E (W - W_{(i)})^2 \cr
&\lesssim& n \E[ n^{k}/p_n\tr\{(X^T X)^{-k}\} \{ \ID(H)- \ID(H_{(i)})\}]^2 \cr
&&+ n \E[ n^{k}/p_n\tr\{(X^T X)^{-k}\} \ID(H_{(i)})- n^{k}/p_n\tr\{(X_{(i)}^T X_{(i)})^{-k}\} \ID(H_{(i)})]^2\cr
&=& \I + \II.
\end{eqnarray*}
Since $X^TX \succeq X_{(i)}^TX_{(i)}$, we know $\|(X^TX)^{-1}\|_2 \leq \| (X_{(i)}^TX_{(i)})^{-1}\|_2$ and hence $H \supseteq H_{(i)}$. Then, $\ID(H)- \ID(H_{(i)}) = \ID(H \cap \bar H_{(i)}) = \ID(H) \ID( \bar H_{(i)})$. From Lemma \ref{Lemma-1}, $$\I \leq n ( n^{k}/p_n) ^2 (p_n n^{-k})^2   \pr (\bar H_{(i)}) = O(1/n).$$

Next, given $H_{(1)}$, we will show that $$n^{2k+1}/p_n^2 \E[ \tr\{(X^T X)^{-k}\} - \tr\{(X_{(1)}^T X_{(1)})^{-k}\} ]^2 = O(1/n).$$
From \eqref{A.1},
we have
\begin{eqnarray*}
(X^T X)^{-k}           =(X_{(1)}^T X_{(1)})^{-k} + \Delta,
\end{eqnarray*}
where $\Delta$ is a sum of $2^k - 1$ terms, each of which can be expressed as $A_1 \times A_2 \times\cdots \times A_k$ with $A_i = (X_{(1)}^T X_{(1)})^{-1}$ or $A_i =B$ $(i=1,\ldots, k)$ where $$B = - (X_{(1)}^T X_{(1)})^{-1} \bfX_1 \bfX_1^T (X_{(1)}^T X_{(1)})^{-1} /\{1 + \bfX_1^T (X_{(1)}^T X_{(1)})^{-1}\bfX_1\},$$ and at least one of $A_1,\ldots, A_k$ is $B$. It suffices to show that for each of the $2^k-1$ terms in $\Delta$, $ \E\{ \tr(A_1A_2 \cdots A_k)\}^2 =O(p_n^2 n^{-2k - 2})$. Without loss of generality, if $A_1 =B$, then from Lemmas \ref{Lemma-1} and S.1, given event $H_{(1)}$,
\begin{eqnarray*}
 \E\{ \tr(A_1A_2 \cdots A_k)\}^2 \leq \E \{ \bfX_1^T (X_{(1)}^T X_{(1)})^{-1}  A_2 \cdots A_k (X_{(1)}^T X_{(1)})^{-1} \bfX_1 \}^2 =O(p_n^2 n^{-2k - 2}).
\end{eqnarray*}

Next, we prove the second result of this lemma. Without loss of generality, assume $\|\bbeta_0\|_2 = 1$, and we will use \eqref{A.12} again with $W =  n^{k } \bbeta_0^T(X^T X)^{-k}\bbeta_0 \ID(H)$ and $W_{(i)} =  n^{k} \bbeta_0^T(X_{(i)}^T X_{(i)})^{-k}\bbeta_0 \ID(H_{(i)}) $ to show that, for each of the $2^k-1$ terms in $\Delta$, $$n^{2k+1} \E\{ \bbeta_0^T A_1A_2 \cdots A_k \bbeta_0 \ID(H_{(i)})\}^2 = O(1/n).$$

We only give the proof of a special case that $A_1= A_2 = B$, and $A_3= \cdots= A_k=(X_{(1)}^T X_{(1)})^{-1}$. From Lemma S.1,
\begin{eqnarray*}
&& \E\{ \bbeta_0^T A_1A_2 \cdots A_k \bbeta_0 \ID(H_{(1)})\}^2 \cr
&\leq& \E \{ \bbeta_0^T(X_{(1)}^T X_{(1)})^{-1} \bfX_1 \bfX_1^T (X_{(1)}^T X_{(1)})^{-2}  \bfX_1
 \bfX_1^T (X_{(1)}^T X_{(1)})^{-1}(X_{(1)}^T X_{(1)})^{-k+2} \bbeta_0 \ID(H_{(1)}) \}^2 \cr
&=& \E \{( \bbeta_0^T(X_{(1)}^T X_{(1)})^{-1} \bfX_1)^2 ( \bfX_1^T (X_{(1)}^T X_{(1)})^{-2} \bfX_1)^2
 ( \bfX_1^T (X_{(1)}^T X_{(1)})^{-k+1} \bbeta_0 )^2 \ID(H_{(1)}) \} \cr
&\leq&  [\E\{( \bbeta_0^T(X_{(1)}^T X_{(1)})^{-1} \bfX_1)^4\ID(H_{(1)})\}]^{1/2}
[ \E \{( \bfX_1^T (X_{(1)}^T X_{(1)})^{-2} \bfX_1)^8\ID(H_{(1)})\}]^{1/4} \cr
&&\cdot[ \E\{( \bfX_1^T (X_{(1)}^T X_{(1)})^{-k+1} \bbeta_0 )^8\ID(H_{(1)}) \}]^{1/4} \cr
&\lesssim&   \{\E\| \bbeta_0^T(X_{(1)}^T X_{(1)})^{-1} \ID(H_{(1)})\|_2^4\}^{1/2}[ \E \{ \bfX_1^T (X_{(1)}^T X_{(1)})^{-2} \bfX_1\ID(H_{(1)}) \}^8]^{1/4} \cr
&&\cdot \{ \E\| (X_{(1)}^T X_{(1)})^{-k+1} \bbeta_0\ID(H_{(1)}) \|_2^8 \}^{1/4} \cr
&\lesssim& n^{-2} n^{-2} n^{-2k + 2}= O(n^{-2k-2}).
\end{eqnarray*}
The proofs for the other terms are similar. We complete the proof.
\endpf

\ni\textbf{Proof of Theorem \ref{Theorem3}: } First, we consider $p_n \to \infty$.
Following the proof of Theorem \ref{Theorem1}, if $c_1 = 1$, $c_2 =  -\E\tr\{(X^TX)^{-1}\ID(K)\} $, $c_3 = 2$, we have $ \hat \zeta_n^2 - \zeta_n^2 = o_{\pr}(t_n^2)$ using the results in Lemmas S.12--S.14
and Proposition \ref{Proposition-1}, where $t_n$ is defined in the proof of Theorem \ref{Theorem1}.

For fixed $p_n$, we first show that $\hat\bbeta^T (X^TX)^{-1} \hat\bbeta/  B \conP 1$, where $B = \bbeta_0^T \E\{(X^TX)^{-1} \ID(K)\} \bbeta_0 = \Omega(\|\bbeta_0\|_2^2/n)$. Note
\begin{eqnarray*}
 \hat\bbeta^T (X^TX)^{-1} \hat\bbeta
  =\bbeta_0^T (X^TX)^{-1} \bbeta_0 + 2 \bbeta_0^T (X^TX)^{-2} X^T\beps + \beps^T X(X^TX)^{-3} X^T\beps.
\end{eqnarray*}
From Lemma \ref{Lemma-2},  $\bbeta_0^T (X^TX)^{-1} \bbeta_0 /  B \conP 1$. Since $2 \bbeta_0^T (X^TX)^{-2} X^T\beps  = O_{\pr}(  \sigma_{\epsilon} \|\bbeta_0\|_2 n^{-3/2}) = o_{\pr}(B)$ and $\beps^T X(X^TX)^{-3} X^T\beps = O_{\pr}(\sigma_{\epsilon} ^2 n^{-2} )= o_{\pr}(B)$, we claim that $\hat\bbeta^T (X^TX)^{-1} \hat\bbeta/  B \conP 1$.

Recall  $\hat \zeta_n^2 =  4 \hat \sigma_{\epsilon}^2  \hat\bbeta^T (X^TX)^{-1} \hat\bbeta - 2 \hat \sigma_{\epsilon}^4    \tr\{ (X^TX)^{-2}\} +  2 \hat \sigma_{\epsilon}^4   [\tr\{(X^TX)^{-1}\}]^2/(n-p_n).$
Then $2 \hat \sigma_{\epsilon}^4    \tr\{ (X^TX)^{-2}\} = O_{\pr}(\sigma_{\epsilon}^4  n^{-2} ) = o_{\pr}(\sigma_{\epsilon}^2 B)$
and  $2 \hat \sigma_{\epsilon}^4   [\tr\{(X^TX)^{-1}\}]^2/(n-p_n) =O_{\pr}(\sigma_{\epsilon}^4  n^{-3} ) = o_{\pr}(\sigma_{\epsilon}^2 B ) $.
Therefore,  from Proposition \ref{Proposition-1}, we have $\hat \zeta_n^2 / (4 \sigma_{\epsilon}^2  B) \conP 1$.

Following the proof of Theorem  \ref{Theorem1}, we have $ \zeta_n^2 / (4 \sigma_{\epsilon}^2  B) \conP 1$, which implies that $\hat \zeta_n^2 /  \zeta_n^2 \conP 1$.
We complete the proof.
\endpf

\bibliography{BIB-MD}

\newpage

\begin{center}
\bf \Large Supplementary Material
\end{center}

\setcounter{equation}{0}
\renewcommand{\theequation} {S.\arabic{section}.\arabic{equation}}
\setcounter{theorem}{0}
\renewcommand{\thetheorem} {S.\arabic{theorem}}
\setcounter{lemma}{0}
\renewcommand{\thelemma} {S.\arabic{lemma}}
\setcounter{proposition}{0}
\renewcommand{\theproposition} {S.\arabic{proposition}}
\setcounter{section}{0}
\renewcommand{\thesection} {S.\arabic{section}}

\textsl{}

 \vspace{.4cm}

\vspace{.4cm}

The supplementary material includes the following technical results:
\begin{itemize}
\item Section \ref{Section-S.1} illustrates that the conditions in \cite{Kelejian_Prucha_2001} are not satisfied in our paper;
\item Section \ref{Section-S.2} briefly explains the possibility of extending our results to centralized data;
\item Section \ref{Section-S.3} conducts inference for $\rho_0$   with large $\SNR$ using the conventional tests;
\item Section \ref{Section-S.4} conducts two-sample inferences;
\item Section \ref{Section-S.5} includes the proofs of Theorems 2, 4 and 5,   Corollary 1, and Propositions 1 and 2 in the paper;
\item Section \ref{Section-S.6} presents the technical lemmas and their proofs which are needed in the proofs of the main results.
\end{itemize}
In the following, we only consider the case that $\sigma_{\epsilon}^2 = \Omega(1)$  in the proofs. For random variables $\{V_n\}_{n\geq 1}$, denote $V_n = o_{L^2}(1)$ if $\E(V_n^2) = o(1)$.  Denote $a\vee b = \max(a,b)$ and $a\wedge b = \min(a,b)$ for $a,b \in \mathbb{R}$.

\section{Violation of conditions in \cite{Kelejian_Prucha_2001} } \label{Section-S.1}
We use a simple example $(n-p_n)\hat \sigma_{\epsilon}^2 = \beps^T \{\bfI_n - X (X^TX)^{-1}X^T\} \beps$ to clarify that the conditions in \cite{Kelejian_Prucha_2001} are not satisfied in our setup. To demonstrate the asymptotic normality of $(n-p_n)\hat \sigma_{\epsilon}^2$ using Theorem 1 in \cite{Kelejian_Prucha_2001}, it's required that $\var\{(n-p_n)\hat \sigma_{\epsilon}^2\}  = \Omega(n)$ and $\|\bfI_n - X (X^TX)^{-1}X^T\|_1\leq C<\infty$ almost surely.
Proposition 1 implies that $\var\{(n-p_n)\hat \sigma_{\epsilon}^2\} = \Omega(n)$. But, we will provide an example that $\|\bfI_n - X (X^TX)^{-1}X^T\|_1$ diverges to infinity in probability when $p_n \to \infty$.

Assuming $\bfX_i\sim N(\bfzero_{p_n}, \bfI_{p_n})$
and noting that $\|\bfI_n - X (X^TX)^{-1}X^T\|_1\geq \|X (X^TX)^{-1}X^T\|_1-1$,
then
\begin{eqnarray*}
&&\|X (X^TX)^{-1}X^T\|_1 \geq \sum_{i=1}^n|\bfX_1^T (X^TX)^{-1}\bfX_i|\cr
&\geq& \sum_{i=h_n}^n|\bfX_1^T (X^TX)^{-1}\bfX_i|
= \sum_{i=h_n}^n|\bfX_1^T (X_{(1)}^TX_{(1)})^{-1}\bfX_i|/\{1 + \bfX_1^T (X_{(1)}^TX_{(1)})^{-1} \bfX_1\} \cr
& =& \sum_{i=h_n}^n\frac{|\bfX_1^T (X_{(1,i)}^TX_{(1,i)})^{-1}\bfX_i|}{\{1 + \bfX_1^T (X_{(1)}^TX_{(1)})^{-1} \bfX_1\}
 \{1 + \bfX_i^T (X_{(1,i)}^TX_{(1,i)})^{-1} \bfX_i\}},
\end{eqnarray*}
  where $h_n = n- d_n $ with $d_n=\lfloor n/p_n^{1/4} \rfloor$. For $i=h_n,\ldots, n$, since $\bfX_i^T (X_{(1,i)}^TX_{(1,i)})^{-1} \bfX_i \leq \bfX_i^T (X_{(1,h_n,\ldots, n)}^TX_{(1,h_n,\ldots, n)})^{-1} \bfX_i$, we have, with probability tending to 1,
  \begin{eqnarray*}
  \max_{i=h_n,\ldots, n }\bfX_i^T (X_{(1,i)}^TX_{(1,i)})^{-1} \bfX_i &\leq& \|(\bfX_{h_n},\ldots, \bfX_n)^T (X_{(1,h_n,\ldots, n)}^TX_{(1,h_n,\ldots, n)})^{-1} (\bfX_{h_n},\ldots, \bfX_n)\|_2\cr
  &\leq&\|(\bfX_{h_n},\ldots, \bfX_n)\|_2^2 \|(X_{(1,h_n,\ldots, n)}^TX_{(1,h_n,\ldots, n)})^{-1}\|_2\cr
  &\leq&C \max(p_n, d_n) n^{-1} \leq C<\infty.
  \end{eqnarray*}
Hence, with probability tending to 1, $$\|X (X^TX)^{-1}X^T\|_1 \gtrsim\sum_{i=h_n}^n|\bfX_1^T (X_{(1,i)}^TX_{(1,i)})^{-1}\bfX_i|.$$

Given $i=h_n,\ldots,n$,
let $D_1 = \sum_{k=2}^{h_n-1} \bfX_k \bfX_k^T$ and $D_2 = \sum_{k=h_n, k\neq i}^{n} \bfX_k \bfX_k^T $, and then $X_{(1,i)}^TX_{(1,i)} = D_1 + D_2$ for $i=h_n,\ldots,n$.
 From the fact that
 \begin{eqnarray*}
 (X_{(1,i)}^TX_{(1,i)})^{-1} = D_1^{-1}   -  D_1^{-1}D_2
 (X_{(1,i)}^TX_{(1,i)})^{-1},
 \end{eqnarray*}
 we have $ |\bfX_1^T (X_{(1,i)}^TX_{(1,i)})^{-1}\bfX_i|\geq
 |\bfX_1^T D_1^{-1}\bfX_i| - |\bfX_1^T D_1^{-1}D_2
 (X_{(1,i)}^TX_{(1,i)})^{-1}\bfX_i|$. Conditioning on $\|(D_1/n)^{-1}\|_2\leq C<\infty$ and $\| (X_{(1,i)}^TX_{(1,i)}/n)^{-1}\|_2\leq C<\infty$, which happens with probability tending to 1, then
 \begin{eqnarray*}
 && \E |\bfX_1^T D_1^{-1}D_2  (X_{(1,i)}^TX_{(1,i)})^{-1}\bfX_i|^2
 =\E\{\bfX_1^T D_1^{-1}D_2  (X_{(1,i)}^TX_{(1,i)})^{-2} D_2 D_1^{-1}\bfX_1\}\cr
 &=&\E \tr\{D_1^{-1}D_2  (X_{(1,i)}^TX_{(1,i)})^{-2} D_2 D_1^{-1}\} \lesssim n^{-4}\E \tr(D_2^2)=n^{-4}\E (\|D_2\|_F^2)\lesssim  n^{-4}(p_n^2 d_n + d_n^2p_n),
 \end{eqnarray*}
 where $\|\cdot\|_F$ denotes the Frobenius norm.

 Hence,
 \begin{eqnarray*}
 \|X (X^TX)^{-1}X^T\|_1 &\gtrsim&\sum_{i=h_n}^n|\bfX_1^T D_1^{-1}\bfX_i| - \sum_{i=h_n}^n|\bfX_1^T D_1^{-1}D_2
 (X_{(1,i)}^TX_{(1,i)})^{-1}\bfX_i| \cr
 &=& \sum_{i=h_n}^n|\bfX_1^T D_1^{-1}\bfX_i|  - O_{\pr}(d_n n^{-2}(p_n d_n^{1/2} + d_n p_n^{1/2}) ).
 \end{eqnarray*}

 Conditioning on $\bfX_1,\ldots, \bfX_{h_n-1}$, then $\{\bfX_1^T D_1^{-1}\bfX_i\}_{i=h_n}^n$ are $\iid$ and  $\sum_{i=h_n}^n|\bfX_1^T D_1^{-1}\bfX_i| = \Omega_{\pr}( (n-h_n)  (\bfX_1^T D_1^{-2}\bfX_1)^{1/2}) = \Omega_{\pr}(d_n p_n^{1/2}/n) $.

 Noting that $d_n p_n^{1/2}/n \to \infty $ and $d_n n^{-2}(p_n d_n^{1/2} + d_n p_n^{1/2})  = o (d_n p_n^{1/2}/n)$, we claim that $\|X (X^TX)^{-1}X^T\|_1\conP \infty$.

\section{Brief explanation for extension to centralized data} \label{Section-S.2}

In the following, we give a brief explanation for extending our results to centralized data for model $Y_i =\alpha_0+ \bfX_i^{T} \bbeta_0 + \epsilon_i$ as discussed in Remark 2 in the paper.  Let $Y_i^* = Y_i - \bar Y$, $\bfX_i^* = \bfX_i - \bar \bfX$ and $\epsilon_i^* = \epsilon_i - \bar\epsilon$ be the centralized responses, predictors and errors, where $\bar \epsilon = n^{-1} \sum_{i=1}^n \epsilon_i$.
The intuition for the possibility of extension is that the proposed test statistics and estimators, i.e., $\|\hat\bbeta\|_2^2$, $\hat\sigma_{\epsilon}^2$, $\tr\{(X^TX)^{-k}\}$, $\hat\zeta_n^2$, $\hat\zeta_0^2$, $\hat\nu_4$, $\hat\eta$ and $\hat\rho$ using non-centralized data ($\bfY$ and $X$) are close to those calculated using the centralized data ($\bfY^*$ and $X^*$), denoted by $\|\hat\bbeta^*\|_2^2$, $\hat\sigma_{\epsilon}^{*2}$, $\tr\{(X^{*T}X^*)^{-k}\}$, $\hat\zeta_n^{*2}$, $\hat\zeta_0^{*2}$, $\hat\nu_4^*$, $\hat\eta^*$ and $\hat\rho^*$, where $k$ is a generic fixed positive integer.

From
$Y_i^* = \bfX_i^{*T} \bbeta_0 + \epsilon_i^*$ and $\hat\bbeta^* = (X^{*T}X^*)^{-1} X^{*T} \bfY^* $, we have that
 \begin{eqnarray*}
\hat\bbeta^* &=& \bbeta_0 + (X^{*T}X^*)^{-1} X^{*T} \beps^* = \bbeta_0 + (X^{*T}X^*)^{-1} (X^T\beps - n \bar \bfX \bar\epsilon)\cr
\hat\sigma_{\epsilon}^{*2} &=& \|\bfY^* - X^{*}\hat\bbeta^*\|_2^2/ (n-p_n) = \beps^{*T}\{\bfI_n - X^*(X^{*T}X^*)^{-1}X^{*T}\} \beps^* / (n-p_n).
 \end{eqnarray*}
 In the following, we briefly show that $\hat\sigma_{\epsilon}^{*2} - \hat\sigma_{\epsilon}^{2}  = O_{\pr}(\sigma_{\epsilon}^{2} / n)$. Without loss of generality, we can assume that $\alpha_0 = 0$ and $\E(\bfX_i) = \bfzero_{p_n}$.

Since $X^TX = X^{*T}X^* + n \bar\bfX \bar\bfX^T$, from the
Sherman-Morrison formula, we have
\begin{eqnarray*}
(X^TX)^{-1} = (X^{*T}X^*)^{-1} - \frac{n (X^{*T}X^*)^{-1}\bar\bfX \bar\bfX^T(X^{*T}X^*)^{-1}}{1 + n\bar\bfX^T(X^{*T}X^*)^{-1}\bar\bfX}.
\end{eqnarray*}
Following the proof of Lemma 1 and the result of \cite{Bai_Silverstein_2010}, the eigenvalues of  $X^{*T}X^* / n$ are bounded away from 0 and $\infty$ with probability tending to 1. In the following, the proof is conditioning on the event that the eigenvalues of  $X^{T}X /n$  and $X^{*T}X^* /n$ are uniformly bounded away from 0 and $\infty$. Then,
\begin{eqnarray*}
&& \hat\sigma_{\epsilon}^{*2} - \hat\sigma_{\epsilon}^{2} \cr
&=& \frac{\beps^{*T}\{\bfI_n - X^*(X^{*T}X^*)^{-1}X^{*T}\} \beps^* - \beps^{T}\{\bfI_n - X(X^{T}X)^{-1}X^{T}\} \beps}{n-p_n} \cr
&=& \frac{ \beps^{*T} \beps^{*} - \beps^{T} \beps -\beps^{*T}  X^*(X^{*T}X^*)^{-1}X^{*T}  \beps^* + \beps^{T} X(X^{T}X)^{-1}X^{T}  \beps}{n-p_n} \cr
&=& \frac{ (\beps^{*T} \beps^{*} - \beps^{T} \beps) -(X^T \beps- n \bar\bfX \bar\epsilon)^T(X^{*T}X^*)^{-1}(X^T \beps- n \bar\bfX \bar\epsilon) + \beps^{T} X(X^{T}X)^{-1}X^{T}  \beps}{n-p_n} \cr
&=& \frac{ -n\bar\epsilon^2
+2 \bar\epsilon n \bar\bfX ^T(X^{*T}X^*)^{-1}X^T \beps
-  n^2 \bar\epsilon^2     \bar\bfX ^T(X^{*T}X^*)^{-1}  \bar\bfX}{n-p_n} \cr
&& + \frac{ \beps^{T} X\{(X^{T}X)^{-1}-(X^{*T}X^*)^{-1}\}X^{T}  \beps}{n-p_n} \cr
&=&\I_1 + \I_2 + \I_3 + \I_4.
\end{eqnarray*}
We have that $\E(|\I_1|) = \sigma_{\epsilon}^2 / (n-p_n)$, $\E(|\I_3 |) =  \sigma_{\epsilon}^2 n/(n-p_n)   \E\{  \bar\bfX ^T(X^{*T}X^*)^{-1}  \bar\bfX\} \lesssim \sigma_{\epsilon}^2 /(n-p_n) \E(  \bar\bfX ^T  \bar\bfX)  \lesssim \sigma_{\epsilon}^2 p_n/\{n(n-p_n) \} $,
\begin{eqnarray*}
\E(|\I_2|)
& \lesssim& n/(n-p_n) [\E( \bar\epsilon^2) \E\{ \bar\bfX ^T(X^{*T}X^*)^{-1}X^T \beps\}^2]^{1/2} \cr
&=&\sigma_{\epsilon}^2 \sqrt n/(n-p_n) [ \E \{ \bar\bfX ^T(X^{*T}X^*)^{-1}X^T X (X^{*T}X^*)^{-1} \bar\bfX\}]^{1/2} \cr
&\lesssim&\sigma_{\epsilon}^2 /(n-p_n) [ \E \{ \bar\bfX ^T\bar\bfX\}]^{1/2} \cr
&\lesssim&\sigma_{\epsilon}^2 \sqrt{p_n} /\{\sqrt n(n-p_n)\},
\end{eqnarray*}
and
\begin{eqnarray*}
\E(|\I_4|)
& =& \E\Big[ \frac{n \beps^{T} X (X^{*T}X^*)^{-1}\bar\bfX \bar\bfX^T(X^{*T}X^*)^{-1} X^{T}  \beps}{\{1 + n\bar\bfX^T(X^{*T}X^*)^{-1}\bar\bfX\}(n-p_n)} \Big]\cr
&\leq&\E\{n \beps^{T} X (X^{*T}X^*)^{-1}\bar\bfX \bar\bfX^T(X^{*T}X^*)^{-1} X^{T}  \beps\}/(n-p_n) \cr
&=& n \sigma_{\epsilon}^2 \E\{   \bar\bfX^T(X^{*T}X^*)^{-1} X^{T} X (X^{*T}X^*)^{-1}\bar\bfX   \}/(n-p_n)\cr
&\lesssim& \sigma_{\epsilon}^2 /(n-p_n)\E(   \bar\bfX^T\bar\bfX   )\cr
&\lesssim&\sigma_{\epsilon}^2 p_n/\{n(n-p_n)\}.
\end{eqnarray*}
Therefore, combining the results for $\I_1, \I_2, \I_3, \I_4$, $\hat\sigma_{\epsilon}^{*2} - \hat\sigma_{\epsilon}^{2} = O_{\pr}(\sigma_{\epsilon}^2 / n)$.

Using similar arguments, i.e.,  the
Sherman-Morrison formula, leave-one-observation-out method and random matrix theory,  we can show that $\|\hat\bbeta\|_2^2$, $\tr\{(X^TX)^{-k}\}$, $\hat\zeta_n^2$, $\hat\zeta_0^2$, $\hat\nu_4$, $\hat\eta$ and $\hat\rho$ are sufficiently approximated by $\|\hat\bbeta^*\|_2^2$, $\tr\{(X^{*T}X^*)^{-k}\}$, $\hat\zeta_n^{*2}$, $\hat\zeta_0^{*2}$, $\hat\nu_4^*$, $\hat\eta^*$ and $\hat\rho^*$, which enables us to extend our results to the centralized data.

\section{Inference for $\rho_0$  with large $\SNR$} \label{Section-S.3}

We first conduct inference for $\eta_0$ using the conventional estimator $\hat\bbeta^T (X^TX/n)\hat \bbeta$ and study its consistency and asymptotic normality.
\begin{proposition}  \label{Proposition-S1}
Assume  $ \tau \in [0,1)$ and Condition  $\mathrm{A}$ for (1.1).
 If and only if $p_n/n = o(\SNR )$, then
 $$\frac{ \hat\bbeta^T (X^TX/n)\hat \bbeta }{ \eta_0 }\conP 1. $$
If and only if $\min(p_n^2/n, p_n/\sqrt n) = o(\SNR )$, we have
$$    \frac{\hat\bbeta^T (X^TX/n)\hat \bbeta - \eta_0} { \sigma_{\eta_0} } \conD N(0, 1),$$
 where $ \sigma^2_{\eta_0}=  [\var \{(\bbeta_0^T \bfX_1)^2\} + 4 \sigma_{\epsilon}^2 \bbeta_0^T \Sigma \bbeta_0]/n.$
\end{proposition}
\pf
Note that
\begin{eqnarray*}
  \hat\bbeta^T (X^TX/n)\hat \bbeta
&=& \{\bbeta_0 + (X^TX)^{-1} X^T \beps\}^T (X^TX/n)  \{\bbeta_0 + (X^TX)^{-1} X^T \beps\} \cr
&=& \bbeta_0^T (X^TX/n) \bbeta_0 + 2 \bbeta_0^T (X^TX/n) (X^TX)^{-1} X^T \beps \cr
&&\quad + \beps^T X (X^TX)^{-1} (X^TX/n)(X^TX)^{-1} X^T \beps  \cr
&=& \bbeta_0^T (X^TX/n) \bbeta_0 + 2 \bbeta_0^T  X^T \beps /n + \beps^T X (X^TX)^{-1}  X^T \beps / n.
 \end{eqnarray*}
From CLT, $\sqrt n \sigma_1^{-1} \{\bbeta_0^T (X^TX/n) \bbeta_0 - \bbeta_0^T\Sigma \bbeta_0\} \conD N(0, 1)$, where $ \sigma_1^2 = \var \{(\bbeta_0^T \bfX_1)^2\}= \Omega( \|\bbeta_0\|_2^4)$. Hence, $\bbeta_0^T (X^TX/n) \bbeta_0 - \bbeta_0^T\Sigma \bbeta_0 = O_{\pr}( \|\bbeta_0\|_2^2/\sqrt n)$.

Also, $\sigma_{2}^{-1} \bbeta_0^T  X^T \beps  /\sqrt n \conD N(0, 1)$ with
 $\sigma_{2}^2 =   \sigma_{\epsilon}^2 \bbeta_0^T \Sigma \bbeta_0$ which indicates that $\bbeta_0^T  X^T \beps  /n = O_{\pr} ( \sigma_{\epsilon}\| \bbeta_0\|_2 /\sqrt n)$.

 It's easy to see $\E\{\beps^T X (X^TX)^{-1}  X^T \beps /   n\} = \sigma_{\epsilon}^2 n^{-1} p_n$.

 Since $\|\bbeta_0\|_2^2 = \Omega(\eta_0)$, we know $ \hat\bbeta^T (X^TX/n)\hat \bbeta - \eta_0 = o_{\pr} ( \eta_0 )$ if and only if $\sigma_{\epsilon}\| \bbeta_0\|_2 /\sqrt n= o( \|\bbeta_0\|_2^2 )$ and $\sigma_{\epsilon}^2 n^{-1} p_n = o(\|\bbeta_0\|_2^2)$ which is equivalent to $ p_n/n = o( \|\bbeta_0\|_2^2 / \sigma_{\epsilon}^2)$.

 From the proof of Theorem 5, it's not hard to show that
 $$\sqrt n \sigma_3^{-1} \{\bbeta_0^T (X^TX/n) \bbeta_0 + 2 \bbeta_0^T  X^T \beps /n -\eta_0 \} \conD N(0, 1),$$
 where $\sigma_3^2 = \sigma_1^2+4\sigma_2^2$, which implies that $\bbeta_0^T (X^TX/n) \bbeta_0 + 2 \bbeta_0^T  X^T \beps /n -\eta_0 = O_{\pr} (\|\bbeta_0\|_2^2/\sqrt n + \sigma_{\epsilon} \|\bbeta_0\|_2/\sqrt n )$

 If and only if $\min(p_n^2/n, p_n/\sqrt n) = o(\|\bfbeta_0\|_2^2/ \sigma_{\epsilon}^2 )$, we have $\beps^T X (X^TX)^{-1}  X^T \beps / n$ is dominated by $ \bbeta_0^T (X^TX/n) \bbeta_0 + 2 \bbeta_0^T  X^T \beps /n$ and $  \sigma^{-1}_{3} \sqrt{n}  \{\hat\bbeta^T (X^TX/n)\hat \bbeta - \eta_0 \} \conD N(0, 1)$. Following the proof of Theorem 4, we complete the proof.
\endpf

From Proposition \ref{Proposition-S1}, it's  straightforward to derive the following asymptotical distribution result for
$ \tilde \rho =   {\hat\bbeta^T (X^TX/n)\hat \bbeta}/\{ \hat\bbeta^T (X^TX/n)\hat \bbeta + \hat\sigma_{\epsilon}^2 \}. $
\begin{theorem}  \label{TheoremS.1}
Assume Condition $\mathrm{A}$ for (1.1), $p_n^2/n = o(1 )$ and $\SNR =\Omega(1)$. Then
$$  \sigma^{-1}_{\tilde\rho} (\tilde \rho - \rho_0 ) \conD N(0,1), $$
where $\sigma^{2}_{\tilde\rho } = n^{-1}(\eta_0 + \sigma_{\epsilon}^2)^{-4} [\sigma_{\epsilon}^4 \{ \E(\sum_{i=1}^n Y_i^4/n) - \nu_4  \} -2 \sigma_{\epsilon}^6 \eta_0 + \eta_0^2(\nu_4 - 2 \sigma_{\epsilon}^4 ) ]$.
 \end{theorem}
\pf
From the proofs of Theorem 1 and Proposition \ref{Proposition-S1} and Proposition 1, we can show that $(\hat\bbeta^T (X^TX/n)\hat \bbeta , \hat \sigma_{\epsilon}^2)$ are asymptotically independent and jointly normal. Using delta method, we complete the proof.
\endpf

Therefore, the conventional test statistic for $H_0$ in (3.3) is
$$\mathbb{U}_0 = \hat\sigma^{-1}_{\tilde\rho }(\tilde\rho - \rho_0^{\nullH0}), $$
where $\hat \sigma^{2}_{\tilde\rho } = n^{-1}\{\hat\bbeta^T (X^TX/n)\hat \bbeta + \hat\sigma_{\epsilon}^2\}^{-4} [\hat\sigma_{\epsilon}^4 ( \sum_{i=1}^n Y_i^4/n - \sum_{i=1}^n \hat\epsilon_i^4 /n ) -2 \hat \sigma_{\epsilon}^6 \hat\bbeta^T $ $ (X^TX/n)\hat \bbeta + \{\hat\bbeta^T (X^TX/n)\hat \bbeta\}^2(\sum_{i=1}^n \hat\epsilon_i^4/n - 2 \hat\sigma_{\epsilon}^4 ) ]$, which is asymptotically standard normal under strong signal and $H_0$ in (3.3).

\section{Two-sample statistical inferences} \label{Section-S.4}
In this section, we study statistical inference problems for two linear regression models. Besides (1.1), consider another linear model
\begin{eqnarray} \label{S.4.1}
\bfW = V \bgamma_0 + \bdelta,
\end{eqnarray}
where $\bfW \in \mathbb{R}^{n'}$ is the vector of responses, $\bgamma_0 =(\gamma_{0,1}, \ldots, \gamma_{0,p_n})^T\in \mathbb{R}^{p_n}$ is the unknown parameter, $V = (\bfV_1, \ldots, \bfV_{n'})^T$ is the random design matrix and $\bdelta\in \mathbb{R}^{n'}$ is the error. Denote by $$\hat\bgamma = (V^TV)^{-1} V^T \bfW$$ the $\OLS$ estimator of $\bgamma_0$. Note that the dimension of $\bbeta_0$ and $\bgamma_0$ are the same, but the sample sizes in the two models could be different.

We assume model \eqref{S.4.1} fulfills the same conditions as (1.1), in the sense that Condition $\mathrm{A}$ is still satisfied if we replace $\bbeta_0$, $X$, $\beps$, $\sigma_{\epsilon}^2$, $\Sigma$ and $\bfZ_i$ therein by $\bgamma_0$, $V$, $\bdelta$, $\sigma_{\delta}^2$, $\Sigma'$ and $\bfZ'_i$, respectively, where $\bfV_i = \Sigma'^{1/2} \bfZ'_i$ and $\E(\delta_i^2)=\sigma_{\delta}^2 $.  Denote $\tau' = \lim_{n'\to\infty} p_n/n'$ and $\SNR' = {\var(\bfV_i^T \bgamma_0)} / {\var(\delta_i)} = \bgamma_0^T \Sigma' \bgamma_0 / \sigma_{\delta}^2 = \Omega (\|\bgamma_0\|_2^2 / \sigma_{\delta}^2 )$. We focus on  the  event below
 \begin{eqnarray*}
&& L := K \cap K', \, \, \text{where}\, \,  K'= H'\cap J', \cr
&& H' =\{\|(V^T V /n' )^{-1}\|_2 < 1/x_2' \}\, \, \text{and} \, \,  J' = \{\|V^T V /n' \|_2 <x_1' \} \, \, \text{with} \cr
&& x_1' =  4(1 + \sqrt{\tau'})^2\|\Sigma'\|_2\, \, \, \,\text{and}\, \, \, \, x_2' = (1 - \sqrt{\tau'})^2/(4\|\Sigma'^{-1}\|_2 ).
\end{eqnarray*}

The first problem is concerned with two-sample testing, i.e.,
\begin{eqnarray} \label{S.4.2}
H_0: \|\bbeta_0 - \bgamma_0\|_2=0 \qquad \text{versus} \qquad H_1:  \|\bbeta_0 - \bgamma_0\|_2 \neq 0.
\end{eqnarray}
Similar to the one-sample case, $\|\hat\bbeta - \hat\bgamma\|_2 ^2$ is a biased estimator with $$\E(\|\hat\bbeta - \hat\bgamma\|_2 ^2) - \|\bbeta_0 - \bgamma_0\|_2^2 = \E \tr\{(X^TX)^{-1}\} \sigma_{\epsilon}^2
+ \E \tr\{(V^TV)^{-1}\} \sigma_{\delta}^2 > 0. $$
 Consequently, we consider the following bias-corrected estimator $$\Widehat{\|\bbeta - \bgamma\|_2^2}= \|\hat\bbeta - \hat\bgamma\|_2 ^2 -  \tr\{(X^TX)^{-1}\} \hat\sigma_{\epsilon}^2
-  \tr\{(V^TV)^{-1}\} \hat\sigma_{\delta}^2,$$
and derive its  limiting distribution as follows. Here, $\hat\sigma_{\delta}^2$ is the estimator of $\sigma_{\delta}^2$, constructed in the same way as $\hat\sigma_{\epsilon}^2$.

\begin{theorem} \label{TheoremS.2}
Let $\{\bfX_i,\epsilon_i\}_{i=1}^n $ and $\{\bfV_i,\delta_i\}_{i=1}^{n'} $ be  independent. Assume Condition  $\mathrm{A}$  for models (1.1) and \eqref{S.4.1}, $\lim_{n\to\infty}p_n = \infty$, $ \tau \in [0,1)$, $ \tau' \in [0,1)$. Then,
 $$   \frac{ \Widehat{\|\bbeta - \bgamma\|_2^2} -\|\bbeta_0 - \bgamma_0\|_2^2 } { \sigma_{ \Widehat{\|\bbeta - \bgamma\|_2^2}} }\conD N(0, 1),$$
where
\begin{eqnarray*}
 \sigma^{2}_{ \Widehat{\|\bbeta - \bgamma\|_2^2}}  &=& 2  \sigma_{\epsilon}^4    \big( \E\tr\{ (X^TX)^{-2} \ID(L) \} +   [\E \tr\{(X^TX)^{-1}\ID(L)\}]^2 /(n-p_n)  \big) \cr
&& +2  \sigma_{\delta}^4    \big( \E \tr\{ (V^TV)^{-2}\ID(L)\} +   [\E \tr\{(V^TV)^{-1}\ID(L)\}]^2/(n' -p_n) \big) \cr
&&+ 4   \sigma_{\epsilon}^2 \sigma_{\delta}^2 \E \tr\{(X^TX)^{-1}(V^TV)^{-1}\ID(L)\}  \cr
&& + 4  \sigma_{\epsilon}^2 (\bbeta_0 - \bgamma_0)^T
\E\{(X^TX)^{-1} \ID(L)\}(\bbeta_0 - \bgamma_0) \cr
&& + 4  \sigma_{\delta}^2 (\bbeta_0 - \bgamma_0)^T
\E\{ (V^TV)^{-1}\ID(L)\}(\bbeta_0 - \bgamma_0).
\end{eqnarray*}
\end{theorem}
\pf Note
\begin{eqnarray*}
&&  \|\hat\bbeta - \hat\bgamma\|_2 ^2 -  \tr\{(X^TX)^{-1}\} \hat\sigma_{\epsilon}^2
-  \tr\{(V^TV)^{-1}\} \hat\sigma_{\delta}^2 \cr
&=&   \|\hat\bbeta - \bbeta_0\|_2 ^2 + \| \hat\bgamma- \bgamma_0\|_2 ^2 + \| \bbeta_0 -  \bgamma_0\|_2^2 - 2 (\hat\bbeta - \bbeta_0)^T ( \hat\bgamma- \bgamma_0) \cr
&&+ 2 (\bbeta_0 - \bgamma_0)^T (\hat\bbeta - \bbeta_0)  - 2 (\bbeta_0 - \bgamma_0)^T (\hat\bgamma - \bgamma_0)\cr
&& -  \tr\{(X^TX)^{-1}\} \hat\sigma_{\epsilon}^2
-  \tr\{(V^TV)^{-1}\} \hat\sigma_{\delta}^2  \cr
&=&   \|\hat\bbeta - \bbeta_0\|_2 ^2 -  \tr\{(X^TX)^{-1}\} \hat\sigma_{\epsilon}^2
+ \| \hat\bgamma- \bgamma_0\|_2 ^2-  \tr\{(V^TV)^{-1}\} \hat\sigma_{\delta}^2   \cr
&& - 2 (\hat\bbeta - \bbeta_0)^T ( \hat\bgamma- \bgamma_0)  + 2 (\bbeta_0 - \bgamma_0)^T (\hat\bbeta - \bbeta_0) - 2 (\bbeta_0 - \bgamma_0)^T (\hat\bgamma - \bgamma_0)\cr
&& + \| \bbeta_0 -  \bgamma_0\|_2^2\cr
&=& \I_1 + \I_2 + \I_3 + \I_4 + \I_5 + \| \bbeta_0 -  \bgamma_0\|_2^2.
\end{eqnarray*}
Following the proof of Theorem \ref{TheoremS.4}, we can show that $  \I_1$, $  \I_2$, $  \I_3$, $  \I_4$ and $ \I_5$ are asymptotically independently normal with mean 0 and  asymptotic variance $\sigma^{2}_{ \Widehat{\|\bbeta - \bgamma\|_2^2}}$.
\endpf

We propose a plug-in estimator for $\sigma^{2}_{ \Widehat{\|\bbeta - \bgamma\|_2^2}}$:
 \begin{eqnarray*}
\hat\sigma^{2}_{ \Widehat{\|\bbeta - \bgamma\|_2^2}}
&=&
2 \hat \sigma_{\epsilon}^4    \big(- \tr\{ (X^TX)^{-2}\} +  [\tr\{(X^TX)^{-1}\}]^2  /(n-p_n) \big) \cr
&&+2  \hat\sigma_{\delta}^4    \big(- \tr\{ (V^TV)^{-2}\} +   [\tr\{(V^TV)^{-1}\}]^2 /(n'-p_n)  \big) \cr
&&-4   \hat \sigma_{\epsilon}^2 \hat\sigma_{\delta}^2 \E \tr\{(X^TX)^{-1}(V^TV)^{-1}\}  \cr
&&+ 4  \hat\sigma_{\epsilon}^2 (\hat\bbeta - \hat\bgamma)^T
(X^TX)^{-1}(\hat\bbeta - \hat\bgamma) \cr
&& + 4  \hat \sigma_{\delta}^2 (\hat\bbeta - \hat\bgamma)^T
(V^TV)^{-1}(\hat\bbeta - \hat\bgamma),
\end{eqnarray*}
and study its consistency as follows.
\begin{theorem} \label{TheoremS.3}
Under the conditions of Theorem \ref{TheoremS.2},
$$ {\hat\sigma^{2}_{ \Widehat{\|\bbeta - \bgamma\|_2^2}} } \Big / { \sigma^{2}_{ \Widehat{\|\bbeta - \bgamma\|_2^2}} } \conP 1.$$
\end{theorem}
\pf
The proof is similar to that of Theorem \ref{TheoremS.5}.
\endpf

A natural test statistic for \eqref{S.4.2} is $$ \mathbb{D}_n =   \hat \sigma^{-1}_{ \Widehat{\|\bbeta - \bgamma\|_2^2}} \Widehat{\|\bbeta - \bgamma\|_2^2},$$ whose null limiting distribution is standard normal. The power of $ \mathbb{D}_n$ under the contiguous alternative hypothesis
\begin{eqnarray*}
H_{1n}: \|\bbeta_0 - \bgamma_0\|_2^2 = \delta_n,
\end{eqnarray*}
is  $ \Phi(\Phi^{-1}(\alpha) + \hat \sigma^{-1}_{ \Widehat{\|\bbeta - \bgamma\|_2^2}}  \delta_n)$.
 Hence, the smallest separation rate in the contiguous alternative is $\delta_n^* = \Omega(\sigma_{ \Widehat{\|\bbeta - \bgamma\|_2^2}}) = \Omega(\sigma_{\epsilon}^2 \sqrt{ p_n}/n +  \sigma_{\delta}^2 \sqrt{ p_n}/{n'})$, by noting that $ \sigma^{2}_{ \Widehat{\|\bbeta - \bgamma\|_2^2}} = O( \sigma_{\epsilon}^4 p_n/n^2 +  \sigma_{\delta}^4 p_n/{n'}^2 + \sigma_{\epsilon}^2\|\bbeta_0 - \bgamma_0\|_2^2/n + \sigma_{\delta}^2\|\bbeta_0 - \bgamma_0\|_2^2/n' )$.

We next focus on a normalized co-heritability $\theta_0= {\bgamma_0 ^T \bbeta_0}/({\|\bgamma_0\|_2\|\bbeta_0\|_2})$; see \cite{Guo_etal_2016} for more introductions. Consider the hypothesis
\begin{eqnarray} \label{S.4.3}
H_0: \theta_0 = \theta_0^{\nullH0} \qquad \text{versus} \qquad H_1:  \theta_0 \neq \theta_0^{\nullH0}.
\end{eqnarray}
The conventional estimator of $\theta_0$ is $\tilde\theta= {\hat\bgamma^T \hat\bbeta }/({\|\hat\bbeta\|_2 \|\hat\bgamma\|_2})$ and the corresponding test statistic for \eqref{S.4.3} is
$$\mathbb{C}_0 = \hat \sigma^{-1}_{\tilde\theta} ( \tilde\theta  - \theta_0^{\nullH0}),$$
where $\hat \sigma^{2}_{\tilde\theta}=1 / (\|\hat \bbeta \|_2 \|\hat\bgamma\|_2 )^2  \hat\sigma_{\delta}^2(\hat\bbeta -\hat\bgamma\hat\bgamma^T \hat\bbeta / \|\hat\bgamma\|_2^2) ^T (V^TV)^{-1} (\hat\bbeta -\hat\bgamma \hat\bgamma^T \hat\bbeta / \|\hat\bgamma\|_2^2) +
1 / (\|\hat\bbeta\|_2 \|\hat\bgamma\|_2 )^2 \hat \sigma_{\epsilon}^2(\hat\bgamma -\hat\bbeta\hat\bgamma^T \hat\bbeta / \|\hat\bbeta\|_2^2) ^T  (X^TX)^{-1}  (\hat\bgamma - \hat\bbeta\hat\bgamma ^T\hat \bbeta  / \|\hat\bbeta \|_2^2)$. Under a large $\SNR$ and $\theta_0^{\nullH0}<1$, the null limiting distribution is standard normal; see Proposition \ref{Proposition-S2} below. As shown previously, $\|\hat\bbeta\|_2^2$ is no longer consistent for $\|\bbeta_0\|_2^2$. Rather, we need the bias-corrected estimator $\hat{\|\bbeta\|_2} := (\hat{\|\bbeta\|_2^2})^{1/2} $ and $\hat{\|\bgamma\|_2} :=(\hat{\|\bgamma\|_2^2})^{1/2} $. Hence, we propose an estimator for $\theta_0$ as
$$\hat \theta =\frac{\hat\bgamma^T \hat\bbeta }{ \hat{\|\bbeta\|_2}\hat{\|\bgamma\|_2}},$$ whose null limiting distribution is given below.

\begin{theorem} \label{TheoremS.4}
Assume the conditions in Theorem \ref{TheoremS.2},  $p_n^{1/2}/n = o(\SNR) $,  $p_n^{1/2}/n' = o(\SNR')$ and $\theta_0 \in[C_1,C_2]$ for some $-1<C_1\leq C_2<1$.
Then, under $H_0$ in \eqref{S.4.3}, $$ \frac{\hat\theta-\theta_0}{ \sigma_{\hat\theta}  } \conD N(0, 1),$$
where
\begin{eqnarray*}
 \sigma^{2}_{\hat \theta}  &=&   \sigma_{\epsilon}^2 \sigma_{\delta}^2  \E\tr\{(X^TX)^{-1}(V^TV)^{-1}\ID(L)\}  / (\|\bbeta_0\|_2 \|\bgamma_0\|_2 )^2  \cr
 &&+
 1/ (\|\bbeta_0\|_2 \|\bgamma_0\|_2 )^2  \sigma_{\delta}^2 (\bbeta_0 -\bgamma_0\bgamma_0^T \bbeta_0 / \|\bgamma_0\|_2^2)^T \cr
 && \qquad \cdot \E\{(V^TV)^{-1} \ID(L)\}(\bbeta_0 -\bgamma_0\bgamma_0^T \bbeta_0 / \|\bgamma_0\|_2^2) \cr
 &&+ 1/ (\|\bbeta_0\|_2 \|\bgamma_0\|_2 )^2  \sigma_{\epsilon}^2(\bgamma_0 -\bbeta_0\bgamma_0^T \bbeta_0 / \|\bbeta_0\|_2^2)^T \cr
 && \qquad\cdot \E\{(X^TX)^{-1}\ID(L)\}(\bgamma_0 -\bbeta_0\bgamma_0^T \bbeta_0 / \|\bbeta_0\|_2^2)\cr
 &&+(\bgamma_0^T \bbeta_0)^2  / (2 \|\bbeta_0\|_2 \|\bgamma_0\|_2^3)^2 2  \sigma_{\delta}^4     \big( \E\tr\{ (V^TV)^{-2}\ID(L)\}\cr
  &&\qquad+  1/(n'-p_n) [\E\tr\{(V^TV)^{-1}\ID(L)\}]^2
  \big)  \cr
&&+(\bgamma_0^T \bbeta_0 )^2 / (2 \|\bbeta_0\|_2^3 \|\bgamma_0\|_2) ^2 2  \sigma_{\epsilon}^4   \big(\E\tr\{ (X^TX)^{-2}\ID(L)\} \cr
&&\qquad+  1/(n-p_n) [\E\tr\{(X^TX)^{-1}\ID(L)\}]^2
  \big).
  \end{eqnarray*}
\end{theorem}
\pf We only consider diverging $p_n$ here. Since $p_n^{1/2}/n = o(\SNR) $ and  $p_n^{1/2}/n' = o(\SNR') $, we have
 $\hat{\|\bbeta\|_2^2} - \|\bfbeta_0\|_2^2= o_{\pr}( \|\bfbeta_0\|_2^2)$ and $\hat{\|\bgamma\|_2^2 }- \|\bgamma_0\|_2^2= o_{\pr}( \|\bgamma_0\|_2^2)$. Note
\begin{eqnarray*}
& &  \hat\theta - \theta_0 \cr
&=&  \frac{\hat\bgamma^T \hat\bbeta }{ [ \|\hat\bfbeta\|_2^2 -  \tr\{(X^TX)^{-1}\} \hat\sigma_{\epsilon}^2 ]^{1/2}
[\|\hat\bgamma\|_2^2 -   \tr\{(V^TV)^{-1}\} \hat\sigma_{\delta}^2  ] ^{1/2} }  - \frac{\bgamma_0^T \bbeta_0 }{ \|\bgamma_0\|_2 \|\bbeta_0\|_2}  \cr
&=&  \frac{  (\hat\bgamma^T \hat\bbeta - \bgamma_0^T \bbeta_0) }{ [ \|\hat\bfbeta\|_2^2 -   \tr\{(X^TX)^{-1}\} \hat\sigma_{\epsilon}^2  ]^{1/2}
[\|\hat\bgamma\|_2^2 -    \tr\{(V^TV)^{-1}\} \hat\sigma_{\delta}^2 ] ^{1/2} } \cr
&&+  \bgamma_0^T \bbeta_0 \Big(\frac 1{[ \|\hat\bfbeta\|_2^2 -    \tr\{(X^TX)^{-1}\} \hat\sigma_{\epsilon}^2  ]^{1/2}
[\|\hat\bgamma\|_2^2 -  \tr\{(V^TV)^{-1}\} \hat\sigma_{\delta}^2  ] ^{1/2}} - \frac 1 {\|\bgamma_0\|_2 \|\bbeta_0\|_2} \Big) \cr
 &=&  \frac{  (\hat\bgamma^T \hat\bbeta - \bgamma_0^T \bbeta_0) }{ [ \|\hat\bfbeta\|_2^2 -   \tr\{(X^TX)^{-1}\} \hat\sigma_{\epsilon}^2  ]^{1/2}
[\|\hat\bgamma\|_2^2 - \tr\{(V^TV)^{-1}\} \hat\sigma_{\delta}^2 ] ^{1/2} } \cr
&&+  \bgamma_0^T \bbeta_0 \Big(\frac 1{[ \|\hat\bfbeta\|_2^2 - \tr\{(X^TX)^{-1}\} \hat\sigma_{\epsilon}^2  ]^{1/2}
[\|\hat\bgamma\|_2^2 - \tr\{(V^TV)^{-1}\} \hat\sigma_{\delta}^2  ] ^{1/2}} \cr
&&-\frac 1{[ \|\hat\bfbeta\|_2^2 -   \tr\{(X^TX)^{-1}\} \hat\sigma_{\epsilon}^2  ]^{1/2}
\|\bgamma_0\|_2} \Big)\cr
&&+ \bgamma_0^T \bbeta_0 \Big(\frac 1{[ \|\hat\bfbeta\|_2^2 -    \tr\{(X^TX)^{-1}\} \hat\sigma_{\epsilon}^2  ]^{1/2} \|\bgamma_0\|_2} - \frac 1 {\|\bgamma_0\|_2 \|\bbeta_0\|_2} \Big)\cr
&=& \big(  \{ (\hat \bgamma - \bgamma_0)^T (\hat\bbeta - \bbeta_0) + (\hat \bgamma - \bgamma_0)^T\bbeta_0 + (\hat\bbeta - \bbeta_0)^T \bgamma_0 \} / (\|\bbeta_0\|_2 \|\bgamma_0\|_2 ) \cr
&& - \bgamma_0^T \bbeta_0  / (2 \|\bbeta_0\|_2 \|\bgamma_0\|_2^3)  [\|\hat\bgamma\|_2^2 - \tr\{(V^TV)^{-1}\} \hat\sigma_{\delta}^2   - \|\bgamma_0\|_2^2 ] \cr
&& - \bgamma_0^T \bbeta_0  / (2 \|\bbeta_0\|_2^3 \|\bgamma_0\|_2)  [\|\hat\bfbeta\|_2^2 - \tr\{(X^TX)^{-1}\} \hat\sigma_{\epsilon}^2 - \|\bbeta_0\|_2^2] \big)\{1+o_{\pr}(1)\}\cr
&=&\big(    (\hat \bgamma - \bgamma_0)^T (\hat\bbeta - \bbeta_0)  / (\|\bbeta_0\|_2 \|\bgamma_0\|_2 )  \cr
&&+1/ (\|\bbeta_0\|_2 \|\bgamma_0\|_2 )   (\hat \bgamma - \bgamma_0)^T(\bbeta_0 -\bgamma_0\bgamma_0^T \bbeta_0 / \|\bgamma_0\|_2^2)\cr
&&+1/ (\|\bbeta_0\|_2 \|\bgamma_0\|_2 )   (\hat \bbeta - \bbeta_0)^T(\bgamma_0 -\bbeta_0\bgamma_0^T \bbeta_0 / \|\bbeta_0\|_2^2)  \cr
&& - \bgamma_0^T \bbeta_0  / (2 \|\bbeta_0\|_2 \|\bgamma_0\|_2^3)  [\|\hat\bgamma- \bgamma_0\|_2^2 - \tr\{(V^TV)^{-1}\} \hat\sigma_{\delta}^2  ] \cr
&& - \bgamma_0^T \bbeta_0  / (2 \|\bbeta_0\|_2^3 \|\bgamma_0\|_2)  [\|\hat\bfbeta-\bbeta_0\|_2^2 - \tr\{(X^TX)^{-1}\} \hat\sigma_{\epsilon}^2 ]\big)\{1+o_{\pr}(1)\}.
\end{eqnarray*}

Let $\bfa = \bgamma_0 -\bbeta_0\bgamma_0^T \bbeta_0 / \|\bbeta_0\|_2^2$, $\bfb = \bbeta_0 -\bgamma_0\bgamma_0^T \bbeta_0 / \|\bgamma_0\|_2^2$ and
\begin{eqnarray*}
\I_0  &=&   (\hat\bgamma - \bgamma_0) ^T (\hat\bfbeta - \bfbeta_0); \cr
\I_1 &=& \|\hat\bfbeta - \bfbeta_0\|_2^2  - \tr\{(X^TX)^{-1}\}  \sigma_{\epsilon}^2; \cr
\I_2 &=&    \hat\sigma_{\epsilon}^2 -  \sigma_{\epsilon}^2;\cr
\I_3 &=& \bfa^T (\hat\bfbeta - \bfbeta_0); \cr
\I_1' &=&  \|\hat\bgamma - \bgamma_0\|_2^2  - \tr\{(V^TV)^{-1}\}  \sigma_{\delta}^2; \cr
\I_2' &=&    \hat\sigma_{\delta}^2 -  \sigma_{\delta}^2;\cr
\I_3' &=&  \bfb^T (\hat\bgamma - \bgamma_0).
\end{eqnarray*}
Then
\begin{eqnarray*}
  \hat\theta - \theta_0
 &=& \I_0/ (\|\bbeta_0\|_2 \|\bgamma_0\|_2 ) + \I_3' /  (\|\bbeta_0\|_2 \|\bgamma_0\|_2 ) + \I_3/ (\|\bbeta_0\|_2 \|\bgamma_0\|_2 ) \cr
 &&- \bgamma_0^T \bbeta_0  / (2 \|\bbeta_0\|_2 \|\bgamma_0\|_2^3) \I_1' - \bgamma_0^T \bbeta_0  / (2 \|\bbeta_0\|_2^3 \|\bgamma_0\|_2) \I_1 \cr
 &&+ \tr\{(V^TV)^{-1}\}\bgamma_0^T \bbeta_0  / (2 \|\bbeta_0\|_2 \|\bgamma_0\|_2^3) \I_2' \cr
&& + \tr\{(X^TX)^{-1}\}\bgamma_0^T \bbeta_0  / (2 \|\bbeta_0\|_2^3 \|\bgamma_0\|_2) \I_2\cr
&=& (\|\bbeta_0\|_2 \|\bgamma_0\|_2 )^{-1}[\I_0 + \I_3'  + \I_3 - \bgamma_0^T \bbeta_0  / (2  \|\bgamma_0\|_2^2) \I_1'  - \bgamma_0^T \bbeta_0  / (2 \|\bbeta_0\|_2^2) \I_1 \cr &&
+ \tr\{(V^TV)^{-1}\}\bgamma_0^T \bbeta_0  / (2  \|\bgamma_0\|_2^2) \I_2'
+ \tr\{(X^TX)^{-1}\}\bgamma_0^T \bbeta_0  / (2 \|\bbeta_0\|_2^2) \I_2].
\end{eqnarray*}

The rest of the proofs follow from those of
Theorems 1 and 5. We can show that $(\I_1, \I_2, \I_3)$ and $(\I_1', \I_2', \I_3')$ are asymptotically trivariate normal, respectively. It's also easy to see that $(\I_1, \I_2, \I_3)$ and $(\I_1', \I_2', \I_3')$ are mutually independent. Conditioning on $(V, X, \bdelta)$, $\sigma_*^{-1}\I_0$ is asymptotically standard normal with $\sigma_*^2=\sigma_{\epsilon}^2 \sigma_{\delta}^2  \E\tr\{(X^TX)^{-1}(V^TV)^{-1}\ID(L)\} $.
 Also, conditioning on $(V, X, \bdelta)$, we can show that $\I_0$ is asymptotically independent with $(\I_1, \I_2, \I_3)$. Hence, $\I_0$, $(\I_1, \I_2, \I_3)$ and $(\I_1', \I_2', \I_3')$ are asymptotically mutually independent. Therefore, $\hat\theta - \theta_0$ is asymptotically normal with mean 0 and variance calculated using the results of Theorem 1. We complete the proof.
\endpf

{The asymptotic normal distribution of $\hat\theta$ requires $p_n^{1/2}/n = o(\SNR) $ and  $p_n^{1/2}/n' = o(\SNR'),$ which guarantees the ratio consistency of $\hat{\|\bbeta\|_2^2}$ and $\hat{\|\bgamma\|_2^2}$ due to Theorem 1. Such an $\SNR$ condition is not required in any one-sample inference, but still much weaker than the strong signal condition.}

The plug-in estimator of $\sigma^{2}_{\hat \theta}$ is developed as follows
\begin{eqnarray*}
\hat \sigma^{2}_{\hat \theta}
  &=& -   \hat\sigma_{\epsilon}^2 \hat \sigma_{\delta}^2  \tr\{(X^TX)^{-1}(V^TV)^{-1}\}  / (\hat{\|\bbeta\|_2^2}\hat{ \|\bgamma \|_2^2} )  \cr
&&+ 1/ (\hat{\|\bbeta\|_2^2}\hat{ \|\bgamma \|_2^2} ) \hat \sigma_{\delta}^2 (\hat\bbeta -\hat\bgamma\hat\bgamma^T \hat\bbeta / \hat{\|\bgamma\|_2^2})^T (V^TV)^{-1} (\hat\bbeta -\hat\bgamma\hat\bgamma^T \hat\bbeta / \hat{\|\bgamma\|_2^2})\cr
&& + 1/ (\hat{\|\bbeta\|_2^2}\hat{ \|\bgamma \|_2^2} ) \hat \sigma_{\epsilon}^2(\hat\bgamma -\hat\bbeta\hat\bgamma^T\hat \bbeta / \hat{\|\bbeta\|_2^2})^T(X^TX)^{-1}(\hat\bgamma -\hat\bbeta\hat\bgamma^T\hat \bbeta / \hat{\|\bbeta\|_2^2}) \cr
&& +(\hat\bgamma^T \hat\bbeta)^2  / (2 \hat{\|\bbeta\|_2}\hat{ \|\bgamma\|_2}^3)^2  2 \hat \sigma_{\delta}^4    \big(- \tr\{ (V^TV)^{-2}\}
 + 1/(n'-p_n) [\tr\{(V^TV)^{-1}\}]^2
  \big) \cr
&&+(\hat\bgamma^T \hat\bbeta )^2 / (2 \hat{\|\bbeta\|_2}^3 \hat{\|\bgamma\|_2}) ^2 2  \hat\sigma_{\epsilon}^4    \big(- \tr\{ (X^TX)^{-2}\}
 +  1 /(n-p_n) [\tr\{(X^TX)^{-1}\}]^2
  \big).
  \end{eqnarray*}
  The consistency  of the proposed estimator is given below.
  \begin{theorem} \label{TheoremS.5}
Under the conditions of Theorem \ref{TheoremS.4},
$$\hat \sigma^{2}_{\hat \theta} /  \sigma^{2}_{\hat \theta} \conP 1.$$
\end{theorem}
\pf We only consider diverging $p_n$ here.
From Theorem 1, $\hat{\|\bbeta\|_2} / \|\bbeta_0\|_2 \conP 1$ and $\hat{\|\bgamma\|_2} / \|\bgamma_0\|_2 \conP 1$.
Let $$
d_n^2  = \frac{p_n}{nn'}   + (1-\theta_0^2) \Big(\frac{\|\bbeta_0\|_2^2 } {n'}+ \frac{\|\bgamma_0\|_2^2 }n\Big ) +  \theta_0^2\Big( \frac{p_n\|\bbeta_0\|_2^2}{ {n'}^2 \|\bgamma_0\|_2^2} +  \frac{p_n\|\bgamma_0\|_2^2 }{n^2\|\bbeta_0\|_2^2 }\Big).$$
Then, $\sigma^{2}_{\hat \theta}\|\bbeta_0\|_2^2\|\bgamma_0\|_2^2=\Omega(d_n^2)$. Hence, it suffices to show that $\hat\sigma^{2}_{\hat \theta}\hat{\|\bbeta\|_2^2} \hat{\|\bgamma\|_2^2}- \sigma^{2}_{\hat \theta}\|\bbeta_0\|_2^2\|\bgamma_0\|_2^2=o_{\pr}(d_n^2)$

  We first consider $\| \bgamma_0\|_2 /\sqrt n + \| \bbeta_0\|_2 /\sqrt {n'} + \sqrt{ p_n/(n n')}= o( \bgamma_0^T\bbeta_0)$.
  Note
   \begin{eqnarray*}
  &&\hat\bgamma^T \hat\bbeta -\bgamma_0^T\bbeta_0 \cr
  &=& \bgamma_0^T(X^TX)^{-1}X^T\beps + \bbeta_0^T (V^TV)^{-1}V^T\bdelta + \beps^T X(X^TX)^{-1}(V^TV)^{-1}V^T\bdelta\cr
   &=& O_{\pr} ( \| \bgamma_0\|_2 /\sqrt n + \| \bbeta_0\|_2 /\sqrt {n'} + \sqrt {p_n/(n n')}) =o_{\pr}(\bgamma_0^T\bbeta_0).
    \end{eqnarray*}
    Then,
  $\hat\bgamma^T \hat\bbeta /\bgamma_0^T\bbeta_0 - 1  = o_{\pr}(1)$, which implies that $\hat\theta / \theta_0 - 1 = o_{\pr}(1)$.

Denote
\begin{eqnarray*}
 &&\sigma^{2}_{\hat \theta} (\|\bbeta_0\|_2 \|\bgamma_0\|_2 )^2  \cr
 &=&   \sigma_{\epsilon}^2 \sigma_{\delta}^2  \E\tr\{(X^TX)^{-1}(V^TV)^{-1}\ID(L)\}   \cr
 &&+
   \sigma_{\delta}^2 (\bbeta_0 -\bgamma_0\bgamma_0^T \bbeta_0 / \|\bgamma_0\|_2^2)^T
   \E\{(V^TV)^{-1} \ID(L)\}(\bbeta_0 -\bgamma_0\bgamma_0^T \bbeta_0 / \|\bgamma_0\|_2^2) \cr
 &&+ \sigma_{\epsilon}^2(\bgamma_0 -\bbeta_0\bgamma_0^T \bbeta_0 / \|\bbeta_0\|_2^2)^T
  \E\{(X^TX)^{-1}\ID(L)\}(\bgamma_0 -\bbeta_0\bgamma_0^T \bbeta_0 / \|\bbeta_0\|_2^2)\cr
 &&+(\bgamma_0^T \bbeta_0)^2  / (2  \|\bgamma_0\|_2^2)^2 2  \sigma_{\delta}^4     \big( \E\tr\{ (V^TV)^{-2}\ID(L)\}
 +  1/(n'-p_n) [\E\tr\{(V^TV)^{-1}\ID(L)\}]^2
  \big)  \cr
&&+(\bgamma_0^T \bbeta_0 )^2 / (2 \|\bbeta_0\|_2^2) ^2 2  \sigma_{\epsilon}^4   \big(\E\tr\{ (X^TX)^{-2}\ID(L)\}
+  1/(n-p_n) [\E\tr\{(X^TX)^{-1}\ID(L)\}]^2
  \big)\cr
  &\equiv& A_1 + A_2 + A_3 + A_4 + A_5.
  \end{eqnarray*}
Then, using the results in Theorems 1 and 3 and Lemma 2, we can show that
\begin{eqnarray*}
&&\hat\sigma_{\epsilon}^2 \hat \sigma_{\delta}^2  \tr\{(X^TX)^{-1}(V^TV)^{-1}\} ,\cr
&&(\hat\bgamma^T \hat\bbeta)^2  / (2 \hat{ \|\bgamma\|_2^2})^2  2 \hat \sigma_{\delta}^4    \big( \tr\{ (V^TV)^{-2}\}
 + 1/(n'-p_n) [\tr\{(V^TV)^{-1}\}]^2
  \big),\cr
  &&(\hat\bgamma^T \hat\bbeta )^2 / (2 \hat{\|\bbeta\|_2^2} ) ^2 2  \hat\sigma_{\epsilon}^4    \big( \tr\{ (X^TX)^{-2}\}
 +  1 /(n-p_n) [\tr\{(X^TX)^{-1}\}]^2
  \big),
\end{eqnarray*}
are ratio consistent for $A_1$, $A_4$ and $A_5$, respectively. In the following, we only need to show that
\begin{eqnarray*}
\hat A_2&=&  \hat \sigma_{\delta}^2 (\hat\bbeta -\hat\bgamma\hat\bgamma^T \hat\bbeta / \hat{\|\bgamma\|_2^2})^T (V^TV)^{-1} (\hat\bbeta -\hat\bgamma\hat\bgamma^T \hat\bbeta / \hat{\|\bgamma\|_2^2}) \cr
&&-\hat\sigma_{\epsilon}^2 \hat \sigma_{\delta}^2  \tr\{(X^TX)^{-1}(V^TV)^{-1}\}
- (\hat\bgamma^T \hat\bbeta)^2  / ( \hat{ \|\bgamma\|_2^2})^2   \hat \sigma_{\delta}^4   \tr\{ (V^TV)^{-2}\}
\end{eqnarray*}
and
\begin{eqnarray*}
\hat A_3&=&  \hat \sigma_{\epsilon}^2(\hat\bgamma -\hat\bbeta\hat\bgamma^T\hat \bbeta / \hat{\|\bbeta\|_2^2})^T(X^TX)^{-1}(\hat\bgamma -\hat\bbeta\hat\bgamma^T\hat \bbeta / \hat{\|\bbeta\|_2^2})
 \cr
&&-\hat\sigma_{\epsilon}^2 \hat \sigma_{\delta}^2  \tr\{(X^TX)^{-1}(V^TV)^{-1}\}
-(\hat\bgamma^T \hat\bbeta )^2 / ( \hat{\|\bbeta\|_2^2} ) ^2  \hat\sigma_{\epsilon}^4    \tr\{ (X^TX)^{-2}\}
\end{eqnarray*}
are consistent for $A_2$ and $A_3$ respectively, i.e., $\hat A_2 - A_2 = o_{\pr}(d_n^2) = \hat A_3 - A_3$.

Since $ \hat\bgamma^T \hat\bbeta /  \hat{\|\bgamma\|_2^2} =   (\bgamma_0^T \bbeta_0 /  \|\bgamma_0\|_2^2)\{1+o_{\pr}(1)\}$, we will replace $\hat\bgamma^T \hat\bbeta /  \hat{\|\bgamma\|_2^2}$ by
$\bgamma_0^T \bbeta_0 /  \|\bgamma_0\|_2^2$ below for ease of presentation.
  Following the proof of Lemma \ref{Lemma-S12},
  \begin{eqnarray*}
  &&  (\hat\bbeta -\hat\bgamma \bgamma_0^T \bbeta_0 /  \|\bgamma_0\|_2^2)^T (V^TV)^{-1} (\hat\bbeta -\hat\bgamma \bgamma_0^T \bbeta_0 /  \|\bgamma_0\|_2^2) \cr
&&
 -    \sigma_{\epsilon}^2 \E \tr\{ (V^TV)^{-1}(X^TX)^{-1} \ID(L)\} \cr
 && -   \sigma_{\delta}^2 (  \bgamma_0^T  \bbeta_0 /  \|\bgamma_0\|_2^2 )^2 \E \tr\{(V^TV)^{-2} \ID(L) \} \cr
&& - ( \bbeta_0 - \bgamma_0 \bgamma_0^T \bbeta_0 /  \|\bgamma_0\|_2^2)^T \E\{ (V^TV)^{-1} \ID(L)\} ( \bbeta _0 - \bgamma _0\bgamma_0^T \bbeta_0 /  \|\bgamma_0\|_2^2) \cr
&=&   \beps^T X (X^TX)^{-1} (V^TV)^{-1} (X^TX)^{-1} X^T \beps
-    \sigma_{\epsilon}^2 \E\tr\{ (V^TV)^{-1}(X^TX)^{-1} \ID(L)\} \cr
&& + (\bgamma_0^T \bbeta_0 /  \|\bgamma_0\|_2^2)^2[ \bdelta ^T V (V^TV)^{-3} V^T \bdelta - \sigma_{\delta}^2 \E \tr\{(V^TV)^{-2}\ID(L) \} ] \cr
&&  + ( \bbeta_0 - \bgamma_0 \bgamma_0^T \bbeta_0 /  \|\bgamma_0\|_2^2)^T[(V^TV)^{-1}-\E\{ (V^TV)^{-1}\ID(L)\}]
 ( \bbeta _0 - \bgamma _0\bgamma_0^T \bbeta_0 /  \|\bgamma_0\|_2^2) \cr
&&+2  (\bbeta_0-\bgamma_0\bgamma_0^T \bbeta_0 /  \|\bgamma_0\|_2^2 )^T (V^TV)^{-1} (X^TX)^{-1} X^T\beps\cr
&& - 2 (\bbeta_0-\bgamma_0\bgamma_0^T \bbeta_0 /  \|\bgamma_0\|_2^2 )^T (V^TV)^{-2} V^T \bdelta \bgamma_0^T \bbeta_0 /  \|\bgamma_0\|_2^2 \cr
&& - 2  \beps^T X (X^TX)^{-1} (V^TV)^{-2} V^T \bdelta \bgamma_0^T \bbeta_0 /  \|\bgamma_0\|_2^2 \cr
&=& o_{\pr} (p_n / (nn')) + o_{\pr} (\theta_0^2 (\|\bbeta_0\|_2^2/ \|\bgamma_0\|_2^2 )(p_n / {n'}^2)) \cr
&& + o_{\pr}(\|\bbeta_0 - \bgamma_0 \bgamma_0^T \bbeta_0 /  \|\bgamma_0\|_2^2\|_2^2 / n') \cr
&& + O_{\pr}(\|\bbeta_0 - \bgamma_0 \bgamma_0^T \bbeta_0 /  \|\bgamma_0\|_2^2\|_2 / (n'\sqrt n )) \cr
&&+ O_{\pr}(\|\bbeta_0 - \bgamma_0 \bgamma_0^T \bbeta_0 /  \|\bgamma_0\|_2^2\|_2 (|\bgamma_0^T \bbeta_0| /  \|\bgamma_0\|_2^2)/ \sqrt {{n'}^3} ) \cr
&& +O_{\pr}( \sqrt{p_n / (n{n'}^3)}|\bgamma_0^T \bbeta_0| /  \|\bgamma_0\|_2^2)\cr
&=& o_{\pr} (p_n / (nn')) + o_{\pr} (\theta_0^2 (\|\bbeta_0\|_2^2/ \|\bgamma_0\|_2^2)( p_n / {n'}^2)) \cr
&& + o_{\pr}(\|\bbeta_0\|_2^2 (1 - \theta_0^2) / n' ) \cr
&& + O_{\pr}(\|\bbeta_0\|_2 \sqrt{1 - \theta_0^2}/ (n'\sqrt n ) ) \cr
&&+ O_{\pr}( \sqrt{1-\theta_0^2} |\theta_0| \|\bbeta_0\|_2^2 /  ({\|\bgamma_0\|_2} \sqrt {{n'}^3} )) \cr
&& +O_{\pr}( \sqrt{p_n / (n{n'}^3)} |\theta_0|\|\bbeta_0\|_2 /  {\|\bgamma_0\|_2}) = o_{\pr} (d_n^2).
  \end{eqnarray*}

 Therefore, from Lemma 2 and Proposition 1, $\hat A_2$ ($\hat A_3$) is consistent for $A_2$ ($A_3$), and hence,
 $$\hat \sigma^{2}_{\hat \theta}  -  \sigma^{2}_{\hat \theta}= o_{\pr} (d_n^2(\|\bbeta_0\|_2 \|\bgamma_0\|_2 )^{-2}) = o_{\pr}(\sigma^{2}_{\hat \theta}). $$

 If $\bgamma_0^T\bbeta_0 = O(\| \bgamma_0\|_2 /\sqrt n + \| \bbeta_0\|_2 /\sqrt {n'} + \sqrt{ p_n/(n n')}  )$, then $\theta_0 = o(1)$ and
 \begin{eqnarray*}
 &&(\hat\bgamma^T \hat\bbeta)^2  / (2 \hat{ \|\bgamma\|_2^2})^2  2 \hat \sigma_{\delta}^4    \big(- \tr\{ (V^TV)^{-2}\} +  1/(n'-p_n) [\tr\{(V^TV)^{-1}\}]^2   \big) \cr
 &=& O_{\pr} ( (p_n/(n n') + \|\bbeta_0\|_2^2/ n' + \|\bgamma_0\|_2 ^2/   n) \|\bgamma_0\|_2^{-4} p_n/{n'}^2) = o_{\pr} (  d_n^2).
  \end{eqnarray*}

 Also,
 \begin{eqnarray*}
 &&  \hat\bbeta ^T (V^TV)^{-1} \hat\bgamma\hat\bgamma^T \hat\bbeta / \hat{\|\bgamma\|_2^2} \cr
 &=&  \bbeta_0 ^T (V^TV)^{-1} \bgamma_0\hat\bgamma^T \hat\bbeta / \hat{\|\bgamma\|_2^2}
  + \bbeta_0 ^T (V^TV)^{-2} V^T \bdelta \hat\bgamma^T \hat\bbeta / \hat{\|\bgamma\|_2^2}  \cr
  && + \beps^T X (X^TX)^{-1} (V^TV)^{-1} \bgamma_0\hat\bgamma^T \hat\bbeta / \hat{\|\bgamma\|_2^2}\cr
&&  + \beps^T X (X^TX)^{-1} (V^TV)^{-2} V \bdelta \hat\bgamma^T \hat\bbeta / \hat{\|\bgamma\|_2^2}  \cr
&=& \{O_{\pr} (\|\bbeta_0\|_2 \|\bgamma_0\|_2 / n' ) + O_{\pr} ({n'}^{-3/2} \|\bbeta_0\|_2 ) + O_{\pr} (n^{-1/2} {n'}^{-1} \|\bgamma_0\|_2 )  + O_{\pr} ( \sqrt{p_n /(n{n'}^3 ) }   )\}\cr
&&
O_{\pr} (\{\| \bgamma_0\|_2 /\sqrt n + \| \bbeta_0\|_2 /\sqrt {n'} + \sqrt{ p_n/(n n')}  \} / \|\bgamma_0\|_2^2 ) \cr
&=&  o_{\pr}(d_n^2),
 \end{eqnarray*}
 and
 \begin{eqnarray*}
 &&  \hat\bgamma ^T (V^TV)^{-1} \hat\bgamma (\hat\bgamma^T \hat\bbeta / \hat{\|\bgamma\|_2^2})^2 \cr
 &=&   \{ \bgamma_0 ^T (V^TV)^{-1} \bgamma_0
     + 2\bdelta^T V  (V^TV)^{-2} \bgamma_0
  + \bdelta^T V  (V^TV)^{-3} V \bdelta \} \cr
  &&O_{\pr} ( \{\| \bgamma_0\|_2^2 /  n + \| \bbeta_0\|_2^2 /  {n'} +  { p_n/(n n')}  \} \|\bgamma_0\|_2^{-4} )\cr
&=& \{O_{\pr} ( \|\bgamma_0\|_2 ^2 /n' ) + O_{\pr} ( \|\bgamma_0\|_2 {n'}^{-3/2}) + O_{\pr} (p_n/{n'}^2 ) \}\cr
 && \cdot O_{\pr} ( \{\| \bgamma_0\|_2^2 /  n + \| \bbeta_0\|_2^2 /  {n'} +  { p_n/(n n')}  \} \|\bgamma_0\|_2^{-4} )  \cr
&=&  o_{\pr}(d_n^2).
 \end{eqnarray*}
Therefore,
\begin{eqnarray*}
\hat \sigma^{2}_{\hat \theta}\hat{\|\bbeta\|_2^2}\hat{ \|\bgamma \|_2^2}
 = -   \hat\sigma_{\epsilon}^2 \hat \sigma_{\delta}^2  \tr\{(X^TX)^{-1}(V^TV)^{-1}\} +  \hat \sigma_{\delta}^2  \hat\bbeta^T (V^TV)^{-1} \hat\bbeta
  +  \hat \sigma_{\epsilon}^2 \hat\bgamma^T(X^TX)^{-1} \hat\bgamma +o_{\pr}(d_n^2).
  \end{eqnarray*}
Following arguments similar to Lemma \ref{Lemma-S12}, we have
  \begin{eqnarray*}
\hat \sigma^{2}_{\hat \theta} -\sigma^{2}_{\hat \theta} = o_{\pr} (d_n^2(\|\bbeta_0\|_2 \|\bgamma_0\|_2 )^{-2}) = o_{\pr}(\sigma^{2}_{\hat \theta}).
  \end{eqnarray*}
 \endpf

The proposed test statistic for \eqref{S.4.3} is thus $$\mathbb{C}_n =    \hat \sigma^{-1}_{\hat\theta} (\hat\theta-\theta_0^{\nullH0}).$$
The power under the contiguous alternative hypothesis
 \begin{eqnarray*}
 H_{1n}: |\theta_0-\theta_0^{\nullH0}|\ge\delta_n,
\end{eqnarray*}
is given by $\Phi(\hat \sigma^{-1}_{\hat\theta} \delta_n  + \Phi^{-1}(\alpha / 2 ) ) +  \Phi(-\hat \sigma^{-1}_{\hat\theta}  \delta_n +  \Phi^{-1}(\alpha / 2 ))$.
The smallest separation rate is $\delta_n^*=\sigma_{\hat \theta}$ with
\begin{eqnarray*}
 \sigma^{2}_{\hat \theta}  = \Omega\Big(\frac{p_n}{ {n n'} \SNR  \SNR'} + \frac{1 - \theta_0^2}{n' \SNR'} + \frac{1 - \theta_0^2}{n \SNR} + \frac{\theta_0^2 p_n}{{n'}^2{\SNR'}^2} + \frac{\theta_0^2p_n }{n^2{\SNR}^2}\Big ).
  \end{eqnarray*}
  If both $\SNR$ and $\SNR'$ are $\Omega(1)$, $\tau \in (0,1)$ and $\tau' \in (0,1)$, then $\sigma^{2}_{\hat \theta} = \Omega( 1/n)$, which means the
  alternative can be detected with large probability if it deviates from the null value for at least $\Omega(n^{-1/2})$.

Next, we conduct inference for $\theta_0$ under large $\SNR$.
The following propositions demonstrate the asymptotic distribution of $ {\hat\bgamma^T \hat\bbeta }/({\|\hat\bbeta\|_2 \|\hat\bgamma\|_2})$ provided sufficient $\SNR$ and $\SNR'$ for $\theta_0 < 1$.
\begin{proposition}   \label{Proposition-S2}
Let $\{\bfX_i, \epsilon_i\}_{i=1}^n $ and $\{\bfV_i, \delta_i\}_{i=1}^{n'} $ be   independent. Assume Condition $\mathrm{A}$ for (1.1) and \eqref{S.4.1},  $ \tau \in [0,1)$ and $ \tau' \in [0,1)$.
 If and only if $p_n/n = o(\SNR )$ and $p_n/n' = o(\SNR')$, then for all $\theta_0 \in [0,1]$,
 $$ {\hat\bgamma^T \hat\bbeta }/({\|\hat\bbeta\|_2 \|\hat\bgamma\|_2})   - \theta_0 = o_{\pr} (1 ). $$
 If the following conditions hold
  \begin{itemize}
 \item[$ \mathrm{(i)}$] $p_n/n = o(\SNR (1-\theta_0^2))$ and $p_n/n' = o(\SNR' (1-\theta_0^2))$,
 \item[$\mathrm{(ii)}$] $  \theta_0^2 p_n^2/n' = o(\SNR' (1-\theta_0^2))$ or $  \theta_0^2 p_n^2 n/{n'}^2 = o({\SNR'}^2/\SNR (1-\theta_0^2))$,
     \item[$\mathrm{(iii)}$] $  \theta_0^2 p_n^2/n = o(\SNR (1-\theta_0^2))$ or $  \theta_0^2 p_n^2 n'/n^2 = o(\SNR^2/\SNR' (1-\theta_0^2))$,
\end{itemize}
then
$$   \sigma^{-1}_{\tilde\theta} \{ {\hat\bgamma^T \hat\bbeta }/({\|\hat\bbeta\|_2 \|\hat\bgamma\|_2})   - \theta_0\} \conD N(0,1), $$
where $ \sigma^{2}_{\tilde\theta}=  1 / (\|\bbeta_0\|_2 \|\bgamma_0\|_2 )^2  \sigma_{\delta}^2(\bbeta_0 -\bgamma_0\bgamma_0^T \bbeta_0 / \|\bgamma_0\|_2^2) ^T \E\{(V^TV)^{-1}\ID(L) \}(\bbeta_0 -\bgamma_0\bgamma_0^T \bbeta_0 / \|\bgamma_0\|_2^2) +
1 / (\|\bbeta_0\|_2 \|\bgamma_0\|_2 )^2  \sigma_{\epsilon}^2(\bgamma_0 -\bbeta_0\bgamma_0^T \bbeta_0 / \|\bbeta_0\|_2^2) ^T \E\{(X^TX)^{-1} \ID(L)\}(\bgamma_0 -\bbeta_0\bgamma_0^T \bbeta_0 / \|\bbeta_0\|_2^2) $.
\end{proposition}
\pf
If   $p_n/n = o( \|\bbeta_0\|_2^2 / \sigma_{\epsilon}^2)$, then $\|\hat\bbeta\|_2^2-\|\bbeta_0\|_2^2 =  o_{\pr}(\|\bbeta_0\|_2^2).$ Similarly, if  $p_n/n' = o( \|\bgamma_0\|_2^2 / \sigma_{\delta}^2)$, then $\|\hat\bgamma\|_2^2-\|\bgamma_0\|_2^2 =  o_{\pr}(\|\bgamma_0\|_2^2).$ Note
\begin{eqnarray*}
\hat\bgamma^T \hat\bbeta
&=& \bbeta_0^T \bgamma_0 + \bbeta_0^T (V^TV)^{-1}V^T \bdelta + \bgamma_0^T (X^TX)^{-1}X^T\beps \cr
&&\quad +\bdelta^T V (V^TV)^{-1} (X^TX)^{-1}X^T\beps \cr
&=& \bbeta_0^T \bgamma_0 +O_{\pr} ( \sigma_{\delta} \|\bbeta_0\|_2/\sqrt {n'} + \sigma_{\epsilon} \|\bgamma_0\|_2/\sqrt n +\sigma_{\epsilon} \sigma_{\delta} \sqrt{p_n/(n n')} ).
\end{eqnarray*}
Therefore,
\begin{eqnarray*}
&& {\hat\bgamma^T \hat\bbeta }/({\|\hat\bbeta\|_2 \|\hat\bgamma\|_2})
={\hat\bgamma^T \hat\bbeta }/({\| \bbeta_0\|_2 \| \bgamma_0\|_2}) \{1 + o_{\pr}(1)\} \cr
&=& \{\bbeta_0^T \bgamma_0  +O_{\pr} ( \sigma_{\delta} \|\bbeta_0\|_2/\sqrt {n'} + \sigma_{\epsilon} \|\bgamma_0\|_2/\sqrt n +\sigma_{\epsilon} \sigma_{\delta} \sqrt{p_n/(n n')} ) \}
 /({\| \bbeta_0\|_2 \| \bgamma_0\|_2}) \{1 + o_{\pr}(1)\}\cr
&=& \theta_0 + o_{\pr}(1).
\end{eqnarray*}
If $p_n/n = \Omega( \|\bbeta_0\|_2^2 / \sigma_{\epsilon}^2)$ and $p_n/n' = O( \|\bgamma_0\|_2^2 / \sigma_{\delta}^2)$, then from the proof of Theorem 4, $\|\hat\bbeta\|_2^2/\|\bbeta_0\|_2^2 \gtrsim1$ with probability tending to 1. It's easy to see that ${\hat\bgamma^T \hat\bbeta }/({\|\hat\bbeta\|_2 \|\hat\bgamma\|_2}) $ is not consistent for $\theta_0$.

If   $  \|\bbeta_0\|_2^2 / \sigma_{\epsilon}^2    = o( p_n/n)$, then $\|\hat\bbeta\|_2^2/\|\bbeta_0\|_2^2 \conP \infty$. It's not hard to verify the inconsistency of ${\hat\bgamma^T \hat\bbeta }/({\|\hat\bbeta\|_2 \|\hat\bgamma\|_2}) $.

Similar arguments imply that if $p_n/n \neq o( \|\bbeta_0\|_2^2 / \sigma_{\epsilon}^2)$ or $p_n/n' \neq o( \|\bgamma_0\|_2^2 / \sigma_{\delta}^2)$, then ${\hat\bgamma^T \hat\bbeta }/({\|\hat\bbeta\|_2 \|\hat\bgamma\|_2}) $ is not consistent for $\theta_0$.

Second,
\begin{eqnarray*}
& &   \hat\bgamma^T\hat\bbeta / (\|\hat\bbeta\|_2 \|\hat\bgamma\|_2) - \theta_0
=  \frac{\hat\bgamma^T \hat\bbeta }{ \|\hat\bbeta\|_2
\|\hat\bgamma\|_2 }  - \frac{\bgamma_0^T \bbeta_0 }{ \|\bgamma_0\|_2 \|\bbeta_0\|_2} \cr
&=&  \frac{  (\hat\bgamma^T \hat\bbeta - \bgamma_0^T \bbeta_0) }{ \|\hat\bbeta\|_2
\|\hat\bgamma\|_2 }  +  \bgamma_0^T \bbeta_0 \Big(\frac 1{\|\hat\bbeta\|_2
\|\hat\bgamma\|_2} - \frac 1 {\|\bgamma_0\|_2 \|\bbeta_0\|_2} \Big) \cr
 &=&  \frac{  (\hat\bgamma^T \hat\bbeta - \bgamma_0^T \bbeta_0) }{ \|\hat\bbeta\|_2
\|\hat\bgamma\|_2 }  + \bgamma_0^T \bbeta_0 \Big(\frac 1{\|\hat\bbeta\|_2
\|\hat\bgamma\|_2}  -\frac 1{\|\hat\bbeta\|_2
\|\bgamma_0\|_2} \Big) \cr
&& + \bgamma_0^T \bbeta_0 \Big(\frac 1{\|\hat\bbeta\|_2 \|\bgamma_0\|_2} - \frac 1 {\|\bgamma_0\|_2 \|\bbeta_0\|_2} \Big)\cr
&=& [  \{ (\hat \bgamma - \bgamma_0)^T (\hat\bbeta - \bbeta_0) + (\hat \bgamma - \bgamma_0)^T\bbeta_0 + (\hat\bbeta - \bbeta_0)^T \bgamma_0 \} / (\|\bbeta_0\|_2 \|\bgamma_0\|_2 ) \cr
&& - \bgamma_0^T \bbeta_0  / (2 \|\bbeta_0\|_2 \|\bgamma_0\|_2^3) (\|\hat\bgamma\|_2^2  - \|\bgamma_0\|_2^2 ) \cr
&& - \bgamma_0^T \bbeta_0  / (2 \|\bbeta_0\|_2^3 \|\bgamma_0\|_2)  (\|\hat\bfbeta\|_2^2  - \|\bbeta_0\|_2^2) ]\{1+o_{\pr}(1)\}\cr
&=&\{1+o_{\pr}(1)\}\{  (\hat \bgamma - \bgamma_0)^T (\hat\bbeta - \bbeta_0)  / (\|\bbeta_0\|_2 \|\bgamma_0\|_2 )\cr
&& +1/ (\|\bbeta_0\|_2 \|\bgamma_0\|_2 )   (\hat \bgamma - \bgamma_0)^T(\bbeta_0 -\bgamma_0\bgamma_0^T \bbeta_0 / \|\bgamma_0\|_2^2)\cr
&&+1/ (\|\bbeta_0\|_2 \|\bgamma_0\|_2 )   (\hat \bbeta - \bbeta_0)^T(\bgamma_0 -\bbeta_0\bgamma_0^T \bbeta_0 / \|\bbeta_0\|_2^2)  \cr
&& - \bgamma_0^T \bbeta_0  / (2 \|\bbeta_0\|_2 \|\bgamma_0\|_2^3)   \|\hat\bgamma- \bgamma_0\|_2^2
   - \bgamma_0^T \bbeta_0  / (2 \|\bbeta_0\|_2^3 \|\bgamma_0\|_2)   \|\hat\bfbeta-\bbeta_0\|_2^2 \} \cr
   &=& (\I_1 + \I_2 + \I_3 + \I_4 + \I_5)\{1+o_{\pr}(1)\}.
\end{eqnarray*}
We know
\begin{eqnarray*}
\I_1 &=&   O_{\pr} (\sigma_{\epsilon} \sigma_{\delta} \sqrt{p_n/(n n')} )/ (\|\bbeta_0\|_2 \|\bgamma_0\|_2 ), \cr
\I_2 &=& 1/ (\|\bbeta_0\|_2 \|\bgamma_0\|_2 ) \|\bbeta_0 -\bgamma_0\bgamma_0^T \bbeta_0 / \|\bgamma_0\|_2^2\|_2 O_{\pr}(\sigma_{\delta} /\sqrt{n'}) \cr
&=& (1-\theta_0^2)^{1/2} O_{\pr}(1/\sqrt{n'}\sigma_{\delta} /\|\bgamma_0\|_2),\cr
\I_3 &=& (1-\theta_0^2)^{1/2} O_{\pr}(1/\sqrt n \sigma_{\epsilon} /\|\bbeta_0\|_2), \cr
\I_4 &=& \theta_0/\|\bgamma_0\|_2^2 O_{\pr} (  \sigma_{\delta}^2 p_n/n'),\cr
\I_5 &=& \theta_0/\|\bbeta_0\|_2^2 O_{\pr} ( \sigma_{\epsilon}^2  p_n/ n).
\end{eqnarray*}

If  $p_n/n = o(\SNR (1-\theta_0^2))$ or $p_n/n' = o(\SNR' (1-\theta_0^2))$, then $\I_2+\I_3$ dominates  $\I_1$.
If  $  \theta_0^2 p_n^2/n' = o(\SNR' (1-\theta_0^2))$ or $  \theta_0^2 p_n^2 n/{n'}^2 = o({\SNR'}^2/\SNR (1-\theta_0^2))$, then $\I_2+\I_3$ dominate  $\I_4$.
If  $  \theta_0^2 p_n^2/n = o(\SNR (1-\theta_0^2))$ or $  \theta_0^2 p_n^2 n'/n^2 = o(\SNR^2/\SNR' (1-\theta_0^2))$, then $\I_2+\I_3$ dominate  $\I_5$.

From the central limit theorem for $\I_2$ and $\I_3$, we finish the proof.
\endpf

\section{Proofs of main theoretical results}  \label{Section-S.5}
This section includes the proofs of the remaining main results in the paper.

\ni{\textbf{Proof of Theorem 2:}}
The joint probability density function of $(Y_i, \bfX_i)$ is
$$f(y, \bfx| \bbeta, \sigma_{\epsilon}^2, \Sigma) =   \{({2\pi})^{p_n+1} \sigma_{\epsilon}^2 |\Sigma|\}^{-1/2}
\exp \{-(y - \bfx^T \bbeta)^2 / (2\sigma_{\epsilon}^2)- \bfx^T \Sigma^{-1} \bfx/2\}.$$
 From \cite{LeCam_1973}, for $\bbeta_1 \neq \bbeta_2$,
\begin{eqnarray} \label{S.5.1}
&& \inf_{T} \sup_{\bbeta \in \mathcal{G}_{\bbeta_0}(c) } \E_{(\bfY,X)|(\bbeta, \sigma_{\epsilon}^2, \Sigma)} (T - \|\bbeta\|_2^2)^2 \cr
&\gtrsim& (\|\bbeta_2\|_2^2 - \|\bbeta_1\|_2^2)^2 \cr
&& \int \cdots \int \Big\{ \prod_{i=1}^n f(y_i, \bfx_i| \bbeta_1, \sigma_{\epsilon}^2, \Sigma)\Big\}  \wedge \Big\{ \prod_{i=1}^n f(y_i, \bfx_i| \bbeta_2, \sigma_{\epsilon}^2, \Sigma)\Big\} d y_1\cdots d y_n d\bfx_1 \cdots d \bfx_n. \qquad \ \
\end{eqnarray}
Next, we will show that for two generic probability density functions $g(x)$ and $h(x)$,
$$\int g(x) \wedge h(x) dx \geq \frac 1 2 \exp \{-\KL(G,H)\},$$
where $\KL(G, H)$ denotes the Kullback-Leibler ($\KL$) distance between the probability measures $G$ and $H$ corresponding to densities $g(x)$ and $h(x)$, respectively, and is defined as
$$\KL(G, H) = \int \log \Big\{\frac{g(x)}{h(x)} \Big\} g(x)dx.$$
Note that $\int g(x) \wedge h(x) dx + \int g(x) \vee h(x) dx = 2$. Then,
\begin{eqnarray*}
&& 2 \int g(x) \wedge h(x) dx \geq \Big(2 -\int g(x) \wedge h(x) dx \Big)\int g(x) \wedge h(x) dx \cr
&=& \int g(x) \vee h(x) dx  \int g(x) \wedge h(x) dx
  \geq \Big(\int \sqrt{\{ g(x) \wedge h(x)\}\{ g(x) \vee h(x)\}} dx \Big)^2 \cr
  &=& \Big(\int \sqrt{  g(x)   h(x)} dx \Big)^2
  = \exp \Big\{ 2 \log \Big(\int  \sqrt{ g(x)   h(x)} dx \Big)  \Big\}\cr
  &=& \exp \Big\{ 2 \log \Big(\int g(x) \sqrt{    h(x)/g(x)} dx \Big)  \Big\}
  \geq \exp \Big\{ 2 \Big(\int g(x) \log \sqrt{    h(x)/g(x)} dx \Big)  \Big\} \cr
  &=& \exp \{- \KL(G,H)\},
\end{eqnarray*}
where we use Jensen's inequality in the last inequality. Therefore, from \eqref{S.5.1},
\begin{eqnarray*}
  \inf_{T} \sup_{\bbeta \in \mathcal{G}_{\bbeta_0}(c) } \E_{(\bfY,X)|(\bbeta, \sigma_{\epsilon}^2, \Sigma)} (T - \|\bbeta\|_2^2)^2
 &\gtrsim& (\|\bbeta_2\|_2^2 - \|\bbeta_1\|_2^2)^2 \exp\{- \KL(F_{n, \bbeta_1}, F_{n, \bbeta_2})\} \cr
 &=&(\|\bbeta_2\|_2^2 - \|\bbeta_1\|_2^2)^2 \exp\{- n \KL(F_{  \bbeta_1}, F_{  \bbeta_2})\} \equiv \I,
\end{eqnarray*}
where $F_{n, \bbeta_j}$ denotes the probability measure for density $\prod_{i=1}^n f(y_i, \bfx_i| \bbeta_j, \sigma_{\epsilon}^2, \Sigma)$ and $F_{  \bbeta_j}$ is that for density
$  f(y , \bfx | \bbeta_j, \sigma_{\epsilon}^2, \Sigma)$ for $j=1,2$.

The $\KL$ distance between $f(y, \bfx| \bbeta_1, \sigma_{\epsilon}^2, \Sigma)$ and $f(y, \bfx| \bbeta_2, \sigma_{\epsilon}^2, \Sigma)$ is
\begin{eqnarray*}
&& \KL(F_{  \bbeta_1}, F_{  \bbeta_2}) \cr
&=& \int\int \log \Big(\frac{f(y, \bfx| \bbeta_1, \sigma_{\epsilon}^2, \Sigma) }{ f(y, \bfx| \bbeta_2, \sigma_{\epsilon}^2, \Sigma)}\Big)f(y, \bfx| \bbeta_1, \sigma_{\epsilon}^2, \Sigma)  d y d \bfx\cr
&=& \E_{(Y,\bfX)|(\bbeta_1, \sigma_{\epsilon}^2, \Sigma)}   \{(Y - \bfX^T \bbeta_2)^2 / (2\sigma_{\epsilon}^2) -(Y - \bfX^T \bbeta_1)^2 / (2\sigma_{\epsilon}^2) \} \cr
&=& (2\sigma_{\epsilon}^2)^{-1}\E_{(Y,\bfX)|(\bbeta_1, \sigma_{\epsilon}^2, \Sigma)}   (-2 Y \bfX^T \bbeta_2 + \bbeta_2^T \bfX \bfX^T \bbeta_2 + 2 Y  \bfX^T \bbeta_1 -\bbeta_1^T \bfX \bfX^T \bbeta_1 ) \cr
&=& (2\sigma_{\epsilon}^2)^{-1}   (-2 \bbeta_1^T  \Sigma \bbeta_2 + \bbeta_2^T \Sigma \bbeta_2 + \bbeta_1^T \Sigma \bbeta_1 ) \cr
&=& (\bbeta_2- \bbeta_1)^T \Sigma (\bbeta_2- \bbeta_1)  / (2\sigma_{\epsilon}^2).
\end{eqnarray*}
We first consider that $p_n/n = O(\SNR) $, in which case $\zeta_n^2 = \Omega(\sigma_{\epsilon}^2 \|\bbeta_0\|_2^2/n )$.
If $\bbeta_1 = \bbeta_0$ and $\bbeta_2 =  \bbeta_0 \{ 1+ \sigma_{\epsilon} / (\sqrt n \|\bbeta_0\|_2) \}$, then $\KL(F_{  \bbeta_1}, F_{  \bbeta_2}) = \bbeta_0^T \Sigma \bbeta_0 / (2n \|\bbeta_0\|_2^2)$ and $(\|\bbeta_2\|_2^2 - \|\bbeta_1\|_2^2)^2 = (2 \sigma_{\epsilon} \|\bbeta_0\|_2 /\sqrt n + \sigma_{\epsilon}^2 / n)^2 $, which implies that $\I =\Omega(\sigma_{\epsilon}^2 \|\bbeta_0\|_2^2/n )= \Omega(\zeta_n^2 )$.

Next, consider the case that $\SNR = O(p_n/n) $, under which $\zeta_n^2 = \Omega(\sigma_{\epsilon}^4 p_n/n^2 )$. Let $\bbeta_1 = \sigma_{\epsilon} /\sqrt n \bfone_{p_n} $ and $\bbeta_2 = \sigma_{\epsilon}/\sqrt n (1 - 1/\sqrt{p_n} )  \bfone_{p_n } $. Then,
$\KL(F_{  \bbeta_1}, F_{  \bbeta_2}) = \Omega(n^{-1})$ and $(\|\bbeta_2\|_2^2 - \|\bbeta_1\|_2^2)^2 = \Omega( \sigma_{\epsilon}^4p_n / n^2 )$, which implies that $\I =\Omega(\sigma_{\epsilon}^4p_n / n^2 )= \Omega(\zeta_n^2 )$.

Therefore, we have shown that
$$\inf_{T} \sup_{\bbeta \in \mathcal{G}_{\bbeta_0}(c) } \E_{(\bfY,X)|(\bbeta, \sigma_{\epsilon}^2, \Sigma)} (T - \|\bbeta\|_2^2)^2 \gtrsim  \zeta_n^2.$$
From the proof of Theorem 1, for any $\bbeta \in \mathcal{G}_{\bbeta_0}(c)$,
$$    \E_{(\bfY,X)|(\bbeta, \sigma_{\epsilon}^2, \Sigma)} (\hat{\|\bbeta\|_2^2} - \|\bbeta\|_2^2)^2 \lesssim  \sigma_{\epsilon}^2 \|\bbeta\|_2^2/n + \sigma_{\epsilon}^4 p_n / n^2 = O(\zeta_n^2).$$
The proof is completed.
\endpf

\ni{\textbf{Proof of Theorem 4:}}
 Note,
\begin{eqnarray*}
\|\hat\bbeta\|_2^2 - \|\bfbeta_0\|_2^2 &=&\|\hat\bfbeta - \bfbeta_0\|_2^2   + 2 \bfbeta_0^T (\hat\bfbeta - \bfbeta_0)\cr
&=& \beps^T X(X^TX)^{-2} X^T \beps + 2 \bbeta_0^T (X^TX)^{-1} X^T \beps \equiv \I_1 + 2 \I_2. \end{eqnarray*}
For term $\I_1$,
\begin{eqnarray*}
 \E\{\I_1 \ID(K)\} &=& \sigma_{\epsilon}^2 \E\tr\{X(X^TX)^{-2}X^T\ID(K)\} \cr
 &=&  \sigma_{\epsilon}^2 \E\tr\{(X^TX)^{-1}\ID(K)\} = \Omega( \sigma_{\epsilon}^2 p_n/  n )\equiv   \Omega( \sigma_{1} ),
\end{eqnarray*}
and
\begin{eqnarray} \label{S.5.2}
&&\E\{\I_1 \ID(K)\}^2 \cr
&=& \nu_4 \sum_i \E\{\bfX_i^T (X^TX)^{-2} \bfX_i\ID(K)\}^2 + 2 \sigma_{\epsilon}^4 \sum_{i\neq j} \E\{\bfX_i^T(X^TX)^{-2} \bfX_j\ID(K)\}^2 \cr
&&+ \sigma_{\epsilon}^4 \sum_{i\neq j} \E\{\bfX_i^T(X^TX)^{-2} \bfX_i \bfX_j^T(X^TX)^{-2} \bfX_j\ID(K)\} \cr
&=& (\nu_4 - 3 \sigma_{\epsilon}^4) \sum_i \E\{\bfX_i^T (X^TX)^{-2} \bfX_i\ID(K)\}^2 + 2 \sigma_{\epsilon}^4 \sum_{i, j} \E\{\bfX_i^T(X^TX)^{-2} \bfX_j\ID(K)\}^2 \cr
&&+ \sigma_{\epsilon}^4 \sum_{i, j} \E\{\bfX_i^T(X^TX)^{-2} \bfX_i \bfX_j^T(X^TX)^{-2} \bfX_j\ID(K)\} \cr
&=& (\nu_4 - 3 \sigma_{\epsilon}^4) \sum_i \E\{\bfX_i^T(X^TX)^{-2} \bfX_i\ID(K)\}^2 + 2 \sigma_{\epsilon}^4 \E\tr\{ (X^TX)^{-2}\ID(K)\} \cr
&&+ \sigma_{\epsilon}^4\E [\tr\{ (X^TX)^{-1}\}\ID(K)]^2\cr
& =& O(\sigma_{\epsilon}^4 p_n^2/n^2 ) = O( \sigma_{1}^2 ).
\end{eqnarray}
Hence, $\E \{\sigma_{1}^{-1}  \I_1 \ID(K)\} ^2 \leq C$, which implies that $\sigma_{1}^{-1} \I_1  \ID(K)$ is uniformly integrable.

For term $\I_2$, following the proof of Theorem 1 and using the results in Lemma 2,
$ \zeta_0^{-1}   2 \I_2   \conD N(0, 1).$
Also,
\begin{eqnarray} \label{S.5.3}
\E\{\I_2\ID(K)\}^2 = \Omega (  \sigma_{\epsilon}^2 \|\bbeta_0\|_2^2 /  n) \equiv \Omega ( \sigma_{2}^2).
\end{eqnarray}
Hence, $\sigma_{2}^{-1} \I_2 \ID(K)$ is uniformly integrable.

First, we will study the consistency of $\|\hat\bbeta\|_2^2$.
\begin{itemize}
\item If $p_n/n = o(\|\bbeta_0\|_2^2/\sigma_{\epsilon}^2) $, then $\I_1 \ID(K) = O_{\pr} ( \sigma_{\epsilon}^2 p_n/  n )= o_{\pr}(\|\bbeta_0\|_2^2)$ and $\I_2 \ID(K) = O_{\pr}(  \sigma_{\epsilon} \|\bbeta_0\|_2 / \sqrt n)= o_{\pr}(\|\bbeta_0\|_2^2)$.

\item If $p_n/n = \Omega(\|\bbeta_0\|_2^2/\sigma_{\epsilon}^2) $, then $\E\{\I_1 \ID(K)\} = \Omega (\|\bbeta_0\|_2^2)$ and $\I_2 \ID(K)= o_{\pr}(\|\bbeta_0\|_2^2)$.

\item If $\|\bbeta_0\|_2^2/\sigma_{\epsilon}^2  = o(p_n/n) $, then $\I_1 \ID(K)  / \|\bbeta_0\|_2^2 \conP \infty$.
\end{itemize}
Hence, if and only if $p_n/n = o(\|\bbeta_0\|_2^2/\sigma_{\epsilon}^2) $, we have $ \|\hat\bbeta\|_2^2 - \|\bfbeta_0\|_2^2 = o_{\pr} ( \|\bfbeta_0\|_2^2 ). $

Next, we show that $p_n/n = o(\SNR )$ is equivalent to the ratio consistency of $\hat\zeta_0^2$ for $\zeta_0^2$. From Proposition 1, it suffices to show that $p_n/n = o(\SNR )$ is equivalent to that $\hat\bbeta^T (X^TX)^{-1} \hat\bbeta$ is ratio consistent for $\bbeta_0^T \E\{(X^TX)^{-1} \ID(K)\} \bbeta_0.$
Note that
\begin{eqnarray*}
 \hat\bbeta^T (X^TX)^{-1} \hat\bbeta
  &=&\bbeta_0^T (X^TX)^{-1} \bbeta_0 + 2 \bbeta_0^T (X^TX)^{-2} X^T\beps + \beps^T X(X^TX)^{-3} X^T\beps \cr
  &=& \II_1 + \II_2 + \II_3.
\end{eqnarray*}
We know $\II_2 \ID(K)= \Omega_{\pr}(\sigma_{\epsilon}\|\bbeta_0\|_2n^{-3/2})$ and
$\II_3\ID(K)= \Omega_{\pr} ( \sigma_{\epsilon}^2 p_n/n^2)$. Following arguments  similar  to \eqref{S.5.2} and \eqref{S.5.3}, we can show that $\II_3\ID(K)/ ( \sigma_{\epsilon}^2 p_n/n^2)$ and $\II_2 \ID(K) /(\sigma_{\epsilon}\|\bbeta_0\|_2n^{-3/2})$ are uniformly integrable.
From Lemma 2, we have $\II_1\ID(K) = \bbeta_0^T \E\{(X^TX)^{-1} \ID(K)\} \bbeta_0 \{1 + o_{\pr}(1)\} = \Omega_{\pr}( \| \bbeta_0\|_2^2 / n)$, which means that $\II_1 \ID(K) / [\bbeta_0^T \E\{(X^TX)^{-1} \ID(K)\} \bbeta_0] \conP 1$.
\begin{itemize}
\item If $p_n/n = o(\SNR) $, then $\II_2 \ID(K)= o_{\pr}(\II_1\ID(K))$ and $\II_3 \ID(K)= o_{\pr}(\II_1\ID(K))$.
Hence, $\hat\bbeta^T (X^TX)^{-1} \hat\bbeta / [\bbeta_0^T \E\{(X^TX)^{-1} \ID(K)\} \bbeta_0] \conP 1.$
\item If $p_n/n \gtrsim\SNR $, suppose that $\hat\bbeta^T (X^TX)^{-1} \hat\bbeta / [\bbeta_0^T \E\{(X^TX)^{-1} \ID(K)\} \bbeta_0] \conP 1.$ Then, we have $(\II_2 + \II_3) / [\bbeta_0^T \E\{(X^TX)^{-1} \ID(K)\} \bbeta_0] \conP 0$, which implies that $(\II_2 + \II_3) / (\sigma_{\epsilon}^2 p_n/n^2)\conP 0$ because $ \sigma_{\epsilon}^2 p_n/n^2 \gtrsim \bbeta_0^T \E\{(X^TX)^{-1} \ID(K)\} \bbeta_0$.
    Due to uniform integrability, we have $\E\{\II_2\ID(K) + \II_3\ID(K)\} / (\sigma_{\epsilon}^2 p_n/n^2) \to 0$.
However, because
    $\E\{\II_3 \ID(K)\}/ (\sigma_{\epsilon}^2 p_n/n^2)  \gtrsim 1$ and $\E\{\II_2 \ID(K)\}= 0$, there is a contradiction. Therefore, the initial assumption $\hat\bbeta^T (X^TX)^{-1} \hat\bbeta / [\bbeta_0^T \E\{(X^TX)^{-1} \ID(K)\} \bbeta_0] \conP 1$ doesn't hold.
\end{itemize}

Lastly, we study the asymptotic normal distribution of $\|\hat\bbeta\|_2^2$.
\begin{itemize}
\item If $p_n^2/n = o(\|\bbeta_0\|_2^2/\sigma_{\epsilon}^2) $, then $\I_1 \ID(K)= o_{\pr}(\I_2\ID(K))$.
From  Slutsky's theorem, we have $  \zeta_0^{-1}   (\|\hat\bbeta\|_2^2 - \|\bfbeta_0\|_2^2 ) \conD N(0, 1)$.
\item If $p_n^2/n = \Omega(\|\bbeta_0\|_2^2/\sigma_{\epsilon}^2) $, assuming $  \zeta_0^{-1}  (\|\hat\bbeta\|_2^2 - \|\bfbeta_0\|_2^2 ) \conD N(0, 1)$, we have $\E \{\zeta_0^{-1}    (\|\hat\bbeta\|_2^2 - \|\bfbeta_0\|_2^2 ) \ID(K)\} \to 0$ due to uniform integrability. However, $\E \{\zeta_0^{-1} (\|\hat\bbeta\|_2^2 - \|\bfbeta_0\|_2^2 ) \ID(K)\} = \zeta_0^{-1} \sigma_{\epsilon}^2 \E\tr\{(X^TX)^{-1}\ID(K)\} = \Omega(1)$. Hence, the assumption $  \zeta_0^{-1} (\|\hat\bbeta\|_2^2 - \|\bfbeta_0\|_2^2 ) \conD N(0, 1)$ does not hold.
\item If $\|\bbeta_0\|_2^2/\sigma_{\epsilon}^2 = o(p_n^2/n) $, then, $ \E\{ \zeta_0^{-1} (\|\hat\bbeta\|_2^2 - \|\bfbeta_0\|_2^2 ) \} \to \infty$.
\end{itemize}
Under Condition A3, $\SNR = \Omega(\| \bbeta_0 \|_2^2 / \sigma_{\epsilon}^2 ) = O(p_n)$. Hence $p_n^2/n = o(\SNR) = o(p_n)$ implies that $p_n /n = o(1)$ i.e. $\tau = 0$.
We complete the proof.
\endpf

\ni\textbf{Proof of Theorem 5:} The proof for fixed $p_n$ is straightforward, and we only focus on diverging $p_n$ here.
Note that
\begin{eqnarray*}
 \hat\bbeta^T (X^TX/n)\hat \bbeta
&=& \{\bbeta_0 + (X^TX)^{-1} X^T \beps\}^T (X^TX/n)  \{\bbeta_0 + (X^TX)^{-1} X^T \beps\} \cr
&=& \bbeta_0^T (X^TX/n) \bbeta_0 + 2 \bbeta_0^T (X^TX/n) (X^TX)^{-1} X^T \beps \cr
&& + \beps^T X (X^TX)^{-1} (X^TX/n)(X^TX)^{-1} X^T \beps  \cr
&=& \bbeta_0^T (X^TX/n) \bbeta_0 + 2 \bbeta_0^T  X^T \beps /n + \beps^T X (X^TX)^{-1}  X^T \beps / n.
 \end{eqnarray*}
Next, we aim to find the joint limiting distribution of   $\bbeta_0^T (X^TX/n) \bbeta_0 - \bbeta_0^T \Sigma \bbeta_0 $, $2 \bbeta_0^T  X^T \beps /n$, $\beps^T X (X^TX)^{-1}  X^T \beps / n  - \hat\sigma_{\epsilon}^2 p_n/n$ and $ \hat\sigma_{\epsilon}^2 - \sigma_{\epsilon}^2$.
Let
\begin{eqnarray*}
 \I_1 &=& \sqrt n \sigma_1^{-1} \{\bbeta_0^T (X^TX/n) \bbeta_0 - \bbeta_0^T \Sigma \bbeta_0\}, \cr
 \I_2 &=& \sigma_2^{-1} \bbeta_0^T  X^T \beps /\sqrt n,  \cr
 \I_3 &=& \sqrt n \sigma_3^{-1}\beps^T \{X (X^TX)^{-1}  X^T  - p_n / n\bfI_n\}\beps / (n-p_n),  \cr
 \I_4 &=& \sqrt n \sigma_4^{-1}( \hat\sigma_{\epsilon}^2 - \sigma_{\epsilon}^2 ) = \sqrt n \sigma_4^{-1}[\beps^T \{\bfI_n - X (X^TX)^{-1}  X^T  \}\beps /(n-p_n)- \sigma_{\epsilon}^2],
\end{eqnarray*}
where
\begin{eqnarray*}
 \sigma_1^2 &=& \var\{(\bbeta_0^T \bfX_1)^2\} = \E\Big (\sum_{i=1}^n Y_i^4/n\Big) - \nu_4 - 6 \sigma_{\epsilon}^2 \bbeta_0^T \Sigma \bbeta_0 - (\bbeta_0^T \Sigma \bbeta_0)^2, \cr
 \sigma_2^2 &=& \sigma_{\epsilon}^2 \bbeta_0^T \Sigma \bbeta_0, \cr
 \sigma_3^2 &=& 2 \sigma_{\epsilon}^4 p_n /(n-p_n),  \cr
 \sigma_4^2 &=& \nu_4 + \sigma_{\epsilon}^4 ( 3\tau - 1)/(1-\tau).
\end{eqnarray*}

From the end of the proof of Theorem 1, we can show that, there exists a set $\mathcal G_n \subseteq \mathbb{R}^{n\times p_n}$, such that conditioning on $X=x\in \mathcal G_n$, the joint cumulative distribution function of $(\I_2,\I_3,\I_4)$ converges  to that of a trivariate normal distribution, with $\pr(X\in \mathcal G_n)\geq 1-\varepsilon_n$ where $\varepsilon_n \to 0$; also, $\I_2$ and $(\I_3, \I_4)$ are asymptotically independent.
  It's easy to see that $\I_1$ is asymptotically standard normal using central limit theorem. Then, up to $\varepsilon_n$, which is arbitrarily close to 0, we have
  \begin{eqnarray*}
 && \E e^{it(c_1\I_1 +c_2\I_2 +c_3\I_3 +c_4\I_4  )} = \E[\E \{e^{it(c_1\I_1 +c_2\I_2 +c_3\I_3 +c_4\I_4  )} |X\} ] = \E[e^{itc_1\I_1}\E \{e^{it(c_2\I_2 +c_3\I_3 +c_4\I_4  )} |X\} ] \cr
 &=& \E[e^{itc_1\I_1}e^{-\sigma^{*2}t^2/2} ] +o(1) = e^{-c_1^2t^2 / 2 -\sigma^{*2}t^2/2} +o(1),
   \end{eqnarray*}
   where $\sigma^{*2} = \var(c_2\I_2 +c_3\I_3 +c_4\I_4)$. Hence, $\I_1$ and $(\I_2, \I_3, \I_4)$ are asymptotically joint normal and independent.

   Next, we will calculate the correlation between $\I_3$ and $\I_4$.
Following arguments similar to those in the proof of Proposition 1,
   \begin{eqnarray*}
 \E(\beps^T \beps) &=& n\sigma_{\epsilon}^2 \cr
 \E\{(\beps^T \beps)^2\} &=& n \nu_4 + n(n-1) \sigma_{\epsilon}^4 \cr
 \E\{\beps^T X (X^TX)^{-1} X^T\beps\} &=& p_n \sigma_{\epsilon}^2  \cr
 \E[\{\beps^T X (X^TX)^{-1} X^T\beps\}^2] &=& (\nu_4 - 3 \sigma_{\epsilon}^4) p_n^2/n \{1 + o(1)\} + 2 \sigma_{\epsilon}^4 p_n + \sigma_{\epsilon}^4 p_n ^2 \cr
 \E[\{\beps^T X (X^TX)^{-1} X^T\beps\} (\beps^T \beps)]&=& \nu_4 p_n  +  \sigma_{\epsilon}^4 (n -1)p_n.
\end{eqnarray*}
Therefore,
\begin{eqnarray*}
 \E(\I_3 \I_4)
 &=& n\sigma_3^{-1} \sigma_4^{-1}/ (n-p_n)^{2} \E [\beps^T \{X (X^TX)^{-1}  X^T  - p_n / n\bfI_n\}\beps
 \beps^T \{\bfI_n - X (X^TX)^{-1}  X^T  \}\beps]\cr
&=& n\sigma_3^{-1} \sigma_4^{-1}/ (n-p_n)^{2} [\nu_4 p_n  +  \sigma_{\epsilon}^4 (n -1)p_n -p_n/n \{n \nu_4 + n(n-1) \sigma_{\epsilon}^4\} \cr
&&-\{(\nu_4 - 3 \sigma_{\epsilon}^4) p_n^2/n \{1 + o(1)\}+ 2 \sigma_{\epsilon}^4 p_n + \sigma_{\epsilon}^4 p_n ^2\} + p_n/n\{\nu_4 p_n  +  \sigma_{\epsilon}^4 (n -1)p_n\} ] \cr
&=& n\sigma_3^{-1} \sigma_4^{-1}/ (n-p_n)^{2} \sigma_{\epsilon}^4 (2p_n^2 /n - 2p_n)\{1 + o(1)\}= 2 \sigma_3^{-1} \sigma_4^{-1} \sigma_{\epsilon}^4 p_n/(p_n-n)\{1 + o(1)\}.
 \end{eqnarray*}

 Then,  $\sqrt n\{ \hat\bbeta^T (X^TX/n)\hat \bbeta- \hat\sigma_{\epsilon}^2 p_n/n -  \bbeta_0^T \Sigma \bbeta_0\}$ and $\sqrt n( \hat\sigma_{\epsilon}^2 - \sigma_{\epsilon}^2)$ are asymptotically jointly normal with mean zero and covariance matrix
 \begin{eqnarray*}
 \Sigma_0=\left(
 \begin{matrix}
  \sigma_1^2 + 4 \sigma_2^2 + \sigma_3^2 & 2  \sigma_{\epsilon}^4 p_n/(p_n-n)\cr
 2  \sigma_{\epsilon}^4 p_n/(p_n-n)& \sigma_4^2
 \end{matrix}
 \right).
 \end{eqnarray*}
By delta method,
 \begin{eqnarray*}
  \sigma_{\hat \rho} ^{-1} (\hat \rho - \rho_0)  \conD N(0, 1)
   \end{eqnarray*}
   with
   \begin{eqnarray*}
 \sigma_{\hat \rho}^2 &=& n^{-1} (\bbeta_0^T\Sigma \bbeta_0 + \sigma_{\epsilon}^2)^{-4} (\sigma_{\epsilon}^2, -\bbeta_0^T\Sigma \bbeta_0  )\Sigma_0 (\sigma_{\epsilon}^2, -\bbeta_0^T\Sigma \bbeta_0  )^T \cr
 &=& n^{-1} (\bbeta_0^T\Sigma \bbeta_0 + \sigma_{\epsilon}^2)^{-4} \Big[\sigma_{\epsilon}^4 \Big\{\E\Big (\sum_{i=1}^n Y_i^4/n\Big) - \nu_4 - 2 \sigma_{\epsilon}^2 \bbeta_0^T \Sigma \bbeta_0 - (\bbeta_0^T \Sigma \bbeta_0)^2 \cr
 & & + 2 \sigma_{\epsilon}^4 p_n /(n-p_n) \Big\} + (\bbeta_0^T\Sigma \bbeta_0)^2 \{\nu_4 + \sigma_{\epsilon}^4 ( 3\tau - 1)/(1-\tau)\} \cr
 && - 4  \sigma_{\epsilon}^6 p_n/(p_n-n) (\bbeta_0^T \Sigma \bbeta_0) \Big]\cr
 &=&n^{-1} (\bbeta_0^T\Sigma \bbeta_0 + \sigma_{\epsilon}^2)^{-4} \Big[2 \sigma_{\epsilon}^8 p_n /(n-p_n) -\{2 + 4 p_n/(p_n-n)\} \sigma_{\epsilon}^6 \bbeta_0^T \Sigma \bbeta_0  \cr
 &&+ \sigma_{\epsilon}^4 \Big\{\E\Big(\sum_{i=1}^n Y_i^4/n\Big) - \nu_4   +(\bbeta_0^T\Sigma \bbeta_0)^2  ( 4\tau - 2)/(1-\tau)\Big\}
  + (\bbeta_0^T\Sigma \bbeta_0)^2 \nu_4  \Big].
   \end{eqnarray*}
\endpf

\ni\textbf{Proof of Corollary 1: }
When $\delta_n = O( \sigma_{\epsilon}^2 p_n^{1/2} /n )$, we have that $4 \sigma_{\epsilon}^2  \bbeta_0^T \E\{(X^TX)^{-1} \ID(K)\} \bbeta_0 = O(\sigma_{\epsilon}^4 p_n^{1/2} /n^2 ) = o(\sigma_{\epsilon}^4 p_n  /n^2  ) = o(\zeta_*^2)$. Therefore, $\zeta_n^2 / \zeta_*^2 \to 1$ as $n\to\infty$.
Due to Lemma 2 and Proposition 1, we can show that  $ \hat \zeta_*^2 / \zeta_*^2 \conP 1$.
From Theorem 1,
$
  ({\hat{\|\bbeta\|_2^2} } - \delta_n)/{\hat \zeta_*} \conD N(0,1),
$
which directly implies the results of this corollary.
\endpf

\ni\textbf{Proof of Proposition 1: }
The proof of Proposition 1 follows from that of Theorem 1. We will only calculate $\zeta_{\epsilon}^2$. Denote $M = \bfI_n - X(X^TX)^{-1}X^T$ and hence
\begin{eqnarray*}
&&\E \{n  (\hat \sigma_{\epsilon}^2  - \sigma_{\epsilon}^2 )^2\} = n \E    (\hat \sigma_{\epsilon}^4) - n\sigma_{\epsilon}^4 \cr
&=& n \E \{ \beps^T M \beps \beps^T M \beps\}/(n-p_n)^2 - n\sigma_{\epsilon}^4 \cr
&=& n /(n-p_n)^2 \E\Big[\nu_4 \sum_{i=1}^n \{M(i,i)\}^2 + 2 \sigma_{\epsilon}^4 \sum_{i\neq j}\{M(i,j)\}^2
+ \sigma_{\epsilon}^4 \sum_{i\neq j}M(i,i)M(j,j)\Big ] - n\sigma_{\epsilon}^4 \cr
&=& \frac n {(n-p_n)^2 } \E\Big[(\nu_4 -3\sigma_{\epsilon}^4)\sum_{i=1}^n \{M(i,i)\}^2 + 2 \sigma_{\epsilon}^4 \sum_{i, j}\{M(i,j)\}^2
+ \sigma_{\epsilon}^4 \sum_{i, j}M(i,i)M(j,j)\Big ] - n\sigma_{\epsilon}^4 \cr
&=& \frac n {(n-p_n)^2 } [(\nu_4 -3\sigma_{\epsilon}^4)n(1-p_n/n)^2\{1 + o(1)\} + 2 \sigma_{\epsilon}^4 (n-p_n)
+ \sigma_{\epsilon}^4 (n-p_n)^2  ] - n\sigma_{\epsilon}^4 \cr
&=& (\nu_4 -3\sigma_{\epsilon}^4) \{1 + o(1)\} + 2 \sigma_{\epsilon}^4 n/(n-p_n) \cr
&=& \{ \nu_4 + \sigma_{\epsilon}^4 ( 3\tau - 1)/(1-\tau) \} \{1+o(1)\}.
\end{eqnarray*}
Since $ \zeta_{\epsilon} = O(\sigma_{\epsilon}^2 / \sqrt n)$, the asymptotically normality directly implies that $\hat \sigma_{\epsilon}^2$ is ratio consistent for $ \sigma_{\epsilon}^2$.
\endpf

\ni\textbf{Proof of Proposition 2: }
Since $0< C_1\leq\rho_0\leq C_2<1$, we know $\eta_0 = \Omega(\sigma_{\epsilon}^2)$ and it's easy to see that $ \sigma^{2}_{\hat \rho}  = \Omega(1/n)$.
Following the proof of Theorem 5, Lemma \ref{Lemma-S15} and Proposition 1, we complete the proof.
\endpf

\section{Technical lemmas}  \label{Section-S.6}
This section includes the lemmas that are needed in the proofs of the main theoretical results in the paper.

\begin{lemma} \label{Lemma-S1} Assume Conditions $\mathrm{A1}$ and $\mathrm{A3}$, for any $k \in \mathbb{N}$ and deterministic $\bfa \in \mathbb{R}^{p_n}$,
we have $\E \|\bfX_1\|_2^k \lesssim  p_n^{k/2}$ and
$\E(\bfa^T \bfX_1) ^{2k} \lesssim \|\bfa\|_2^{2k} $.
\end{lemma}
\pf The result is straightforward for $k = 0,1$. Next, we only consider $k\geq 2$.
 First, from H\"{o}lder's inequality, $\E \|\bfX_1\|_2^k \lesssim \E \|\bfZ_1\|_2^k = \E (\sum_{i=1}^{p_n} z_{1i}^2)^{k/2} \leq p_n^{k/2-1}\E (\sum_{i=1}^{p_n} |z_{1i}|^k) \lesssim p_n^{k/2}$. Next, let $\bfa^* =  \Sigma^{1/2} \bfa$. Then,
$\E(\bfa^T \bfX_1) ^{2k} = \E(\bfa^{*T} \bfZ_1) ^{2k} = \E(\sum_{i=1}^{p_n} a^*_i z_{1i}) ^{2k}$.
From Theorem 2.6.3 of \cite{Vershynin_2018}, we have
\begin{eqnarray*}
\pr\Big( \Big|\sum_{i=1}^{p_n} a^*_i z_{1i}\Big| \geq t \Big ) \leq 2 \exp\Big( - \frac{C t^2 }{K_0^2 \| \bfa^* \|_2^2 }\Big),
\end{eqnarray*}
where $K_0\leq 1/\sqrt{c^*}$. Hence,
$ \E(\sum_{i=1}^{p_n} a^*_i z_{1i}) ^{2k} \lesssim \|\bfa^*\|_2^{2k}\lesssim \|\bfa\|_2^{2k} $. We finish the proof.
\endpf
\begin{lemma} \label{Lemma-S2} Assume the conditions of Lemma 2, $\lim_{n\to\infty} p_n = \infty$ and $\bfX_0$ is an $\iid$ copy of $\bfX_1$. For any integer $k >0 $,
$\var \{ n^{k }/p_n \bfX_0^T(X^T X)^{-k}\bfX_0 \ID(H) \} = o(1)$. Also, $\bfX_1^T (X^T X)^{-1} \bfX_1 - p_n/n = o_{\pr}(1).$
\end{lemma}
\pf Let $W =  n^{k }/p_n \bfX_0^T(X^T X)^{-k}\bfX_0 \ID(H) $, $W_{(0)} =  n^{k }/p_n \tr\{(X^T X)^{-k}\Sigma\} \ID(H) $.

First, we will show that $\E(W - W_{(0)})^2 = o(1)$. Let $B = \Sigma^{1/2} (X^TX)^{-k} \Sigma^{1/2}$ and $\bfZ = \Sigma^{-1/2}\bfX_0 \equiv (Z_1, \ldots, Z_{p_n})^T$. Since $\E(W|X) = W_{(0)}$, we have $\E(W - W_{(0)})^2 = \E W^2 - \E W_{(0)}^2 = \E \{n^{k}/p_n\bfZ^T B \bfZ\ID(H)\}^2 - \E \{n^{k}/p_n\tr(B)\ID(H)\}^2$.
Note that
\begin{eqnarray*}
& & \E \{n^{k}/p_n\bfZ^T B \bfZ \ID(H) \}^2 = n^{2k}/p_n^2 \sum_{ijkh} \E (Z_i Z_j Z_k Z_h B_{ij} B_{kh} \ID(H)) \cr
&=& n^{2k}/p_n^2 \sum_i \E(Z_i^4) \E (B_{ii}^2 \ID(H)) +  n^{2k}/p_n^2 \sum_{i \neq k} \E (Z_i ^2 Z_k ^2) \E( B_{ii} B_{kk} \ID(H))\cr
&& +
 2n^{2k}/p_n^2 \sum_{i \neq j} \E (Z_i ^2Z_j ^2)\E( B_{ij} B_{ij} \ID(H))\cr
 &=& O(1/p_n) + n^{2k}/p_n^2 \E\Big\{\sum_{i }  B_{ii} \ID(H) \Big\}^2 + 2n^{2k}/p_n^2 \sum_{i  j} \E\{ B_{ij}^2 \ID(H) \}\cr
 &=& O(1/p_n) +  n^{2k}/p_n^2 \E\{\tr(B) \ID(H)\}^2 + 2n^{2k}/p_n^2  \E \{\tr(B^2) \ID(H)\}\cr
 &=&\E \{n^{k}/p_n \tr(B) \ID(H) \}^2 + O(1/p_n).
\end{eqnarray*}
Hence, $\E(W - W_{(0)})^2 = O(1/p_n)$.
Following the proof of Lemma 2, we can show that $\var(W_{(0)}) = o(1)$. Hence $\var(W) \leq 2 \var(W_{(0)}) + 2 \var(W - W_{(0)}) = o(1)$.

Next, since
 \begin{eqnarray*}
\bfX_1^T (X^T X)^{-1} \bfX_1 = \bfX_1^T (X_{(1)}^T X_{(1)})^{-1}\bfX_1 /\{1 + \bfX_1^T (X_{(1)}^T X_{(1)})^{-1}\bfX_1\},
\end{eqnarray*}
and $ \bfX_1^T (X_{(1)}^T X_{(1)})^{-1}\bfX_1 - \E\{ \bfX_1^T (X_{(1)}^T X_{(1)})^{-1}\bfX_1\} = o_{\pr}(1)$,
 by continuous mapping theorem
 $$\bfX_1^T (X^T X)^{-1} \bfX_1 - \frac{\E\{\bfX_1^T (X_{(1)}^T X_{(1)})^{-1}\bfX_1\} }{1 +\E\{ \bfX_1^T (X_{(1)}^T X_{(1)})^{-1}\bfX_1\}} = o_{\pr}(1).$$
By dominant convergence theorem,
$$\E \{\bfX_1^T (X^T X)^{-1} \bfX_1\} - \frac{ \E\{\bfX_1^T (X_{(1)}^T X_{(1)})^{-1}\bfX_1\} }{1 +\E\{ \bfX_1^T (X_{(1)}^T X_{(1)})^{-1}\bfX_1\}}= o(1).$$
Since $\E \{\bfX_1^T (X^T X)^{-1} \bfX_1\} =\sum_{i=1}^n\E \{\bfX_i^T (X^T X)^{-1} \bfX_i\} /n = \E \tr\{ X (X^T X)^{-1}  X \} /n = p_n/n$, we have
$\bfX_1^T (X^T X)^{-1} \bfX_1 - p_n/n = o_{\pr}(1).$
\endpf
\begin{lemma} \label{Lemma-S3} Under the conditions of Lemma 2,
for any positive integer $k  $,
\begin{eqnarray*}
\E\{\bfX_1^T(X^TX)^{-1} \bfX_2 \ID(H) \}^{2k} &=& O(p_n^kn^{-2k}),\cr
\E\{\bbeta_0^T(X^TX)^{-1} \bfX_2 \ID(H) \}^{2k} &=& O( \|\bbeta_0\|_2^{2k}n^{-2k}),\cr
\E\{\bfX_1^T(X^TX)^{-2} \bfX_2\ID(H)\}^{2k} &=& O(p_n^kn^{-4k}), \cr
\E\{\bfX_1^T(X^TX)^{-2} \Sigma (X^TX)^{-2} \bfX_2\ID(H)\}^{2k} &=& O(p_n^kn^{-8k}), \cr
\E\{\bfX_1^T(X^TX)^{-2} \Sigma (X^TX)^{-1} \bfX_2\ID(H)\}^{2k} &=& O(p_n^kn^{-6k}), \cr
\E\{\bfX_1^T(X^TX)^{-1} \Sigma (X^TX)^{-1} \bfX_2\ID(H)\}^{2k} &=& O(p_n^kn^{-4k}).
\end{eqnarray*}
\end{lemma}
\pf
First, from (A.3) that $$\bfX_1^T(X^T X)^{-1} \bfX_2 = \bfX_1^T(X_{(2)}^T X_{(2)})^{-1}\bfX_2 /\{1 + \bfX_2^T (X_{(2)}^T X_{(2)})^{-1}\bfX_2\},$$
  using Lemmas 1 and \ref{Lemma-S1}, by Cauchy-Schwarz inequality,
\begin{eqnarray*}
&& \E\{\bfX_1^T(X^TX)^{-1} \bfX_2 \ID(H)\}^{2k} \cr
&\lesssim& \E\{\bfX_1^T(X^TX)^{-1} \bfX_2 \ID(H_{(2)})\}^{2k} + \E\{\bfX_1^T(X^TX)^{-1} \bfX_2 \ID(H) \ID(\bar H_{(2)})\}^{2k} \cr &\lesssim& \E\{ \bfX_1^T(X_{(2)}^T X_{(2)})^{-1}\bfX_2 \ID(H_{(2)} ) \}^{2k}  + o(p_n^k n^{-2k}) \cr
&\lesssim&
\E\| \bfX_1^T(X_{(2)}^T X_{(2)})^{-1} \ID(H_{(2)}) \|_2^{2k}  + o(p_n^k n^{-2k}) \cr
&\lesssim&
\{\E\| \bfX_1\|_2^{4k} \E\| (X_{(2)}^T X_{(2)})^{-1} \ID(H_{(2)}) \|_2^{4k}\}^{1/2}  + o(p_n^k n^{-2k}) = O(p_n^k n^{-2k}).
 \end{eqnarray*}

Second, similarly, from (A.3) and Lemma \ref{Lemma-S1},
\begin{eqnarray*}
&& \E\{\bbeta_0^T(X^TX)^{-1} \bfX_2 \ID(H) \}^{2k} \cr
&\lesssim& \E\{ \bbeta_0^T(X_{(2)}^T X_{(2)})^{-1}\bfX_2 \ID(H_{(2)} ) \}^{2k}
+\E\{ \bbeta_0^T(X_{(2)}^T X_{(2)})^{-1}\bfX_2 \ID(H) \ID(\bar H_{(2)}) \}^{2k} \cr
&\lesssim&
\E\| \bbeta_0^T(X_{(2)}^T X_{(2)})^{-1} \ID(H_{(2)} )\|_2^{2k} +o( \|\bbeta_0\|_2^{2k} n^{-2k})= O( \|\bbeta_0\|_2^{2k}n^{-2k}).
\end{eqnarray*}

Third, from (A.4), Lemma \ref{Lemma-S1} and Cauchy-Schwarz inequality,
\begin{eqnarray*}
&&\E \{\bfX_i^T(X^T X)^{-2}\bfX_j \ID(H)\}^{2k} \cr
&\lesssim& \E \{\bfX_i^T(X^T X)^{-2}\bfX_j \ID(H_{(2)} ) \}^{2k}+\E \{\bfX_i^T(X^T X)^{-2}\bfX_j \ID(H) \ID(\bar H_{(2)} )\}^{2k} \cr
&\lesssim& \E\{\bfX_i^T(X_{(j)}^T X_{(j)})^{-2}\bfX_j \ID(H_{(2)} ) \}^{2k} \cr
&&+ \E\{\bfX_i^T(X_{(j)}^T X_{(j)})^{-1} \bfX_j \bfX_j^T (X_{(j)}^T X_{(j)})^{-2} \bfX_j \ID(H_{(2)} ) \}^{2k} + o(n^{-4k})\cr
&\lesssim& \E\{\|\bfX_i^T(X_{(j)}^T X_{(j)})^{-2}\|_2^{2k} \ID(H_{(2)} )\}  +[\E\{\bfX_i^T(X_{(j)}^T X_{(j)})^{-1} \bfX_j \ID(H_{(2)} ) \}^{4k} \cr
&&\cdot \E\{ \bfX_j^T (X_{(j)}^T X_{(j)})^{-2} \bfX_j \ID(H_{(2)} ) \}^{4k}]^{1/2}+ o(n^{-4k})\cr
&=& O(p_n^kn^{-4k}).
\end{eqnarray*}
Fourth,
\begin{eqnarray*}
&& \E\{\bfX_2^T(X^TX)^{-2} \Sigma (X^TX)^{-2} \bfX_1 \ID(H) \}^{2k} \cr
&\lesssim& \E\{\bfX_2^T(X^TX)^{-2} \Sigma (X^TX)^{-2} \bfX_1 \ID(H_{(1)} ) \}^{2k}
 + \E\{\bfX_2^T(X^TX)^{-2} \Sigma (X^TX)^{-2} \bfX_1 \ID(H) \ID(\bar H_{(1)} ) \}^{2k} \cr
 &\lesssim& \E\{\bfX_2^T(X^TX)^{-2} \Sigma (X^TX)^{-2} \bfX_1 \ID(H_{(1)} ) \}^{2k} +o(n^{-8k}).
\end{eqnarray*}
Then, conditioning on $H_{(1)}$, from (A.2), by Lemma \ref{Lemma-S1} and Cauchy-Schwarz inequality,
\begin{eqnarray*}
&&\E\{\bfX_2^T(X^TX)^{-2} \Sigma (X^TX)^{-2} \bfX_1\}^{2k} \cr
&\lesssim& \E\{\bfX_2^T(X_{(1)}^TX_{(1)})^{-2} \Sigma (X_{(1)}^TX_{(1)})^{-2} \bfX_1\}^{2k} \cr
&&+ \E\{\bfX_2^T(X_{(1)}^T X_{(1)})^{-2} \bfX_1 \bfX_1^T (X_{(1)}^T X_{(1)})^{-1}\Sigma (X_{(1)}^TX_{(1)})^{-2} \bfX_1\}^{2k} \cr
&& + \E\{\bfX_2^T(X_{(1)}^T X_{(1)})^{-1} \bfX_1 \bfX_1^T (X_{(1)}^T X_{(1)})^{-2}\Sigma (X_{(1)}^TX_{(1)})^{-2} \bfX_1\}^{2k} \cr
&& + \E\{\bfX_2^T(X_{(1)}^T X_{(1)})^{-1} \bfX_1 \bfX_1^T (X_{(1)}^T X_{(1)})^{-2} \bfX_1
 \bfX_1^T (X_{(1)}^T X_{(1)})^{-1}\Sigma (X_{(1)}^TX_{(1)})^{-2} \bfX_1\}^{2k}\cr
&&+\E\{\bfX_2^T(X_{(1)}^TX_{(1)})^{-2} \Sigma (X_{(1)}^TX_{(1)})^{-1} \bfX_1 \bfX_1^T (X_{(1)}^TX_{(1)})^{-2} \bfX_1\}^{2k} \cr
&&+ \E\{\bfX_2^T(X_{(1)}^T X_{(1)})^{-2} \bfX_1 \bfX_1^T (X_{(1)}^T X_{(1)})^{-1}\Sigma (X_{(1)}^TX_{(1)})^{-1} \bfX_1
 \bfX_1^T (X_{(1)}^TX_{(1)})^{-2} \bfX_1\}^{2k} \cr
&& + \E\{\bfX_2^T(X_{(1)}^T X_{(1)})^{-1} \bfX_1 \bfX_1^T (X_{(1)}^T X_{(1)})^{-2}\Sigma (X_{(1)}^TX_{(1)})^{-1} \bfX_1
\bfX_1^T (X_{(1)}^TX_{(1)})^{-2} \bfX_1\}^{2k} \cr
&& + \E\{\bfX_2^T(X_{(1)}^T X_{(1)})^{-1} \bfX_1 \bfX_1^T (X_{(1)}^T X_{(1)})^{-2} \bfX_1 \cr
&&\quad \cdot \bfX_1^T (X_{(1)}^T X_{(1)})^{-1}\Sigma (X_{(1)}^TX_{(1)})^{-1} \bfX_1  \bfX_1^T (X_{(1)}^TX_{(1)})^{-2} \bfX_1\}^{2k}\cr
&=& O(p_n^kn^{-8k}).
\end{eqnarray*}
Following similar arguments, we can show the other results.
\endpf

\begin{lemma} \label{Lemma-S4}
Under the conditions of Lemma 2,
for any positive integer $k$,
\begin{eqnarray*}
\E[\{\bfX_1^T(X^TX)^{-1} \bfX_1 -\bfX_1^T(X_{(3)}^TX_{(3)})^{-1} \bfX_1 \}^{2k} \ID(H)]&=& O(p_n^{2k}n^{-4k}), \cr
\E[\{\bfX_1^T(X^TX)^{-2} \bfX_1-\bfX_1^T(X_{(3)}^TX_{(3)})^{-2} \bfX_1\}^{2k}  \ID(H)]&=& O(p_n^{2k}n^{-6k}), \cr
\E[\{\bfX_1^T(X^TX)^{-1} \bfX_2 -\bfX_1^T(X_{(3)}^TX_{(3)})^{-1} \bfX_2 \}^{2k}  \ID(H)]&=& O(p_n^{2k}n^{-4k}), \cr
\E[\{\bfX_1^T(X^TX)^{-2} \bfX_2-\bfX_1^T(X_{(3)}^TX_{(3)})^{-2} \bfX_2\}^{2k}  \ID(H)]&=& O(p_n^{2k}n^{-6k}), \cr
\E[\{\bbeta_0^T(X^TX)^{-1} \bfX_1 -\bbeta_0^T(X_{(3)}^TX_{(3)})^{-1} \bfX_1 \}^{2k}  \ID(H)]&=& O(\|\bbeta_0\|_2^{2k}n^{-3k}), \cr
\E[\{\bbeta_0^T(X^TX)^{-2} \bbeta_0 -\bbeta_0^T(X_{(3)}^TX_{(3)})^{-2} \bbeta_0 \}^{2k}  \ID(H)]&=& O(\|\bbeta_0\|_2^{4k}n^{-6k}).
\end{eqnarray*}
\end{lemma}
\pf
From (A.1), we know
\begin{eqnarray*}
&&\bfX_1^T(X^TX)^{-1} \bfX_2 -\bfX_1^T(X_{(3)}^TX_{(3)})^{-1} \bfX_2 \cr
&=&  - \bfX_1^T(X_{(3)}^T X_{(3)})^{-1} \bfX_3 \bfX_3^T (X_{(3)}^T X_{(3)})^{-1} \bfX_2 /\{1 + \bfX_3^T (X_{(3)}^T X_{(3)})^{-1}\bfX_3\},
\end{eqnarray*}
which together with Cauchy-Schwarz inequality implies that,
\begin{eqnarray*}
&&\E[\{\bfX_1^T(X^TX)^{-1} \bfX_2 -\bfX_1^T(X_{(3)}^TX_{(3)})^{-1} \bfX_2\}^{2k} \ID(H)] \cr
&\lesssim & \E[\{\bfX_1^T(X^TX)^{-1} \bfX_2 -\bfX_1^T(X_{(3)}^TX_{(3)})^{-1} \bfX_2\}^{2k} \ID(H_{(3)})]\cr
 &&+  \E[\{\bfX_1^T(X^TX)^{-1} \bfX_2 -\bfX_1^T(X_{(3)}^TX_{(3)})^{-1} \bfX_2\}^{2k} \ID(H) \ID(\bar H_{(3)})] \cr
&\lesssim& [\E\{\bfX_1^T(X_{(3)}^T X_{(3)})^{-1} \bfX_3\ID(H_{(3)})\}^{4k} \E\{\bfX_3^T (X_{(3)}^T X_{(3)})^{-1} \bfX_2\ID(H_{(3)})\}^{4k}]^{1/2}
 + o(p_n^{2k}n^{-4k})\cr
&\lesssim&
[\E\{\bfX_1^T(X_{(3)}^T X_{(3)})^{-2} \bfX_1\ID(H_{(3)})\}^{2k} \E\{\bfX_2^T (X_{(3)}^T X_{(3)})^{-2} \bfX_2\ID(H_{(3)})\}^{2k}]^{1/2}
 +o(p_n^{2k}n^{-4k}) \cr
&=&O(p_n^{2k}n^{-4k}).
\end{eqnarray*}
Following the same arguments, we can show other results.
\endpf
\begin{lemma} \label{Lemma-S5} With the notations and conditions in the proof of part (a) of Theorem 1 with $\lim_{n\to\infty}p_n=\infty$,
$$\sum_{j=1}^n \E \{ U_j  \ID(K)\}^4 = o(t_n^4).$$
\end{lemma}
\pf Recalling $t_n = \|\bbeta_0\|_2/\sqrt n + \sqrt{p_n}/n$, it suffices to show that for any $j=1,\ldots, n$, $k = 1$ or $2$, $c_1=\Omega( 1)$, $c_2 =\Omega( p_n / n)$ and $c_3=\Omega(1)$,
\begin{eqnarray}
\E\{\bfbeta_0^T (X^T X)^{-1} \bfX_j \epsilon_j / \|\bbeta_0\|_2\ID(K) \}^4&\leq& C n^{-4}, \label{S.6.1}\\
\E [\{c_kM_k(j,j) \sqrt{n/p_n}\}^4 (\epsilon_j ^2 - \sigma_{\epsilon}^2)^4 \ID(K) ] &\leq& C n^{-4}, \label{S.6.2}\\
  \E \Big\{\sum_{1 \leq i < j} c_k M_k(i,j) \sqrt{n/p_n} \epsilon_i\epsilon_j \ID(K) \Big\}^4 &\leq& C n^{-4}. \label{S.6.3}
\end{eqnarray}
First, \eqref{S.6.1} follows from Lemma \ref{Lemma-S3}.

 Second, for \eqref{S.6.2}, since $ \E\{(\epsilon_j ^2 - \sigma_{\epsilon}^2)^4\} \leq C<\infty$, from Lemma \ref{Lemma-S3},
\begin{eqnarray*}
&& \E [\{c_1M_1(j,j) \sqrt{n/p_n}\}^4 (\epsilon_j ^2 - \sigma_{\epsilon}^2)^4 \ID(K) ]  \lesssim n^2 /p_n^2 \E \{M_1(j,j)\ID(K)\}^4 \cr &=& n^2 /p_n^2 \E \{\bfX_j (X^TX)^{-2}\bfX_j \ID(K) \}^4
\lesssim n^2 /p_n^2 O( p_n^4 n^{-8})
= O( p_n^2 n^{-6}) = O(n^{-4}),
\end{eqnarray*}
and
\begin{eqnarray*}
&& \E [\{c_2 M_2(j,j)\sqrt{n/p_n}\}^4 (\epsilon_j ^2 - \sigma_{\epsilon}^2)^4 \ID(K) ] \lesssim p_n^2 / n^2 \E \{M_2(j,j)\ID(K)\}^4 \cr &\lesssim& n^{-4} +  \E \{\bfX_j (X^TX)^{-1}\bfX_j /n \ID(K)\}^4
 =O( n^{-4}).
\end{eqnarray*}
Third, for \eqref{S.6.3},
\begin{eqnarray} \label{S.6.4}
 && \E \Big\{\sum_{1 \leq i < j} c_k M_k(i,j) \sqrt{n/p_n} \epsilon_i\epsilon_j \ID(K)\Big\}^4
\lesssim \frac {n^2}{p_n^2}\E \Big[ \Big\{\sum_{1 \leq i < j} c_k M_k(i,j)  \epsilon_i \Big\}^4\ID(K)\Big] \cr
&\lesssim& \frac {n^2}{p_n^2} \E  \Big[\sum_{1 \leq i < j} \{ c_k M_k(i,j)\} ^4 \ID(K)\Big]
 + \frac {n^2}{p_n^2} \E \Big [ \sum_{1 \leq i < j} \{c_k M_k(i,j) \} ^2 \ID(K)\Big]^2 \cr
&\lesssim&  \frac {n^3}{p_n^2} \E \{ c_k M_k(i,j) \ID(K)\} ^4 + \frac {n^2}{p_n^2} \E \Big [ \sum_{1 \leq i < j} \{ c_k M_k(i,j) \} ^2 \ID(K) \Big]^2.
\end{eqnarray}
From Lemma \ref{Lemma-S3}, for the first term in \eqref{S.6.4} with $k=1,2$,
\begin{eqnarray} \label{S.6.5}
n^3 / p_n^2 \E \{ c_1 M_1(i,j) \ID(K) \} ^4 &\lesssim& (n^3 / p_n^2)( p_n^2 n^{-8}) = O( n^{-5}), \cr
n^3 / p_n^2 \E \{ c_2 M_2(i,j) \ID(K) \} ^4& = &(n^3 / p_n^2)( p_n^4 /n^4) [ n^{-4} \E \{ \bfX_i^T(X^T X)^{-1}\bfX_j \ID(K)\} ^4 ] \cr
&\lesssim& (n^3 / p_n^2)( p_n^4 /n^4)  n^{-4} (p_n^2 / n^4)   = O(n^{-5}).
\end{eqnarray}
From Lemma 1 and Cauchy-Schwarz inequality, for the second term in \eqref{S.6.4} with $k=1$,
\begin{eqnarray} \label{S.6.6}
&&\E \Big[ \sum_{1 \leq i < j}  \{c_1 M_1(i,j)\}^2 \ID(K) \Big]^2 \cr
&=& \E \Big[ \sum_{1 \leq i < j}  \{\bfX_i^T(X^TX)^{-2}\bfX_j\}^2 \ID(K) \Big]^2 \cr
 &=&  \E \Big[  \bfX_j^T(X^TX)^{-2}\Big\{\sum_{1 \leq i < j} (\bfX_i   \bfX_i^T)\Big\}(X^TX)^{-2}\bfX_j \ID(K) \Big]^2\cr
&\leq&   \E \Big[  \|\bfX_j\|_2^4 \|(X^TX)^{-2}\|_2^4 \Big\|\sum_{1 \leq i < j} (\bfX_i   \bfX_i^T)\Big\|_2^2 \ID(K) \Big] \cr
&\leq &   [\E \{  \|\bfX_j\|_2^8  \|(X^TX)^{-2}\|_2^8\} \ID(K) ]^{1/2} \Big\{\E \Big\|\sum_{1 \leq i < j} (\bfX_i   \bfX_i^T)\ID(K)\Big\|_2^4 \Big\}^{1/2} \cr
&\leq&    (\E   \|\bfX_j\|_2^{16})^{1/4} \{\E \|(X^TX)^{-2}\ID(K)\|_2^{16}\}^{1/4}
 \Big\{\E \Big\|\sum_{1 \leq i < j} (\bfX_i   \bfX_i^T)\ID(K)\Big\|_2^4 \Big\}^{1/2} \cr
& \lesssim& p_n^2 n^{-8} (\E \|X^TX\|_2^4 \ID(K))^{1/2} \leq C p_n^2 n^{-6},
\end{eqnarray}
and with $k=2$,
\begin{eqnarray} \label{S.6.7}
&&  \E \Big[ \sum_{1 \leq i < j}  \{c_2 M_2(i,j)\}^2 \ID(K) \Big]^2 \cr
&=& \frac{p_n^4} {n^4}\E \Big[ \sum_{1 \leq i < j}  \{\bfX_i^T(X^TX)^{-1}\bfX_j/n\}^2 \ID(K) \Big]^2 \cr
&=&  \frac{p_n^4} {n^4} n^{-4}\E \Big[  \bfX_j^T(X^TX)^{-1}\Big\{\sum_{1 \leq i < j} (\bfX_i   \bfX_i^T) \Big\}(X^TX)^{-1}\bfX_j\ID(K) \Big]^2\cr
&\leq& \frac{p_n^4} {n^4}  n^{-4}\E \Big[  \|\bfX_j\|_2^4 \|(X^TX)^{-1}\|_2^4 \Big\|\sum_{1 \leq i < j} (\bfX_i   \bfX_i^T) \Big\|_2^2 \ID(K)\Big] \cr
&\leq &  \frac{p_n^4} {n^4} n^{-4}[\E \{  \|\bfX_j\|_2^8  \|(X^TX)^{-1}\|_2^8\ID(K)\}]^{1/2}
\Big\{\E \Big\|\sum_{1 \leq i < j} (\bfX_i   \bfX_i^T)\ID(K)\Big\|_2^4 \Big\}^{1/2} \cr
&\leq&  \frac{p_n^4} {n^4}  n^{-4}(\E   \|\bfX_j\|_2^{16})^{1/4} \{\E \|(X^TX)^{-1}\ID(K)\|_2^{16}\}^{1/4}
 \Big\{\E \Big\|\sum_{1 \leq i < j} (\bfX_i   \bfX_i^T)\ID(K)\Big\|_2^4 \Big\}^{1/2} \cr
& \lesssim& \frac{p_n^4} {n^4} n^{-4}p_n^2 n^{-4} \Big\{\E \Big\|\sum_{i=1}^n (\bfX_i   \bfX_i^T)\ID(K)\Big\|_2^4 \Big\}^{1/2} \lesssim \frac{p_n^6} {n^{12}}   (\E \|X^TX\|_2^4 )^{1/2}\cr
& =& O( p_n^6 n^{-10}).
\end{eqnarray}
From \eqref{S.6.4}, \eqref{S.6.5}, \eqref{S.6.6}, \eqref{S.6.7},
we complete the proof.
\endpf
\begin{lemma} \label{Lemma-S6} With the notations and conditions in the proof of part (a) of Theorem 1 with $\lim_{n\to\infty}p_n=\infty$,
$$\var \Big ( \sum_{j = 1}^n  \Big  [\sum_{1 \leq i < j} \{c_1 M_1(i,j) + c_2M_2(i,j)\} \epsilon_i \Big ]^2 \ID(K) \Big ) = o( t_n^4).$$
\end{lemma}
\pf Denote $Q_n = \sum_{j = 1}^n     [\sum_{1 \leq i < j} \{c_1 M_1(i,j) + c_2M_2(i,j)\}  \epsilon_i   ]^2 \ID(K)$. First, we calculate the expectation of $Q_n$,
\begin{eqnarray} \label{S.6.8}
\E(Q_n)&=& \sum_{j = 1}^n  \E\Big[\sum_{1 \leq i < j}\{c_1 M_1(i,j) + c_2M_2(i,j)\} \epsilon_i\ID(K)\Big]^2 \cr
&=& \sigma_{\epsilon}^2  \sum_{j = 1}^n  \sum_{1 \leq i < j} \E[\{c_1 M_1(i,j) + c_2M_2(i,j)\} ^2\ID(K)] \cr
&=& \sigma_{\epsilon}^2 \sum_{j = 1}^n  \sum_{1 \leq i < j} \E[\{c_1 \bfX_i^T(X^TX)^{-2} \bfX_j
- c_2\bfX_i^T(X^TX)^{-1} \bfX_j/(n-p_n)\} ^2 \ID(K)]\cr
&=& \frac 1 2 \sigma_{\epsilon}^2 \sum_{j = 1}^n  \sum_{ i \neq j} \E[\{c_1^2 \bfX_j^T(X^TX)^{-2} \bfX_i \bfX_i^T(X^TX)^{-2} \bfX_j \cr
&&- 2c_1 c_2\bfX_j^T(X^TX)^{-2} \bfX_i \bfX_i^T(X^TX)^{-1} \bfX_j/(n-p_n) \cr
&& + c_2^2 \bfX_j^T(X^TX)^{-1} \bfX_i \bfX_i^T(X^TX)^{-1} \bfX_j/(n-p_n)^2\}\ID(K)] \cr
 &=& \frac 1 2 \sigma_{\epsilon}^2 \sum_{j = 1}^n  \E[\{c_1^2 \bfX_j^T(X^TX)^{-3} \bfX_j  - 2c_1 c_2 \bfX_j^T(X^TX)^{-2}\bfX_j/(n-p_n)\cr
 &&\quad
 + c_2^2  \bfX_j^T(X^TX)^{-1}  \bfX_j/(n-p_n)^2\}\ID(K)] \cr
 && - \frac 1 2 \sigma_{\epsilon}^2 \sum_{j = 1}^n   \E[\{c_1^2 \bfX_j^T(X^TX)^{-2} \bfX_j \bfX_j^T(X^TX)^{-2} \bfX_j \cr
&&\quad- 2c_1 c_2\bfX_j^T(X^TX)^{-2} \bfX_j \bfX_j^T(X^TX)^{-1} \bfX_j/(n-p_n) \cr
&&\quad + c_2^2 \bfX_j^T(X^TX)^{-1} \bfX_j \bfX_j^T(X^TX)^{-1} \bfX_j/(n-p_n)^2\}\ID(K)] \cr
&=& \frac 1 2 \sigma_{\epsilon}^2   \E[ \tr\{c_1^2 (X^TX)^{-2} - 2c_1 c_2 (X^TX)^{-1}/(n-p_n)
 + c_2^2 \bfI_{p_n}/(n-p_n)^2 \}\ID(K)]  \cr
  &&- \frac 1 2 \sigma_{\epsilon}^2 \sum_{j = 1}^n   \E[\{c_1 \bfX_j^T(X^TX)^{-2} \bfX_j -  c_2 \bfX_j^T(X^TX)^{-1} \bfX_j/(n-p_n)\}^2\ID(K)]. \qquad
\end{eqnarray}

Second, we calculate the expectation of $Q_n^2$,
\begin{eqnarray*}
&&\E(Q_n^2)\cr
&=&\E\Big( \sum_{j = 1}^n  \Big[\sum_{1 \leq i < j} \{c_1 M_1(i,j) + c_2M_2(i,j)\} \epsilon_i\Big]^2 \ID(K)\Big)^2 \cr
&=& \sum_{j = 1}^n \sum_{j' = 1}^n \E \Big(\Big[\sum_{1 \leq i < j} \{c_1 M_1(i,j) + c_2M_2(i,j)\} \epsilon_i\Big]^2  \cr && \qquad \Big[\sum_{1 \leq i' < j'} \{c_1 M_1(i',j') + c_2M_2(i',j')\} \epsilon_{i'}\Big]^2\ID(K)\Big) \cr
&=& \sum_{j = 1}^n \sum_{j' = 1}^n \sum_{  i < j} \sum_{  k < j} \sum_{  i' < j'} \sum_{  k' < j'} \E[\{c_1 M_1(i,j) + c_2M_2(i,j)\} \epsilon_i
 \{c_1 M_1(k,j) + c_2M_2(k,j)\} \epsilon_k \cr
&&\qquad\{c_1 M_1(i',j') + c_2M_2(i',j')\} \epsilon_{i'}
 \{c_1 M_1(k',j') + c_2M_2(k',j')\} \epsilon_{k'} \ID(K)] \cr
&=& \sum_{j = 1}^n \sum_{j' = 1}^n \sum_{1 \leq i < j \wedge j'} \E[\{c_1 M_1(i,j) + c_2M_2(i,j)\}^2
 \{c_1 M_1(i,j') + c_2M_2(i,j')\}^2 \E\epsilon_{i}^4 \ID(K)] \cr
&&+ \sigma_{\epsilon}^4 \sum_{j = 1}^n \sum_{j' = 1}^n \sum_{  i < j}  \sum_{  i' < j', i'\neq i} \E[\{c_1 M_1(i,j) + c_2M_2(i,j)\}^2
 \{c_1 M_1(i',j') + c_2M_2(i',j')\}  ^2 \ID(K)] \cr
&&+\sigma_{\epsilon}^4 \sum_{j = 1}^n \sum_{j' = 1}^n \sum_{  i < j\wedge j'} \sum_{  k < j \wedge j', k\neq i}
\E[\{c_1 M_1(i,j) + c_2M_2(i,j)\} \{c_1 M_1(k,j) + c_2M_2(k,j)\} \cr
&&\qquad  \{c_1 M_1(i,j') + c_2M_2(i,j')\}  \{c_1 M_1(k,j') + c_2M_2(k,j')\} \ID(K) ] \cr
&&+ \sigma_{\epsilon}^4 \sum_{j = 1}^n \sum_{j' = 1}^n \sum_{  i < j\wedge j'} \sum_{ k < j\wedge j', k\neq i}
\E[\{c_1 M_1(i,j) + c_2M_2(i,j)\} \{c_1 M_1(k,j) + c_2M_2(k,j)\} \cr
&&\qquad \{c_1 M_1(k,j') + c_2M_2(k,j')\} \{c_1 M_1(i,j') + c_2M_2(i,j')\} \ID(K) ] \cr
&=& \III_1 + \III_2 + \III_3 + \III_4.
\end{eqnarray*}

Following the proof for \eqref{S.6.5}, $\III_1 /t_n^4= n^{3}O(p_n^2/n^{8} + p_n^6 / n^{12}) /t_n^4 = O( 1/n)$. Note
\begin{eqnarray*}
\III_3&=& \sigma_{\epsilon}^4 \sum_{j = 1}^n \sum_{  i < j} \sum_{  k < j , k\neq i}
\E[\{c_1 M_1(i,j) + c_2M_2(i,j)\} \{c_1 M_1(k,j) + c_2M_2(k,j)\}\cr
 && \qquad\{c_1 M_1(i,j) + c_2M_2(i,j)\}  \{c_1 M_1(k,j) + c_2M_2(k,j)\} \ID(K) ] \cr
 && +\sigma_{\epsilon}^4 \sum_{j = 1}^n \sum_{j' = 1, j' \neq j}^n \sum_{ i < j\wedge j'} \sum_{  k < j \wedge j', k\neq i} \E[\{c_1 M_1(i,j) + c_2M_2(i,j)\}  \cr
   &&\qquad \{c_1 M_1(k,j) + c_2M_2(k,j)\} \{c_1 M_1(i,j') + c_2M_2(i,j')\}  \cr
   &&\qquad \{c_1 M_1(k,j') + c_2M_2(k,j')\} \ID(K) ]\cr
   &=& \III_{3,1} + \III_{3,2}.
 \end{eqnarray*}
Similar to  term  $\III_1$, it's easy to show that $\III_{3,1}/t_n^4 = n^{3}O(p_n^2/n^{8} + p_n^6 / n^{12})O(n^4 / p_n^2) = O( 1/n)$.

Noting that $c_1 = \Omega(1)$ and
\begin{eqnarray*}
 \III_{3,2}
 &\lesssim& n^{4} \E[\{c_1 M_1(1,3) + c_2M_2(1,3)\} \{c_1 M_1(2,3) + c_2M_2(2,3)\}  \cr
&& \{c_1 M_1(2,4) + c_2M_2(2,4)\} \{c_1 M_1(1,4) + c_2M_2(1,4)\}  \ID(K)],
 \end{eqnarray*}
 we first calculate the rate of
 \begin{eqnarray*}
 &&n^{4}\E\{  M_1(1,3)  M_1(2,3)  M_1(2,4)   M_1(1,4)    \ID(K)\}\cr
 &=&n^{4} \E\{  \bfX_1^T (X^TX)^{-2}\bfX_3
 \bfX_2^T (X^TX)^{-2}\bfX_3   \bfX_2^T (X^TX)^{-2}\bfX_4   \bfX_1^T (X^TX)^{-2}\bfX_4 \ID(K)\}
 \end{eqnarray*}
  by replacing $(X^TX)^{-2}$ with $(X_{(1,2,3,4)}^TX_{(1,2,3,4)})^{-2}$, where $X_{(1,2,3,4)}$ is the design matrix without the first 4 observations.

 From Lemmas \ref{Lemma-S3} and \ref{Lemma-S4},
 \begin{eqnarray*}
&& n^{4} |\E[\{  M_1(1,3)  -   \bfX_1^T (X_{(2)}^TX_{(2)})^{-2}\bfX_3 \} M_1(2,3)  M_1(2,4)   M_1(1,4)   \ID(K)]| \cr
&\leq& n^{4}[ \E\{  M_1(1,3)  -   \bfX_1^T (X_{(2)}^TX_{(2)})^{-2}\bfX_3 \}^2 \E\{ M_1(2,3)  M_1(2,4)   M_1(1,4)   \ID(K)\}^2]^{1/2} \cr
&\lesssim& n^{4} \{p_n^2 / n^6 (p_n/n^4)^3\}^{1/2} = p_n^{2.5} / n^5.
 \end{eqnarray*}
By similar arguments, we can show that
\begin{eqnarray} \label{S.6.9}
&&n^{4}\E\{  M_1(1,3)  M_1(2,3)  M_1(2,4)   M_1(1,4)    \ID(K)\} \cr
&=& n^{4} \E\{  \bfX_1^T (X_{(2,4)}^TX_{(2,4)})^{-2}\bfX_3
 \bfX_2^T (X_{(1,4)}^TX_{(1,4)})^{-2}\bfX_3  \cr
&&  \bfX_2^T (X_{(1,3)}^TX_{(1,3)})^{-2}\bfX_4   \bfX_1^T (X_{(2,3)}^TX_{(2,3)})^{-2}\bfX_4 \ID(K)\} +O(p_n^{2.5} / n^5) \cr
&\equiv& \IV + O(p_n^{2.5} / n^5).
 \end{eqnarray}
 From Lemma \ref{Lemma-S2}, $$\bfX_1^T (X_{(1,2,4)}^T X_{(1,2,4)})^{-1}\bfX_1-
 \E\{\bfX_1^T (X_{(1,2,4)}^T X_{(1,2,4)})^{-1}\bfX_1 \} = o_{\pr}(p_n/n).$$
For any positive integer $k$, $n^k/p_n^k \E[\bfX_1^T (X_{(1,2,4)}^T X_{(1,2,4)})^{-1}\bfX_1-
 \E\{\bfX_1^T (X_{(1,2,4)}^T X_{(1,2,4)})^{-1}\bfX_1 \}]^k =O(1)$ given $K_{(1,2,4)}$ (which is event $K$ calculated using $X_{(1,2,4)}$), and hence, due to uniform integrability,
\begin{eqnarray*}
\E[\bfX_1^T (X_{(1,2,4)}^T X_{(1,2,4)})^{-1}\bfX_1-
 \E\{\bfX_1^T (X_{(1,2,4)}^T X_{(1,2,4)})^{-1}\bfX_1 \}]^k = o(p_n^k/n^k).
 \end{eqnarray*}
 Similarly, given $K_{(1,2,4)}$,
\begin{eqnarray*}
\E[\bfX_1^T (X_{(1,2,4)}^T X_{(1,2,4)})^{-2}\bfX_1-
 \E\{\bfX_1^T (X_{(1,2,4)}^T X_{(1,2,4)})^{-2}\bfX_1 \}]^k = o(p_n^k/n^{2k}).
 \end{eqnarray*}
 Therefore, for any random variable $R$ satisfying $\E(R^2)=O(1)$ and any positive integer $k$,
 \begin{eqnarray} \label{S.6.added1}
&&\E|[n^k/p_n\bfX_1^T (X_{(1,2,4)}^T X_{(1,2,4)})^{-k}\bfX_1 -
 n^k/p_n\E\{\bfX_1^T (X_{(1,2,4)}^T X_{(1,2,4)})^{-k}\bfX_1 \}] R|^2 \cr
 &\leq& \E[n^k/p_n\bfX_1^T (X_{(1,2,4)}^T X_{(1,2,4)})^{-k}\bfX_1-
 n^k/p_n\E\{\bfX_1^T (X_{(1,2,4)}^T X_{(1,2,4)})^{-k}\bfX_1 \}]^2 \E (R^2)\cr
 &=&o(1).
 \end{eqnarray}
From \eqref{S.6.added1}, both
$n^k/p_n\E\{\bfX_1^T (X_{(1,2,4)}^T X_{(1,2,4)})^{-k}\bfX_1 R\}$ and $n^k/p_n\E\{\bfX_1^T (X_{(1,2,4)}^T X_{(1,2,4)})^{-k}\bfX_1 \} \E(R)$ are  $O(1)$ but their difference is $o(1)$. That is, we can calculate  $n^k/p_n\E\{\bfX_1^T (X_{(1,2,4)}^T X_{(1,2,4)})^{-k}\bfX_1 R\}$ by replacing $\bfX_1^T (X_{(1,2,4)}^T X_{(1,2,4)})^{-k}\bfX_1$ with $\E\{\bfX_1^T (X_{(1,2,4)}^T X_{(1,2,4)})^{-k}\bfX_1\}$ without influencing its magnitude.

By (A.4),
 \begin{eqnarray} \label{S.6.10}
\IV&=&n^{4} \E\Big(\Big[
 \frac{\bfX_1^T (X_{(1,2,4)}^T X_{(1,2,4)})^{-2} \bfX_3 } { 1 + \bfX_1^T (X_{(1,2,4)}^T X_{(1,2,4)})^{-1}\bfX_1 } \cr
 &&
 - \frac{ \bfX_1^T (X_{(1,2,4)}^T X_{(1,2,4)})^{-1}  \bfX_3 \bfX_1^T (X_{(1,2,4)}^T X_{(1,2,4)})^{-2} \bfX_1 }{\{1 + \bfX_1^T (X_{(1,2,4)}^T X_{(1,2,4)})^{-1}\bfX_1\}^2} \Big]
 \bfX_2^T (X_{(1,4)}^TX_{(1,4)})^{-2}\bfX_3\cr
 &&
   \bfX_2^T (X_{(1,3)}^TX_{(1,3)})^{-2}\bfX_4  \bfX_1^T (X_{(2,3)}^TX_{(2,3)})^{-2}\bfX_4 \ID(K_{(1,2,4)})\Big)  + o(p_n^2 / n^4)\cr
   &=&
n^{4} C' \E\{ \bfX_1^T (X_{(1,2,4)}^T X_{(1,2,4)})^{-2} \bfX_3 \bfX_2^T (X_{(1,4)}^TX_{(1,4)})^{-2}\bfX_3 \cr
&&   \bfX_2^T (X_{(1,3)}^TX_{(1,3)})^{-2}\bfX_4   \bfX_1^T (X_{(2,3)}^TX_{(2,3)})^{-2}\bfX_4 \ID(K_{(1,2,4)})\} \cr
&&
   + n^{3} C'' \E\{ \bfX_1^T (X_{(1,2,4)}^T X_{(1,2,4)})^{-1} \bfX_3 \bfX_2^T (X_{(1,4)}^TX_{(1,4)})^{-2}\bfX_3 \cr
&&   \bfX_2^T (X_{(1,3)}^TX_{(1,3)})^{-2}\bfX_4   \bfX_1^T (X_{(2,3)}^TX_{(2,3)})^{-2}\bfX_4 \ID(K_{(1,2,4)})\} + o(p_n^2 / n^4),
 \end{eqnarray}
 where
 \begin{eqnarray*}
 C' &=& [1 + \E\{\bfX_1^T (X_{(1,2,4)}^T X_{(1,2,4)})^{-1}\bfX_1 \}]^{-1}\cr
 C'' &=& - { n\E\{ \bfX_1^T (X_{(1,2,4)}^T X_{(1,2,4)})^{-2} \bfX_1\} } / {[1 + \E\{\bfX_1^T (X_{(1,2,4)}^T X_{(1,2,4)})^{-1}\bfX_1\}]^2}
   \end{eqnarray*}
   are constants bounded away from infinity. The intuition for the last equation in  \eqref{S.6.10} is that we can replace $\bfX_1^T (X_{(1,2,4)}^T X_{(1,2,4)})^{-k}\bfX_1$ by $
 \E\{\bfX_1^T (X_{(1,2,4)}^T X_{(1,2,4)})^{-k}\bfX_1 \}$ ($k=1,2$) without influencing the magnitude of $\IV$ due to \eqref{S.6.added1} and the fact that $\{\E(\IV^2) \}^{1/2}= O(p_n^2/n^4)$ using Lemma \ref{Lemma-S3}.

 Similarly, we can repeat the procedure in \eqref{S.6.10} and obtain that
 \begin{eqnarray*}
\IV   &=&  \sum_{i_1 = 1}^2 \sum_{i_2 = 1}^2 \sum_{i_3 = 1}^2 \sum_{i_4 = 1}^2
c_{i_1,i_2,i_3,i_4}
n^{i_1+i_2+i_3+i_4-4} \E\{  \bfX_1^T (X_{(1,2,3,4)}^TX_{(1,2,3,4)})^{-i_1}\bfX_3   \cr
&&\bfX_2^T (X_{(1,2,3,4)}^TX_{(1,2,3,4)})^{-i_2}\bfX_3    \bfX_2^T (X_{(1,2,3,4)}^TX_{(1,2,3,4)})^{-i_3}\bfX_4 \cr
&&  \bfX_1^T (X_{(1,2,3,4)}^TX_{(1,2,3,4)})^{-i_4}\bfX_4 \ID(K_{(1,2,3,4)})\}  + o(p_n^2 / n^4)\cr
 &=&  \sum_{i_1 = 1}^2 \sum_{i_2 = 1}^2 \sum_{i_3 = 1}^2 \sum_{i_4 = 1}^2
c_{i_1,i_2,i_3,i_4}
n^{i_1+i_2+i_3+i_4-4} \E\{  \bfX_1^T (X_{(1,2,3,4)}^TX_{(1,2,3,4)})^{-i_1}\Sigma   \cr
&& (X_{(1,2,3,4)}^TX_{(1,2,3,4)})^{-i_2} \Sigma (X_{(1,2,3,4)}^TX_{(1,2,3,4)})^{-i_3}\Sigma \cr
&&    (X_{(1,2,3,4)}^TX_{(1,2,3,4)})^{-i_4}\bfX_1 \ID(K_{(1,2,3,4)})\}  + o(p_n^2 / n^4)\cr
&=&  \sum_{i_1 = 1}^2 \sum_{i_2 = 1}^2 \sum_{i_3 = 1}^2 \sum_{i_4 = 1}^2
c_{i_1,i_2,i_3,i_4}
n^{i_1+i_2+i_3+i_4-4} O(p_n/n^{i_1+i_2+i_3+i_4})  + o(p_n^2 / n^4)\cr
&=& o(p_n^2 / n^4),
 \end{eqnarray*}
 where $c_{i_1,i_2,i_3,i_4} $ are constants that are bounded away from infinity and  $ K_{(1,2,3,4)}$ is defined in  the same way as $K$ by excluding $\bfX_1,\ldots, \bfX_4$ in $X$.
  Therefore, from \eqref{S.6.9},
 \begin{eqnarray*}
 n^{4}\E\{  M_1(1,3)  M_1(2,3)  M_1(2,4)   M_1(1,4)    \ID(K)\} = o(p_n^2 / n^4)+ O(p_n^{2.5} / n^5) = o(p_n^2 / n^4).
 \end{eqnarray*}
Following similar arguments for the other terms, we can show that $\III_{3,2} = o(p_n^2/n^4)=o(t_n^4)$.
Hence, $\III_3 /t_n^4 = o(1).$
Similarly, $\III_4/t_n^4 = o(1).$

Therefore, $\E(Q_n^2) - \III_2 = o(t_n^4)$. Next, it suffices to prove that $ \III_2 - \{\E(Q_n)\}^2= o(t_n^4)$.
Up to $O(p_n^2/n^{5} + p_n^6 / n^9)$, we have
\begin{eqnarray} \label{S.6.11}
\III_2&=&\sigma_{\epsilon}^4 \sum_{j = 1}^n \sum_{j' = 1}^n \sum_{1 \leq i < j}  \sum_{1 \leq i' < j'}
\E[\{c_1 M_1(i,j) + c_2M_2(i,j)\}^2  \{c_1 M_1(i',j') + c_2M_2(i',j')\}  ^2 \ID(K)] \cr
&=&\sigma_{\epsilon}^4   \E \Big [\sum_{j = 1}^n \sum_{1 \leq i < j} \{c_1 M_1(i,j) + c_2M_2(i,j)\}^2
 \sum_{j' = 1}^n  \sum_{1 \leq i' < j'} \{c_1 M_1(i',j') + c_2M_2(i',j')\}  ^2 \ID(K) \Big ] \cr
&=&\sigma_{\epsilon}^4    \E \Big( \Big [\sum_{j = 1}^n \sum_{1 \leq i < j} \{c_1 M_1(i,j) + c_2M_2(i,j)\}^2 \Big ] ^2 \ID(K) \Big ) \cr
&=&\sigma_{\epsilon}^4 \E \Big( \Big [\sum_{j = 1}^n \sum_{1 \leq i < j} \{c_1^2 \bfX_j^T (X^TX)^{-2} \bfX_i \bfX_i^T (X^TX)^{-2} \bfX_j\cr
&&
- 2c_1 c_2 \bfX_j^T (X^TX)^{-2} \bfX_i \bfX_i^T (X^TX)^{-1} \bfX_j/(n-p_n)\cr
&&
+ c_2^2 \bfX_j^T (X^TX)^{-1} \bfX_i \bfX_i^T (X^TX)^{-1} \bfX_j  /(n-p_n)^2 \} \Big ] ^2 \ID(K) \Big ) \cr
&=&\frac 1 4 \sigma_{\epsilon}^4  \E \Big( \Big [\sum_{j = 1}^n \sum_{i \neq j} \{c_1^2 \bfX_j^T (X^TX)^{-2} \bfX_i \bfX_i^T (X^TX)^{-2} \bfX_j\cr
&&
- 2c_1 c_2 \bfX_j^T (X^TX)^{-2} \bfX_i\bfX_i^T (X^TX)^{-1} \bfX_j/(n-p_n) \cr
&&
+ c_2^2 \bfX_j^T (X^TX)^{-1} \bfX_i \bfX_i^T (X^TX)^{-1} \bfX_j  /(n-p_n)^2 \} \Big ] ^2 \ID(K) \Big ) \cr
&=&\frac 1 4 \sigma_{\epsilon}^4   \E \Big( \Big [\sum_{j = 1}^n \{c_1^2 \bfX_j^T (X^TX)^{-3}  \bfX_j
- 2c_1 c_2 \bfX_j^T (X^TX)^{-2}  \bfX_j/(n-p_n)\cr
&&
+ c_2^2 \bfX_j^T (X^TX)^{-1} \bfX_j  /(n-p_n)^2 \} \cr
 &&-\sum_{j = 1}^n \{c_1^2 \bfX_j^T (X^TX)^{-2} \bfX_j \bfX_j^T (X^TX)^{-2} \bfX_j \cr
&&- 2c_1 c_2 \bfX_j^T (X^TX)^{-2} \bfX_j\bfX_j^T (X^TX)^{-1} \bfX_j/(n-p_n) \cr
&&+ c_2^2 \bfX_j^T (X^TX)^{-1} \bfX_j \bfX_j^T (X^TX)^{-1} \bfX_j  /(n-p_n)^2 \} \Big ] ^2 \ID(K) \Big ) \cr
&=&\frac 1 4 \sigma_{\epsilon}^4   \E \Big( \Big [\tr \{c_1^2  (X^TX)^{-2}
- 2c_1 c_2  (X^TX)^{-1}  /(n-p_n)
+ c_2^2 \bfI_{p_n} /(n-p_n)^2 \} \cr
&& -\sum_{j = 1}^n \{c_1 \bfX_j^T (X^TX)^{-2} \bfX_j - c_2 \bfX_j^T (X^TX)^{-1} \bfX_j  /(n-p_n) \}^2 \Big ] ^2 \ID(K) \Big ).
\end{eqnarray}
From Lemma 2 and the proof of Lemma \ref{Lemma-S8}, the variances of $ \tr \{  (X^TX)^{-2} \ID(K) \} /t_n^2 $, $ \tr \{  (X^TX)^{-1} $ $ \ID(K) \} /(n t_n^2)$ and $ \sum_{j = 1}^n \{c_1 \bfX_j^T (X^TX)^{-2} \bfX_j - c_2 \bfX_j^T (X^TX)^{-1} \bfX_j  /(n-p_n) \}^2 /t_n^2\ID(K) $ are all $o(1)$, combining \eqref{S.6.8} and \eqref{S.6.11} we can prove that $ \III_2 - \{\E(Q_n)\}^2= o(t_n^4)$ and hence $\E(Q_n^2/t_n^4) - \{\E(Q_n)/t_n^2\}^2 = o(1)$.
\endpf

\begin{lemma}\label{Lemma-S7}
With the notations and conditions in the proof of part (a) of Theorem 1 with $\lim_{n\to\infty}p_n=\infty$,
$$\var \Big [ \sum_{j}  \{\bfbeta_0^T (X^T X)^{-1} \bfX_j\}^2 /t_n^2 \ID(K) \Big] = o(1).$$
\end{lemma}
\pf
Since $ \sum_{j}  \{\bfbeta_0^T (X^T X)^{-1} \bfX_j\}^2  = \bfbeta_0^T (X^T X)^{-1} \bfbeta_0  $, the proof follows from Lemma 2.
\endpf
\begin{lemma}\label{Lemma-S8} With the notations and conditions in the proof of part (a) of Theorem 1 with $\lim_{n\to\infty}p_n=\infty$,
$$ \var \Big [n\sum_j \{c_1 M_1(j,j) + c_2M_2(j,j)\}^2 (n/p_n) ^2 \ID(K) \Big]=o(1).$$
\end{lemma}
\pf Denoting $ \III_{j,j} = c_1 M_1(j,j) + c_2M_2(j,j)$, then
\begin{eqnarray*}
&&\var \Big [n\sum_j \{c_1 M_1(j,j) + c_2M_2(j,j)\}^2 (n/p_n) ^2 \ID(K) \Big ] \cr
&\leq& n^4 \var [\{c_1 M_1(j,j) + c_2M_2(j,j)\}^2(n/p_n) ^2  \ID(K) ] \cr
&=& n^4 \E \{ \III_{j,j}^4(n/p_n) ^4 \ID(K)\} - n^4[\E\{\III_{j,j}^2 (n/p_n) ^2 \ID(K)\}]^2.
\end{eqnarray*}

From (A.3) and (A.5),
\begin{eqnarray*}
\bfX_1(X^T X)^{-1} \bfX_1 &=& \bfX_1^T (X_{(1)}^T X_{(1)})^{-1}\bfX_1 /\{1 + \bfX_1^T (X_{(1)}^T X_{(1)})^{-1}\bfX_1\}\cr
\bfX_1^T(X^T X)^{-2} \bfX_1&=& \bfX_1^T(X_{(1)}^T X_{(1)})^{-2}\bfX_1 /\{1 + \bfX_1^T (X_{(1)}^T X_{(1)})^{-1}\bfX_1\}^2,
\end{eqnarray*}
which together with Lemma \ref{Lemma-S2} implies
\begin{eqnarray*}
  n  \III_{j,j}(n/p_n) \ID(K)  - n \E\{ \III_{j,j}(n/p_n)  \ID(K)\} = o_{\pr}(1).
\end{eqnarray*}
Since $\E \{n  \III_{j,j} (n/p_n) \ID(K) \}^8 = O (1)$,
\begin{eqnarray*}
n^2  \III_{j,j}^2(n/p_n)^2\ID(K)  - n^2 [ \E\{\III_{j,j}(n/p_n)\ID(K)\}]^2 &=& o_{\pr}(1), \cr
n^4 \III_{j,j}^4(n/p_n)^4\ID(K) - n^4 [ \E\{\III_{j,j}(n/p_n)\ID(K)\}]^4 &=& o_{\pr}(1).
 \end{eqnarray*}
Then, by uniformly integrability,
\begin{eqnarray*}
 n^2 \E[\III_{j,j}(n/p_n)\ID(K)]^2  - n^2 [ \E\{\III_{j,j}(n/p_n)\ID(K)\}]^2 &=& o(1),\cr
 n^4 \E[\III_{j,j}(n/p_n)\ID(K)]^4  - n^4 [ \E\{\III_{j,j}(n/p_n)\ID(K)\}]^4 &=& o(1).
\end{eqnarray*}
  We complete the proof.
\endpf
\begin{lemma}\label{Lemma-S9} With the notations and conditions in the proof of part (a) of Theorem 1 with $\lim_{n\to\infty}p_n=\infty$,
 $$ \sum_{j=1}^n\II_{4,j} /t_n^2  \ID(H) = o_{L^2}(1).$$
\end{lemma}
\pf
First, we show
\begin{eqnarray}\label{S.6.12}
\sum_{j=1}^n \sum_{i=1}^{j-1} M_1(i,j) M_1(j,j) \epsilon_i /t_n^2 \ID(H) = o_{L^2}(1).
 \end{eqnarray}
 It's easy to see that  $\E\{\sum_{j=1}^n \sum_{i=1}^{j-1} M_1(i,j) M_1(j,j) \epsilon_i /t_n^2\ID(H) \} =0$. Note that
\begin{eqnarray*}
 \sum_{j=1}^n \sum_{i=1}^{j-1} M_1(i,j) M_1(j,j) \epsilon_i = \sum_{i=1}^{n-1} \epsilon_i \sum_{j=i+1}^{n} M_1(i,j) M_1(j,j).
\end{eqnarray*}
We only need to show that $ \sum_{i=1}^{n-1} \E\{ \epsilon_i \sum_{j=i+1}^{n} M_1(i,j) M_1(j,j)/t_n^2\ID(H) \}^2 =o(1)$. By (A.4) and Cauchy-Schwarz inequality, using the result that $\E\{\bfX_1^T (X_{(1)}^T X_{(1)})^{-2} \bfX_1\} = O(p_nn^{-2})$,
\begin{eqnarray*}
&&  \sum_{i=1}^{n-1} \E\Big \{ \epsilon_i \sum_{j=i+1}^{n} M_1(i,j) M_1(j,j)\ID(H) \Big\}^2 = \sigma_{\epsilon}^2 \sum_{i=1}^{n-1} \E\Big\{  \sum_{j=i+1}^{n} M_1(i,j) M_1(j,j) \ID(H) \Big\}^2 \cr
&=&\sigma_{\epsilon}^2  \sum_{i=1}^{n-1} \sum_{j=i+1}^{n} \sum_{k =i+1}^{n} \E\{ M_1(i,j) M_1(j,j)  M_1(i,k) M_1(k,k) \ID(H) \} \cr
&{=}& C'  n^2 \E\{ M_1(1,2)^2 M_1(2,2)^2 \ID(H)\}
 + C'' n^{3} \E\{ M_1(1,2) M_1(2,2)  M_1(1,3) M_1(3,3) \ID(H) \}\cr
&\stackrel{*}=&  n^2 O( p_n/n^4 p_n^2/n^4)+ n^{3}O(p_n^2/n^4) \E\{ M_1(1,2)  M_1(1,3)\ID(H)\}  +o({p_n}^{3}/n^5) \cr
&=&p_n^2 /n \E\{ \bfX_1 ^T(X^TX)^{-2} \bfX_2\bfX_1 ^T(X^TX)^{-2} \bfX_3 \ID(H) \}+o({p_n}^{3}/n^5)\cr
&=&p_n^2 /n \E\{ \bfX_1 ^T(X^TX)^{-2} \bfX_2\bfX_1 ^T(X^TX)^{-2} \bfX_3 \ID(H_{(1)}) \}+o({p_n}^{3}/n^5)\cr
&\lesssim& p_n^2 /n |\E\{ \bfX_2 ^T(X_{(1)}^TX_{(1)})^{-2} \Sigma (X_{(1)}^TX_{(1)})^{-2} \bfX_3 \ID(H_{(1)}) \}| \cr
&&+ p_n^{2} /n^2| \E\{ \bfX_2 ^T(X_{(1)}^TX_{(1)})^{-2} \Sigma (X_{(1)}^TX_{(1)})^{-1} \bfX_3\ID(H_{(1)})\}| \cr
&&+ p_n^{2} /n^2|\E\{ \bfX_2 ^T(X_{(1)}^TX_{(1)})^{-1} \Sigma (X_{(1)}^TX_{(1)})^{-2} \bfX_3\ID(H_{(1)})\}| \cr
&&+ p_n^2/n^3 |\E\{ \bfX_2 ^T(X_{(1)}^TX_{(1)})^{-1} \Sigma (X_{(1)}^TX_{(1)})^{-1} \bfX_3\ID(H_{(1)})\}|  + o({p_n}^{3}/n^5)\cr
&=& O(np_n^2 \sqrt {p_n}/n^6 +p_n^{2}\sqrt {p_n}/n^5 +   p_n^2/n \sqrt {p_n}/n^4 ) + o({p_n}^{3}/n^5) \cr
 &=& o({p_n}^{3}/n^5 ).
\end{eqnarray*}
The intuition for equation $*$ above is: to calculate $\E\{ M_1(1,2) M_1(2,2)  M_1(1,3) M_1(3,3) \ID(H) \}$, for $$M_1(k,k)=\bfX_k^T(X^T X)^{-2} \bfX_k= {\bfX_k^T(X_{(k)}^T X_{(k)})^{-2}\bfX_k }/{\{1 + \bfX_k^T (X_{(k)}^T X_{(k)})^{-1}\bfX_k\}^2}$$ we can replace $\bfX_k^T(X_{(k)}^T X_{(k)})^{-i}\bfX_k$ by $\E\{\bfX_k^T(X_{(k)}^T X_{(k)})^{-i}\bfX_k\}$ ($k=2,3$, $i=1,2$) without influencing the magnitude of $\E\{ M_1(1,2) M_1(2,2)  M_1(1,3) M_1(3,3) \ID(H) \}$, due to \eqref{S.6.added1}. The intuition for inequality ``$\lesssim$'' above is due to (A.4) and that we can replace $ \bfX_1^T (X_{(1)}^T X_{(1)})^{-2} \bfX_1   $ by $\E\{\bfX_1^T (X_{(1)}^T X_{(1)})^{-2} \bfX_1\}  $ and also $\E\{\bfX_1^T (X_{(1)}^T X_{(1)})^{-2} \bfX_1\} = O(p_nn^{-2})=O(n^{-1})$.
Since $({p_n}^{3}/n^5 ) / t_n^4 =O( 1)  $, we finish the proof for \eqref{S.6.12}.
The proof for the other terms are similarly.
\endpf
\begin{lemma}\label{Lemma-S10} With the notations and conditions in the proof of part (a) of Theorem 1 with $\lim_{n\to\infty}p_n=\infty$,
$$  \sum_{j=1}^n\II_{5,j}  /t_n^2 \ID(H) = o_{L^2}(1).$$
\end{lemma}
\pf First,
\begin{eqnarray*}
&& t_n^{-4}  \E\Big\{ \sum_{j = 1}^n  \sum_{i=1}^{j-1} M_1(i,j) \epsilon_i  v_j\ID(H) \Big\}^2
 = t_n^{-4}  \E\Big\{ \sum_{i=1}^{n-1} \epsilon_i \sum_{j = i+1}^n   M_1(i,j)   v_j \ID(H) \Big\}^2  \cr
&=&  t_n^{-4}  \sigma_{\epsilon}^2 \sum_{i=1}^{n-1} \E \Big\{\sum_{j = i+1}^n   M_1(i,j)   v_j \ID(H)\Big\}^2 \cr
&=&  t_n^{-4}  \sigma_{\epsilon}^2 \sum_{i=1}^{n-1}\sum_{j = i+1}^n\sum_{k = i+1}^n \E \{   M_1(i,j)   v_j  M_1(i,k)   v_k \ID(H) \} \cr
&\lesssim&  t_n^{-4}n^3  \E \{   M_1(1,2)   v_2  M_1(1,3)   v_3 \ID(H) \} +  t_n^{-4}n^2  \E \{   M_1(1,2)   v_2 \ID(H) \}^2 \cr
&=&   t_n^{-4}n^3  \E \{  \bfX_1^T (X^T X)^{-2}\bfX_2 \bbeta_0^T (X^T X)^{-1}\bfX_2
\bfX_1^T (X^T X)^{-2}\bfX_3 \bbeta_0^T (X^T X)^{-1}\bfX_3 \ID(H) \}\cr
&&  +  t_n^{-4}n^2  \E \{  \bfX_1^T (X^T X)^{-2}\bfX_2 \bbeta_0^T (X^T X)^{-1}\bfX_2  \ID(H) \}^2\cr
&=&\IV_1 + \IV_2.
\end{eqnarray*}

First, from Lemma \ref{Lemma-S3} and Cauchy-Schwarz inequality,  $\IV_1 =  t_n^{-4}n^3 O\{(p_n^{1/2} / n^2)^2 (\|\bbeta_0\|_2 / n)^2\} =t_n^{-4} O( p_n \|\bbeta_0\|_2^2 / n^3)=t_n^{-4} O( p_n^2 / n^4 +  \|\bbeta_0\|_2^4 / n^2) = O(1)$.
From Lemmas \ref{Lemma-S3} and \ref{Lemma-S4}, we have that
\begin{eqnarray*}
&&|\IV_1  -  t_n^{-4}n^3  \E \{  \bfX_1^T (X^T X)^{-2}\bfX_2 \bbeta_0^T (X_{(1)}^T X_{(1)})^{-1}\bfX_2 \cr
&& \qquad \cdot\bfX_1^T (X^T X)^{-2}\bfX_3 \bbeta_0^T (X^T X)^{-1}\bfX_3 \ID(H) \}| \cr
&=& t_n^{-4}n^3  |\E [  \bfX_1^T (X^T X)^{-2}\bfX_2 \bbeta_0^T \{(X^T X)^{-1} - (X_{(1)}^T X_{(1)})^{-1}\}\bfX_2 \cr
&&\qquad \cdot \bfX_1^T (X^T X)^{-2}\bfX_3 \bbeta_0^T (X^T X)^{-1}\bfX_3 \ID(H) ]|\cr
&=& t_n^{-4}n^3 O(p_n^{1/2} / n^2) o(\|\bbeta_0\|_2 / n)  O(p_n^{1/2} / n^2) O(\|\bbeta_0\|_2 / n) \cr
&=&o(1).
\end{eqnarray*}
By similar arguments and (A.4),  we have,
\begin{eqnarray} \label{S.6.13}
  \IV_1
&=&   t_n^{-4}n^3  \E \{  \bfX_1^T (X^T X)^{-2}\bfX_2 \bbeta_0^T (X_{(1)}^T X_{(1)})^{-1}\bfX_2  \cr
&&\quad \cdot \bfX_1^T (X^T X)^{-2}\bfX_3 \bbeta_0^T (X_{(1)}^T X_{(1)})^{-1}\bfX_3\ID(H_{(1)})\} + o(1)\cr
&\lesssim&
  t_n^{-4}n^3 |\E[ \bbeta_0^T (X_{(1)}^T X_{(1)})^{-1}\bfX_2 \bbeta_0^T (X_{(1)}^T X_{(1)})^{-1}\bfX_3 \cr
  &&\quad \cdot \bfX_2 ^T(X_{(1)}^TX_{(1)})^{-2} \Sigma (X_{(1)}^TX_{(1)})^{-2} \bfX_3 \ID(H_{(1)}) ]| \cr
&&+ t_n^{-4}n^2 |\E[ \bbeta_0^T (X_{(1)}^T X_{(1)})^{-1}\bfX_2 \bbeta_0^T (X_{(1)}^T X_{(1)})^{-1}\bfX_3 \cr
  &&\quad \cdot \bfX_2 ^T(X_{(1)}^TX_{(1)})^{-2} \Sigma (X_{(1)}^TX_{(1)})^{-1} \bfX_3 \ID(H_{(1)}) ]|\cr
&&+ t_n^{-4}n^2 |\E[ \bbeta_0^T (X_{(1)}^T X_{(1)})^{-1}\bfX_2 \bbeta_0^T (X_{(1)}^T X_{(1)})^{-1}\bfX_3  \cr
  &&\quad \cdot \bfX_2 ^T(X_{(1)}^TX_{(1)})^{-1} \Sigma (X_{(1)}^TX_{(1)})^{-2} \bfX_3\ID(H_{(1)})]|\cr
&&+ t_n^{-4}n |\E[ \bbeta_0^T (X_{(1)}^T X_{(1)})^{-1}\bfX_2 \bbeta_0^T (X_{(1)}^T X_{(1)})^{-1}\bfX_3  \cr
  &&\quad \cdot \bfX_2 ^T(X_{(1)}^TX_{(1)})^{-1} \Sigma (X_{(1)}^TX_{(1)})^{-1} \bfX_3\ID(H_{(1)})]| + o(1) \cr
&=& t_n^{-4} \|\bbeta_0\|_2^2 O(\sqrt {p_n}/n^3 ) +o(1)= o(1).
\end{eqnarray}
The intuition for inequality ``$ \lesssim$'' above is due to \eqref{S.6.added1} that in
\begin{eqnarray*}
\bfX_2^T(X^T X)^{-2} \bfX_1 = \frac{\bfX_2^T(X_{(1)}^T X_{(1)})^{-2} \bfX_1 } { 1 + \bfX_1^T (X_{(1)}^T X_{(1)})^{-1}\bfX_1 }
 - \frac{\bfX_2^T (X_{(1)}^T X_{(1)})^{-1}  \bfX_1 \bfX_1^T (X_{(1)}^T X_{(1)})^{-2} \bfX_1 }{\{1 + \bfX_1^T (X_{(1)}^T X_{(1)})^{-1}\bfX_1\}^2},
 \end{eqnarray*}
 we can replace $ \bfX_1^T (X_{(1)}^T X_{(1)})^{-2} \bfX_1   $ by $\E\{\bfX_1^T (X_{(1)}^T X_{(1)})^{-2} \bfX_1\}  $ and also $\E\{\bfX_1^T (X_{(1)}^T X_{(1)})^{-2} \bfX_1\} = O(p_nn^{-2})=O(n^{-1})$.

Also,
\begin{eqnarray*}
  \IV_2
&\leq&t_n^{-4} n^2 [\E \{  \bfX_1^T (X^T X)^{-2}\bfX_2\ID(H)\}^4 \E\{\bbeta_0^T (X^T X)^{-1}\bfX_2  \ID(H) \}^4]^{1/2}\cr
&=& t_n^{-4} O(n^2 \sqrt{p_n^2n^{-8} \|\bbeta_0\|_2^4n^{-4}}) = O(t_n^{-4}\|\bbeta_0\|_2^2 p_n n^{-4}) = o(1).
\end{eqnarray*}
Similarly, we can show that $ \E\{ \sum_{j = 1}^n  \sum_{i=1}^{j-1} c_2 M_2(i,j) \epsilon_i  v_j / t_n^2\}^2 = o(1)$. We complete the proof.
\endpf
\begin{lemma}\label{Lemma-S11} With the notations and conditions in the proof of part (a) of Theorem 1 with $\lim_{n\to\infty}p_n=\infty$,
$$ \sum_{j = 1}^n \II_{6,j} /t_n^2 \ID(H) =o_{L^2}(1).$$
\end{lemma}
\pf
It suffices to show  $$ \sum_{j=1}^n c_1 M_1(j,j) v_j/t_n^2\ID(H) = o_{L^2}(1),\quad \sum_{j=1}^nc_2 M_2(j,j) v_j /t_n^2\ID(H)= o_{L^2}(1).$$ First, we will show $\E\{ \sum_{j=1}^n M_1(j,j) v_j/t_n^2\ID(H)\}^2 = o(1)$. From Lemmas \ref{Lemma-S2} and \ref{Lemma-S4}, following arguments similar to those in \eqref{S.6.13},
\begin{eqnarray*}
&&\E \Big\{  \sum_{j=1}^n M_1(j,j) v_j /t_n^2 \ID(H) \Big\}^2 \cr
&=& t_n^{-4} \E\Big\{ \sum_{j=1}^n \bfX_j^T (X^TX)^{-2} \bfX_j  \bbeta_0^T (X^TX)^{-1} \bfX_j \ID(H)\Big\}^2\cr
&=& t_n^{-4} \sum_{j=1}^n\sum_{i=1}^n \E\{ \bbeta_0^T (X^TX)^{-1} \bfX_j \bfX_j^T (X^TX)^{-2} \bfX_j
 \bbeta_0^T (X^TX)^{-1} \bfX_i \bfX_i^T (X^TX)^{-2} \bfX_i \ID(H)\} \cr
&=& t_n^{-4}n \E[\{ \bbeta_0^T (X^TX)^{-1} \bfX_1\}^2 \{ \bfX_1^T (X^TX)^{-2} \bfX_1  \}^2 \ID(H)] \cr
&& + t_n^{-4}n(n-1)  \E\{ \bbeta_0^T (X^TX)^{-1} \bfX_1 \bfX_1^T (X^TX)^{-2} \bfX_1
 \bbeta_0^T (X^TX)^{-1} \bfX_2 \bfX_2^T (X^TX)^{-2} \bfX_2 \ID(H)\}   \cr
&\lesssim& t_n^{-4} O(\|\bbeta_0\|_2^2 p_n^2n^{-5})
+ t_n^{-4} p_n^2/n^2  \E\{ \bbeta_0^T (X^TX)^{-1} \bfX_1 \bbeta_0^T (X^TX)^{-1} \bfX_2 \ID(H) \} + o(1) \cr
&\lesssim& o(1) + t_n^{-4} p_n^2/n^2 \E\{ \bbeta_0^T (X_{(2)}^TX_{(2)})^{-1} \bfX_1 \bbeta_0^T (X_{(1)}^TX_{(1)})^{-1} \bfX_2 \ID(H) \} \cr
&=& o(1) + t_n^{-4} p_n^2/n^2 \E\Big\{\frac{ \bbeta_0^T (X_{(1,2)}^TX_{(1,2)})^{-1} \bfX_1 \bbeta_0^T (X_{(1)}^TX_{(1)})^{-1} \bfX_2 \ID(H_{(1)})}{ 1+\bfX_1^T (X_{(1,2)}^TX_{(1,2)})^{-1} \bfX_1} \Big\} \cr
&\lesssim&o(1) + t_n^{-4} p_n^2/n^2 \frac{\E \{ \bbeta_0^T (X_{(1,2)}^TX_{(1,2)})^{-1} \bfX_1 \bbeta_0^T (X_{(1)}^TX_{(1)})^{-1} \bfX_2 \ID(H_{(1)})\} }{ 1+\E\{\bfX_1^T (X_{(1,2)}^TX_{(1,2)})^{-1} \bfX_1\}} \cr
&=&o(1).
\end{eqnarray*}
The intuition for the last inequality is due to \eqref{S.6.added1} that we can replace $\bfX_1^T (X_{(1,2)}^TX_{(1,2)})^{-1} \bfX_1$
with $\E\{\bfX_1^T (X_{(1,2)}^TX_{(1,2)})^{-1} \bfX_1\}$.

Similarly arguments can be applied to show that $ \sum_{j=1}^n c_2M_2(j,j) v_j /t_n^2\ID(H) = o_{L^2}(1)$.
\endpf

\begin{lemma}\label{Lemma-S12} With the notations and conditions in the proof of part (a) of Theorem 1 with $\lim_{n\to\infty}p_n=\infty$, then $\sum_j \E\{ v_j^2 \ID(K)\} =\bbeta_0^T \E\{(X^TX)^{-1} \ID(K)\} \bbeta_0 =\Omega( \|\bbeta_0\|_2^2 / n)$ and
$$  \hat\bbeta^T (X^TX)^{-1} \hat\bbeta - \sigma_{\epsilon}^2  \E\tr\{ (X^TX)^{-2} \ID(K)\} -  \sum_j \E\{ v_j^2 \ID(K)\}   = o_{\pr}(t_n^2).$$
Consequently,
$$   \hat\bbeta^T (X^TX)^{-1} \hat\bbeta -  \hat\sigma_{\epsilon}^2   \tr\{ (X^TX)^{-2}  \} -   \bfbeta_0^T  \E\{ (X^T X)^{-1} \ID(K)\} \bfbeta_0  = o_{\pr}(t_n^2).$$
\end{lemma}
\pf
Note
$ \sum_j  v_j^2 =  \sum_{j}  \{\bfbeta_0^T (X^T X)^{-1} \bfX_j\}^2 =  \bfbeta_0^T (X^T X)^{-1} \bfbeta_0$,
and
\begin{eqnarray*}
 \hat\bbeta^T (X^TX)^{-1} \hat\bbeta
  =\bbeta_0^T (X^TX)^{-1} \bbeta_0 + 2 \bbeta_0^T (X^TX)^{-2} X^T\beps + \beps^T X(X^TX)^{-3} X^T\beps.
\end{eqnarray*}
First, we show that $\bbeta_0^T (X^TX)^{-2} X^T\beps \ID(K)/t_n^2= o_{L^2}(1)$.  Given event $H$,
\begin{eqnarray*}
&&\E \{\bbeta_0^T (X^TX)^{-2} X^T\beps\}^2 =  \E \{\bbeta_0^T (X^TX)^{-2} X^T\beps \beps^T X (X^TX)^{-2}\bbeta_0\} \cr
&=& \sigma_{\epsilon}^2 \E \{\bbeta_0^T (X^TX)^{-2} X^T X (X^TX)^{-2}\bbeta_0\} = \sigma_{\epsilon}^2 \E \{\bbeta_0^T (X^TX)^{-3} \bbeta_0\} \cr
&=& O(\|\bbeta_0\|_2^2/n^3) = o(t_n^4),
\end{eqnarray*}
by noting that $t_n^4 \gtrsim \|\bbeta_0\|_2^4 /n^2 + p_n^2 / n^4 \gtrsim \|\bbeta_0\|_2^2 p_n /n^3$.

Next, we show that $\var\{\beps^T X(X^TX)^{-3} X^T\beps \ID(K) /t_n^2\} = o(1)$. Note
\begin{eqnarray*}
\E \{\beps^T X(X^TX)^{-3} X^T\beps/t_n^2\} = 1/t_n^2 \sigma_{\epsilon}^2 \E\tr\{ X(X^TX)^{-3} X^T\} = 1/t_n^2 \sigma_{\epsilon}^2 \E\tr\{ (X^TX)^{-2} \}.
\end{eqnarray*}
Denoting $B = X(X^TX)^{-3} X^T$, given event $H$,
\begin{eqnarray*}
&&\E \{ \beps^T X(X^TX)^{-3} X^T\beps/t_n^2\}^2 = t_n^{-4} \sum_{ijkh}\E(\epsilon_i\epsilon_j\epsilon_k\epsilon_h B_{ij} B_{kh} ) \cr
&=& t_n^{-4} \sum_{i}\E(\epsilon_i^4 B_{ii} ^2 ) + t_n^{-4} \sum_{i\neq k}\E(\epsilon_i^2\epsilon_k^2 B_{ii} B_{kk} ) +2 t_n^{-4} \sum_{i\neq j}\E(\epsilon_i^2\epsilon_j^2 B_{ij} ^2 )\cr
&=& t_n^{-4}O( n p_n^2 /n^6)
+ t_n^{-4}\sigma_{\epsilon}^4 \E \Big(\sum_{i } B_{ii} \Big)^2 +2 t_n^{-4}\sigma_{\epsilon}^4 \E\{\tr( B ^2 )\} \cr
&=&   t_n^{-4} \sigma_{\epsilon}^4 \E[ \tr\{(X^TX)^{-2}\} ]^2 + o(1).
\end{eqnarray*}

By Lemma 2, $\var\{\beps^T X(X^TX)^{-3} X^T\beps \ID(K)t_n^{-2} \} = o(1)$, implying that $\beps^T X(X^TX)^{-3} X^T\beps - \sigma_{\epsilon}^2 \E\tr\{ (X^TX)^{-2} \ID(K)\} = o_{\pr}(t_n^{2})$. From Lemma 2, we have $\var\{n\bbeta_0^T (X^TX)^{-1} \bbeta_0 / \|\bbeta_0\|_2^2\ID(K)\} = o(1)$.
Hence,
\begin{eqnarray*}&&\hat\bbeta^T (X^TX)^{-1} \hat\bbeta - \bbeta_0^T \E\{(X^TX)^{-1} \ID(K)\} \bbeta_0-  \sigma_{\epsilon}^2 \E\tr\{ (X^TX)^{-2} \ID(K) \}= o_{\pr}(t_n^{2}).
\end{eqnarray*}
\endpf
\begin{lemma}\label{Lemma-S13} With the notations and conditions in the proof of part (a) of Theorem 1 with $\lim_{n\to\infty}p_n=\infty$, $ \E [\sum_j  \{c_1 M_1(j,j) + c_2M_2(j,j)\}^2 \ID(K) ]=O( p_n^2/n^3)$. Additionally, if $c_1 = 1$ and $c_2 = - \E\tr\{(X^TX)^{-1}\ID(K)\} $, then $$ \E\Big[\sum_j  \{c_1 M_1(j,j) + c_2M_2(j,j)\}^2 \ID(K)\Big] = o(t_n^2).$$
\end{lemma}
\pf
 Note
\begin{eqnarray*}
&& \E\Big[\sum_j  \{c_1 M_1(j,j) + c_2M_2(j,j)\}^2  \ID(K)\Big]
=  n  \E[\{c_1 M_1(1,1) + c_2M_2(1,1)\}^2\ID(K)] \cr
&=& n  \E [c_1 \bfX_1^T (X^TX)^{-2} \bfX_1 \ID(K)+ c_2 \{1 - \bfX_1^T(X^TX)^{-1} \bfX_1\}/(n-p_n)\ID(K)]^2 \cr
&\equiv& n\E( \III^2).
\end{eqnarray*}
It's easy to see that $n\E( \III^2)=O( p_n^2/n^3) $.
Next, it suffices to show that for $c_1 = 1$ and $c_2 = - \E\tr\{(X^TX)^{-1}\ID(K)\} $,
\begin{eqnarray}
 n t_n^{-2} \var (\III) &=& o(1),\qquad \label{S.6.14} \\
\sqrt n t_n^{-1} \E (\III) &=& o(1). \quad\label{S.6.15}
\end{eqnarray}
From
\begin{eqnarray*}
\bfX_1^T(X^T X)^{-1}\bfX_1 &=& \bfX_1^T(X_{(1)}^T X_{(1)})^{-1}\bfX_1 /\{1 + \bfX_1^T (X_{(1)}^T X_{(1)})^{-1}\bfX_1\}\cr
 \bfX_1^T(X^T X)^{-2}\bfX_1 &=& \bfX_1^T(X_{(1)}^T X_{(1)})^{-2}\bfX_1 /\{1 + \bfX_1^T (X_{(1)}^T X_{(1)})^{-1}\bfX_1\}^2,
\end{eqnarray*}
and Lemma \ref{Lemma-S2}, we have
\begin{eqnarray*}
\var\{n/p_n\bfX_1^T(X^T X)^{-1}\bfX_1\ID(K)\} &=& o(1),\cr
 \var\{n^2/p_n\bfX_1^T(X^T X)^{-2}\bfX_1\ID(K)\} &=& o(1).
\end{eqnarray*}
Then, \eqref{S.6.14} is true.
Next, for \eqref{S.6.15}, we first observe that
\begin{eqnarray*}
&&n \E [c_1 \bfX_1^T (X^TX)^{-2} \bfX_1 \ID(K) + c_2 \{1 - \bfX_1^T(X^TX)^{-1} \bfX_1\}/(n-p_n)  ] \cr
&=& \sum_{j=1}^n \E [c_1 \bfX_j^T (X^TX)^{-2} \bfX_j \ID(K)- c_2  \bfX_j^T(X^TX)^{-1} \bfX_j /(n-p_n) ]  +c_2 n/(n-p_n) \cr
&=& \E \tr \Big\{ c_1  (X^TX)^{-2}\sum_{j=1}^n \bfX_j \bfX_j^T \ID(K) - c_2  (X^TX)^{-1} \sum_{j=1}^n \bfX_j \bfX_j^T/(n-p_n)  \Big\}
 + c_2 n/ (n-p_n)   \cr
&=& \E \tr \Big\{c_1  (X^TX)^{-1}\ID(K)- c_2  \bfI_{p_n}/(n-p_n)  \Big\} + c_2 n/(n-p_n)   \cr
&=& c_1 \E \tr\{ (X^TX)^{-1} \ID(K)\} -c_2 p_n  /(n-p_n) + c_2n/(n-p_n)   \cr
&=& \E \tr\{ (X^TX)^{-1}\ID(K)\} +  c_2 = 0.
\end{eqnarray*}
Hence, from  Lemma 1 that $\pr(\bar K) = o(n^{-\ell})$ for any $\ell  \in \mathbb{N}$,  we have
\begin{eqnarray*}
 |\E (\III)|
&=&  |\E [ c_2 \{1 - \bfX_1^T(X^TX)^{-1} \bfX_1\}/(n-p_n) \ID(\bar K) ]  | \cr
&\leq&  \big( \E [ c_2 \{1 - \bfX_1^T(X^TX)^{-1} \bfX_1\}/(n-p_n)]^2 \pr(\bar K)  \big)^{1/2} \cr
&=& o(n^{-\ell/2}) = o(t_n / \sqrt n),
\end{eqnarray*}
for $\ell$ large enough,
which together with \eqref{S.6.14} completes the proof.
\endpf
\begin{lemma}\label{Lemma-S14} With the notations and conditions in the proof of part (a) of Theorem 1 with $\lim_{n\to\infty}p_n=\infty$,
$$ \sum_{j = 1}^n  \E\Big [\sum_{1 \leq i < j}\{c_1 M_1(i,j) + c_2M_2(i,j)\} \epsilon_i\ID(K) \Big]^2 =O(p_n/n^2).$$ Additionally, if $c_1 = 1$ and $c_2 = -  \E\tr\{(X^TX)^{-1}\ID(K)\}  $, then
\begin{eqnarray*}
&& 2  \sigma_{\epsilon}^4    \big(\tr\{ (X^TX)^{-2}\} +  1/(n-p_n) [\tr\{(X^TX)^{-1}\}]^2   \big) \cr
&&\qquad - \sum_{j = 1}^n  \E\Big [\sum_{1 \leq i < j} 2\{c_1 M_1(i,j) + c_2M_2(i,j)\} \epsilon_i\ID(K) \Big]^2 \sigma_{\epsilon}^2 = o_{\pr}(t_n^2).
 \end{eqnarray*}
\end{lemma}
\pf
We know
\begin{eqnarray*}
&& \sum_{j = 1}^n  \E\Big[\sum_{1 \leq i < j}\{c_1 M_1(i,j) + c_2M_2(i,j)\} \epsilon_i\ID(K)\Big]^2 \cr
&=& \sigma_{\epsilon}^2  \sum_{j = 1}^n  \sum_{1 \leq i < j} \E[\{c_1 M_1(i,j) + c_2M_2(i,j)\} ^2\ID(K)] \cr
&=& \sigma_{\epsilon}^2 \sum_{j = 1}^n  \sum_{1 \leq i < j} \E[\{c_1 \bfX_i^T(X^TX)^{-2} \bfX_j  - c_2\bfX_i^T(X^TX)^{-1} \bfX_j/(n-p_n)\} ^2 \ID(K)]\cr
&=& \frac 1 2 \sigma_{\epsilon}^2 \sum_{j = 1}^n  \sum_{ i \neq j} \E[\{c_1^2 \bfX_j^T(X^TX)^{-2} \bfX_i \bfX_i^T(X^TX)^{-2} \bfX_j \cr
&&- 2c_1 c_2\bfX_j^T(X^TX)^{-2} \bfX_i \bfX_i^T(X^TX)^{-1} \bfX_j/(n-p_n) \cr
&& + c_2^2 \bfX_j^T(X^TX)^{-1} \bfX_i \bfX_i^T(X^TX)^{-1} \bfX_j/(n-p_n)^2\}\ID(K)] \cr
 &=& \frac 1 2 \sigma_{\epsilon}^2 \sum_{j = 1}^n  \E[\{c_1^2 \bfX_j^T(X^TX)^{-3} \bfX_j - 2c_1 c_2 \bfX_j^T(X^TX)^{-2}\bfX_j/(n-p_n)\cr
 &&
 + c_2^2  \bfX_j^T(X^TX)^{-1}  \bfX_j/(n-p_n)^2\}\ID(K)] \cr
 && - \frac 1 2 \sigma_{\epsilon}^2 \sum_{j = 1}^n   \E\{c_1^2 \bfX_j^T(X^TX)^{-2} \bfX_j \bfX_j^T(X^TX)^{-2} \bfX_j \cr
&&- 2c_1 c_2\bfX_j^T(X^TX)^{-2} \bfX_j \bfX_j^T(X^TX)^{-1} \bfX_j/(n-p_n) \cr
&& + c_2^2 \bfX_j^T(X^TX)^{-1} \bfX_j \bfX_j^T(X^TX)^{-1} \bfX_j/(n-p_n)^2\}\ID(K)] \cr
&=& \frac 1 2 \sigma_{\epsilon}^2   \E[ \tr\{c_1^2 (X^TX)^{-2} - 2c_1 c_2 (X^TX)^{-1}/(n-p_n) + c_2^2 \bfI_{p_n}/(n-p_n)^2 \}\ID(K)] \cr
 && - \frac 1 2 \sigma_{\epsilon}^2 \sum_{j = 1}^n   \E[\{c_1 \bfX_j^T(X^TX)^{-2} \bfX_j -  c_2 \bfX_j^T(X^TX)^{-1} \bfX_j/(n-p_n)\}^2\ID(K)]\cr
&\equiv& \frac 1 2 \sigma_{\epsilon}^2 \III_1 - \frac 1 2 \sigma_{\epsilon}^2  \III_2 =O( p_n/n^2).
\end{eqnarray*}
From Lemma \ref{Lemma-S13} and \eqref{S.6.15}, for $c_1 = 1$ and $c_2 = - \E\tr\{(X^TX)^{-1}\ID(K)\} $,
\begin{eqnarray*}
&& \sum_{j = 1}^n   \E \big([c_1 \bfX_j^T(X^TX)^{-2} \bfX_j  + c_2 \{1 - \bfX_j^T(X^TX)^{-1} \bfX_j\}/(n-p_n) ]^2\ID(K) \big) = o(t_n^2), \cr
&& \sum_{j = 1}^n \E [c_1 \bfX_j^T (X^TX)^{-2} \bfX_j \ID(K)   + c_2 \{1 - \bfX_j^T(X^TX)^{-1} \bfX_j\}/(n-p_n) ] = 0.
\end{eqnarray*}
Then,
\begin{eqnarray*}
  \III_2
&=&  \sum_{j = 1}^n \E [ \{c_1 \bfX_j^T(X^TX)^{-2} \bfX_j+ c_2/(n-p_n)   -  c_2 \bfX_j^T(X^TX)^{-1} \bfX_j/(n-p_n)\}^2 \ID(K) ] \cr
&&- 2 c_2/(n-p_n) \sum_{j = 1}^n  \E [\{c_1 \bfX_j^T(X^TX)^{-2} \bfX_j+ c_2/(n-p_n) \cr
&&-  c_2 \bfX_j^T(X^TX)^{-1} \bfX_j/(n-p_n)\} \ID(K)] + c_2^2n/(n-p_n)^2 \pr(K)\cr
&=& o (t_n^2) + c_2^2n/(n-p_n)^2.
 \end{eqnarray*}
 From Lemma 2,
 \begin{eqnarray*}
&& \frac 1 2 \sigma_{\epsilon}^2 \III_1 - \frac 1 2 \sigma_{\epsilon}^2  \III_2 \cr
&=& 1 / 2   \sigma_{\epsilon}^2    \E[ \tr\{c_1^2 (X^TX)^{-2} - 2c_1 c_2 (X^TX)^{-1}/(n-p_n) \cr
&& + c_2^2 \bfI_{p_n}/(n-p_n)^2 \} \ID(K)]
 -  1/ 2   \sigma_{\epsilon}^2 c_2^2 n/(n-p_n)^2+ o(t_n^2) \cr
 &=& 1 / 2   \sigma_{\epsilon}^2    [ \E \tr\{ (X^TX)^{-2}\ID(K)\} - 2 c_2 \E\tr\{(X^TX)^{-1}\ID(K)\}/(n-p_n) \cr
&& + c_2^2 {p_n}/(n-p_n)^2  ]-  1/ 2   \sigma_{\epsilon}^2 c_2^2 n/(n-p_n)^2 + o(t_n^2) \cr
 &=& 1 / 2   \sigma_{\epsilon}^2    \big( \E\tr\{ (X^TX)^{-2}\ID(K)\} + 2 /(n-p_n) [\E\tr\{(X^TX)^{-1}\ID(K)\}]^2 \cr
&& - 1/(n-p_n) [\E\tr\{(X^TX)^{-1}\ID(K)\}]^2   \big) + o(t_n^2) \cr
 &=& 1 / 2   \sigma_{\epsilon}^2    \big(\E\tr\{ (X^TX)^{-2}\ID(K)\} +  1/(n-p_n) [\E\tr\{(X^TX)^{-1}\ID(K)\}]^2
  \big)  + o (t_n^2)\cr
  &=& 1 / 2   \sigma_{\epsilon}^2    \big(\tr\{ (X^TX)^{-2}\} +  1/(n-p_n) [\tr\{(X^TX)^{-1}\}]^2
  \big)  + o _{\pr}(t_n^2).
 \end{eqnarray*}
 The proof is completed.
\endpf
\begin{lemma}\label{Lemma-S15} Under the conditions of part (a) of Theorem 1, we have
$ \var(1/n \sum_{i = 1}^n  \hat\epsilon_i^4) = o(1) $ and $(1-p_n/n)^{-4} \{1/n \sum_{i = 1}^n  \hat\epsilon_i^4-3 \hat \sigma_{\epsilon}^4 (p_n/n) (1-p_n/n)^2(2-p_n/n)\} - \nu_4 = o_{\pr}(1)$.
\end{lemma}
\pf We know
$ 1/n \sum_{i = 1}^n  \hat\epsilon_i^4 = 1/n \sum_{i = 1}^n (\epsilon_i - \bfX_i^T (X^TX)^{-1} X^T\beps)^4 =  1/n \sum_{i = 1}^n (\epsilon_i - \sum_j N_{ij} \epsilon_j)^4$ where $N = X(X^TX)^{-1} X^T$. Then, $ 1/n \sum_{i = 1}^n  \E(\hat\epsilon_i^4) =\E(\hat\epsilon_1^4) = \E(\epsilon_1 - \sum_j N_{1j} \epsilon_j)^4$. We have,
\begin{eqnarray*}
& &  \E \Big (\epsilon_1 - \sum_j N_{1j} \epsilon_j \Big)^4 \cr
&=&  \E(\epsilon_1^4) - 4 \E\Big\{ \epsilon_1 \Big(\sum_j N_{1j} \epsilon_j\Big)^3 \Big\} + 6 \E \Big\{\epsilon_1^2 \Big(\sum_j N_{1j} \epsilon_j\Big)^2\Big\}
- 4 \E\Big\{ \epsilon_1^3 \Big(\sum_j N_{1j} \epsilon_j\Big) \Big\} +\E \Big(\sum_j N_{1j} \epsilon_j\Big) ^4,\qquad
\end{eqnarray*}
with $\E(\epsilon_1^4) = \nu_4$,
\begin{eqnarray*}
 \E\Big\{ \epsilon_1 \Big(\sum_j N_{1j} \epsilon_j\Big)^3 \Big\} &=& \sum_j \sum_k \sum_h \E\{  N_{1j} N_{1k} N_{1h}\epsilon_1 \epsilon_j\epsilon_k \epsilon_h \}
 = \E ( N_{11}^3) \nu_4 + 3 \sigma_{\epsilon}^4  \sum_{j, j\neq 1} \E(  N_{11} N_{1j}^2 )\cr
&=& \E ( N_{11}^3) \nu_4 + 3 \sigma_{\epsilon}^4  (n-1) \E(  N_{11} N_{12}^2 ),\cr
 \E \Big\{\epsilon_1^2 \Big(\sum_j N_{1j} \epsilon_j\Big)^2\Big\}
&=& \sum_j \sum_k \E (N_{1j} N_{1k} \epsilon_1^2  \epsilon_j \epsilon_k)
 =  \E (N_{11}^2) \nu_4 +  \sigma_{\epsilon} ^4 \sum_{j, j\neq 1} \E (N_{1j}^2)\cr
&=& \E (N_{11}^2) \nu_4 +  \sigma_{\epsilon} ^4 (n-1) \E (N_{12}^2),\cr
\E\Big\{ \epsilon_1^3 \Big(\sum_j N_{1j} \epsilon_j\Big) \Big\} &=& \E( \epsilon_1^4  N_{11}) = \E(N_{11}) \nu_4,\cr
 \E \Big(\sum_j N_{1j} \epsilon_j\Big) ^4
&=&\sum_j\sum_k\sum_h\sum_{\ell}\E ( N_{1j}N_{1k}N_{1h}N_{1\ell} \epsilon_j \epsilon_k \epsilon_h \epsilon_{\ell}) \cr
&=& \sum_j \E ( N_{1j}^4) \nu_4 + 3 \sigma_{\epsilon}^4 \sum_j\sum_{k, k\neq j}\E ( N_{1j}^2N_{1k}^2) \cr
&=& \E ( N_{11}^4) \nu_4 + \sum_{j\neq 1} \E ( N_{1j}^4) \nu_4+ 6 \sigma_{\epsilon}^4 \sum_{k, k\neq 1}\E ( N_{11}^2N_{1k}^2) \cr
&&\quad  + 3 \sigma_{\epsilon}^4 \sum_{j,j\neq 1}\sum_{k, k\neq j,k\neq 1}\E ( N_{1j}^2N_{1k}^2)\cr
&=& \E ( N_{11}^4) \nu_4 + 6 \sigma_{\epsilon}^4 (n-1)\E ( N_{11}^2N_{12}^2) \cr
&&\quad + 3 \sigma_{\epsilon}^4(n-1)(n-2)\E ( N_{12}^2N_{13}^2) +o(1).
\end{eqnarray*}
From Lemma \ref{Lemma-S2},  $N_{11} \conP \tau$, and $\E(N_{11}^k) = O(1)$, which implies that $\E(N_{11}^k) \to \tau^k$ due to uniform integrability, for any $k\in \mathbb{N}$. Also,
\begin{eqnarray*}
&& (n-1) \E (N_{12}^2) = \sum_{j=2}^n \E (N_{1j}^2) \cr
&=& \sum_{j=1}^n \E \{\bfX_1^T (X^TX)^{-1} \bfX_j\bfX_j^T (X^TX)^{-1}\bfX_1 \} - \E [\{\bfX_1^T (X^TX)^{-1} \bfX_1 \}^2 ]\cr
&=& \E \{\bfX_1^T   (X^TX)^{-1}\bfX_1 \}- \E [\{\bfX_1^T (X^TX)^{-1} \bfX_1 \}^2 ]= (1-\tau) \tau +o(1),
 \end{eqnarray*}
 and
 \begin{eqnarray*}
&& (n-1)(n-2) \E (N_{12}^2 N_{13}^2 ) = \sum_{i=2}^n \sum_{j=2, j\neq i}^n \E (N_{1i}^2 N_{1j}^2) \cr
&=& \sum_{i=1}^n \sum_{j=1}^n  \E \{\bfX_1^T (X^TX)^{-1} \bfX_j\bfX_j^T (X^TX)^{-1}\bfX_1 \bfX_1^T (X^TX)^{-1} \bfX_i\bfX_i^T (X^TX)^{-1}\bfX_1\} \cr
&& - \sum_{j=1}^n \E \{\bfX_1^T (X^TX)^{-1} \bfX_1\bfX_1^T (X^TX)^{-1}\bfX_1
 \bfX_1^T (X^TX)^{-1} \bfX_j\bfX_j^T (X^TX)^{-1}\bfX_1\} \cr
&& - \sum_{i=1}^n \E \{\bfX_1^T (X^TX)^{-1} \bfX_i\bfX_i^T (X^TX)^{-1}\bfX_1
\bfX_1^T (X^TX)^{-1} \bfX_1\bfX_1^T (X^TX)^{-1}\bfX_1\} \cr
&& - \sum_{i=1}^n \E \{\bfX_1^T (X^TX)^{-1} \bfX_i\bfX_i^T (X^TX)^{-1}\bfX_1
\bfX_1^T (X^TX)^{-1} \bfX_i\bfX_i^T (X^TX)^{-1}\bfX_1\}  \cr
&&+2 \E [\{\bfX_1^T (X^TX)^{-1} \bfX_1 \}^4 ]\cr
&=& \E [\{\bfX_1^T (X^TX)^{-1} \bfX_1 \}^2 ] - 2\E [\{\bfX_1^T (X^TX)^{-1} \bfX_1 \}^3 ] \cr
&&-\E [\{\bfX_1^T (X^TX)^{-1} \bfX_1 \}^4 ] + O(p_n^2 / n^3) +2 \E [\{\bfX_1^T (X^TX)^{-1} \bfX_1 \}^4 ]\cr
&=& \tau^2 - 2\tau^3 + \tau^4 + o(1) = \tau^2 (1- \tau)^2 + o(1).
 \end{eqnarray*}
Therefore, $ 1/n \sum_{i = 1}^n  \E(\hat\epsilon_i^4) = \nu_4 - 4 [\tau^3 \nu_4 + 3 \sigma_{\epsilon}^4 \tau^2  (1-\tau)] + 6[\tau^2 \nu_4 +  \sigma_{\epsilon} ^4  \tau (1-\tau) ] - 4 \tau \nu_4 + \tau^4 \nu_4 + 6 \sigma_{\epsilon}^4 \tau^3(1-\tau)  + 3 \sigma_{\epsilon}^4\tau^2(1-\tau)^2 +o(1) = (1  - \tau )^4\nu_4 +3\sigma_{\epsilon} ^4 \tau (1-\tau)^2 (2-\tau) +o(1)$.

By similar arguments and using Condition A2, we have $\E\{(\epsilon_i - \bfX_i^T (X^TX)^{-1} X^T\beps)^8\} < \infty$, and hence,
 \begin{eqnarray*}
  \E\Big ( 1/n \sum_{i = 1}^n  \hat\epsilon_i^4\Big)^2
 &=& 1/n^2 \sum_{i = 1}^n \sum_{k = 1}^n \E\{(\epsilon_i - \bfX_i^T (X^TX)^{-1} X^T\beps)^4 (\epsilon_k - \bfX_k^T (X^TX)^{-1} X^T\beps)^4\} \cr
 &=& \E\{(\epsilon_1 - \bfX_1^T (X^TX)^{-1} X^T\beps)^4 (\epsilon_2 - \bfX_2^T (X^TX)^{-1} X^T\beps)^4\} + o(1).
 \end{eqnarray*}
 It suffices to show that the covariance of $(\epsilon_1 - \bfX_1^T (X^TX)^{-1} X^T\beps)^4 $ and $(\epsilon_2 - \bfX_2^T (X^TX)^{-1} X^T\beps)^4$ is o(1).
Note
\begin{eqnarray}
&&\{\epsilon_1 - \bfX_1^T (X^TX)^{-1} X^T\beps\}^4  \cr
&=& \epsilon_1^4 - 4 \epsilon_1^3 \Big(\sum_j N_{1j} \epsilon_j\Big) + 6 \epsilon_1^2 \Big(\sum_j N_{1j} \epsilon_j\Big)^2 - 4  \epsilon_1 \Big(\sum_j N_{1j} \epsilon_j\Big)^3 +\Big(\sum_j N_{1j} \epsilon_j\Big) ^4, \label{S.6.16}\\
&&\{\epsilon_2 - \bfX_2^T (X^TX)^{-1} X^T\beps\}^4  \cr
&=& \epsilon_2^4 - 4 \epsilon_2^3 \Big(\sum_j N_{2j} \epsilon_j\Big) + 6 \epsilon_2^2 \Big(\sum_j N_{2j} \epsilon_j\Big)^2 - 4  \epsilon_2 \Big(\sum_j N_{2j} \epsilon_j\Big)^3 +\Big(\sum_j N_{2j} \epsilon_j\Big) ^4.\label{S.6.17}
\end{eqnarray}
In the following, we aim to show that the covariance of each term in \eqref{S.6.16} and each term in \eqref{S.6.17} is $o(1)$. We mainly use the results in Lemmas \ref{Lemma-S2} and \ref{Lemma-S3}.
For $\epsilon_1^4 $ and $\epsilon_2^4$, we know
 $\E(\epsilon_1^4 \epsilon_2^4) = \E(\epsilon_1^4) \E(\epsilon_2^4)$.
 For $\epsilon_2^4 $ and $\epsilon_1^3  (\sum_j N_{1j} \epsilon_j )$, we have
\begin{eqnarray*}
 \E \Big\{\epsilon_2^4 \epsilon_1^3 \Big(\sum_j N_{1j} \epsilon_j\Big)\Big\} = \E(\epsilon_2^4 \epsilon_1^4 N_{11} ) + \E(\epsilon_2^5 \epsilon_1^3 N_{12})
 =\E( \epsilon_2^4) \E\Big\{ \epsilon_1^3 \Big(\sum_j N_{1j} \epsilon_j\Big)\Big\} +o(1).
\end{eqnarray*}
For $\epsilon_2^4$ and $ \epsilon_1^2  (\sum_j N_{1j} \epsilon_j )^2$,
\begin{eqnarray*}
&&\E\Big\{\epsilon_2^4 \epsilon_1^2 \Big(\sum_j N_{1j} \epsilon_j\Big)^2\Big\} = \E\Big(\epsilon_2^4 \epsilon_1^2 \sum_j N_{1j}^2 \epsilon_j^2\Big) + 2 \E(\epsilon_2^5 \epsilon_1^3  N_{11} N_{12} )\cr
&=&\E(\epsilon_2^4 \epsilon_1^4  N_{11}^2 ) +\E(\epsilon_2^6 \epsilon_1^2  N_{12}^2 ) + \E\Big(\epsilon_2^4 \epsilon_1^2 \sum_{j=3}^n N_{1j}^2 \epsilon_j^2\Big) +o(1)\cr
&=&\nu_4^2 \E(N_{11}^2 ) + \nu_4\sigma_{\epsilon}^4 (n - 2) \E(N_{12}^2) +o(1) \cr
&=& \E(\epsilon_2^4 )\E\Big\{\epsilon_1^2 \Big(\sum_j N_{1j} \epsilon_j\Big)^2\Big\}+o(1).
\end{eqnarray*}
For $\epsilon_2^4 $ and $\epsilon_1  (\sum_j N_{1j} \epsilon_j )^3$,
\begin{eqnarray*}
&&\E\Big\{\epsilon_2^4 \epsilon_1 \Big(\sum_j N_{1j} \epsilon_j\Big)^3\Big\} \cr
& =& \E(\epsilon_2^4 \epsilon_1^4  N_{11} ^3)
+3 \E(\epsilon_2^5 \epsilon_1^3  N_{11}^2 N_{12})+3 \E(\epsilon_2^6 \epsilon_1^2  N_{11} N_{12}^2)
+3 \E\Big\{\epsilon_2^4 \epsilon_1^2  N_{11}\Big(\sum_{j=3}^n N_{1j}^2 \epsilon_j^2\Big)\Big\} \cr
&=& \nu_4^2 \E(N_{11} ^3)  + 3(n-2) \nu_4 \sigma_{\epsilon}^4 \E(N_{11} N_{12}^2) +o(1) \cr
&=& \E(\epsilon_2^4 )\E\Big\{\epsilon_1 \Big(\sum_j N_{1j} \epsilon_j\Big)^3\Big\} +o(1).
\end{eqnarray*}
For $\epsilon_2^4$ and $ (\sum_j N_{1j} \epsilon_j ) ^4$,
\begin{eqnarray*}
&&\E\Big\{\epsilon_2^4\Big(\sum_j N_{1j} \epsilon_j\Big) ^4 \Big\} \cr
&=& \E\Big\{\epsilon_2^4\Big(\sum_j N_{1j}^4 \epsilon_j^4\Big)  \Big\} + 4 \E \Big\{\epsilon_2^5 N_{12}\Big(\sum_j N_{1j}^3 \epsilon_j^3\Big) \Big\}
+ 3 \E\Big\{\epsilon_2^4\Big(\sum_j\sum_{k\neq j} N_{1j}^2  N_{1k}^2 \epsilon_j^2\epsilon_k^2\Big) \Big\}\cr
&=& \nu_4^2 \E(N_{11}^4)+\E(\epsilon_2^8) \E(N_{12}^4)+\nu_4^2 (n-2) \E(N_{12}^4)
+4 \E (\epsilon_2^5\epsilon_1^3 N_{12} N_{11}^3 )
\cr
&&+4 \E (\epsilon_2^8 N_{12}^4 )+4 (n-2)\E (\epsilon_2^5 \epsilon_3^3 N_{12} N_{13}^3 )
+6 \E\Big\{\epsilon_2^4 N_{12}^2 \epsilon_2^2 \Big(\sum_{k\neq 2} N_{1k}^2 \epsilon_k^2\Big) \Big\}\cr
&&+3 \E\Big\{\epsilon_2^4\Big(\sum_{j\neq 2}\sum_{k\neq j, k\neq 2} N_{1j}^2  N_{1k}^2 \epsilon_j^2\epsilon_k^2\Big) \Big\} \cr
&=& \nu_4^2 \E(N_{11}^4)+6 \E\Big\{\epsilon_2^4N_{11}^2 \epsilon_1^2\Big(\sum_{k=3}^n   N_{1k}^2 \epsilon_k^2\Big) \Big\}
+3 \E\Big\{\epsilon_2^4\Big(\sum_{j=3}^n\sum_{k=3, k\neq j}^n N_{1j}^2  N_{1k}^2 \epsilon_j^2\epsilon_k^2\Big) \Big\} + o(1)\cr
&=&\nu_4^2 \E(N_{11}^4)+6\nu_4 \sigma_{\epsilon}^4 (n-2)\E(N_{11}^2   N_{13}^2)
+3 \nu_4 \sigma_{\epsilon}^4 (n-2)(n-3) \E(N_{13}^2  N_{14}^2 ) +o(1)\cr
&=& \E(\epsilon_2^4) \E\Big\{\Big(\sum_j N_{1j} \epsilon_j\Big) ^4 \Big\} + o(1).
\end{eqnarray*}
For $\epsilon_2^3  (\sum_j N_{2j} \epsilon_j ) $ and $\epsilon_1^3  (\sum_k N_{1k} \epsilon_k )$,
\begin{eqnarray*}
&&\E\Big\{\epsilon_2^3 \Big(\sum_j N_{2j} \epsilon_j\Big) \epsilon_1^3 \Big(\sum_k N_{1k} \epsilon_k\Big)\Big\} \cr
&=& \E\{\epsilon_2^4 \epsilon_1^4 ( N_{12}^2 + N_{22}N_{11})\}  + \E\Big\{\epsilon_2^3 \epsilon_1^3 \Big(\sum_j N_{1j} N_{2j} \epsilon_j^2\Big) \Big\} \cr
&=& \nu_4^2 \E(N_{22}N_{11}) + \E(\epsilon_2^3 \epsilon_1^3 \epsilon_3^2) (n-2) \E(N_{13} N_{23}) +o(1) \cr
&=&\nu_4^2 \E(N_{22}N_{11}) +o(1) = \E\Big\{\epsilon_2^3 \Big(\sum_j N_{2j} \epsilon_j\Big)\Big\}\E\Big\{ \epsilon_1^3 \Big(\sum_k N_{1k} \epsilon_k\Big)\Big\} + o(1).
\end{eqnarray*}
For $\epsilon_2^3  (\sum_j N_{2j} \epsilon_j ) $ and $ \epsilon_1^2 (\sum_k N_{1k} \epsilon_k )^2$,
\begin{eqnarray*}
&&\E\Big\{\epsilon_2^3 \Big(\sum_j N_{2j} \epsilon_j\Big)  \epsilon_1^2 \Big(\sum_k N_{1k} \epsilon_k\Big)^2\Big\}\cr
 &=& \E\Big\{\epsilon_2^3\epsilon_1^2  \Big(\sum_j N_{2j} N_{1j}^2  \epsilon_j^3\Big) \Big\} + \E\{\epsilon_2^5 \epsilon_1^3 (2N_{22} N_{12} N_{11} + N_{12}^3)\} \cr
&&+\E\{\epsilon_2^4 \epsilon_1^4 (N_{11}^2 N_{22} + 2 N_{11} N_{12}^2)\} +\E\Big\{\epsilon_2^3  \epsilon_1^3 N_{12} \Big(\sum_{k=3}^n N_{1k}^2 \epsilon_k^2\Big)\Big\} \cr
&&+2 \E\Big\{\epsilon_2^3  \epsilon_1^3 N_{11} \Big(\sum_{k=3}^n N_{1k}N_{2k} \epsilon_k^2\Big)\Big\}\cr
&&+\E\Big\{\epsilon_2^4  \epsilon_1^2 N_{22} \Big(\sum_{k=3}^n N_{1k}^2 \epsilon_k^2\Big)\Big\} +2 \E\Big\{\epsilon_2^4  \epsilon_1^2 N_{12} \Big(\sum_{k=3}^n N_{1k}N_{2k} \epsilon_k^2\Big)\Big\}\cr
&=& \nu_4^2 \E(N_{11}^2 N_{22}) + \nu_4 \sigma_{\epsilon}^4 (n-2) \E(N_{22} N_{13}^2)+o(1)\cr
&= &\E\Big\{\epsilon_2^3 \Big(\sum_j N_{2j} \epsilon_j\Big)\Big\}\E\Big\{  \epsilon_1^2 \Big(\sum_k N_{1k} \epsilon_k\Big)^2\Big\}+o(1).
\end{eqnarray*}
For $\epsilon_2^3  (\sum_j N_{2j} \epsilon_j) $ and $ \epsilon_1 (\sum_k N_{1k} \epsilon_k )^3$,
\begin{eqnarray*}
&&\E\Big\{\epsilon_2^3 \Big(\sum_j N_{2j} \epsilon_j\Big)  \epsilon_1 \Big(\sum_k N_{1k} \epsilon_k\Big)^3\Big\} \cr&=&\E(\epsilon_2^3  \epsilon_1 ^ 5 N_{21}N_{11}^3) + \E\{\epsilon_2^4 \epsilon_1^4 ( N_{22}  N_{11}^3 + 3 N_{12}^2 N_{11}^2 ) \}  \cr
&& + \E\Big\{\epsilon_2^3 \epsilon_1^2  N_{21}   \Big(\sum_{k=2}^n N_{1k}^3 \epsilon_k^3\Big)\Big\} +
3 \E\Big\{\epsilon_2^3 \epsilon_1^2 N_{11} \Big(\sum_{j=2}^n N_{2j} N_{1j}^2 \epsilon_j^3\Big)  \Big\} \cr
&&+3 \E\Big\{\epsilon_2^4 \epsilon_1^2  N_{22} N_{11}   \Big(\sum_{k=3}^n N_{1k}^2 \epsilon_k^2\Big)\Big\} +
3 \E\Big\{\epsilon_2^4 \epsilon_1^2 N_{12} ^2 \Big(\sum_{k=3}^n N_{1k}^2 \epsilon_k^2\Big)\Big\}\cr
&&+ 6 \E\Big\{\epsilon_2^4 \epsilon_1^2 N_{11} N_{12}  \Big(\sum_{k=3}^n N_{1k} N_{2k} \epsilon_k^2\Big)\Big\} \cr
&& + 3 \E\Big\{\epsilon_2^3 \epsilon_1^3  N_{21} N_{11}    \Big(\sum_{k=2}^n N_{1k}^2 \epsilon_k^2\Big)\Big\}
+ 3 \E\Big\{\epsilon_2^3 \epsilon_1^3  N_{11}^2  \Big(\sum_{k=2}^n N_{1k} N_{2k} \epsilon_k^2\Big)\Big\}\cr
&=& \nu_4^2 \E (N_{22} N_{11}^3) + 3 \nu_4 \sigma_{\epsilon}^4 (n-2)\E(  N_{22} N_{11}   N_{13}^2 ) +o(1) \cr
&= &\E\Big\{\epsilon_2^3 \Big(\sum_j N_{2j} \epsilon_j\Big) \Big\}\E\Big\{ \epsilon_1 \Big(\sum_k N_{1k} \epsilon_k\Big)^3\Big\}+o(1).
\end{eqnarray*}
For $\epsilon_2^3  (\sum_j N_{2j} \epsilon_j ) $ and $  (\sum_k N_{1k} \epsilon_k ) ^4$,
\begin{eqnarray*}
&&\E\Big\{\epsilon_2^3 \Big(\sum_j N_{2j} \epsilon_j\Big) \Big (\sum_k N_{1k} \epsilon_k\Big) ^4 \Big\} \cr&=&
\E\Big\{\epsilon_2^3 \Big(\sum_j N_{2j} N_{1j}^4 \epsilon_j^5\Big)  \Big \}
+ \E\Big\{\epsilon_2^4 N_{22}  \Big(\sum_{k\neq 2} N_{1k}^4 \epsilon_k^4\Big)  \Big\}
+
4 \E\Big\{\epsilon_2^4 N_{12} \Big(\sum_{j\neq 2} N_{2j} N_{1j}^3 \epsilon_j^4\Big) \Big \} \cr
&& + \E\Big\{\epsilon_2^3 \Big(\sum_j\sum_{k\neq j} (6N_{2j}N_{1j}^2  N_{1k}^2 + 4 N_{1j}^3N_{1k}N_{2k}) \epsilon_j^3 \epsilon_k^2\Big) \Big \} \cr
&& + 3 \E\Big\{\epsilon_2^4  N_{22} \Big(\sum_{j\neq 2}\sum_{k\neq j, 2}  N_{1j}^2 N_{1k}^2 \epsilon_j^2  \epsilon_k^2\Big)\Big\}
+ 12 \E\Big\{\epsilon_2^4  N_{12} \Big(\sum_{j\neq 2}\sum_{k\neq j, 2}  N_{1j}^2 N_{1k}N_{2k} \epsilon_j^2  \epsilon_k^2\Big)\Big\} \cr
&=& \nu_4^2 \E(N_{11} ^4N_{22})+ 6 \nu_4 \sigma_{\epsilon}^4 \E\Big\{  N_{22} N_{11}^2 \Big(\sum_{k=3}^n  N_{1k}^2\Big )\Big\} \cr
&&+ 3 \nu_4 \sigma_{\epsilon}^4 (n-2)(n-3) \E(  N_{22}    N_{13}^2N_{14}^2 ) + o(1)\cr
&=&\E\Big\{\epsilon_2^3 \Big(\sum_j N_{2j} \epsilon_j\Big)\Big\}\E\Big\{  \Big(\sum_k N_{1k} \epsilon_k\Big) ^4 \Big\} +o(1).
\end{eqnarray*}
For $\epsilon_2^2  (\sum_j N_{2j} \epsilon_j )^2 $ and $\epsilon_1^2  (\sum_k N_{1k} \epsilon_k )^2$,
\begin{eqnarray*}
&&\E\Big\{\epsilon_2^2 \Big(\sum_j N_{2j} \epsilon_j\Big)^2 \epsilon_1^2 \Big(\sum_k N_{1k} \epsilon_k\Big)^2\Big \} \cr
&=&
\E\Big\{\epsilon_2^2 \epsilon_1^2 \Big(\sum_j N_{2j}^2 N_{1j}^2 \epsilon_j^4\Big) \Big\}
+ 2 \E\Big\{\epsilon_2^3 \epsilon_1^2 N_{22} \Big(\sum_{j \neq 2} N_{2j} N_{1j}^2\epsilon_j^3\Big)\Big \}\cr
&&+ 2 \E\Big\{\epsilon_2^3 \epsilon_1^2 N_{12} \Big(\sum_{j \neq 2} N_{2j}^2 N_{1j}\epsilon_j^3\Big)\Big \} \cr
&&+ 2 \E\Big\{\epsilon_2^2 \epsilon_1^3 N_{12} \Big(\sum_{j \neq 1} N_{2j} N_{1j}^2\epsilon_j^3\Big)\Big \}
+ 2 \E\Big\{\epsilon_2^2 \epsilon_1^3 N_{11} \Big(\sum_{j \neq 1} N_{2j}^2 N_{1j}\epsilon_j^3\Big) \Big\} \cr
&&+\E\Big\{\epsilon_2^3 \epsilon_1^3 N_{21} N_{22}   \Big(\sum_{k=3}^n N_{1k}^2 \epsilon_k^2\Big)\Big \}
+\E\Big\{\epsilon_2^3 \epsilon_1^3 N_{11} N_{12}  \Big (\sum_{k=3}^n N_{2k}^2 \epsilon_k^2\Big)\Big \}\cr
&&+4 \E\Big\{\epsilon_2^3 \epsilon_1^3 N_{11} N_{22}  \Big (\sum_{k=3}^n N_{1k} N_{2k} \epsilon_k^2\Big)\Big \}+4 \E\Big\{\epsilon_2^3 \epsilon_1^3 N_{12}^2  \Big (\sum_{k=3}^n N_{1k} N_{2k} \epsilon_k^2\Big)\Big \} \cr
&&+\E\Big\{\epsilon_2^2 \epsilon_1^2 \Big(\sum_j\sum_{k\neq j} N_{2j}^2  N_{1k}^2 \epsilon_j^2\epsilon_k^2\Big)\Big \}
+2\E\Big\{\epsilon_2^2 \epsilon_1^2 \Big(\sum_j\sum_{k\neq j} N_{1j}N_{2j}  N_{1k}N_{2k}\epsilon_j^2\epsilon_k^2\Big)\Big \} \cr
&=& \nu_4^2\E(  N_{22}^2  N_{11}^2 ) +
 2 \nu_4 \sigma_{\epsilon}^4 (n-2) \E( N_{11}^2  N_{23}^2  ) \cr
 &&\quad +\sigma_{\epsilon}^8 (n-2) (n-3) \E( N_{23}^2  N_{14}^2 ) +o(1)\cr
 &=&\E\Big\{\epsilon_2^2 \Big(\sum_j N_{2j} \epsilon_j\Big)^2\Big \} \E\Big\{\epsilon_1^2 \Big(\sum_k N_{1k} \epsilon_k\Big)^2 \Big\}+o(1).
 \end{eqnarray*}
 For $\epsilon_2^2  (\sum_j N_{2j} \epsilon_j )^2$ and $   \epsilon_1 (\sum_k N_{1k} \epsilon_k )^3$,
 \begin{eqnarray*}
&&\E\Big\{ \epsilon_2^2 \Big(\sum_j N_{2j} \epsilon_j\Big)^2   \epsilon_1 \Big(\sum_k N_{1k} \epsilon_k\Big)^3\Big\} \cr
&=&\E\Big\{ \epsilon_2^2 \epsilon_1^4 N_{11}^3 \Big(\sum_{j=2}^n N_{2j}^2 \epsilon_j^2\Big)\Big\}
+ 3 \E\Big\{ \epsilon_2^2 \epsilon_1 ^2 N_{11} \Big(\sum_{j=2}^n\sum_{k=2, k\neq j}^n N_{2j}^2 N_{1k} ^2\epsilon_j^2\epsilon_k^2\Big)\Big\}+o(1) \cr
&=& \nu_4^2 \E(N_{11}^3 N_{22}^2) + \nu_4 \sigma_{\epsilon}^4 (n-2)\E(N_{11}^3 N_{23}^2) + 3 \nu_4 \sigma_{\epsilon}^4 (n-2)\E(N_{11} N_{22}^2 N_{13}^2)\cr
&& + 3 \nu_4 \sigma_{\epsilon}^4 (n-2)\E(N_{11} N_{12}^2 N_{23}^2) + 3 \sigma_{\epsilon}^8 (n-2)(n-3)\E(N_{11} N_{23}^2 N_{14}^2) +o(1)\cr
&=&\E\Big\{ \epsilon_2^2 \Big(\sum_j N_{2j} \epsilon_j\Big)^2 \Big\} \E\Big\{  \epsilon_1 \Big(\sum_k N_{1k} \epsilon_k\Big)^3\Big\} + o(1).
\end{eqnarray*}
 For $\epsilon_2^2  (\sum_j N_{2j} \epsilon_j )^2 $ and $ (\sum_k N_{1k} \epsilon_k ) ^4$,
 \begin{eqnarray*}
&&\E\Big\{\epsilon_2^2 \Big(\sum_j N_{2j} \epsilon_j\Big)^2 \Big(\sum_k N_{1k} \epsilon_k\Big) ^4\Big\}\cr
&=&\E\Big\{\epsilon_2^2 \Big(\sum_j\sum_{k \neq j} N_{2j}^2N_{1k}^4 \epsilon_j^2  \epsilon_k ^4\Big)\Big\}
+3 \E\Big\{\epsilon_2^2 \Big(\sum_j \sum_{k \neq j} \sum_{h \neq k,j} N_{2j}^2 N_{1k}^2 N_{1h}^2 \epsilon_j^2\epsilon_k ^2 \epsilon_h ^2\Big)\Big\} \cr
&=& \nu_4^2 \E(N_{22}^2N_{11}^4) + \nu_4 \sigma_{\epsilon}^4 (n-2)\E(N_{11}^4 N_{23}^2  )+6\nu_4 \sigma_{\epsilon}^4 (n-2)\E(N_{22}^2 N_{11}^2N_{13}^2 ) \cr
&& + 3 \nu_4 \sigma_{\epsilon}^4 n^2\E( N_{22}^2  N_{13}^2 N_{14}^2 ) +6\sigma_{\epsilon}^8 n^2\E(N_{11}^2 N_{23}^2N_{14}^2 )
+3\sigma_{\epsilon}^8 n^3\E(N_{23}^2 N_{14}^2N_{15}^2 ) +o(1)\cr
&=&\E\Big\{\epsilon_2^2 \Big(\sum_j N_{2j} \epsilon_j\Big)^2\Big\}\E\Big\{ \Big(\sum_k N_{1k} \epsilon_k\Big) ^4\Big\} + o(1).
\end{eqnarray*}
For $\epsilon_2 (\sum_j N_{2j} \epsilon_j )^3$ and $ \epsilon_1  (\sum_k N_{1k} \epsilon_k )^3$,
\begin{eqnarray*}
&&\E\Big\{\epsilon_2 \Big(\sum_j N_{2j} \epsilon_j\Big)^3 \epsilon_1 \Big(\sum_k N_{1k} \epsilon_k\Big)^3 \Big\} \cr
&=&\E(\epsilon_2^4 \epsilon_1^4 N_{11}^3 N_{22}^3)
+3\E\Big\{\epsilon_2^2  \epsilon_1 ^4 N_{22}N_{11}^3 \Big(\sum_{j=3}^n N_{2j}^2 \epsilon_j^2\Big)\Big \}
+3\E\Big\{\epsilon_1^2  \epsilon_2 ^4  N_{11}N_{22}^3 \Big(\sum_{j=3}^n N_{1j}^2 \epsilon_j^2\Big) \Big\}\cr
&&+9\E\Big\{\epsilon_2^2 \epsilon_1^2 N_{11}N_{22} \Big(\sum_{j=3}^n\sum_{k=3,k\neq j}^n N_{1j}^2N_{2k}^2  \epsilon_j^2 \epsilon_k^2  \Big)\Big\} +o(1)\cr
&=& \nu_4^2 \E( N_{11}^3 N_{22}^3) + 3 \nu_4 \sigma_{\epsilon}^4 (n-2) \E(N_{22}N_{11}^3 N_{23}^2)\cr
&&+3\nu_4 \sigma_{\epsilon}^4 (n-2) \E( N_{11}N_{22}^3  N_{13}^2)
 + 9 \sigma_{\epsilon}^8 n^2 \E(N_{11}N_{22} N_{13}^2 N_{24}^2) + o(1)\cr
&=& \E\Big\{\epsilon_2 \Big(\sum_j N_{2j} \epsilon_j\Big)^3 \Big\}\E\Big\{\epsilon_1 \Big(\sum_k N_{1k} \epsilon_k\Big)^3 \Big\}+ o(1).
\end{eqnarray*}
For $\epsilon_2  (\sum_j N_{2j} \epsilon_j )^3 $ and $ (\sum_k N_{1k} \epsilon_k )^4$,
\begin{eqnarray*}
&& \E\Big\{\epsilon_2 \Big(\sum_j N_{2j} \epsilon_j\Big)^3 \Big(\sum_k N_{1k} \epsilon_k\Big)^4 \Big\}\cr
&=&\E\Big\{\epsilon_2^4 N_{22}^3 \Big(\sum_{j\neq 2} N_{1j}^4 \epsilon_j^4\Big) \Big\}
+3 \E\Big[\epsilon_2^2  N_{22} \Big\{\sum_{j\neq 2}\sum_{k\neq 2,j} N_{1j}^4 N_{2k}^2 \epsilon_j^4 \epsilon_k^2 \Big\}\Big]  \cr
&&+3\E\Big[\epsilon_2^4  N_{22}^3 \Big\{\sum_{j\neq 2}\sum_{k\neq 2,j} N_{1j}^2N_{1k}^2 \epsilon_j^2\epsilon_k^2 \Big\}\Big]\cr
&&+9\E\Big[\epsilon_2^2  N_{22} \Big\{\sum_{j\neq 2}\sum_{k\neq 2,j}\sum_{h\neq 2,j,k} N_{1j}^2N_{1k}^2 N_{2h}^2  \epsilon_j^2\epsilon_k^2\epsilon_h^2 \Big\}\Big] + o(1)\cr
&=& \nu_4^2 \E(N_{11}^4 N_{22}^3) + 3\nu_4 \sigma_{\epsilon}^4 (n-2) \E(N_{22} N_{11}^4 N_{23}^2)
+ 3 \nu_4 \sigma_{\epsilon}^4 \E\{N_{22}^3(2nN_{11}^2 N_{13}^2 +n^2 N_{13}^2 N_{14}^2 )\} \cr
&& + 9 \sigma_{\epsilon}^8 \E\{N_{22}(2 n^2 N_{11}^2 N_{13}^2 N_{24}^2 + n^3 N_{15}^2 N_{13}^2 N_{24}^2)\} + o(1)\cr
&=&\E\Big\{\epsilon_2 \Big(\sum_j N_{2j} \epsilon_j\Big)^3\Big\}\E\Big\{ \Big(\sum_k N_{1k} \epsilon_k\Big)^4 \Big\} + o(1).
\end{eqnarray*}
Lastly, for $ (\sum_j N_{2j} \epsilon_j )^4 $ and $ (\sum_k N_{1k} \epsilon_k )^4$,
\begin{eqnarray*}
&& \E\Big\{ \Big(\sum_j N_{2j} \epsilon_j\Big)^4 \Big(\sum_k N_{1k} \epsilon_k\Big)^4 \Big\}\cr
&=&\E \Big\{\sum_j\sum_{k\neq j} N_{1j}^4 N_{2k}^4  \epsilon_j^4  \epsilon_k^4\Big\}
+ 3 \E \Big\{\sum_j\sum_{k\neq j}\sum_{h\neq j,k} (N_{1j}^4 N_{2k}^2N_{2h}^2 +N_{2j}^4 N_{1k}^2N_{1h}^2) \epsilon_j^4  \epsilon_k^2\epsilon_h^2\Big\} \cr
&&+ 9\E \Big\{\sum_j\sum_{k\neq j}\sum_{h\neq j,k}\sum_{\ell\neq j,k,h} N_{1j}^2 N_{1k}^2N_{2h}^2 N_{2\ell}^2\epsilon_j^2  \epsilon_k^2\epsilon_h^2\epsilon_{\ell}^2\Big\} \cr
&=& \nu_4^2 \E ( N_{11}^4 N_{22}^4 ) + 6 \nu_4 \sigma_{\epsilon}^4 \E\{N_{11}^4 (n^2N_{23}^2 N_{24}^2 + 2n N_{22}^2N_{23}^2 ) \}\cr
&& + 9 \sigma_{\epsilon}^8 \E (4 n^2N_{11}^2 N_{22}^2 N_{13}^2 N_{24}^2   + 2 n^3 N_{11}^2 N_{13}^2 N_{24}^2 N_{25}^2 + 2 n^3 N_{22}^2 N_{13}^2 N_{14}^2 N_{25}^2\cr
&& + n^4 N_{13}^2 N_{14}^2 N_{25}^2 N_{26}^2) + o(1) \cr
&=& \E\Big\{ \Big(\sum_j N_{2j} \epsilon_j\Big)^4\Big\} \E\Big\{ \Big(\sum_k N_{1k} \epsilon_k\Big)^4 \Big\} + o(1).
\end{eqnarray*}
We have finished the proof.
\endpf
\end{document}